\documentclass[11pt,a4paper,reqno]{amsart}

\usepackage{amsmath,amssymb,amsthm}
\usepackage[utf8]{inputenc}
\usepackage[T1]{fontenc}
\usepackage{geometry}
\usepackage{hyperref}
\usepackage{xurl}
\hypersetup{
    colorlinks=true,
    linkcolor=blue,
    citecolor=red,
    urlcolor=blue,
}
\usepackage{cleveref}
\usepackage{tikz}
\usetikzlibrary{decorations.pathmorphing,arrows.meta}

\theoremstyle{plain}
\newtheorem{theorem}{Theorem}[section]
\newtheorem{lemma}[theorem]{Lemma}
\newtheorem{proposition}[theorem]{Proposition}
\newtheorem{corollary}[theorem]{Corollary}

\newtheorem{mainthm}{Main Theorem}

\theoremstyle{definition}
\newtheorem{definition}[theorem]{Definition}
\newtheorem{remark}[theorem]{Remark}
\newtheorem{assumption}[theorem]{Assumption}
\newtheorem{problem}[theorem]{Problem}

\numberwithin{equation}{section}

\ifdefined\RobinMacrosLoaded \fi
\let\RobinMacrosLoaded\relax

\newcommand{\Rn}{\mathbb{R}^n}
\newcommand{\Rnm}{\mathbb{R}^{n-1}}          %
\newcommand{\Rnp}{\mathbb{R}^n_+}            %
\newcommand{\xpar}{x'}                       %
\newcommand{\xipar}{\xi'}                    %
\newcommand{\Rhalf}{\mathbb{R}_+}            %
\newcommand{\Sm}{S^{n-1}_+}                  %
\newcommand{\Ssurf}{S^{n-2}}                 %
\newcommand{\Ltwo}[1]{L^2(#1)}
\newcommand{\norm}[1]{\left\lVert #1 \right\rVert}

\newcommand{\Bsp}{\mathcal{B}}
\newcommand{\Bst}{\mathcal{B}^*}
\newcommand{\Bstz}{\mathcal{B}^*_0}
\newcommand{\Bnorm}[1]{\left\lVert #1 \right\rVert_{\Bsp}}
\newcommand{\Bstnorm}[1]{\left\lVert #1 \right\rVert_{\Bst}}
\newcommand{\Bspd}[1]{\mathcal{B}(#1)}       %
\newcommand{\Bstd}[1]{\mathcal{B}^*(#1)}     %
\newcommand{\shell}[1]{\Omega_{#1}}          %
\newcommand{\shellp}[1]{\Omega_{#1}^+}       %

\newcommand{\Halpha}{H_\alpha}               %
\newcommand{\Ralpha}[1]{R_\alpha(#1)}        %
\newcommand{\DN}{\Lambda}                    %
\newcommand{\Dom}{\mathcal{D}}               %
\newcommand{\hform}{\mathfrak{h}_\alpha}     %
\newcommand{\hfib}{h}                        %
\newcommand{\Lop}[2]{\mathcal{L}(#1,#2)}     %
\DeclareMathOperator{\Ran}{Ran}

\newcommand{\sac}{\sigma_{\mathrm{ac}}}
\newcommand{\spp}{\sigma_{\mathrm{pp}}}
\newcommand{\ssc}{\sigma_{\mathrm{sc}}}
\newcommand{\sess}{\sigma_{\mathrm{ess}}}
\newcommand{\Thr}{\mathcal{T}}               %
\newcommand{\Ereg}{\mathcal{E}}              %

\newcommand{\kbulk}{\sqrt{E}}                %
\newcommand{\ksurf}{\sqrt{E+\alpha^2}}       %
\newcommand{\Echan}[1]{\mathcal{H}(#1)}      %
\newcommand{\Hdirect}{\mathcal{H}}           %

\newcommand{\Psib}{\Psi_{\mathrm{bulk}}}     %
\newcommand{\Psis}{\Psi_{\mathrm{surf}}}     %
\newcommand{\Fzero}{\mathcal{F}_0}           %
\newcommand{\Fbulk}{\mathcal{F}_{\mathrm{bulk}}}
\newcommand{\Fsurf}{\mathcal{F}_{\mathrm{surf}}}
\newcommand{\Adens}[1]{A(#1)}                %
\newcommand{\cph}[1]{c(#1)}                  %
\newcommand{\zext}[1]{#1^{\circ}}            %

\newcommand{\Smat}{S}

\newcommand{\refl}[1]{\mathcal{R}(#1)}   %

\newcommand{\usurf}{u_{\mathrm{surf}}}

\newcommand{\psizero}{\psi_0}                %

\newcommand{\Jbulk}{J_{\mathrm{bulk}}}
\newcommand{\Jsurf}{J_{\mathrm{surf}}}
\newcommand{\Jchan}{J}
\newcommand{\yout}{y^{(+)}}                  %
\newcommand{\yin}{y^{(-)}}                   %
\newcommand{\gambulk}{\gamma_{\mathrm{bulk}}}
\newcommand{\gamsurf}{\gamma_{\mathrm{surf}}}

\newcommand{\Ball}[1]{B_{#1}^{+}}

\newcommand{\Nsol}[1]{\mathcal{N}(#1)}

\newcommand{\Wb}[2]{W^{(#1)}_{\mathrm{bulk}}[#2]}
\newcommand{\Ws}[2]{W^{(#1)}_{\mathrm{surf}}[#2]}

\newcommand{\eps}{\varepsilon}
\newcommand{\dist}{\operatorname{dist}}
\DeclareMathOperator{\re}{Re}                %
\DeclareMathOperator{\im}{Im}                %

\ifdefined\RobinMacrosPartTwoLoaded \fi
\let\RobinMacrosPartTwoLoaded\relax

\newcommand{\Om}{\Omega}                     %
\newcommand{\Gpert}{\partial\Om}             %
\newcommand{\Rzero}{R_0}                     %
\newcommand{\Jzero}{J_0}                     %
\newcommand{\Bfull}[1]{B_{#1}}               %
\newcommand{\Omin}[1]{\Om^{\mathrm{in}}_{#1}}%
\newcommand{\Omex}[1]{\Om^{\mathrm{ex}}_{#1}}%
\newcommand{\Nvec}{N}                        %

\newcommand{\Hpert}{H_{\mathrm{pert}}}       %
\newcommand{\Rpert}[1]{R_{\mathrm{pert}}(#1)}%
\newcommand{\hpert}{\mathfrak{h}_{\mathrm{pert}}} %

\newcommand{\Bstzd}[1]{\mathcal{B}^*_0(#1)}          %
\newcommand{\BnormO}[1]{\left\lVert #1 \right\rVert_{\Bspd{\Om}}}
\newcommand{\BstnormO}[1]{\left\lVert #1 \right\rVert_{\Bstd{\Om}}}
\newcommand{\Bnormsh}[1]{{\left\lVert #1 \right\rVert_{\Bspd{\Om}}^{\mathrm{sh}}}}
\newcommand{\Bstnormsh}[1]{{\left\lVert #1 \right\rVert_{\Bstd{\Om}}^{\mathrm{sh}}}}
\newcommand{\Csh}{C_{\mathrm{sh}}}           %
\newcommand{\Qsc}[1]{Q(#1)}                  %

\newcommand{\Eband}{\mathcal{E}_{\mathrm{b}}}  %
\newcommand{\wsig}{w_{\sigma}}
\newcommand{\Cyl}[1]{C_{#1}}                 %
\newcommand{\tsub}{T}                        %

\newcommand{\cut}{\chi}                      %

\newcommand{\Cdiv}{C_{\mathrm{div}}}         %
\newcommand{\Sexc}{\mathcal{S}}

\newcommand{\Phipert}[2]{\Phi^{(#1)}_{\mathrm{pert}}(#2)}

\newcommand{\Fpert}[2]{\mathcal{F}^{(#1)}_{\mathrm{pert}}(#2)}
\newcommand{\Fpertg}[1]{\mathcal{F}^{(#1)}_{\mathrm{pert}}}
\newcommand{\Fbulkpert}[2]{\mathcal{F}^{(#1)}_{\mathrm{bulk},\mathrm{pert}}(#2)}
\newcommand{\Fsurfpert}[2]{\mathcal{F}^{(#1)}_{\mathrm{surf},\mathrm{pert}}(#2)}
\newcommand{\Hac}{\mathcal{H}_{\mathrm{ac}}}
\newcommand{\Adenspert}[1]{A_{\mathrm{pert}}(#1)}
\newcommand{\Ampl}[1]{\mathcal{A}_{\mathrm{pert}}(#1)}

\newcommand{\Nsolpert}[1]{\mathcal{N}_{\mathrm{pert}}(#1)}
\newcommand{\Phase}[1]{\Theta(#1)}
\newcommand{\Smatpert}[1]{S_{\mathrm{pert}}(#1)}
\newcommand{\Smatrel}[1]{S_{\mathrm{rel}}(#1)}
\newcommand{\Sbb}{S_{\mathrm{bb}}}
\newcommand{\Sbs}{S_{\mathrm{bs}}}
\newcommand{\Ssb}{S_{\mathrm{sb}}}
\newcommand{\Sss}{S_{\mathrm{ss}}}

\newcommand{\Sch}[1]{\mathbb{S}_{#1}}
\newcommand{\Psich}[1]{\Psi_{#1}}
\newcommand{\kapch}[2]{\varkappa_{#1}(#2)}

\newcommand{\Phidir}[2]{\Phi^{(#1)}_{#2}}
\newcommand{\Akern}[4]{A_{#1}(#2;#3,#4)}

\newcommand{\Sart}[1]{\Sigma_{#1}}
\newcommand{\Cdata}{\mathcal{C}}

\newcommand{\Hpertj}[1]{H^{(#1)}_{\mathrm{pert}}}
\newcommand{\Rpertj}[2]{R^{(#1)}_{\mathrm{pert}}(#2)}
\newcommand{\Amplj}[2]{\mathcal{A}^{(#1)}_{\mathrm{pert}}(#2)}
\newcommand{\Nsolpertj}[2]{\mathcal{N}^{(#1)}_{\mathrm{pert}}(#2)}
\newcommand{\Sexcj}[1]{\mathcal{S}^{(#1)}}
\newcommand{\Phipertj}[3]{\Phi^{(#2)}_{\mathrm{pert},#1}(#3)}

\newcommand{\Odiff}{\mathcal{O}}
\newcommand{\Gacc}{\Gamma_{\mathrm{acc}}}
\newcommand{\Gfar}{\Gamma_{\mathrm{far}}}
\newcommand{\Omcap}{\Om_{12}}
\newcommand{\Vmix}[1]{V_E(#1)}
\newcommand{\Phidirj}[3]{\Phi^{(#2)}_{#3,#1}}

\newcommand{\Dfat}{D^{\mathrm{1}}}

\begin{document}

\title[Scattering theory for the impedance half-space, II]%
{Scattering theory for bulk and surface waves in an impedance half-space, II:
 local perturbations, the scattering matrix and an inverse problem}

\author{Michiyuki Watanabe}
\address{Department of Applied Mathematics,
Okayama University of Science,
1--1 Ridaicho, Kita--ku, Okayama 700--0005, Japan}
\email{mic-watanabe@ous.ac.jp}
\thanks{ORCID: 0000-0003-0612-8411}

\begin{abstract}
We study the Robin, or impedance, Laplacian on a half-space of dimension
$n \ge 3$ which is perturbed inside a fixed ball. The boundary may be deformed,
and obstacles may be present in the interior or attached to the boundary, in
any combination and with no smallness assumption. Outside the ball the
configuration is the flat impedance half-space of Part~I, with constant
impedance $\alpha > 0$, where a bulk channel and a surface channel propagate at
every positive energy and only the surface channel propagates on the band
$(-\alpha^2,0)$. We prove a limiting absorption principle, valid at every
positive energy and outside a discrete set of embedded eigenvalues in the band.
We construct the generalized Fourier transform and prove its completeness by a
stationary argument which constructs no wave operator and reduces completeness
to the uniqueness of the scattering solution at each energy. We then obtain the
scattering matrix, its unitarity and a two-channel optical theorem. The main
result is a uniqueness theorem for the inverse problem. Under two geometric
conditions on the pair, the scattering amplitude at a \emph{single} positive
energy determines the perturbation completely, and the argument does not
distinguish a deformation of the boundary from an obstacle. It follows that no
nontrivial perturbation satisfying the same conditions is invisible at any
positive energy. On the surface band the corresponding statement is reduced to
one explicitly stated question about the exterior Robin problem.
\end{abstract}

\maketitle

\tableofcontents

\section{Introduction}
\label{sec:intro}
\subsection{Background and Aim}
\label{subsec:intro-perturbed}

Part~I of this work developed the spectral and scattering theory of the Robin,
or impedance, Laplacian on the flat half-space
$\Rnp = \{x = (\xpar,x_n) \in \Rn : x_n > 0\}$, $n \ge 3$. That is the
self-adjoint operator
\begin{equation}
  \label{eq:intro2-flat-operator}
  \Halpha = -\Delta \ \text{ on } \Rnp,
  \qquad
  \partial_\nu u - \alpha\,u = 0 \ \text{ on } \Gamma = \partial\Rnp,
  \qquad \alpha > 0 ,
\end{equation}
where $\nu = -e_n$ is the outward unit normal, so that in the coordinate $x_n$
the boundary condition reads $\partial_{x_n}u + \alpha u = 0$. This is a Robin,
or impedance, condition $\partial_\nu u + \lambda u = 0$ with a negative
impedance $\lambda = -\alpha$, and the sign is decisive. It is the sign for
which the transverse operator $-d^2/dx_n^2$ on the half-line, carrying this
same boundary condition, has an eigenvalue. The boundary carries the bound mode
\begin{equation}
  \label{eq:intro2-psizero}
  \psizero(x_n) = \sqrt{2\alpha}\,e^{-\alpha x_n},
  \qquad\text{of eigenvalue } {-\alpha^2} ,
\end{equation}
confined to a layer of thickness $1/\alpha$ along a boundary which is
unbounded. For the opposite sign there is no such mode. The transverse spectrum
is then purely continuous, and the theory is a variant of the Dirichlet or
Neumann theory in the half-space. From \S\ref{sec:setting} onwards the boundary
condition is written in the form $\gamma_1 u - \alpha\gamma_0 u = 0$ used in
Part~I, with $\gamma_0$ the Dirichlet trace and $\gamma_1 = \partial_\nu$ the
conormal derivative.

The spectrum of $\Halpha$ is accordingly $[-\alpha^2,\infty)$, and a fixed
energy $E$ is carried in two ways at once. The continuum of the transverse
fibres carries it on the ball $|\xipar| \le \kbulk$ of tangential frequencies,
and the bound mode carries it on the single sphere $|\xipar| = \ksurf$. Two
channels are open when $E > 0$. On the surface band $\Eband = (-\alpha^2,0)$
only the second is. At infinity the two appear as spherical waves of different
dimensions, $r^{-(n-1)/2}e^{i\kbulk r}$ into the interior and
$\rho^{-(n-2)/2}e^{i\ksurf\rho}\psizero(x_n)$ along the boundary, with
$r = |x|$ and $\rho = |\xpar|$, and these
rates are each exactly critical for the Agmon--H\"ormander space $\Bst$. Both
waves lie in it with finite nonzero norm, and neither lies in $\Bstz$. That is
why the theory is carried out in the pair $\Bsp$--$\Bst$. Within that framework
Part~I established a limiting absorption principle, the generalized Fourier
transform $\Fzero$ with its two-channel representation space $\Echan{E}$, an
intrinsic radiation condition, a Rellich-type uniqueness theorem, and the
scattering matrix
\begin{equation}
  \label{eq:intro2-flat-smatrix}
  \Smat(E) = \refl{\kbulk\theta_n} \oplus \mathrm{Id}_{\mathrm{surf}} ,
\end{equation}
in which $\theta \in \Sm$ is the direction of the bulk wave and
$\refl{k} = (ik-\alpha)/(ik+\alpha)$ is the reflection coefficient of the flat
boundary. Its surface block is the identity at every energy and for every
$\alpha$.

\begin{figure}[htbp]
\centering
\begin{tikzpicture}[scale=0.92, >=Stealth, line join=round]

\def\bdry{
  (-6.4,0) -- (-3.4,0)
  .. controls (-3.05,0) and (-2.95,1.25) .. (-2.60,1.25)
  .. controls (-2.25,1.25) and (-2.15,0) .. (-1.80,0)
  -- (1.10,0)
  .. controls (1.45,0) and (1.55,-0.90) .. (2.00,-0.90)
  .. controls (2.45,-0.90) and (2.55,0) .. (2.90,0)
  -- (6.4,0)
}

\begin{scope}
  \clip (-6.4,-1.30) rectangle (6.4,4.45);
  \fill[black!5] \bdry -- (6.4,4.45) -- (-6.4,4.45) -- cycle;
\end{scope}

\fill[black!22] (0.35,2.05) ellipse [x radius=0.80, y radius=0.58];
\draw[very thick] (0.35,2.05) ellipse [x radius=0.80, y radius=0.58];
\begin{scope}
  \clip (-6.4,0) rectangle (6.4,4.45);
  \fill[black!22] (-4.25,0) circle (0.58);
  \draw[very thick] (-4.25,0) circle (0.58);
\end{scope}

\draw[very thick] \bdry;

\draw[dashed] (0,0) circle (5.0);
\draw[<->] (0,0) -- node[fill=white,inner sep=1.3pt,pos=0.62,font=\small]
  {$\Rzero$} (-4.10,-2.87);
\node[font=\small,anchor=west] at (4.62,-4.30) {$\partial\Bfull{\Rzero}$};
\draw[thin] (4.50,-4.22) -- (3.410,-3.657);
\fill (0,0) circle (1.3pt);
\node[below,font=\small,inner sep=3pt] at (0.02,-0.06) {$0$};

\node[font=\small] at (-2.55,2.80) {$\Omin{\Rzero}$};
\node[font=\small] at (5.05,2.20) {$\Omex{\Rzero}$};
\node at (2.95,3.35) {$\Om$};

\node[font=\small,anchor=north] at (-5.85,-0.10) {$\Gamma$};
\node[font=\small,anchor=north] at (5.85,-0.10) {$\Gamma$};
\node[font=\small,anchor=north] at (2.00,-1.00) {$\Gpert$};
\draw[thin] (-2.966,0.718) -- (-3.85,1.45);
\node[font=\small,anchor=east] at (-3.90,1.47) {$\Gpert$};

\draw[->,thick] (4.35,0) -- (4.35,-0.70);
\node[font=\small,anchor=west] at (4.42,-0.42) {$\nu$};
\draw[->,thick] (-2.60,1.25) -- (-2.60,0.70);
\node[font=\small,anchor=north] at (-2.60,0.60) {$\nu$};

\draw[->,very thick] (-5.75,4.25) -- (-3.30,1.85);
\node[anchor=west,font=\small] at (-5.30,4.05) {bulk, $\kbulk$};
\draw[->,very thick] (-2.45,1.35) -- (1.35,4.20);
\node[anchor=west,font=\small] at (1.45,4.12) {scattered};

\draw[->,very thick,decorate,
      decoration={snake,amplitude=0.62mm,segment length=3.0mm,post length=2.2mm}]
  (3.25,0.26) -- (6.30,0.26);
\node[anchor=south,font=\small] at (4.80,0.46) {surface, $\ksurf$};

\end{tikzpicture}
\caption{The perturbed half-space $\Om$ of
Assumption~\ref{assump:geometry}. Outside the ball of radius $\Rzero$ the
configuration is the flat impedance half-space of Part~I. The two propagating
channels are those of \eqref{eq:surface-band}, a bulk wave of wavenumber
$\kbulk$ and a surface wave of wavenumber $\ksurf$ travelling along $\Gamma$
with the transverse profile $\psizero$ of \eqref{eq:psizero}. Here
$\Omin{\Rzero} = \Om \cap \Bfull{\Rzero}$ and
$\Omex{\Rzero} = \Rnp \setminus \overline{\Bfull{\Rzero}}$.
Inside the ball the hypothesis allows a deformation of the boundary, obstacles
detached from it and obstacles touching it, in any combination and with no
smallness assumption. What such a perturbation does to the scattering matrix is
the subject of \S\ref{sec:smatrix}.}
\label{fig:geometry}
\end{figure}
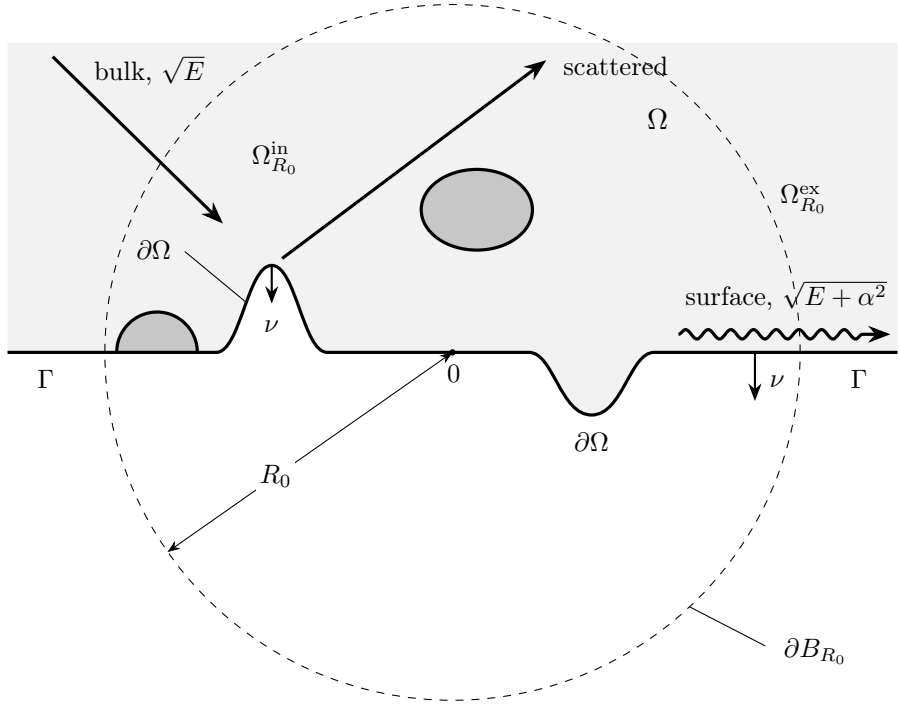

The present Part~II treats a local perturbation of that configuration. The
hypotheses are stated in Assumption~\ref{assump:geometry} and drawn in
Figure~\ref{fig:geometry}. The set $\Om \subset \Rn$ is open and connected, its
boundary is a $C^2$ hypersurface, and it coincides with $\Rnp$ outside a ball
of radius $\Rzero$. The operator $\Hpert$ is
\eqref{eq:intro2-flat-operator} with $\Rnp$ replaced by $\Om$. Inside the ball
nothing further is assumed. A single hypothesis therefore covers an arbitrary
deformation of the boundary, obstacles removed from the interior, obstacles
touching the boundary, and combinations of these, with no smallness,
star-shapedness or topological restriction.

The aim of this paper is to show that the two-channel scattering theory of
Part~I survives a compactly supported geometric perturbation of the half-space,
and that the perturbation is uniquely determined by its scattering amplitude at
a single positive energy. More precisely, we answer the following four
questions.
\begin{enumerate}
  \item[\textup{(Q1)}] How does the perturbation change the spectrum? In
  particular, can it trap waves, that is, create eigenvalues embedded in the
  continuous spectrum?
  \item[\textup{(Q2)}] How is the scattering by the perturbation described?
  \item[\textup{(Q3)}] Does the optical theorem hold, that is, is the total
  energy scattered by the perturbation, into the bulk and along the boundary
  together, given by the imaginary part of the forward scattering amplitude?
  \item[\textup{(Q4)}] Is the perturbation uniquely determined by the
  scattered waves?
\end{enumerate}
The answers are Main Theorems~\ref{mainthm:lap-pert} and
\ref{mainthm:completeness} for (Q1), parts (i) and (ii) of Main
Theorem~\ref{mainthm:smatrix} for (Q2) and (Q3), and Main
Theorem~\ref{mainthm:uniqueness} for (Q4).

Outside $\Bfull{\Rzero}$ the configuration is literally the flat one. The whole
apparatus of Part~I is therefore available there verbatim, with no
perturbative error to control, and that includes the channel space
$\Echan{E}$, the model waves, and the far-field data $\yin$ and $\yout$. This
is what makes a comparison possible, and it is why hypothesis~(G2) is stated
as an exact coincidence rather than as a decay condition.

\subsection{Main Results}
\label{subsec:intro-main}

All results assume $n \ge 3$, $\alpha > 0$ and
Assumption~\ref{assump:geometry}. We write $\Thr = \{-\alpha^2,0\}$ for the
two thresholds, $\Ereg$ for the non-threshold energies of the essential
spectrum, and $\Bspd{\Om}$, $\Bstd{\Om}$ for the Agmon--H\"ormander spaces of
the perturbed domain (\S\ref{subsec:AH-pert}), which coincide with the spaces
obtained by restricting the dyadic shells of $\Rn$ to $\Om$
(Lemma~\ref{lem:AH-equivalence}). The starting point is that $\Hpert$ is
self-adjoint on a domain contained in $H^2(\Om)$, that its essential spectrum
is again $[-\alpha^2,\infty)$, and that below $-\alpha^2$ its spectrum is
discrete (Theorem~\ref{thm:pert-sa-spec}). Below the essential spectrum the
perturbation may create eigenvalues, and we make no claim either way.

\begin{mainthm}[Limiting absorption principle;
Theorem~\ref{thm:lap-pert} and Proposition~\ref{prop:exceptional-discrete}]
\label{mainthm:lap-pert}
The set
\begin{equation}
  \label{eq:intro2-exceptional}
  \Sexc = \spp(\Hpert) \cap \Eband
\end{equation}
is a closed discrete subset of the band $\Eband$, each of whose points is an
eigenvalue of $\Hpert$ of finite multiplicity. For every
$E \in \Ereg \setminus \Sexc$ the boundary values
\[
  \Rpert{E \pm i0}
  = \lim_{\eps \downarrow 0}\,(\Hpert - E \mp i\eps)^{-1}
\]
exist as bounded operators $\Bspd{\Om} \to \Bstd{\Om}$, with a bound uniform on
compact subsets of $\Ereg \setminus \Sexc$, and the pairings
$E \mapsto (\Rpert{E \pm i0}f, g)$ are continuous there. In particular no
exceptional energy occurs at positive energy, and the principle holds without
restriction on the whole of $(0,\infty)$.
\end{mainthm}

At positive energy $\Hpert$ has no eigenvalue at all
(Corollary~\ref{cor:global-uniqueness}(a)), so no wave is trapped there. On the
band, eigenvalues embedded in the continuous spectrum are not excluded, and no
hypothesis is imposed to exclude them. What is proved instead is that they
form a discrete set, and that suffices for everything that follows.

The limiting absorption principle yields an intrinsic radiation condition on
$\Om$ (Definition~\ref{def:radiation-pert}, Theorem~\ref{thm:radiation-equivalence})
and with it the scattering solution. For every incident wave $\phi$, that is,
for every element of the space $\Nsol{E}$ of $\Bst$-solutions of the
\emph{unperturbed} problem, there is exactly one solution of the perturbed
problem whose scattered part is outgoing
(Theorem~\ref{thm:scattering-solution}). Writing $\Phipert{\pm}{E}\phi$ for it,
one obtains the perturbed generalized Fourier transform and, with it, the
following.

\begin{mainthm}[Spectral representation and completeness;
Theorem~\ref{thm:spectral-rep-pert} and Proposition~\ref{prop:ac-subspace-pert}]
\label{mainthm:completeness}
For each choice of sign there is a bounded operator
\[
  \Fpertg{\pm} : \Ltwo{\Om} \longrightarrow
  \Hdirect = \int_{\Ereg}^{\oplus} \Echan{E}\,dE ,
\]
given on $\Bspd{\Om}$ by $f \mapsto \bigl( E \mapsto \Fpert{\pm}{E}f \bigr)$,
which is a partial isometry with initial space $\Hac(\Hpert)$ and final space
\emph{all} of $\Hdirect$, and which diagonalizes $\Hpert$. In particular the
absolutely continuous part of $\Hpert$ is unitarily equivalent to $\Halpha$.
Moreover $\ssc(\Hpert) = \varnothing$ and
$\spp(\Hpert) \cap \Ereg = \Sexc$.
\end{mainthm}

Main Theorems~\ref{mainthm:lap-pert} and \ref{mainthm:completeness} together
answer (Q1). The perturbation leaves the absolutely continuous spectrum and its
multiplicity unchanged (Corollary~\ref{cor:ac-spectrum-pert}) and creates no
singular continuous spectrum. It traps no wave at positive energy, and in the
band the trapped waves are confined to the discrete set $\Sexc$, whose
emptiness is not decided here. The substantive half of
Main Theorem~\ref{mainthm:completeness} is the surjectivity, and it is proved
without constructing a wave operator (\S\ref{subsec:intro-novelty}).

\begin{mainthm}[The scattering matrix and the optical theorem;
Theorems~\ref{thm:smatrix-far-field} and \ref{thm:optical}]
\label{mainthm:smatrix}
Let $E \in \Ereg \setminus \Sexc$.
\begin{enumerate}
  \item[\textup{(i)}] There is a unitary operator $\Smatpert{E}$ on
  $\Echan{E}$ such that every $\Bst$-solution $u$ of $(\Hpert - E)u = 0$
  satisfies $\yout[u] = \Smatpert{E}\,\Jchan\,\yin[u]$, where $\Jchan$ is the
  involution of \cite[Definition 2.57]{PartI}, and
  \begin{equation}
    \label{eq:intro2-smatrix}
    \Smatpert{E} = \Smat(E) + 2\pi i\,\Phase{E}\,\Ampl{E}\,\Phase{E} ,
  \end{equation}
  where $\Smat(E)$ is the unperturbed matrix \eqref{eq:intro2-flat-smatrix},
  $\Ampl{E}$ is the scattering amplitude of
  Definition~\ref{def:scattering-amplitude}, and $\Phase{E}$ is the multiplier
  with $\Phase{E}^2 = \Smat(E)$.
  \item[\textup{(ii)}] The amplitude obeys
  \begin{equation}
    \label{eq:intro2-optical}
    \tfrac{1}{2i}\bigl( \Ampl{E} - \Ampl{E}^{*} \bigr)
      = \pi\,\Ampl{E}^{*}\Ampl{E} .
  \end{equation}
  Equivalently, for an incident wave of channel datum $g$ whose scattered
  field has outgoing amplitudes $A$ and $C$,
  \[
    \im\bigl( \Ampl{E}g,\, g \bigr)_{\Echan{E}}
    = \kbulk\,\norm{A}_{\Ltwo{\Sm}}^{2} + \ksurf\,\norm{C}_{\Ltwo{\Ssurf}}^{2} ,
  \]
  the first term on the right being absent when $E \in \Eband$.
\end{enumerate}
\end{mainthm}

Part~(i) answers (Q2). Identity \eqref{eq:intro2-smatrix} is the analogue of
the classical relation between the scattering matrix and the scattering
amplitude, with the unperturbed matrix in the place of the identity, and it
says that $\Smatpert{E} = \Smat(E)$ precisely when $\Ampl{E} = 0$. Nothing
beyond this equivalence is asserted about the individual blocks. Part~(ii)
answers (Q3). The weights on the right are the channel wavenumbers, and the
identity accounts for the whole of the energy removed from the incident wave,
whether it departs as a bulk wave or travels away along the boundary.

For (Q4) the datum is the scattering amplitude at a single energy, which by
\eqref{eq:intro2-smatrix} is interchangeable with the scattering matrix. Two
configurations are to be compared, and besides
Assumption~\ref{assump:geometry} for each, two conditions on the pair are
needed.

\begin{mainthm}[Uniqueness of the perturbation;
Theorem~\ref{thm:main-uniqueness}]
\label{mainthm:uniqueness}
Let $\Om_1$ and $\Om_2$ both satisfy Assumption~\ref{assump:geometry}, with the
same radius $\Rzero$, the same impedance $\alpha > 0$ and the same dimension
$n \ge 3$, and assume that
\begin{enumerate}
  \item[\textup{(U1)}] $\Om_1 \cap \Om_2$ is connected, and
  \item[\textup{(U2)}] every connected component of
  $\Om_1 \setminus \overline{\Om_2}$ and of $\Om_2 \setminus \overline{\Om_1}$
  has Lipschitz boundary.
\end{enumerate}
If $\Amplj{1}{E} = \Amplj{2}{E}$ for one single energy $E > 0$, then
$\Om_1 = \Om_2$.
\end{mainthm}

The data are those of one positive energy and the conclusion is the whole
perturbation. No smallness, star-shapedness or graph condition is imposed, and
the complements need not be connected. Above all, \emph{a deformation of the
boundary and an obstacle in the interior are not distinguished}. One might
expect a uniqueness theorem for obstacles above a perturbed boundary to assume
the boundary already known, and the proof shows that this is unnecessary. Both
(U1) and (U2) are conditions on the pair, and Remark~\ref{rem:hypotheses}
explains where each is used. Taking one of the two configurations to be $\Rnp$,
whose amplitude vanishes, gives Corollary~\ref{cor:no-invisible}. If
$\Om \ne \Rnp$, then $\Ampl{E} \ne 0$, and hence $\Smatpert{E} \ne \Smat(E)$,
at \emph{every} positive energy. In that case (U2) asks that each piece added
to or removed from the half-space have Lipschitz boundary, and (U1) that the
remainder of the half-space not be disconnected.

On the band the same argument runs unchanged except at a single point, the
Rellich-type theorem. What results there is a conditional theorem
(Theorem~\ref{thm:band-conditional}) together with one explicitly stated open
problem (Problem~\ref{prob:band-rellich}), to which \S\ref{subsec:intro-open}
returns.

\subsection{Method}
\label{subsec:intro-method}

The chain of results is the one used in Part~I and, before it, in the
stationary multichannel theory from which the template is taken
\cite{Isozaki2001,IsozakiKurylev2014}. It runs from a uniqueness theorem near
infinity to a limiting absorption principle, an eigenfunction expansion, a
radiation condition and a scattering matrix, each link supplying the
hypotheses of the next. The perturbed configuration is compared with the flat
one through a cut-off $\eta$ equal to $1$ near infinity. Comparison through a
cut-off is standard for geometric perturbations, see
\cite{IsozakiKadowakiWatanabe2025JOS} and the references there.

\subsubsection*{The Rellich-type theorem as the axis}

One theorem carries most of the weight. The Rellich-type theorem of Part~I
imposes its hypothesis only near infinity, so it applies verbatim in
$\Omex{\Rzero}$, where the configuration is exactly the flat one
(Theorem~\ref{thm:rellich-positive}). Combined with unique continuation it
becomes a statement about all of $\Om$
(Corollary~\ref{cor:global-uniqueness}(a)). Unique continuation is cheap here,
since the equation has constant coefficients in the interior and a solution is
therefore real analytic there by the theorem of Morrey and Nirenberg
\cite{MorreyNirenberg1957} (Proposition~\ref{prop:ucp}). The resulting
uniqueness statement is then used three times, at the three decisive points of
the paper.
\begin{enumerate}
  \item[\textup{(a)}] It gives the uniqueness of outgoing solutions
  (Theorem~\ref{thm:outgoing-uniqueness}), which is the one nonelementary input
  of the compactness argument of \cite[Ch.~2, \S2.5]{IsozakiKurylev2014} and hence of
  the limiting absorption principle.
  \item[\textup{(b)}] It gives the uniqueness of the scattering solution, that
  is, the injectivity of $\Phipert{\pm}{E}$ at each energy
  (Corollary~\ref{cor:scattering-parametrization}). This is what
  completeness is reduced to.
  \item[\textup{(c)}] It carries the equality of the scattering amplitudes over
  to the equality of the scattered fields outside the perturbation
  (Theorem~\ref{thm:exterior-equality}), which is where the inverse problem
  begins.
\end{enumerate}
Completeness, the scattering matrix and the inverse problem thus rest on one
uniqueness statement, which at positive energy is the Rellich-type theorem.
That theorem was named in \cite{IsozakiKadowakiWatanabe2025JOS} as
the first step that an inverse problem for the elastic half-space would
require.

\subsubsection*{The limiting absorption principle and the surface band}

Uniqueness of outgoing solutions turns on one identity
(Lemma~\ref{lem:zero-flux}). Because $\alpha$ is real, the boundary integrand
$u\,\overline{\partial_\nu u} - \bar u\,\partial_\nu u$ vanishes identically on
$\Gpert$, so a solution on all of $\Om$ has spherical flux with vanishing
imaginary part. Part~I's Stone formula and shell-average identity then force
the channel datum to vanish, whence the solution lies in $\Bstzd{\Om}$. Nothing
in this passes through the positivity of the energy. On the band, however, the
Rellich-type theorem is unavailable, and a $\Bstzd{\Om}$-solution at a band
energy is only shown to decay exponentially, and hence to be an
$L^2$-eigenfunction (Theorem~\ref{thm:band-equivalence}). The decay comes from
an Agmon-type estimate with a slowly varying tangential weight. It is uniform
on compact subintervals, and the compactness it confers is what makes $\Sexc$
discrete (Proposition~\ref{prop:exceptional-discrete}). A parametrix and the
analytic Fredholm theorem would not give this discreteness
(Remark~\ref{rem:no-fredholm}).

\subsubsection*{Transforms, completeness and the scattering matrix}

An incident wave is a function on $\Rnp$, and $\Om$ need not be contained in
$\Rnp$. The total field is therefore built as $u = \eta\phi + v$, and the
relation $\simeq$ of asymptotic equality, which does not see a compactly
supported difference, makes all far-field statements independent of $\eta$.
The perturbed transform is defined through its adjoint,
$\Fpert{\pm}{E}^{*} = \Phipert{\pm}{E}\circ\Fzero(E)^{*}$, which quotes
$\Fzero(E)^{*}$ of Part~I unchanged. Everything downstream rests on an explicit
formula expressing $\Fpert{\mp}{E}f$ through the free transform of the
transplanted source $\eta f + [-\Delta,\eta]\Rpert{E \pm i0}f$
(Theorem~\ref{thm:Fpert-explicit}), and on the antisymmetry of the commutator
$[-\Delta,\eta]$ (Lemma~\ref{lem:commutator-antisym}), which uses neither the
equation nor the boundary condition.

The isometry of $\Fpertg{\pm}$ follows from Stone's formula. For the
surjectivity, the range of $\Fpertg{\pm}$ is closed and invariant under
multiplication by bounded functions of the energy. A vector orthogonal to it
is therefore orthogonal, at almost every energy, to the range of the
transform at that energy. It vanishes as soon as that range is dense. That is
a statement at a single energy, and by duality it is (b) above.
The scattering matrix is defined by its action on
the far-field data, and \eqref{eq:intro2-smatrix} is proved by comparing the
two normalizations directly. The optical theorem is then the unitarity of
$\Smatpert{E}$ written out.

\subsubsection*{The inverse problem}

The inverse problem is Schiffer's argument, adapted twice. The boundary
condition is of impedance type, so the overdetermined problem on the region
where the two configurations differ is a mixed one, carrying $\alpha$ with
opposite signs on the two parts of its boundary. The sign reversal is harmless.
The boundary term is subordinate to the Dirichlet integral whatever its sign,
so the solution space is still finite dimensional (Lemma~\ref{lem:finite-dim}).
The second adaptation concerns the supply of infinitely many solutions. In the
classical argument it comes from a density statement, which does not transfer.
What replaces it is the linear independence of the directional eigenfunctions
(Lemma~\ref{lem:independence}), and here an obstruction of the direct problem
becomes the mechanism. The generalized eigenfunctions
$\Psich{\mu}(\cdot,\theta)$ do not lie in $\Bst$, so a nontrivial finite
combination of directional eigenfunctions cannot lie in $\Bstd{\Om}$, which it
would have to do if it vanished on an open set. The surface directions alone
already supply infinitely many, since $\Ssurf$ is infinite for $n \ge 3$.

\subsection{Novelty and Relation to Earlier Work}
\label{subsec:intro-novelty}

\subsubsection*{Scattering on the perturbed half-space}

No novelty is claimed for the Agmon--H\"ormander framework
\cite{Agmon1975,AgmonHormander1976}, which Part~I already adopted, nor for the
individual techniques of \S\ref{sec:lap-pert}. For the half-space carrying an
impedance condition, the outgoing problem was analysed by Dur\'an, Muga and
N\'ed\'elec \cite{DuranMugaNedelec2009} in weighted $L^2$ spaces, for the flat
boundary at positive energy. Perturbed stratified media and open waveguides
have a substantial literature
\cite{Weder1991,Zhang1994,Zhang2000,BonnetBenDhiaDakhia2009,Hazard2015,%
AlHumaikani2025,IsozakiKadowakiWatanabe2020}. From it we use the dichotomy of
Hazard \cite{Hazard2015} between open and closed waveguides. On the band the
bulk medium lies below cutoff and $\psizero$ is the only propagating mode, so
the configuration is closed in that sense, and the argument by which trapped
modes are excluded in open waveguides is unavailable there
(Remark~\ref{rem:band-closed-waveguide}). The spectral theory of the
Robin Laplacian on unbounded domains is highly developed
\cite{LevitinParnovski2008,BehrndtLanger2007,%
LotoreichikRohleder2012,BehrndtSchlosser2018,CossettiKrejcirik2020}. It
concerns eigenvalue asymptotics and the dependence of the spectrum on the
geometry rather than scattering at a fixed energy.

The closest analogue of the whole picture is the elastic half-space of Isozaki,
Kadowaki and Watanabe \cite{IsozakiKadowakiWatanabe2023,%
IsozakiKadowakiWatanabe2025,IsozakiKadowakiWatanabe2025JOS}, where body waves
and Rayleigh surface waves occur together and the boundary is itself locally
perturbed. There the asymptotic expansion at infinity, a radiation condition
adapted to the two wave types and the uniqueness of the scattering solutions
are established. The extension of the Rellich-type theorem to the elastic
half-space is recorded in \cite[\S5]{IsozakiKadowakiWatanabe2025JOS} as an open
problem, and as the first step that an inverse problem would require. The
present paper takes that step in the corresponding scalar situation
(Theorem~\ref{thm:rellich-positive}), and \S\ref{sec:inverse} is what having it
makes possible.

\subsubsection*{Completeness without wave operators}

The chain of results follows the stationary multichannel template of
\cite{Isozaki2001,IsozakiKurylev2014} with one exception, the proof of
completeness. In \cite{Isozaki2001} completeness is obtained from the
asymptotic completeness of time-dependent channel wave operators, and in
\cite[Ch.~2, \S4]{IsozakiKurylev2014} from stationary wave operators built from the
resolvents of two operators acting on one Hilbert space with a common domain.
Main Theorem~\ref{mainthm:completeness} is proved by neither route, and no wave
operator is constructed anywhere in this paper. Since the transform
diagonalizes the operator, a vector orthogonal to its range is orthogonal, at
almost every energy, to the range of the transform at that energy. The range
is therefore all of the direct integral as soon as it is dense at almost every
energy. That is a statement at a single energy, where it follows by duality
from the uniqueness of the scattering solution, and so, at positive energy,
from the Rellich-type theorem.
Completeness thus requires nothing beyond the fixed-energy theory which the
scattering matrix needs in any case. The argument uses no identification of
the two dynamics, which the time-dependent theory
\cite{Kato1976,ReedSimonIV,Yafaev} would require, and no hypothesis excluding
embedded eigenvalues in the band. In its structure the argument is the
stationary one of Agmon \cite[Theorem 6.2]{Agmon1975}, and Weder
\cite[Ch.~2, Lemma 7.1]{Weder1991} uses it for perturbed stratified media. The
surface channel enters only through the fixed-energy theory, in particular
through the characterization of solutions of $\Bst$ growth
\cite[Theorem 2.68]{PartI}. We are not aware of a proof of this kind for a
geometric perturbation in which a surface channel propagates alongside the
bulk one.

\subsubsection*{The inverse problem and the reflection principle}

Inverse problems in perturbed half-spaces were considered by Eskin and Ralston
\cite{EskinRalston1999} for the coefficients of a second-order elliptic operator
under a Dirichlet condition. There the perturbed boundary, a graph tending to
zero at infinity, is known, and a coefficient is recovered from the
Dirichlet-to-Neumann map at a fixed energy. Here the boundary itself is the
unknown, and the data are given at infinity.
For a locally perturbed half-space carrying a Dirichlet or a Neumann condition,
the corresponding uniqueness theorem is due to Zhao, Hu and Yan
\cite{ZhaoHuYan2021}. They recover a locally perturbed Lipschitz graph in three
dimensions from the far-field patterns of infinitely many incident plane waves,
and their argument rests on the mixed reciprocity relation together with the
\emph{reflection principle}. A solution of the Helmholtz equation vanishing on
a flat Dirichlet boundary, or having vanishing normal derivative on a flat
Neumann boundary, extends to a solution on $\Rn$ by odd, respectively even,
reflection.

That device is not available here. Zhao, Hu and Yan remark
\cite[Remark~2.1]{ZhaoHuYan2021} that the impedance case can be treated with
the corresponding reflection principle. That principle reflects
$v = \partial_{x_n}u + \alpha u$, which vanishes on $\Gamma$, and recovers $u$
below $\Gamma$ by integrating $\partial_{x_n}u + \alpha u = v$ along the
normals. The operator $\partial_{x_n} + \alpha$ annihilates
$g(\xpar)\psizero(x_n)$, so the surface wave is invisible to $v$, and the
integration reintroduces it as a term growing like $e^{\alpha|x_n|}$ below
$\Gamma$. The extended function is then not a radiating solution on $\Rn$. For
the opposite sign there is no bound mode, and this obstruction does not arise.
The second channel and the failure of the reflection principle are two
descriptions of one fact. What takes the place
of the reflection principle is the Rellich-type theorem transplanted in
Theorem~\ref{thm:rellich-positive}. The mixed reciprocity relation is not
needed for the uniqueness argument, and its role in producing enough solutions
is taken over by Lemma~\ref{lem:independence}. It is nevertheless available
here (\S\ref{subsec:inverse-kernels}), where it follows from the antisymmetry
of the commutator without an asymptotic expansion of a fundamental solution
(Remark~\ref{rem:kernel-scope}).

Three further differences should be noted. The data are, as in
\cite{ZhaoHuYan2021}, those of a single energy, but no graph hypothesis is
imposed on the perturbed boundary. No distinction is made between a deformation of the
boundary and an obstacle, which in the obstacle scattering literature are
ordinarily the subject of separate theorems. The argument is Schiffer's
\cite[\S5.1]{ColtonKress2013}, \cite[Ch.~6]{Isakov2017}, and what is new in the
adaptation is the mixed problem with reversed signs and the replacement of the
density argument by the linear independence of the directional
eigenfunctions.

The limitations should be set beside these. The hypotheses (U1) and (U2) are
conditions on the pair rather than on either domain, and (U2) refers to a set
built from the two unknowns. The boundary is assumed $C^2$ where
\cite{ZhaoHuYan2021} assumes only a Lipschitz graph. The impedance $\alpha$ is
constant and is assumed known. Finally, the theorem says nothing on the surface
band, and \S\ref{subsec:intro-open} is devoted to that.

\subsection{The Surface Band and Open Problems}
\label{subsec:intro-open}

The incompleteness of this paper is located on the surface band, and sharply.
Everything in \S\S\ref{sec:setting}--\ref{sec:smatrix} holds on
$\Eband \setminus \Sexc$ as it holds at positive energy, the bulk terms simply
dropping out. In \S\ref{sec:inverse} every step except one is proved for every
$E \in \Ereg$. Positivity of the energy is used only in
Theorem~\ref{thm:exterior-equality}, and there only through the Rellich-type
theorem. Supply that one implication and Main Theorem~\ref{mainthm:uniqueness}
extends to the band. That is the content of Theorem~\ref{thm:band-conditional},
and the missing implication is written out as Problem~\ref{prob:band-rellich}.

A positive answer would determine the entire perturbation, boundary deformation
and obstacles alike, from the surface-to-surface block of the scattering matrix
at a single band energy, that is, from data collected without launching a wave
into the bulk. It would also reverse the role of the surface channel. By
\eqref{eq:intro2-flat-smatrix} that channel is rigid on the flat half-space and
carries no information about $\alpha$, whereas it would then carry the whole of
the geometry. Section~\ref{subsec:surface-only} sets out what we know about the
question. The exceptional set is a second, more elementary open question. We
prove that $\Sexc$ is discrete but do not decide whether it is ever nonempty. A
resonance theory for the surface channel, that is, a meromorphic continuation
of $\Rpert{z}$ across $\Eband$, would presumably settle both this and the
discreteness at once. Part~I develops a Mourre estimate for $\Halpha$ which is
a natural starting point for such a continuation. It is used nowhere in the
present paper, and we do not claim that it survives the perturbation.

Two directions outside the band remain open. The first is a variable
impedance $\alpha = \alpha(\xpar)$, constant outside a compact set, for which
the exact fibration on which Part~I rests would be lost. The inverse problem in
that setting is the natural common generalization of the theorem of Part~I,
which recovers a constant $\alpha$ from the bulk block, and of Main
Theorem~\ref{mainthm:uniqueness}, which recovers the geometry with $\alpha$
known. The second is the elastic half-space discussed above. The mechanisms
used here are not specific to the scalar setting, although the elastic problem
has one feature with no scalar analogue, the singularities carried by the
critical-angle cone between the two body-wave speeds.

\subsection{Organization of the Paper}
\label{subsec:intro-plan}

Section~\ref{sec:setting} fixes the notation, states the geometric hypothesis
and its three configurations (Assumption~\ref{assump:geometry},
Figure~\ref{fig:geometry}), defines $\Hpert$ and establishes its
self-adjointness and the location of its essential spectrum
(Theorem~\ref{thm:pert-sa-spec}).
Section~\ref{sec:lap-pert} is the analytic core.
\S\ref{subsec:AH-pert} defines the Agmon--H\"ormander spaces of $\Om$ and
identifies them with the restricted ones (Lemma~\ref{lem:AH-equivalence}).
\S\ref{subsec:rellich-pert} transplants the Rellich-type theorem
(Theorem~\ref{thm:rellich-positive}) and proves the band reduction
(Theorem~\ref{thm:band-equivalence}), from which the discreteness of $\Sexc$
follows (Proposition~\ref{prop:exceptional-discrete}).
\S\ref{subsec:outgoing} proves the uniqueness theorem for outgoing solutions,
and \S\ref{subsec:lap-pert-main} the limiting absorption principle
(Main Theorem~\ref{mainthm:lap-pert}).
Section~\ref{sec:scattering-solution} converts the limiting absorption
principle into a radiation condition
(Definition~\ref{def:radiation-pert}, Theorem~\ref{thm:radiation-equivalence}).
\S\ref{subsec:incident} identifies the incident waves as the $\Bst$-solutions
of the unperturbed problem, and \S\ref{subsec:scattering-sol} constructs the
scattering solution (Theorem~\ref{thm:scattering-solution}).
Section~\ref{sec:gft-pert} defines the perturbed generalized Fourier transform,
proves the transplanted source formula
(Theorem~\ref{thm:Fpert-explicit}), and establishes the spectral
representation and completeness (Main Theorem~\ref{mainthm:completeness})
together with the far field and the scattering amplitude.
Section~\ref{sec:smatrix} constructs the scattering matrix, relates it to the
amplitude and proves unitarity and the optical theorem
(Main Theorem~\ref{mainthm:smatrix}).
Section~\ref{sec:inverse} treats the inverse problem.
\S\ref{subsec:inverse-kernels} constructs the kernels and proves the mixed
reciprocity relation, and \S\ref{subsec:inverse-exterior} recovers the field
outside the perturbation from the data
(Theorem~\ref{thm:exterior-equality}). \S\ref{subsec:schiffer} proves
Main Theorem~\ref{mainthm:uniqueness}, and \S\ref{subsec:surface-only}
examines what the surface channel determines on its own.
Appendix~\ref{app:regularity} contains the elliptic regularity and trace
statements used for the perturbed domain, and Appendix~\ref{app:ucp} the unique
continuation statement.

\section{The Perturbed Half-Space and the Robin Laplacian}
\label{sec:setting}
In this section we fix the notation, describe the class of perturbed domains
that the paper treats, define the perturbed Robin Laplacian, and establish its
self-adjointness and the location of its essential spectrum. Everything here is
preparatory: the analysis proper begins in \S\ref{sec:lap-pert}. At this stage
the only difference from Part~I lies in the lower bound of the quadratic form,
discussed in Remark~\ref{rem:bottom}.

\subsection{Notation and Function Spaces}
\label{subsec:pert-notation}

The unperturbed configuration is the one studied in Part~I. We write
\begin{equation}
  \label{eq:halfspace}
  \Rnp = \{x = (\xpar, x_n) \in \Rn : x_n > 0\},
  \qquad n \ge 3,
  \qquad
  \Gamma = \partial\Rnp \cong \Rnm,
\end{equation}
and a point of $\Rn$ is always split as $x = (\xpar, x_n)$ with
$\xpar = (x_1,\dots,x_{n-1})$ the tangential variable and $x_n$ the normal
variable. The dual
variable is split in the same way, $\xi = (\xipar, \xi_n)$. Accordingly
\begin{equation*}
  \nabla = (\nabla', \partial_{x_n}),
  \qquad
  \nabla' = (\partial_{x_1},\dots,\partial_{x_{n-1}}),
  \qquad
  \Delta = \Delta' + \partial_{x_n}^2 .
\end{equation*}
On $\Gamma$ the outward unit normal of $\Rnp$ is $\nu = -e_n$, so that
$\partial_\nu = -\partial_{x_n}$ there. We write $\Bfull{R} = \{x \in \Rn :
|x| < R\}$ for the open ball of radius $R$ centred at the origin, and $\Ball{R}
= \Bfull{R} \cap \Rnp$ for the corresponding half-ball.

The letter $\alpha > 0$ is the impedance, or Robin, parameter, a constant fixed
once and for all. Since $\alpha$ is thus unavailable, multi-indices are
written $\beta$. The spectral parameter is denoted $E$ when real and $z$ when
complex. The hat is reserved for the Fourier transform. Directions on the unit
spheres are written $\theta \in \Sm$ for the bulk hemisphere and
$\omega \in \Ssurf$ for the surface sphere. Finally $C$ denotes a generic
positive constant whose value may change from line to line. Such constants
never depend on the functions involved, and their dependence on the energy, on
$\alpha$ and on the geometry is displayed whenever it matters.

Because the half-space is invariant under translations in $\xpar$, the Fourier
transform we use is the partial one,
\begin{equation}
  \label{eq:fourier-partial}
  (\mathcal{F}' u)(\xipar, x_n) = \widehat{u}(\xipar, x_n)
    = \left(\frac{1}{2\pi}\right)^{(n-1)/2}
       \int_{\Rnm} e^{-i \xpar \cdot \xipar}\, u(\xpar, x_n)\, d\xpar ,
\end{equation}
normalized so as to be unitary on $\Ltwo{\Rnp} = L^2(\Rnm_{\xpar};
\Ltwo{\Rhalf})$. Once the domain is perturbed, this transform no longer
diagonalizes the operator. That is the central difficulty of Part~II. The
transform remains available in the exterior region, where it is used
throughout.

\begin{remark}[Standing assumption on the dimension]
\label{rem:standing-dimension}
Throughout this paper $n \ge 3$. Wherever a result holds for every $n \ge 2$ we
state it in that generality and say so. This is the case for the whole of the
present section. The restriction is imposed by the Rellich-type theorem of
Part~I, which is quoted in \S\ref{sec:lap-pert} and \S\ref{sec:scattering-solution}
and which requires $n - 1 \ge 2$: its argument applies the classical Rellich
theorem in dimension $n-1$ to the surface component and then propagates the
vanishing through the exterior region $\{\xpar \in \Rnm : |\xpar| > \Rzero\}$,
and for $n = 2$ that region is a pair of half-lines, along which there is
nothing to propagate. A second, milder reason is that the set of surface
directions $\Ssurf$ degenerates to the two-point set $S^0$ when $n = 2$.
\end{remark}

All function spaces are complex, and derivatives are taken in the weak sense.
For an open set $D \subset \Rn$ we use $\Ltwo{D}$ with the inner product
$(u,v)_{\Ltwo{D}} = \int_D u \overline{v}\, dx$, and the Sobolev spaces
$H^m(D)$, $m \in \mathbb{N}$, with the norm
$\norm{u}_{H^m(D)} = \bigl(\sum_{|\beta| \le m}
\norm{\partial^\beta u}_{\Ltwo{D}}^2\bigr)^{1/2}$. On a boundary of class $C^2$
the trace operators are available in the usual form: $\gamma_0 u = u|_{\partial
D}$ extends to a bounded surjection $H^1(D) \to H^{1/2}(\partial D)$, and
$\gamma_1 u = (\nu \cdot \nabla u)|_{\partial D}$ is bounded from $H^2(D)$ to
$H^{1/2}(\partial D)$, where $\nu$ is the outward unit normal
\cite[Ch.~1]{Grisvard1985}. We use the same letters $\gamma_0, \gamma_1$ for
the traces on $\Gamma$ and on the perturbed boundary introduced below. Which is
meant is always clear from the domain.

The Agmon--H\"ormander spaces $\Bsp$ and $\Bst$ in which the limiting
absorption principle is formulated are Besov-type Banach spaces, not weighted
Sobolev spaces, and their definition on the perturbed domain requires a
dyadic decomposition adapted to $\Om$. We postpone it to
\S\ref{sec:lap-pert}, where it is first needed.

\subsection{The Perturbed Half-Space}
\label{subsec:pert-domain}

The geometry treated in this paper is a perturbation of $\Rnp$ that is
localized in a ball and exact outside it.

\begin{assumption}[The perturbed half-space]
\label{assump:geometry}
$\Om \subset \Rn$ is a connected open set with the following two properties.
\begin{enumerate}
  \item[\textup{(G1)}] The boundary $\Gpert$ is a hypersurface of class $C^2$.
  \item[\textup{(G2)}] There is a radius $\Rzero > 0$ such that $\Om$ coincides
  with the half-space outside the ball of that radius:
  \begin{equation}
    \label{eq:exterior-coincidence}
    \Om \setminus \overline{\Bfull{\Rzero}}
      \;=\; \Rnp \setminus \overline{\Bfull{\Rzero}} .
  \end{equation}
\end{enumerate}
We fix such a radius $\Rzero$ once and for all, and abbreviate
\begin{equation*}
  \Omin{R} = \Om \cap \Bfull{R},
  \qquad
  \Omex{R} = \Om \setminus \overline{\Bfull{R}}
  \qquad (R \ge \Rzero),
\end{equation*}
so that $\Omex{R} = \Rnp \setminus \overline{\Bfull{R}}$ is a subset of the
unperturbed half-space for every $R \ge \Rzero$.
\end{assumption}

Assumption~\ref{assump:geometry} covers simultaneously the three configurations
of interest, and no further hypothesis distinguishes them. The first is an
arbitrary $C^2$ deformation of the flat boundary inside $\Bfull{\Rzero}$. The
second is a finite family of disjoint compact obstacles with $C^2$ boundary,
removed from $\Bfull{\Rzero} \cap \Rnp$, whether or not they touch the
boundary. The third is any combination of the two. In particular neither $\Gpert$ nor $\Rn \setminus \Om$
need be connected, and the perturbation is not required to be small in any
sense. Connectedness of $\Om$ itself is used only in \S\ref{sec:lap-pert} and
\S\ref{sec:inverse}, where a solution is continued across the perturbed region
and where a set that is both open and closed is identified with the whole
domain.

Hypothesis~(G1) is what places the domain of the operator inside $H^2(\Om)$,
through the elliptic regularity of Appendix~\ref{app:regularity} and the trace
estimate of Lemma~\ref{lem:trace-pert}. Hypothesis~(G2) is stated as an exact
coincidence rather than an asymptotic one, and this is a deliberate
restriction: all the objects constructed in Part~I are then available in
$\Omex{\Rzero}$ verbatim, with no perturbative error to control, and the whole
difficulty is confined to the compact region $\Omin{\Rzero}$. Perturbations
decaying at infinity, and a non-constant impedance $\alpha = \alpha(x)$, are
natural next steps and are not pursued here.

\subsection{The Perturbed Robin Laplacian}
\label{subsec:pert-operator}

Let $\nu$ denote the outward unit normal of $\Om$ along $\Gpert$. By (G1) it is
a $C^1$ vector field, and by (G2) it equals $-e_n$ at every boundary point
outside $\Bfull{\Rzero}$. The impedance boundary condition is imposed in the
same form as in Part~I,
\begin{equation}
  \label{eq:robin-pert}
  \gamma_1 u - \alpha\, \gamma_0 u = 0
  \qquad \text{on } \Gpert,
  \qquad \alpha > 0 .
\end{equation}

\begin{definition}[The perturbed Robin Laplacian]
\label{def:Hpert}
The \emph{perturbed Robin Laplacian} is the operator in $\Ltwo{\Om}$ given by
\begin{equation}
  \label{eq:Hpert-def}
  \Hpert u = -\Delta u,
  \qquad
  \Dom(\Hpert)
    = \bigl\{\, u \in H^2(\Om) \;\bigm|\;
       \gamma_1 u - \alpha\, \gamma_0 u = 0 \text{ on } \Gpert \,\bigr\}.
\end{equation}
Its natural sesquilinear form is
\begin{equation}
  \label{eq:form-pert}
  \hpert[u,v]
    = (\nabla u, \nabla v)_{\Ltwo{\Om}}
      - \alpha\,(\gamma_0 u, \gamma_0 v)_{\Ltwo{\Gpert}},
  \qquad u, v \in H^1(\Om),
\end{equation}
and we write $\hpert[u] := \hpert[u,u]$.
\end{definition}

\begin{remark}[The sign of the impedance]
\label{rem:sign}
On the flat part of the boundary, where $\nu = -e_n$, \eqref{eq:robin-pert}
reads $\partial_{x_n} u + \alpha u = 0$. This is the impedance condition
$\partial_\nu u + \lambda u = 0$ with $\lambda = -\alpha < 0$, the sign
discussed in the introduction. It makes the boundary term in
\eqref{eq:form-pert} negative, and that term produces both the spectrum below
$0$ and the surface channel. For the opposite sign one would have
$\Hpert \ge 0$ and no surface channel.
\end{remark}

For comparison we keep the notation of Part~I for the unperturbed operator:
$\Halpha$ is the half-space Robin Laplacian, that is, \eqref{eq:Hpert-def} with
$\Om$ replaced by $\Rnp$ and $\Gpert$ by $\Gamma$. We shall use repeatedly the
following two facts about it, both established in Part~I. First, its spectrum
is a half-line and is purely absolutely continuous
\cite[Theorem~2.2(b),(c)]{PartI},
\begin{equation}
  \label{eq:PartI-spectrum}
  \sigma(\Halpha) = \sac(\Halpha) = [-\alpha^2, \infty),
  \qquad
  \spp(\Halpha) = \ssc(\Halpha) = \varnothing .
\end{equation}
Secondly, by \cite[Proposition~2.4]{PartI} the partial Fourier transform
\eqref{eq:fourier-partial} diagonalizes $\Halpha$ into the direct integral of
the half-line Robin operators $\hfib(\xipar) = -d^2/dx_n^2 + |\xipar|^2$ on $\Ltwo{\Rhalf}$ with
boundary condition $v'(0) + \alpha v(0) = 0$, each of which has the simple
eigenvalue
\begin{equation}
  \label{eq:surface-band}
  \lambda_{\mathrm{surf}}(\xipar) = |\xipar|^2 - \alpha^2
\end{equation}
below a continuum $[\,|\xipar|^2, \infty)$, with the $\xipar$-independent
normalized eigenfunction
\begin{equation}
  \label{eq:psizero}
  \psizero(x_n) = \sqrt{2\alpha}\, e^{-\alpha x_n}.
\end{equation}
The two branches are the two scattering channels: the surface channel, of
internal energy $-\alpha^2$ and transverse profile $\psizero$, open for
$E > -\alpha^2$, and the bulk channel, of internal energy $0$, open for
$E > 0$. The thresholds are $\Thr = \{-\alpha^2, 0\}$.

\subsection{Self-Adjointness and the Spectrum}
\label{subsec:pert-sa-spec}

Semiboundedness of the form \eqref{eq:form-pert} is not immediate, since the
boundary term is negative. On the flat half-space it follows from a trace
inequality whose constant is optimal and returns $-\alpha^2$ exactly. On $\Om$
the same argument yields semiboundedness, with a constant that depends on the
geometry through a single vector field.

\begin{lemma}[Trace inequality on the perturbed domain]
\label{lem:trace-pert}
Let $n \ge 2$ and let $\Om$ satisfy Assumption~\ref{assump:geometry}. There is
a constant $\Cdiv = \Cdiv(\Om) \ge 0$ such that
\begin{equation}
  \label{eq:trace-ineq}
  \norm{\gamma_0 u}_{\Ltwo{\Gpert}}^2
    \;\le\; \varepsilon \norm{\nabla u}_{\Ltwo{\Om}}^2
      + \bigl(\varepsilon^{-1} + \Cdiv\bigr)
        \norm{u}_{\Ltwo{\Om}}^2
  \qquad
  \bigl(u \in H^1(\Om),\ \varepsilon > 0\bigr).
\end{equation}
On the unperturbed half-space one may take $\Cdiv = 0$, and the constant
$\varepsilon^{-1}$ is then optimal for each fixed $\varepsilon$, although it is
not attained.
\end{lemma}

\begin{proof}
The whole inequality comes from a single vector field. We claim that there is
$\Nvec \in C^1(\overline{\Om}; \Rn)$, bounded together with its first
derivatives, such that
\begin{equation}
  \label{eq:vector-field}
  \Nvec = -e_n \ \text{ on } \Om \setminus \Bfull{\Rzero+1},
  \qquad
  \Nvec = \nu \ \text{ on } \Gpert,
  \qquad
  |\Nvec| \le 1 \ \text{ on } \overline{\Om} .
\end{equation}
Indeed, the portion $\Gpert \cap \overline{\Bfull{\Rzero+2}}$ of the
boundary is compact and of class $C^2$ by (G1). Hence there are
$\delta \in (0,\tfrac12)$ and a neighbourhood
$V = \{x : |x| < \Rzero + \tfrac32,\ \operatorname{dist}(x,\Gpert) < \delta\}$
on which the distance to $\Gpert$, signed so as to be positive outside $\Om$,
is a function $d$ of class $C^2$. Its gradient $\nabla d$ has length $1$ and
equals $\nu$ on $\Gpert$. By (G2) the boundary is flat at every point with
$|x| > \Rzero$. A point $x \in V$ with $|x| > \Rzero + 1$ therefore has its
nearest boundary point on the flat part, so that $d(x) = -x_n$ and
$\nabla d(x) = -e_n$. Choose $\chi \in C^\infty(\Rn)$ with
$0 \le \chi \le 1$, supported in $V$ and equal to $1$ near
$\Gpert \cap \overline{\Bfull{\Rzero+1}}$, and put
\begin{equation*}
  \Nvec = \chi\,\nabla d - (1 - \chi)\,e_n .
\end{equation*}
At each point $\Nvec$ is a convex combination of two vectors of length at most
$1$, so $|\Nvec| \le 1$. On $\Gpert \cap \overline{\Bfull{\Rzero+1}}$ one has
$\chi = 1$ and hence $\Nvec = \nu$. On the rest of $\Gpert$ both $\nabla d$
and $-e_n$ equal $\nu = -e_n$, so again $\Nvec = \nu$. Finally $\Nvec$ differs
from $-e_n$ only where $\chi \ne 0$ and $\nabla d \ne -e_n$, that is, inside
$\Bfull{\Rzero+1}$. Since $\chi$ and $\nabla d$ are of class $C^1$ with bounded
derivatives on the support of $\chi$, this proves
\eqref{eq:vector-field}.

Let now $u \in H^1(\Om)$. Applying the divergence theorem to the vector field
$|u|^2 \Nvec$, which is admissible because $\Nvec$ is bounded and Lipschitz and
$|u|^2 \in W^{1,1}(\Om)$, and recalling that $\nu$ is the outward normal,
\begin{equation*}
  \int_{\Gpert} |\gamma_0 u|^2\, (\Nvec \cdot \nu)\, dS
    = \int_{\Om} |u|^2 \operatorname{div}\Nvec \, dx
      + 2\operatorname{Re}\int_{\Om} \overline{u}\, (\Nvec \cdot \nabla u)\, dx .
\end{equation*}
By \eqref{eq:vector-field} we have $\Nvec \cdot \nu = |\nu|^2 = 1$ on $\Gpert$,
so the left-hand side is exactly
$\norm{\gamma_0 u}_{\Ltwo{\Gpert}}^2$, while on the right
$|\Nvec| \le 1$ and $\operatorname{div}\Nvec$ is supported in
$\Omin{\Rzero+1}$ and bounded there, say by $\Cdiv$. Hence
\begin{equation*}
  \norm{\gamma_0 u}_{\Ltwo{\Gpert}}^2
    \le \Cdiv \norm{u}_{\Ltwo{\Om}}^2
      + 2 \norm{u}_{\Ltwo{\Om}} \norm{\nabla u}_{\Ltwo{\Om}},
\end{equation*}
and \eqref{eq:trace-ineq} follows from Young's inequality in the form
$2\norm{u}\norm{\nabla u} \le \varepsilon \norm{\nabla u}^2
+ \varepsilon^{-1}\norm{u}^2$, valid for every $\varepsilon > 0$.

On $\Rnp$ one may take $\Nvec \equiv -e_n$, whose divergence vanishes, so that
$\Cdiv = 0$. The resulting inequality is the classical one obtained by
integrating $\partial_{x_n}|u|^2$ over the normal variable. Its constant cannot
be improved, but neither is it attained. Indeed, for
$u_\delta(\xpar,x_n) = w(\delta\xpar)\,e^{-x_n/\varepsilon}$ with
$w \in C_0^\infty(\Rnm) \setminus \{0\}$ the two sides of
\eqref{eq:trace-ineq} differ by exactly
$\tfrac{\varepsilon^2}{2}\norm{\nabla'(w(\delta\,\cdot))}_{\Ltwo{\Rnm}}^2$,
which is positive for every $\delta > 0$ and is $o(1)$ relative to
$\norm{\gamma_0 u_\delta}^2$ as $\delta \to 0$.
\end{proof}

\begin{theorem}[Self-adjointness and the essential spectrum]
\label{thm:pert-sa-spec}
Let $n \ge 2$, let $\alpha > 0$, and let $\Om$ satisfy
Assumption~\ref{assump:geometry}.
\begin{enumerate}
  \item[\textup{(a)}] The form $\hpert$ of \eqref{eq:form-pert} is closed,
  symmetric and semibounded on the form domain $H^1(\Om)$, with
  \begin{equation}
    \label{eq:lower-bound-pert}
    \hpert[u]
      \ge \tfrac12 \norm{\nabla u}_{\Ltwo{\Om}}^2
        - \bigl(2\alpha^2 + \alpha\,\Cdiv\bigr)
          \norm{u}_{\Ltwo{\Om}}^2
    \qquad \bigl(u \in H^1(\Om)\bigr),
  \end{equation}
  and the self-adjoint operator associated with it is exactly the operator
  $\Hpert$ of \eqref{eq:Hpert-def}. In particular $\Hpert$ is self-adjoint and
  bounded below, and there is $C > 0$ with
  \begin{equation}
    \label{eq:apriori-H2}
    \norm{u}_{H^2(\Om)}
      \le C\bigl( \norm{\Hpert u}_{\Ltwo{\Om}} + \norm{u}_{\Ltwo{\Om}} \bigr)
    \qquad \bigl(u \in \Dom(\Hpert)\bigr).
  \end{equation}
  \item[\textup{(b)}] The essential spectrum is that of the unperturbed
  operator:
  \begin{equation}
    \label{eq:sess}
    \sess(\Hpert) = \sess(\Halpha) = [-\alpha^2, \infty).
  \end{equation}
  \item[\textup{(c)}] Consequently $\sigma(\Hpert) \setminus [-\alpha^2,\infty)$
  consists of at most countably many eigenvalues of finite multiplicity, lying
  in the interval $[\,\inf\sigma(\Hpert), -\alpha^2)$ and having no accumulation
  point other than $-\alpha^2$.
\end{enumerate}
\end{theorem}

Part~(a) is standard once Lemma~\ref{lem:trace-pert} is available. Inserting
\eqref{eq:trace-ineq} into \eqref{eq:form-pert} gives
$\hpert[u] \ge (1 - \alpha\varepsilon)\norm{\nabla u}^2
- \alpha(\varepsilon^{-1} + \Cdiv)\norm{u}^2$ for every $\varepsilon > 0$, and
\eqref{eq:lower-bound-pert} is the choice $\varepsilon = 1/(2\alpha)$. The
identification of the operator associated with the form, together with the
elliptic regularity giving \eqref{eq:apriori-H2}, is carried out in
Appendix~\ref{app:regularity}. Part~(c) is the general spectral theorem for a
semibounded self-adjoint operator, given (b). We prove (b) here, since the
argument is short.

\begin{proof}[Proof of Theorem~\textup{\ref{thm:pert-sa-spec}(b)}]
\emph{The inclusion $[-\alpha^2,\infty) \subset \sess(\Hpert)$.} Let
$E \ge -\alpha^2$. By \eqref{eq:surface-band} the sphere
$\{|\xipar| = \ksurf\}$, of radius the surface wavenumber
$\ksurf \ge 0$, carries the energy $E$ in the surface band. Fix
$\xipar_0 \in \Rnm$ with $|\xipar_0| = \ksurf$, and fix
$\varphi \in C_0^\infty(\Rnm)$ with $\norm{\varphi}_{\Ltwo{\Rnm}} = 1$ and
$\operatorname{supp}\varphi \subset \{|\xpar| < \rho\}$. Let $a_k \in \Rnm$ be
any sequence with
\begin{equation}
  \label{eq:weyl-centres}
  |a_k| = \Rzero + (\rho + 1)\,k
  \qquad (k \in \mathbb{N}),
\end{equation}
and set
\begin{equation}
  \label{eq:weyl-sequence}
  u_k(x)
    = k^{-(n-1)/2}\, \varphi\Bigl(\frac{\xpar - a_k}{k}\Bigr)\,
      e^{i \xpar \cdot \xipar_0}\, \psizero(x_n).
\end{equation}
Each $u_k$ is smooth and, since $\psizero'(0) + \alpha\psizero(0) = 0$, it
satisfies the boundary condition \eqref{eq:robin-pert} on the flat part of the
boundary. Its support is contained in $\{|\xpar - a_k| \le \rho k\}$, hence in
$\{|\xpar| \ge \Rzero + k\}$ by \eqref{eq:weyl-centres}. Note that the support
grows linearly in $k$, so that a condition of the form $|a_k| \to \infty$
alone would not suffice here. In particular the support of $u_k$ lies in
$\Omex{\Rzero}$, where $\Om$ and $\Rnp$ agree by
\eqref{eq:exterior-coincidence}. Hence $u_k \in \Dom(\Hpert)$ and
$u_k \in \Dom(\Halpha)$, and the two operators act on it identically.

The normalization of \eqref{eq:weyl-sequence} gives
$\norm{u_k}_{\Ltwo{\Om}} = 1$ for every $k$, since $\psizero$ is normalized in
$\Ltwo{\Rhalf}$. Moreover $u_k \rightharpoonup 0$ weakly, because
$\operatorname{supp} u_k \subset \{|\xpar| \ge \Rzero + k\}$ and these sets
escape to infinity. Finally, using $\psizero'' = \alpha^2 \psizero$ and
$|\xipar_0|^2 - \alpha^2 = E$, a direct computation gives
\begin{equation*}
  (\Hpert - E)u_k
    = -k^{-(n-1)/2}
      \Bigl[ (\Delta' \varphi_k) + 2i\,\xipar_0 \cdot (\nabla'\varphi_k) \Bigr]
      e^{i\xpar\cdot\xipar_0} \psizero,
  \qquad
  \varphi_k(\xpar) = \varphi\Bigl(\frac{\xpar - a_k}{k}\Bigr),
\end{equation*}
whose $\Ltwo{\Om}$ norm is bounded by $C(1 + |\xipar_0|)/k \to 0$. Thus
$(u_k)$ is a singular sequence for $\Hpert$ at $E$, and $E \in \sess(\Hpert)$
by Weyl's criterion.

\emph{The inclusion $\sess(\Hpert) \subset [-\alpha^2,\infty)$.} Let
$E \in \sess(\Hpert)$ and let $(u_k) \subset \Dom(\Hpert)$ be a singular
sequence, so that $\norm{u_k}_{\Ltwo{\Om}} = 1$, $u_k \rightharpoonup 0$ and
$f_k := (\Hpert - E)u_k \to 0$ in $\Ltwo{\Om}$. By \eqref{eq:apriori-H2} the
sequence is bounded in $H^2(\Om)$.

Choose a radial cut-off $\cut \in C^\infty(\Rn)$ with $\cut = 0$ on
$\Bfull{\Rzero + 1}$ and $\cut = 1$ outside $\Bfull{\Rzero + 2}$, and set
$K = \Om \cap \{\Rzero + 1 \le |x| \le \Rzero + 2\}$, on which $\nabla\cut$ and
$\Delta\cut$ are supported. By \eqref{eq:exterior-coincidence} this is a
bounded subset of the region where $\Om$ and $\Rnp$ agree. Since $\cut$
vanishes near $\Bfull{\Rzero}$, the function $\cut u_k$, extended by zero,
belongs to $H^2(\Rnp)$. It also satisfies the boundary condition on $\Gamma$.
Indeed $\cut$ is radial, so that $\partial_{x_n}\cut$ vanishes on
$\{x_n = 0\}$ and therefore $\gamma_1(\cut u_k) = \cut\, \gamma_1 u_k$. Hence
$\cut u_k \in \Dom(\Halpha)$, and
\begin{equation*}
  (\Halpha - E)(\cut u_k)
    = \cut f_k - [\Delta, \cut] u_k
    = \cut f_k - 2\nabla\cut \cdot \nabla u_k - (\Delta \cut) u_k .
\end{equation*}
Now $\Omin{\Rzero+2}$ is a bounded Lipschitz set, so the embedding
$H^2(\Omin{\Rzero+2}) \hookrightarrow H^1(\Omin{\Rzero+2})$ is compact by the
Rellich--Kondrachov theorem. Since $(u_k)$ is bounded in $H^2(\Om)$, every
subsequence therefore has a further subsequence converging in
$H^1(\Omin{\Rzero+2})$, and the limit is $0$ because $u_k \rightharpoonup 0$ in
$\Ltwo{\Om}$. Hence the full sequence converges to $0$ in
$H^1(\Omin{\Rzero+2})$. Both commutator terms are supported in
$K \subset \Omin{\Rzero+2}$ and are controlled by the $H^1(K)$ norm of $u_k$.
Therefore $\norm{(\Halpha - E)(\cut u_k)}_{\Ltwo{\Rnp}} \to 0$. The same
convergence gives $\norm{(1 - \cut)u_k}_{\Ltwo{\Om}} \to 0$, whence
$\norm{\cut u_k}_{\Ltwo{\Rnp}} \to 1$, and $\cut u_k \rightharpoonup 0$.
After normalizing, $(\cut u_k)$ is a singular sequence for $\Halpha$ at $E$, so
$E \in \sess(\Halpha) = [-\alpha^2,\infty)$ by \eqref{eq:PartI-spectrum}.
\end{proof}

\begin{remark}[The bottom of the spectrum]
\label{rem:bottom}
The lower bound \eqref{eq:lower-bound-pert} is weaker than its counterpart
$\hform[u] \ge -\alpha^2\norm{u}^2$ in Part~I, for two separate reasons. The
first is bookkeeping. Estimate \eqref{eq:lower-bound-pert} keeps half of the
gradient term, as form boundedness requires. Taking instead
$\varepsilon = 1/\alpha$ in \eqref{eq:trace-ineq} gives
$\hpert[u] \ge -(\alpha^2 + \alpha\Cdiv)\norm{u}^2$, which on the flat
half-space is the sharp bound of Part~I. The second is the term $\alpha\Cdiv$,
which comes from the divergence of the field $\Nvec$ in the proof of
Lemma~\ref{lem:trace-pert} and vanishes on the flat half-space, where
$\Nvec = -e_n$. Whether $\Hpert$ has spectrum below $-\alpha^2$ depends on the
geometry, and we make no claim either way. By part~(c) of
Theorem~\ref{thm:pert-sa-spec} any such spectrum consists of eigenvalues of
finite multiplicity, and the scattering theory of the sequel, which takes place
on $[-\alpha^2,\infty)$, is unaffected by it.
\end{remark}

Eigenvalues embedded in the continuous spectrum are treated in
\S\ref{sec:lap-pert}. Above the bulk threshold there are none. The
Rellich-type theorem of Part~I transplants to the perturbed domain and excludes
them. On the band $(-\alpha^2,0)$ they are not excluded. They form a discrete
subset of the band, and each has finite multiplicity. That distinction
governs where the limiting absorption principle holds.

\section{The Limiting Absorption Principle}
\label{sec:lap-pert}
The limiting absorption principle for $\Hpert$ is the analytic core of this
paper, and it is established here in the same spaces in which it was
established in Part~I for the flat half-space: the Agmon--H\"ormander spaces
$\Bsp$ and $\Bst$. The reason for that choice was explained in
\S\ref{sec:intro} and is unchanged by the perturbation. The two propagating
channels carry different rates and live in different dimensions, and it is only
at the Agmon--H\"ormander endpoint that both are resolved with finite nonzero
norm. A weighted space $L^{2,s}$ either contains both without seeing their
amplitudes, or excludes both.

What does change is that the spaces themselves have to be defined on $\Om$, and
that the resolvent is no longer diagonalized by the partial Fourier transform,
so that the exact decomposition used in Part~I is unavailable. The route taken
here is the classical one: the free resolvent is transplanted from $\Rnp$
through the exterior region by a cut-off, the error this produces is supported
in a compact set, and a compactness argument reduces the resolvent estimate to
a uniqueness statement, which is in turn supplied by a Rellich-type
theorem. The first subsection sets up the spaces, in a form that makes the
transplantation in the later subsections literal rather than approximate.

\subsection{The Agmon--H\"ormander Spaces on the Perturbed Domain}
\label{subsec:AH-pert}

\subsubsection*{Notation, and a convention on subscripts}

We recall the dyadic decomposition of Part~I
\cite[\S2.3.1]{PartI}. In $\Rn$ one puts
\begin{equation}
  \label{eq:shells-pert}
  \shell{0} = \{x \in \Rn : |x| < 1\},
  \qquad
  \shell{j} = \{x \in \Rn : 2^{j-1} \le |x| < 2^{j}\}
  \quad (j \ge 1),
\end{equation}
and writes $\shellp{j} = \shell{j} \cap \Rnp$. The Agmon--H\"ormander norms on
the flat half-space are
\begin{equation}
  \label{eq:B-norms-flat}
  \Bnorm{v} = \sum_{j \ge 0} 2^{j/2}\,\norm{v}_{\Ltwo{\shellp{j}}},
  \qquad
  \Bstnorm{v} = \sup_{j \ge 0} 2^{-j/2}\,\norm{v}_{\Ltwo{\shellp{j}}},
\end{equation}
with $\Bsp$ and $\Bst$ the spaces on which these are finite, and
\begin{equation}
  \label{eq:Bzero-flat}
  \Bstz = \Bigl\{\, v \in \Bst \;\Bigm|\;
    \lim_{j \to \infty} 2^{-j/2}\norm{v}_{\Ltwo{\shellp{j}}} = 0 \,\Bigr\}.
\end{equation}
Throughout this section, $\Bnorm{\cdot}$ and $\Bstnorm{\cdot}$ without further
decoration always denote these norms on $\Rnp$.

Two families of sets indexed by integers or radii now occur together: the shells
\eqref{eq:shells-pert}, which are subsets of $\Rn$ and have nothing to do with
the perturbation, and the truncations $\Omin{R}$, $\Omex{R}$ of $\Om$
introduced in Assumption~\ref{assump:geometry}. We distinguish them
typographically: a bare $\Omega$ with a numerical subscript is always a dyadic
shell, whereas the parts of the perturbed domain always carry the superscript
$\mathrm{in}$ or $\mathrm{ex}$.

\subsubsection*{Normalization of the radius}

Assumption~\ref{assump:geometry} is not affected by enlarging $\Rzero$: if
\eqref{eq:exterior-coincidence} holds for one radius it holds for every larger
one. We may therefore fix $\Rzero$ to be a power of two,
\begin{equation}
  \label{eq:R0-dyadic}
  \Rzero = 2^{\Jzero},
  \qquad \Jzero \in \mathbb{N},
\end{equation}
and we do so from now on, taking $\Jzero \ge 2$, so that $\Rzero \ge 4$. The
point of the dyadic normalization is that the decomposition then splits exactly
at $\Rzero$:
\begin{equation}
  \label{eq:shells-split}
  \bigcup_{j = 0}^{\Jzero} \shell{j} = \Bfull{\Rzero},
  \qquad
  \shell{j} \cap \Om = \shellp{j}
  \quad (j > \Jzero),
\end{equation}
the second identity holding up to the null set
$\shell{j} \cap \{|x| = \Rzero\}$ and being nothing but
\eqref{eq:exterior-coincidence}, since $\shell{j} \subset \{|x| \ge \Rzero\}$
for $j > \Jzero$. All the quantities below are $L^2$ norms, so the null set is
immaterial. Thus the shells of index $j > \Jzero$ are \emph{the same sets} as in
Part~I, and only the finitely many shells inside $\Bfull{\Rzero}$ see the
perturbation. We group those into a single block. The lower bound
$\Rzero \ge 4$ is imposed for a separate reason, and is used only once, in the
proof of Lemma~\ref{lem:incident-transplant}: it puts every admissible cut-off
annulus inside the region where the pointwise bound of
\cite[Lemma 2.75]{PartI} is available.

\begin{definition}[Agmon--H\"ormander spaces on $\Om$]
\label{def:AH-pert}
For $u \in L^2_{\mathrm{loc}}(\overline{\Om})$ set
\begin{equation}
  \label{eq:B-norms-pert}
  \begin{aligned}
    \BnormO{u}
      &= \Rzero^{1/2}\,\norm{u}_{\Ltwo{\Omin{\Rzero}}}
         + \sum_{j > \Jzero} 2^{j/2}\,\norm{u}_{\Ltwo{\shellp{j}}},
    \\[1mm]
    \BstnormO{u}
      &= \max\Bigl\{\,
           \Rzero^{-1/2}\,\norm{u}_{\Ltwo{\Omin{\Rzero}}},\;
           \sup_{j > \Jzero} 2^{-j/2}\,\norm{u}_{\Ltwo{\shellp{j}}}
         \,\Bigr\},
  \end{aligned}
\end{equation}
and let $\Bspd{\Om}$ and $\Bstd{\Om}$ be the spaces of those $u$ for which the
respective quantity is finite. We further set
\begin{equation}
  \label{eq:Bzero-pert}
  \Bstzd{\Om}
    = \Bigl\{\, u \in \Bstd{\Om} \;\Bigm|\;
        \lim_{j \to \infty} 2^{-j/2}\norm{u}_{\Ltwo{\shellp{j}}} = 0
      \,\Bigr\}.
\end{equation}
\end{definition}

The block $\Omin{\Rzero}$ is thus treated as a single shell sitting at radius
$\Rzero$, and it carries the weights $\Rzero^{\pm 1/2}$ that a shell at that
radius would carry in \eqref{eq:B-norms-flat}. This is what keeps the two norms
in exact duality with those of Part~I on functions supported outside
$\Bfull{\Rzero}$. The two consequences are recorded in
Lemma~\ref{lem:AH-exterior} and Proposition~\ref{prop:AH-elementary}(a). Note
also that the defining condition \eqref{eq:Bzero-pert} of $\Bstzd{\Om}$ involves
only the limit $j \to \infty$, and so does not refer to the block at all.

\subsubsection*{Comparison with the shell definition}

There is an alternative convention, namely to intersect every shell with the
domain. It is the one used in Part~I, and by Isozaki, Kadowaki and Watanabe for
an exterior domain \cite{IsozakiKadowakiWatanabe2025JOS}. Write
\begin{equation}
  \label{eq:B-norms-shell}
  \Bnormsh{u} = \sum_{j \ge 0} 2^{j/2}\,\norm{u}_{\Ltwo{\shell{j} \cap \Om}},
  \qquad
  \Bstnormsh{u} = \sup_{j \ge 0} 2^{-j/2}\,\norm{u}_{\Ltwo{\shell{j} \cap \Om}},
\end{equation}
and let
\begin{equation}
  \label{eq:Q-scale}
  \Qsc{u} = \sup_{R > 1} \frac{1}{R}
    \int_{\Om \cap \Bfull{R}} |u(x)|^2\,dx
\end{equation}
be the scale-invariant quantity by which $\Bst$ is often defined directly. The
next lemma says that nothing is at stake in the choice.

\begin{lemma}[Equivalence of the three conventions]
\label{lem:AH-equivalence}
Put $\Csh = \Rzero^{1/2}\bigl(1 + \log_2 \Rzero\bigr)^{1/2}$. For every
$u \in L^2_{\mathrm{loc}}(\overline{\Om})$,
\begin{equation}
  \label{eq:AH-equivalence}
  \Csh^{-1}\,\Bnormsh{u} \le \BnormO{u} \le \Csh\,\Bnormsh{u},
  \qquad
  \Csh^{-1}\,\Bstnormsh{u} \le \BstnormO{u} \le \Csh\,\Bstnormsh{u},
\end{equation}
and moreover
\begin{equation}
  \label{eq:AH-scale-invariant}
  \tfrac{1}{\sqrt2}\, \Bstnormsh{u}
    \;\le\; \Qsc{u}^{1/2}
    \;\le\; 2\, \Bstnormsh{u}.
\end{equation}
Consequently the three conventions define the same two spaces with equivalent
norms, and $\Bstzd{\Om}$ is exactly the set of $u \in \Bstd{\Om}$ for which
$R^{-1}\int_{\Om \cap \Bfull{R}}|u|^2 \to 0$ as $R \to \infty$. All constants
depend on $\Rzero$ alone.
\end{lemma}

\begin{proof}
By \eqref{eq:shells-split} the two definitions \eqref{eq:B-norms-pert} and
\eqref{eq:B-norms-shell} have identical contributions from the indices
$j > \Jzero$, so it suffices to compare the blocks. Write
$a_j = \norm{u}_{\Ltwo{\shell{j} \cap \Om}}$ for $0 \le j \le \Jzero$ and
$A = \norm{u}_{\Ltwo{\Omin{\Rzero}}}$, so that
$A = \bigl(\sum_{j \le \Jzero} a_j^2\bigr)^{1/2}$ by the first identity in
\eqref{eq:shells-split}.

For $\Bsp$: since $2^{j/2} \ge 1$ and the $\ell^1$ norm dominates the $\ell^2$
norm,
\begin{equation*}
  A \;\le\; \sum_{j \le \Jzero} a_j
    \;\le\; \sum_{j \le \Jzero} 2^{j/2} a_j ,
\end{equation*}
so that $\Rzero^{1/2} A \le \Rzero^{1/2}\sum_{j \le \Jzero}2^{j/2}a_j$. By
the Cauchy--Schwarz inequality in the $\Jzero + 1$ indices,
\begin{equation*}
  \sum_{j \le \Jzero} 2^{j/2} a_j
    \;\le\; 2^{\Jzero/2} \sum_{j \le \Jzero} a_j
    \;\le\; 2^{\Jzero/2} (\Jzero + 1)^{1/2} A
    \;=\; (1 + \log_2 \Rzero)^{1/2}\, \bigl(\Rzero^{1/2} A\bigr) .
\end{equation*}
The two displays are the first half of \eqref{eq:AH-equivalence}, both constants
being at most $\Csh$.

For $\Bst$: in the same notation the block terms are $\Rzero^{-1/2}A$ and
$\max_{j \le \Jzero} 2^{-j/2}a_j$, and
\begin{equation*}
  \max_{j \le \Jzero} 2^{-j/2} a_j
    \;\le\; \max_{j \le \Jzero} a_j
    \;\le\; A
    \;=\; \Rzero^{1/2}\bigl(\Rzero^{-1/2}A\bigr),
\end{equation*}
while conversely $\max_{j \le \Jzero} a_j \ge (\Jzero+1)^{-1/2} A$ gives
\begin{equation*}
  \max_{j \le \Jzero} 2^{-j/2} a_j
    \;\ge\; 2^{-\Jzero/2} \max_{j \le \Jzero} a_j
    \;\ge\; (1 + \log_2 \Rzero)^{-1/2} \bigl(\Rzero^{-1/2} A\bigr).
\end{equation*}
Taking the maximum with the common tail preserves both inequalities, and the
second half of \eqref{eq:AH-equivalence} follows.

For \eqref{eq:AH-scale-invariant}, set $c_j = \norm{u}_{\Ltwo{\shell{j}\cap\Om}}$
for every $j \ge 0$ and $M = \sup_{j \ge 0} 2^{-j/2}c_j = \Bstnormsh{u}$. Given
$R > 1$, choose $J \ge 1$ with $2^{J-1} \le R < 2^{J}$. Then
$\Om \cap \Bfull{R} \subset \bigcup_{j \le J}(\shell{j}\cap\Om)$ and
\begin{equation*}
  \int_{\Om \cap \Bfull{R}} |u|^2
    \;\le\; \sum_{j \le J} c_j^2
    \;\le\; M^2 \Bigl( 1 + \sum_{1 \le j \le J} 2^{j} \Bigr)
    \;\le\; 2^{J+1} M^2
    \;\le\; 4 R M^2 ,
\end{equation*}
which gives $\Qsc{u} \le 4M^2$, the upper bound in
\eqref{eq:AH-scale-invariant}. Conversely
$c_j^2 \le \int_{\Om \cap \Bfull{2^j}}|u|^2 \le 2^{j}\Qsc{u}$ for $j \ge 1$, so
that $2^{-j}c_j^2 \le \Qsc{u}$, while
$c_0^2 \le \int_{\Om\cap\Bfull{2}}|u|^2 \le 2\Qsc{u}$. Hence
$M^2 \le 2\Qsc{u}$.

For the last statement, apply the two displays of the preceding paragraph to
$u$ restricted to $\Om \setminus \Bfull{\rho}$, $\rho > 0$. The first gives
\begin{equation*}
  \frac{1}{R}\int_{\Om\cap\Bfull{R}}|u|^2
    \;\le\; \frac{1}{R}\int_{\Om\cap\Bfull{\rho}}|u|^2
      + 4 \sup_{2^{j} > \rho} 2^{-j}c_j^2
  \qquad (R > \rho),
\end{equation*}
so that if $2^{-j/2}c_j \to 0$ then the left-hand side tends to $0$ as
$R \to \infty$, the first term on the right doing so for fixed $\rho$. The
second gives $2^{-j}c_j^2 \le \sup_{R > 2^{j-1}} R^{-1}\int_{\Om\cap\Bfull{R}}
|u|^2$ for $j \ge 1$, so that the converse implication holds as well. Finally
$2^{-j/2}c_j$ and $2^{-j/2}\norm{u}_{\Ltwo{\shellp{j}}}$ agree for
$j > \Jzero$ by \eqref{eq:shells-split}, which identifies the condition just
obtained with \eqref{eq:Bzero-pert}.
\end{proof}

By Lemma~\ref{lem:AH-equivalence} we may quote results stated in the shell
convention, in particular
those of \cite{IsozakiKadowakiWatanabe2025JOS}, without further comment. What
Definition~\ref{def:AH-pert} adds is the exact identity of
Lemma~\ref{lem:AH-exterior} below, which lets the estimates of Part~I be
transplanted with their constants unchanged.

\subsubsection*{Exact agreement in the exterior}

\begin{lemma}[Functions supported outside the perturbation]
\label{lem:AH-exterior}
Let $u \in L^2_{\mathrm{loc}}(\overline{\Om})$ vanish almost everywhere on
$\Omin{\Rzero}$, and let $\zext{u}$ denote its extension by zero to $\Rnp$.
Then
\begin{equation}
  \label{eq:AH-exterior}
  \BnormO{u} = \Bnorm{\zext{u}},
  \qquad
  \BstnormO{u} = \Bstnorm{\zext{u}},
\end{equation}
the right-hand norms being those of \eqref{eq:B-norms-flat} on $\Rnp$. Moreover
$u \in \Bstzd{\Om}$ if and only if $\zext{u} \in \Bstz$.
\end{lemma}

\begin{proof}
The extension is well defined because $\Omex{\Rzero} = \Rnp \setminus
\overline{\Bfull{\Rzero}}$ by \eqref{eq:exterior-coincidence}. The block term in
\eqref{eq:B-norms-pert} vanishes by hypothesis, and so do the terms of
\eqref{eq:B-norms-flat} with $j \le \Jzero$, since $\shellp{j} \subset
\Bfull{\Rzero}$ for those indices by \eqref{eq:shells-split}. The remaining
terms coincide one by one, again by \eqref{eq:shells-split}. The statement about
$\Bstz$ is immediate from \eqref{eq:Bzero-pert} and \eqref{eq:Bzero-flat}.
\end{proof}

The use to which this will be put is the following. Every estimate of Part~I
whose two sides involve only functions supported in $\Omex{\Rzero}$ is, by
Lemma~\ref{lem:AH-exterior}, an estimate on $\Om$ with the \emph{same constant}.
In particular the Rellich-type theorem of Part~I applies on $\Omex{\Rzero}$
without any modification of its constants or of its hypotheses, which is the
starting point of \S\ref{sec:lap-pert} and the reason the perturbation can be
confined to a compact region throughout.

\subsubsection*{Elementary properties}

\begin{proposition}
\label{prop:AH-elementary}
Let $\Om$ satisfy Assumption~\ref{assump:geometry} with $\Rzero$ normalized as
in \eqref{eq:R0-dyadic}.
\begin{enumerate}
  \item[\textup{(a)}] For $u \in \Bstd{\Om}$ and $g \in \Bspd{\Om}$ the pairing
  $(u,g) = \int_{\Om} u\overline{g}\,dx$ converges absolutely, with
  \begin{equation}
    \label{eq:AH-pairing}
    |(u,g)| \le \BstnormO{u}\,\BnormO{g}.
  \end{equation}
  \item[\textup{(b)}] With $L^{2,s}(\Om)$ denoting the weighted space with norm
  $\norm{(1+|x|)^s u}_{\Ltwo{\Om}}$ one has the chain of continuous inclusions
  \begin{equation}
    \label{eq:AH-inclusions-pert}
    L^{2,s}(\Om) \subset \Bspd{\Om} \subset L^{2,1/2}(\Om)
      \subset \Ltwo{\Om} \subset L^{2,-1/2}(\Om)
      \subset \Bstzd{\Om} \subset \Bstd{\Om} \subset L^{2,-s}(\Om)
    \qquad (s > 1/2),
  \end{equation}
  the inclusions $\Bspd{\Om} \subset \Ltwo{\Om}$ and
  $\Ltwo{\Om} \subset \Bstd{\Om}$ having norm at most $1$.
  \item[\textup{(c)}] Every operator bounded on $\Ltwo{\Om}$ is bounded from
  $\Bspd{\Om}$ to $\Bstd{\Om}$, with the same bound.
\end{enumerate}
\end{proposition}

\begin{proof}
(a) Split the integral over $\Omin{\Rzero}$ and the shells $\shellp{j}$,
$j > \Jzero$, apply the Cauchy--Schwarz inequality on each piece, and insert the
reciprocal weights:
\begin{equation*}
  |(u,g)|
    \le \bigl(\Rzero^{-1/2}\norm{u}_{\Ltwo{\Omin{\Rzero}}}\bigr)
        \bigl(\Rzero^{1/2}\norm{g}_{\Ltwo{\Omin{\Rzero}}}\bigr)
      + \sum_{j > \Jzero}
        \bigl(2^{-j/2}\norm{u}_{\Ltwo{\shellp{j}}}\bigr)
        \bigl(2^{j/2}\norm{g}_{\Ltwo{\shellp{j}}}\bigr).
\end{equation*}
Bounding every factor of the first kind by $\BstnormO{u}$ and summing the
factors of the second kind gives \eqref{eq:AH-pairing}. This is the elementary
$\ell^1$--$\ell^\infty$ pairing of Part~I, transported to $\Om$.

(b) The inclusion $\Ltwo{\Om} \subset \Bstd{\Om}$ holds because each of the
quantities $\Rzero^{-1/2}\norm{u}_{\Ltwo{\Omin{\Rzero}}}$ and
$2^{-j/2}\norm{u}_{\Ltwo{\shellp{j}}}$ is at most $\norm{u}_{\Ltwo{\Om}}$, the
weights being at most $1$ since $\Rzero \ge 1$ and $j \ge 0$. The inclusion
$\Bspd{\Om} \subset \Ltwo{\Om}$ holds because the weights in $\BnormO{\cdot}$
are at least $1$ and the $\ell^2$ norm of a sequence is at most its $\ell^1$
norm:
\begin{equation*}
  \norm{u}_{\Ltwo{\Om}}^2
    = \norm{u}_{\Ltwo{\Omin{\Rzero}}}^2
      + \sum_{j > \Jzero} \norm{u}_{\Ltwo{\shellp{j}}}^2
    \le \BnormO{u}^2 .
\end{equation*}
For the remaining inclusions of \eqref{eq:AH-inclusions-pert} it suffices, by
Lemma~\ref{lem:AH-equivalence}, to argue with the shell norms
\eqref{eq:B-norms-shell}. For those the proof given on $\Rnp$ in Part~I
\cite[\S2.3.1]{PartI} applies word for word. It uses only the comparison
of $(1+|x|)^{\pm s}$ with $2^{\pm js}$ on $\shell{j}$, and no property of the
underlying domain.

(c) Immediate from the two inclusions
$\Bspd{\Om} \subset \Ltwo{\Om} \subset \Bstd{\Om}$ established at the start of
(b), both of norm at most $1$.
\end{proof}

\subsection{The Rellich-Type Theorem on the Perturbed Domain}
\label{subsec:rellich-pert}

The uniqueness statements that drive \S\ref{sec:lap-pert} concern, in the first
instance, solutions defined only near infinity. We fix the terminology first.

\begin{definition}[The exterior Robin problem]
\label{def:exterior-problem}
Let $E \in \mathbb{R}$ and $R \ge \Rzero$. A function $u$ \emph{solves the
exterior Robin problem at energy $E$ outside $\Bfull{R}$} if
$u \in H^2_{\mathrm{loc}}(\overline{\Omex{R}})$ and
\begin{equation}
  \label{eq:exterior-robin}
  (-\Delta - E)u = 0 \ \text{ in } \Omex{R},
  \qquad
  \gamma_1 u - \alpha\,\gamma_0 u = 0
  \ \text{ on } \Gpert \cap \{|x| > R\}.
\end{equation}
\end{definition}

By \eqref{eq:exterior-coincidence} the set $\Omex{R}$ is a subset of $\Rnp$ and
$\Gpert \cap \{|x|>R\}$ is a subset of $\Gamma$, so \eqref{eq:exterior-robin} is
a problem for the \emph{flat} half-space. The perturbation enters only through
the requirement $R \ge \Rzero$. We record two elementary consequences of
membership in $\Bstd{\Om}$, both of which are used repeatedly below.

\begin{lemma}[Local energy growth and weighted integrability]
\label{lem:local-growth}
Let $u \in \Bstd{\Om}$.
\begin{enumerate}
  \item[\textup{(a)}] For every $R > 1$,
  \begin{equation}
    \label{eq:local-growth}
    \int_{\Om \cap \Bfull{R}} |u|^2\,dx
      \;\le\; 4\,\Csh^2\, R\, \BstnormO{u}^2 .
  \end{equation}
  \item[\textup{(b)}] Let $s > 1/2$, $\sigma > 0$, and write
  $\wsig(\xpar) = \langle \sigma\xpar\rangle^{-s}$. For every $b > 0$,
  \begin{equation}
    \label{eq:slab-weighted}
    \int_{\Om \cap \{0 < x_n < b\}} \wsig(\xpar)^2\,|u|^2\,dx
      \;\le\; C\,\BstnormO{u}^2,
    \qquad C = C(n,s,\sigma,b).
  \end{equation}
\end{enumerate}
\end{lemma}

\begin{proof}
(a) By Lemma~\ref{lem:AH-equivalence}, $\Qsc{u} \le 4\Bstnormsh{u}^2
\le 4\Csh^2\BstnormO{u}^2$, and \eqref{eq:local-growth} is the definition
\eqref{eq:Q-scale} of $\Qsc{u}$.

(b) Decompose the tangential variable dyadically. The slab piece over
$\{|\xpar| < 1\}$ is contained in $\Om \cap \Bfull{1+b}$ and contributes at most
$4\Csh^2(1+b)\BstnormO{u}^2$ by (a). For $j \ge 1$ the piece over
$\{2^{j-1} \le |\xpar| < 2^{j}\}$ is contained in $\Om \cap \Bfull{2^{j}+b}$,
where $\wsig \le (\sigma 2^{j-1})^{-s}$, so by (a) it contributes at most
$4\Csh^2 (\sigma 2^{j-1})^{-2s}(2^{j}+b)\BstnormO{u}^2$. Since $2s > 1$ the
series $\sum_{j \ge 1} 2^{-2s(j-1)}(2^{j}+b)$ converges, which is
\eqref{eq:slab-weighted}.
\end{proof}

The gradient of a solution is controlled by the solution itself, at a cost that
does not depend on where in the exterior the estimate is made. This uniformity
is what allows the next two propositions to sum a dyadic decomposition, and we
isolate it once.

\begin{lemma}[Weighted Caccioppoli inequality in the exterior]
\label{lem:caccioppoli}
Let $E \in \mathbb{R}$ and $R \ge \Rzero$, and let $u$ solve the exterior Robin
problem at energy $E$ outside $\Bfull{R}$. Let $s > 1/2$ and $0 < \sigma \le 1$.
Let $D \subset \Om$ be open with $\dist(D, \Bfull{R}) \ge 1$, and let $\Dfat$ be
the set of points of $\Om$ at distance less than $1$ from $D$. Then
\begin{equation}
  \label{eq:caccioppoli}
  \int_{D} \wsig^2\,|\nabla u|^2\,dx
    \;\le\; C \int_{\Dfat} \wsig^2\,|u|^2\,dx ,
  \qquad C = 4\,(E + 2\alpha^2)_+ + 32 + 10\,s^2 ,
\end{equation}
where $t_+ = \max\{t, 0\}$. In particular the constant does not depend on $D$,
on $R$ or on $\sigma$.
\end{lemma}

\begin{proof}
It suffices to treat bounded $D$. Indeed, for $\rho > 0$ the bounded open set
$D \cap \Bfull{\rho}$ satisfies the same hypothesis, and the points of $\Om$ at
distance less than $1$ from it lie in $\Dfat$. Applying the bounded case to it
bounds the integral of $\wsig^2|\nabla u|^2$ over $D \cap \Bfull{\rho}$ by the
right-hand side of \eqref{eq:caccioppoli}, and monotone convergence as
$\rho \to \infty$ gives the general case.

Let $D$ be bounded. Fix $\chi \in C^\infty(\Rn)$ with $0 \le \chi \le 1$,
$\chi = 1$ on $D$, $|\nabla\chi| \le 2$ and
$\operatorname{supp}\chi \subset \{x \in \Rn : \dist(x, D) < 1\}$, and put
$\varphi = \chi\wsig$. The support of $\chi$ is compact. It does not meet
$\overline{\Bfull{R}}$, because $\dist(x, \Bfull{R}) \ge \dist(D, \Bfull{R})
- \dist(x, D) \ge 1 - \dist(x, D) > 0$ for each of its points. Hence
$\varphi u$ vanishes near
$\partial\Bfull{R}$, its extension by zero to $\Rnp$ lies in $H^1(\Rnp)$ and
has compact support, and every integral below is finite. Moreover $\varphi$ and
$\nabla\varphi$ vanish on $\Rnp \setminus \Dfat$, since
$\Rnp \setminus \overline{\Bfull{R}} = \Omex{R} \subset \Om$ by
\eqref{eq:exterior-coincidence}.

Multiply $(-\Delta - E)u = 0$ by $\varphi^2\bar u$, integrate over $\Omex{R}$ and
take real parts. The only boundary contribution is on $\Gamma$, where
$\gamma_1 u = \alpha\gamma_0 u$ and
$\gamma_0(\varphi^2 u) = \varphi^2\gamma_0 u$. Using the localization identity
$\re\int \nabla u \cdot \overline{\nabla(\varphi^2 u)}
= \int |\nabla(\varphi u)|^2 - \int |\nabla\varphi|^2 |u|^2$, valid for real
$\varphi \in W^{1,\infty}$, we obtain
\begin{equation*}
  \int_{\Rnp} |\nabla(\varphi u)|^2
    - \alpha \int_{\Gamma} |\gamma_0(\varphi u)|^2
  \;=\; E \int_{\Rnp} \varphi^2 |u|^2
    + \int_{\Rnp} |\nabla\varphi|^2 |u|^2 .
\end{equation*}
Apply Lemma~\ref{lem:trace-pert} on $\Rnp$, where $\Cdiv = 0$, with
$\varepsilon = 1/(2\alpha)$. This bounds the boundary term by
$\tfrac12\norm{\nabla(\varphi u)}^2 + 2\alpha^2\norm{\varphi u}^2$, so that
\begin{equation*}
  \tfrac12 \int_{\Rnp} |\nabla(\varphi u)|^2
  \;\le\; (E + 2\alpha^2) \int_{\Rnp} \varphi^2 |u|^2
    + \int_{\Rnp} |\nabla\varphi|^2 |u|^2 .
\end{equation*}
We bound the right-hand side by an integral over $\Dfat$. Since
$\varphi^2 \le \wsig^2$ and $\varphi$ vanishes outside $\Dfat$, the first term
is at most $(E + 2\alpha^2)_+ \int_{\Dfat} \wsig^2|u|^2$. The positive part
is needed for $E < -2\alpha^2$, where the first term is negative and is
dropped. For the second term, differentiation gives
$\nabla'\wsig = -s\sigma^2\xpar\langle\sigma\xpar\rangle^{-2}\wsig$, hence
$|\nabla'\wsig| \le s\sigma\wsig \le s\wsig$ because $\sigma \le 1$, and
\begin{equation*}
  |\nabla\varphi|^2
    \le 2\wsig^2|\nabla\chi|^2 + 2\chi^2|\nabla'\wsig|^2
    \le (8 + 2s^2)\,\wsig^2 .
\end{equation*}
As $\nabla\varphi$ also vanishes outside $\Dfat$, we arrive at
\begin{equation*}
  \int_{\Rnp} |\nabla(\varphi u)|^2
  \;\le\; 2\bigl( (E + 2\alpha^2)_+ + 8 + 2s^2 \bigr)
    \int_{\Dfat} \wsig^2 |u|^2 .
\end{equation*}
We now return to $D$. The set $D$ is open and $\chi = 1$ on it, so
$\nabla\chi = 0$, $\varphi = \wsig$ and $\nabla\varphi = \nabla\wsig$ there,
with $|\nabla\wsig| = |\nabla'\wsig|$ because $\wsig$ does not depend on $x_n$.
Writing $\varphi\nabla u = \nabla(\varphi u) - u\nabla\varphi$ and using
$|a - b|^2 \le 2|a|^2 + 2|b|^2$, we obtain, pointwise on $D$,
\begin{equation*}
  \wsig^2|\nabla u|^2
    \le 2|\nabla(\varphi u)|^2 + 2|\nabla'\wsig|^2|u|^2
    \le 2|\nabla(\varphi u)|^2 + 2s^2\wsig^2|u|^2 .
\end{equation*}
This bound holds on $D$ only, since $\chi$ may be smaller than $1$ on
$\Dfat \setminus D$. Integrate it over $D$. Enlarge the domain of the first
integral to $\Rnp$ and that of the second to $\Dfat$, and insert the preceding
display. This gives \eqref{eq:caccioppoli} with
$C = 4(E + 2\alpha^2)_+ + 32 + 8s^2 + 2s^2$, which is the stated constant. It
involves neither $D$, nor $R$, nor $\sigma$.
\end{proof}

\subsubsection*{Energies above the bulk threshold}

For $E > 0$ the theorem of Part~I applies without modification.

\begin{theorem}[Rellich-type theorem, $E > 0$]
\label{thm:rellich-positive}
Let $n \ge 3$, $\alpha > 0$, $E > 0$ and $R \ge \Rzero$. If $u$ solves the
exterior Robin problem at energy $E$ outside $\Bfull{R}$ and
$u \in \Bstzd{\Om}$, then there is $R_1 \ge R$ with $u = 0$ in $\Omex{R_1}$.
\end{theorem}

\begin{proof}
Let $u^{\mathrm{e}}$ denote the function on $\Rnp$ equal to $u$ on
$\Omex{R} = \Rnp \setminus \overline{\Bfull{R}}$ and to $0$ on
$\Rnp \cap \Bfull{R}$. Comparing \eqref{eq:B-norms-flat} with
\eqref{eq:B-norms-shell} term by term, and using
$\shellp{j} \cap \{|x| > R\} \subset \shell{j} \cap \Om$, gives
$\Bstnorm{u^{\mathrm{e}}} \le \Bstnormsh{u}$, and likewise
$u \in \Bstzd{\Om}$ forces $u^{\mathrm{e}} \in \Bstz$. Since $u$ solves
\eqref{eq:exterior-robin} and $\Omex{R}$, $\Gpert \cap \{|x|>R\}$ are subsets of
$\Rnp$, $\Gamma$ respectively, $u$ satisfies on $\Rnp \setminus
\overline{\Bfull{R}}$ exactly the hypotheses of the Rellich-type theorem of
Part~I \cite[Theorem~2.76]{PartI}, whose conclusion is the assertion.
\end{proof}

\subsubsection*{Energies in the surface band: decay away from the boundary}

Below the bulk threshold the situation is different, and the remainder of this
subsection is devoted to it. Write
\begin{equation}
  \label{eq:band}
  \Eband = (-\alpha^2, 0)
\end{equation}
for the surface band. The first observation is that $\Omex{\Rzero}$ contains a
full slab. Indeed, if $x_n > \Rzero$ then $|x| \ge x_n > \Rzero$, so
\begin{equation}
  \label{eq:slab-inclusion}
  \{x \in \Rn : x_n > R\} \subset \Omex{R}
  \qquad (R \ge \Rzero),
\end{equation}
and on that slab $u$ solves $(-\Delta - E)u = 0$ with no boundary condition at
all. For $E \in \Eband$ the operator $-\Delta - E = -\Delta + |E|$ is uniformly
positive there, and an Agmon-type estimate gives exponential decay in $x_n$. The
rate $\sqrt{|E|}$ is the one predicted by the transverse decomposition. On the
Fourier side the slab equation reads $-\partial_n^2\widehat{u} + \mu(\xipar)^2
\widehat{u} = 0$ with $\mu(\xipar) = (|\xipar|^2 + |E|)^{1/2} \ge |E|^{1/2}$,
uniformly in $\xipar$. This last inequality fails for $E > 0$. There $\mu$
vanishes on the sphere $|\xipar| = \sqrt{E}$ and becomes imaginary inside it.

\begin{proposition}[Exponential decay in the normal variable]
\label{prop:band-normal-decay}
Let $E \in \Eband$, let $R \ge \Rzero$, and let $u$ solve the exterior Robin
problem at energy $E$ outside $\Bfull{R}$ with $u \in \Bstd{\Om}$. Let
$0 < \theta < |E|^{1/2}$ and $s > 1/2$. Then there is
$\sigma_0 = \sigma_0(E,s,\theta) > 0$ such that for every
$0 < \sigma \le \sigma_0$,
\begin{equation}
  \label{eq:band-normal-decay}
  \int_{\{x_n > R+1\}} \wsig(\xpar)^2\, e^{2\theta x_n}
    \bigl( |u|^2 + |\nabla u|^2 \bigr) dx
  \;\le\; C\, \BstnormO{u}^2,
\end{equation}
with $C = C(n,E,s,\theta,\sigma,R)$. Only $u \in \Bstd{\Om}$ is used. The
hypothesis $u \in \Bstzd{\Om}$ is not needed.
\end{proposition}

\begin{proof}
Fix $\lambda$ with $\theta^2/|E| < \lambda < 1$ and put
$\varepsilon = (1-\lambda)/2$, so that
\begin{equation}
  \label{eq:agmon-constants}
  \kappa_1 := 1 - \lambda - \varepsilon > 0,
  \qquad
  \kappa_2 := |E| - \theta^2/\lambda > 0 .
\end{equation}
Choose $\zeta \in C^\infty(\mathbb{R})$ with $0 \le \zeta \le 1$, $\zeta = 0$
for $x_n \le R$ and $\zeta = 1$ for $x_n \ge R+1$, and for $T > R+1$ set
$\phi_T(x_n) = \theta\min(x_n, T)$, a Lipschitz function with
$|\phi_T'| \le \theta$. Fix also a nonincreasing $\chi \in C^\infty(\mathbb{R})$
with $\chi = 1$ on $(-\infty, 1]$ and $\chi = 0$ on $[2, \infty)$, and for
$\varrho > R+1$ put $\chi_\varrho(x) = \chi(|x|/\varrho)$. Write
$\Psi = \wsig \zeta \chi_\varrho$. Then $\Psi \in W^{1,\infty}$ has compact
support, with
\begin{equation}
  \label{eq:cutoff-gradient}
  |\nabla \Psi|^2
    \le 2 s^2\sigma^2 \Psi^2 + 2\,\wsig^2\,|\nabla(\zeta\chi_\varrho)|^2 ,
\end{equation}
because $|\nabla'\wsig| \le s\sigma\wsig$.

Multiply $(-\Delta + |E|)u = 0$ by $\overline{\Psi^2 e^{2\phi_T} u}$, integrate
over the slab $\{x_n > R\}$, which lies in $\Omex{R}$ by
\eqref{eq:slab-inclusion}, and take real parts. There is no boundary
contribution, since $\Psi$ vanishes for $x_n \le R$ and for
$|x| \ge 2\varrho$. Every integral is finite, because $\Psi$ has compact support
in $\overline{\Omex{R}}$ and $u \in H^2_{\mathrm{loc}}(\overline{\Omex{R}})$.
One obtains
\begin{equation*}
  \int \Psi^2 e^{2\phi_T}|\nabla u|^2
  + 2\operatorname{Re}\!\int \Psi e^{2\phi_T}\bar u\,\nabla\Psi\cdot\nabla u
  + 2\operatorname{Re}\!\int \Psi^2 e^{2\phi_T}\bar u\,\nabla\phi_T\cdot\nabla u
  + |E|\int \Psi^2 e^{2\phi_T}|u|^2 = 0 .
\end{equation*}
By Young's inequality the third term is bounded in modulus by
\begin{equation*}
  \lambda \int \Psi^2 e^{2\phi_T}|\nabla u|^2
    + \frac{\theta^2}{\lambda} \int \Psi^2 e^{2\phi_T}|u|^2 ,
\end{equation*}
and the second by
\begin{equation*}
  \varepsilon \int \Psi^2 e^{2\phi_T}|\nabla u|^2
    + \frac{1}{\varepsilon} \int |\nabla\Psi|^2 e^{2\phi_T}|u|^2 .
\end{equation*}
With \eqref{eq:agmon-constants} and \eqref{eq:cutoff-gradient} this gives
\begin{equation*}
  \kappa_1 \int \Psi^2 e^{2\phi_T}|\nabla u|^2
    + \Bigl( \kappa_2 - \frac{2s^2\sigma^2}{\varepsilon} \Bigr)
      \int \Psi^2 e^{2\phi_T}|u|^2
  \;\le\; \frac{2}{\varepsilon}
    \int \wsig^2 |\nabla(\zeta\chi_\varrho)|^2 e^{2\phi_T}|u|^2 .
\end{equation*}
Let $\sigma_0$ be defined by $2s^2\sigma_0^2/\varepsilon = \kappa_2/2$. For
$\sigma \le \sigma_0$ the second coefficient on the left is at least
$\kappa_2/2$. On the right we use
$|\nabla(\zeta\chi_\varrho)|^2 \le 2|\zeta'|^2 + 2|\nabla\chi_\varrho|^2$. The
function $\zeta'$ is supported in $\{R < x_n < R+1\}$, where
$e^{2\phi_T} \le e^{2\theta(R+1)}$, so its part is at most $C\BstnormO{u}^2$ by
Lemma~\ref{lem:local-growth}(b). The function $\nabla\chi_\varrho$ is supported
in $\{\varrho \le |x| \le 2\varrho\}$, where
$|\nabla\chi_\varrho| \le \norm{\chi'}_\infty/\varrho$, $\wsig \le 1$ and
$e^{2\phi_T} \le e^{2\theta T}$. By Lemma~\ref{lem:local-growth}(a) its part is
at most $32\,\varepsilon^{-1}\norm{\chi'}_\infty^2\,\Csh^2\,e^{2\theta T}
\varrho^{-1}\BstnormO{u}^2$, which tends to $0$ as $\varrho \to \infty$. As
$\varrho$ increases, $\Psi^2$ increases to $\wsig^2\zeta^2$. Monotone
convergence as $\varrho \to \infty$ therefore bounds the left-hand side with
$\Psi = \wsig\zeta$ by $C\BstnormO{u}^2$, and this bound is independent of $T$.
Since $\Psi = \wsig$ and $\phi_T = \theta x_n$ on $\{R + 1 < x_n < T\}$,
monotone convergence as $T \to \infty$ yields \eqref{eq:band-normal-decay}.
\end{proof}

\begin{remark}[The scope of the normal decay]
\label{rem:band-vs-positive}
Proposition~\ref{prop:band-normal-decay} is the first point at which the two
energy regimes differ. The estimate uses only that $-\Delta - E$ is uniformly
positive on the slab, that is, that $E$ lies below the bottom $0$ of the bulk
continuum. It uses neither the boundary condition, nor the value of $\alpha$,
nor the geometry inside $\Bfull{\Rzero}$, nor $u \in \Bstzd{\Om}$.
Correspondingly it says nothing about the surface channel, whose transverse
profile $\psizero$ decays at the rate $\alpha$. This is compatible with
\eqref{eq:band-normal-decay}, since $\mu(\xipar) = (\ksurf^2 + |E|)^{1/2} =
\alpha$ for $|\xipar| = \ksurf$. The surface channel is treated in
Proposition~\ref{prop:band-tangential-decay}.
\end{remark}

\begin{proposition}[Exponential decay in the tangential variables]
\label{prop:band-tangential-decay}
Let $E \in \Eband$, let $R \ge \Rzero$, and let $u$ solve the exterior Robin
problem at energy $E$ outside $\Bfull{R}$ with $u \in \Bstd{\Om}$. Assume
that for some $\rho_1 \ge R+1$ the surface component vanishes outside
$\Bfull{\rho_1}$, that is,
\begin{equation}
  \label{eq:usurf-vanishes}
  \usurf(\xpar) := \int_0^\infty u(\xpar, x_n)\,\psizero(x_n)\,dx_n = 0
  \qquad (|\xpar| > \rho_1).
\end{equation}
Let $s > 1/2$ and let $\sigma > 0$ satisfy $\tau := s\sigma|E|^{-1/2} < 1$.
Then, with $\Cyl{\rho} = \{\xpar \in \Rnm : |\xpar| > \rho\} \times (0,\infty)$,
the quantity
\begin{equation}
  \label{eq:G-def}
  G(\rho) = \int_{\Cyl{\rho}} \wsig^2
    \bigl( |\nabla' u|^2 + |E|\,|u|^2 \bigr) dx
\end{equation}
is finite for $\rho \ge \rho_1$ and satisfies
\begin{equation}
  \label{eq:G-decay}
  G(\rho) \le G(\rho_1)\, e^{-2\beta(\rho - \rho_1)},
  \qquad \beta := (1-\tau)\,|E|^{1/2} > 0 .
\end{equation}
Consequently $u \in \Ltwo{\Cyl{\rho_1}}$, and for every $\beta' < \beta$,
\begin{equation}
  \label{eq:tangential-L2-decay}
  \int_{\Cyl{\rho_1}} e^{2\beta'|\xpar|}\,|u|^2\,dx
    \;\le\; C\, G(\rho_1),
  \qquad C = C(n,E,s,\sigma,\beta',\rho_1).
\end{equation}
Finally $G(\rho_1) \le C'\,\BstnormO{u}^2$ with $C' = C'(n,E,s,\sigma,\theta,R)$,
for every rate $\theta$ admissible in
Proposition~\ref{prop:band-normal-decay}. In particular $C'$ does not depend on
$\rho_1$.
\end{proposition}

\begin{proof}
\emph{Finiteness.} Since $|x| \ge |\xpar|$, every point of $\Cyl{R+1}$ lies
in $\Omex{R}$ and has distance at least $1$ from $\Bfull{R}$. Split
$\Cyl{R+1}$ at $x_n = R+1$. On
$\{x_n > R+1\}$ the integrand is dominated by the left-hand side of
\eqref{eq:band-normal-decay}. On $\{0 < x_n < R+1\}$ the term $|E||u|^2$ is
controlled by Lemma~\ref{lem:local-growth}(b) and the term $|\nabla' u|^2$ by
Lemma~\ref{lem:caccioppoli}, whose hypothesis is the distance just noted. Both
bounds are of the form $C\,\BstnormO{u}^2$ with $C$ depending on $R$ but not on
$\rho_1$. Hence $G(\rho_1) \le G(R+1) \le C'\,\BstnormO{u}^2 < \infty$, which
is the last assertion of the proposition. The same decomposition shows, by
Fubini's theorem, that $u(\xpar, \cdot) \in H^1(0,\infty)$ for almost every
$\xpar$ with $|\xpar| > \rho_1$. In particular the trace
$u(\xpar, 0)$ and the integral in \eqref{eq:usurf-vanishes} are defined for
almost every such $\xpar$, the latter because
$\norm{e^{\theta \cdot}u(\xpar,\cdot)}_{\Ltwo{\Rhalf}}
\norm{e^{-\theta\cdot}\psizero}_{\Ltwo{\Rhalf}} < \infty$ by
Proposition~\ref{prop:band-normal-decay}.

\emph{The energy identity.} Fix $\rho \ge \rho_1$. Multiplying
$(-\Delta - E)u = 0$ by $\overline{\wsig^2 u}$, integrating over $\Cyl{\rho}$
and taking real parts gives, since $\gamma_1 u = \alpha\gamma_0 u$ on
$\{x_n = 0,\ |\xpar| > \rho\}$,
\begin{equation}
  \label{eq:cylinder-identity}
  \begin{aligned}
    \int_{\Cyl{\rho}} \wsig^2 |\nabla u|^2
      &- \alpha \int_{|\xpar| > \rho} \wsig^2 |u(\xpar,0)|^2 \,d\xpar
       + |E| \int_{\Cyl{\rho}} \wsig^2 |u|^2
       + 2\operatorname{Re}\!\int_{\Cyl{\rho}}
           \wsig\, \bar u\, \nabla'\wsig \cdot \nabla' u
    \\
      &= -\operatorname{Re} \int_{|\xpar| = \rho}
           \wsig^2\, \bar u \,\partial_r u \, dS ,
  \end{aligned}
\end{equation}
where $\partial_r$ is the derivative in the radial direction of $\xpar$ and
$dS$ is the surface measure on $\{|\xpar| = \rho\}\times(0,\infty)$. The
divergence theorem is applied on the truncated cylinders
$\{\rho < |\xpar| < \varrho\} \times (0,\tsub)$ and the limits
$\varrho \to \infty$, $\tsub \to \infty$ are taken. The two families of
boundary terms vanish along suitable sequences, because
$G(\rho_1) < \infty$ and by \eqref{eq:band-normal-decay} respectively. Note
that $\wsig$ is radial in $\xpar$, so it is constant on
$\{|\xpar| = \rho\}$ and does not interfere with the boundary term.

\emph{Transverse positivity.} The half-line Robin form is nonnegative on the
orthogonal complement of its ground state $\psizero$. Explicitly, for
$v \in H^1(0,\infty)$ with $\int_0^\infty v\,\psizero \,dx_n = 0$,
\begin{equation}
  \label{eq:transverse-positivity}
  \int_0^\infty |v'(x_n)|^2 dx_n - \alpha\,|v(0)|^2 \;\ge\; 0 ,
\end{equation}
since the spectrum of $\hfib(0)$ on $\{\psizero\}^\perp$ is $[0,\infty)$.
Applying \eqref{eq:transverse-positivity} to $v = u(\xpar,\cdot)$, which is
admissible for almost every $|\xpar| > \rho \ge \rho_1$ by
\eqref{eq:usurf-vanishes}, multiplying by $\wsig(\xpar)^2 \ge 0$ and
integrating in $\xpar$, we find that the first two terms of
\eqref{eq:cylinder-identity} together dominate
$\int_{\Cyl{\rho}}\wsig^2|\nabla' u|^2$. It is here, and only here, that the
hypothesis \eqref{eq:usurf-vanishes} is used.

\emph{The weight error.} By $|\nabla'\wsig| \le s\sigma\wsig$ and Young's
inequality with parameter $|E|^{-1/2}$,
\begin{equation*}
  \Bigl| 2\operatorname{Re}\!\int_{\Cyl{\rho}}
    \wsig \bar u\, \nabla'\wsig\cdot\nabla' u \Bigr|
  \;\le\; \frac{s\sigma}{|E|^{1/2}}
    \Bigl( \int_{\Cyl{\rho}} \wsig^2|\nabla' u|^2
         + |E| \int_{\Cyl{\rho}} \wsig^2 |u|^2 \Bigr)
  \;=\; \tau\, G(\rho).
\end{equation*}
Combining the last three paragraphs,
\begin{equation}
  \label{eq:G-vs-flux}
  (1-\tau)\, G(\rho)
    \;\le\; -\operatorname{Re}\int_{|\xpar|=\rho}
       \wsig^2\, \bar u\, \partial_r u\, dS .
\end{equation}

\emph{Gronwall.} Write $\xpar = r\omega$ with $r = |\xpar|$ and
$\omega \in \Ssurf$, and keep $x_n$ as it is, so that
$dx = r^{n-2}\,dr\,d\omega\,dx_n$. The set $\{|\xpar| = r\} \times (0,\infty)$
is a hypersurface of $\Rn$, not a sphere in $\Rnm$. Its surface measure, the
measure $dS$ of \eqref{eq:cylinder-identity}, is
$dS = r^{n-2}\,d\omega\,dx_n$. Tonelli's theorem therefore gives
$G(\rho) = \int_\rho^\infty g(r)\,dr$, where
\begin{equation*}
  g(r) = \int_{|\xpar| = r} \wsig^2
    \bigl( |\nabla' u|^2 + |E| |u|^2 \bigr) dS
  = r^{n-2} \int_0^\infty \!\! \int_{\Ssurf} \wsig^2
    \bigl( |\nabla' u|^2 + |E|\,|u|^2 \bigr)(r\omega, x_n)\,
    d\omega\, dx_n .
\end{equation*}
The function $g$ lies in $L^1(\rho_1,\infty)$ because $G(\rho_1) < \infty$. In
particular $g(r) < \infty$ for almost every $r$, $G$ is absolutely continuous on
$[\rho_1,\infty)$ and $G'(\rho) = -g(\rho)$ almost everywhere. We use
\eqref{eq:G-vs-flux} only at radii $\rho$ with $g(\rho) < \infty$, where its
right-hand side is finite. Since $\partial_r = (\xpar/|\xpar|) \cdot \nabla'$,
one has $|\partial_r u| \le |\nabla' u|$. Young's inequality with the parameter
$|E|^{-1/2}$ therefore bounds the right-hand side of \eqref{eq:G-vs-flux} by
\begin{equation*}
  \frac{|E|^{1/2}}{2}\int_{|\xpar|=\rho}\wsig^2|u|^2
    + \frac{1}{2|E|^{1/2}}\int_{|\xpar|=\rho}\wsig^2|\nabla' u|^2
  \;=\; \frac{1}{2|E|^{1/2}}\, g(\rho),
\end{equation*}
the two coefficients having been matched by that choice of parameter. Hence
\begin{equation*}
  (1-\tau)\,G(\rho) \;\le\; -\frac{G'(\rho)}{2|E|^{1/2}} ,
  \qquad\text{that is,}\qquad
  G'(\rho) \le -2\beta\,G(\rho) ,
\end{equation*}
for almost every $\rho \ge \rho_1$. The function $\rho \mapsto e^{2\beta\rho}G(\rho)$
is absolutely continuous, and its derivative
$e^{2\beta\rho}\bigl( G'(\rho) + 2\beta G(\rho) \bigr)$ is nonpositive almost
everywhere. It is therefore nonincreasing on $[\rho_1,\infty)$, which is
\eqref{eq:G-decay}.

\emph{Removing the weight.} There is $c_\sigma > 0$ such that
$\wsig(\xpar)^2 \ge c_\sigma (1+\rho)^{-2s}$ holds whenever
$\rho < |\xpar| < \rho+1$. Hence \eqref{eq:G-decay} gives
\begin{equation*}
  \int_{\Cyl{\rho}\setminus\Cyl{\rho+1}} |E|\,|u|^2\,dx
    \;\le\; c_\sigma^{-1}(1+\rho)^{2s}\, G(\rho_1)\,
       e^{-2\beta(\rho-\rho_1)} .
\end{equation*}
Multiplying by $e^{2\beta'(\rho+1)}$ and summing over
$\rho \in \rho_1 + \mathbb{N}$, the polynomial factor being absorbed by
$e^{-2(\beta-\beta')\rho}$, gives $u \in \Ltwo{\Cyl{\rho_1}}$ and the bound
\eqref{eq:tangential-L2-decay}, the constant collecting $|E|^{-1}$,
$c_\sigma^{-1}$, $e^{2\beta'\rho_1}$ and the value of the convergent series.
\end{proof}

\subsubsection*{Scope and limits of the band results}

Propositions~\ref{prop:band-normal-decay} and
\ref{prop:band-tangential-decay} together turn a solution of the exterior
problem into a square-integrable one. Combined with the surface step, which is
insensitive to the sign of $E$, they identify the Rellich-type conclusion on the
band with a statement about the point spectrum.

\begin{theorem}[The band: reduction to an eigenvalue problem]
\label{thm:band-equivalence}
Let $n \ge 3$, $\alpha > 0$, $E \in \Eband$ and $R \ge \Rzero$, and let $u$
solve the exterior Robin problem at energy $E$ outside $\Bfull{R}$ with
$u \in \Bstzd{\Om}$. Then there are $a \ge R$ and $\beta > 0$ with
\begin{equation}
  \label{eq:band-L2}
  \int_{\Omex{a}} e^{2\beta |x|}\,|u|^2\,dx < \infty .
\end{equation}
If, in addition, the surface component \eqref{eq:usurf-vanishes} vanishes
already outside $\Bfull{R}$ itself, then one may take $a = R+1$, and there are
$\beta > 0$ and $C$, depending only on $n$, $\alpha$, $R$ and a compact
subinterval of $\Eband$ containing $E$, such that
\begin{equation}
  \label{eq:band-L2-quantitative}
  \int_{\Omex{R+1}} e^{2\beta |x|}\,|u|^2\,dx
    \;\le\; C\,\BstnormO{u}^2 .
\end{equation}
Suppose now that $u$ solves \eqref{eq:exterior-robin} on all of $\Om$. By this
we mean that $u \in H^2_{\mathrm{loc}}(\overline{\Om})$, that
$(-\Delta - E)u = 0$ in $\Om$, and that $\gamma_1 u = \alpha \gamma_0 u$ on
$\Gpert$. Then $u \in \Dom(\Hpert)$ and $\Hpert u = Eu$.
\end{theorem}

\begin{proof}
The surface component $\usurf$ of \eqref{eq:usurf-vanishes} satisfies
$(-\Delta' - \ksurf^2)\usurf = 0$ on $\{|\xpar| > R\}$, obtained by pairing
\eqref{eq:exterior-robin} with $\psizero$ and integrating by parts in $x_n$,
the boundary terms cancelling because $\psizero'(0) + \alpha\psizero(0) = 0$.
Moreover $\usurf \in \Bstz(\Rnm)$ because $u \in \Bstzd{\Om}$. Since
$\ksurf = (E+\alpha^2)^{1/2} > 0$ for $E \in \Eband$, the classical Rellich
theorem in dimension $n-1$ gives $\usurf = 0$ outside some ball
$\Bfull{\rho_1} \subset \Rnm$. This is Steps~1--3 of
\cite[Theorem~2.76]{PartI}, and they are valid for every $E > -\alpha^2$. It is
only the subsequent treatment of the bulk component that requires $E > 0$.

Now fix $s > 1/2$ and $\sigma$ small enough for both
Proposition~\ref{prop:band-normal-decay} and
Proposition~\ref{prop:band-tangential-decay}, and enlarge $\rho_1$ so that
$\rho_1 \ge R+1$. The vanishing \eqref{eq:usurf-vanishes} persists under this
enlargement. Set $a = \rho_1$ and split $\Omex{a\sqrt2}$ according to
which of $|\xpar|$ and $x_n$ is the larger. Since
$|x| \le \sqrt2\,\max(|\xpar|, x_n)$, on each piece the larger of the two
variables exceeds $a$.

On the piece where $|\xpar|$ is the larger we have $|\xpar| > a$, and
$e^{2\beta|x|} \le e^{2\sqrt2\beta|\xpar|}$, so
\eqref{eq:tangential-L2-decay} applies as soon as $\sqrt2\beta \le \beta'$. On
the piece where $x_n$ is the larger we have $x_n > a \ge R+1$, and
$e^{2\beta|x|} \le e^{2\sqrt2\beta x_n}$. There the factor $\wsig(\xpar)^{-2}$
needed to pass from \eqref{eq:band-normal-decay} to an unweighted integral is at
most $\langle\sigma x_n\rangle^{2s}$, since $|\xpar| < x_n$, and a polynomial in
$x_n$ is absorbed by any gap between $2\sqrt2\beta$ and $2\theta$. Choose
$\beta$ with $\sqrt2\beta < \min(\beta', \theta)$. It remains to treat the
bounded set $\Omex{a} \setminus \Omex{a\sqrt2}$. There
$e^{2\beta|x|} \le e^{2\sqrt2\beta a}$, and Lemma~\ref{lem:local-growth}(a)
bounds its contribution by
$4\sqrt2\,\Csh^2\,a\,e^{2\sqrt2\beta a}\BstnormO{u}^2$. Together these give
\eqref{eq:band-L2}.

If \eqref{eq:usurf-vanishes} holds already with $\rho_1 = R$, then the
enlargement gives $\rho_1 = a = R+1$, and every constant entering the last
paragraph is quantitative. The right-hand side of
\eqref{eq:band-normal-decay} is $C\BstnormO{u}^2$, and so is the right-hand
side of \eqref{eq:tangential-L2-decay}, since $G(\rho_1) \le C\BstnormO{u}^2$
by Proposition~\ref{prop:band-tangential-decay}. For $E$ ranging over a compact
subinterval of $\Eband$ one may fix $s$, $\theta$, $\sigma$ and $\beta'$ once
and for all, since $|E|$ is then bounded away from $0$. This lower bound is the
only property of the subinterval that enters, and the resulting constants may
therefore be taken uniform in $E$. This is
\eqref{eq:band-L2-quantitative}.

For the second assertion, \eqref{eq:band-L2} gives $u \in \Ltwo{\Omex{a}}$,
while $u \in H^2(\Omin{a})$ by hypothesis and (G1). On $\Omex{a}$ the domain is
the flat half-space and the problem is invariant under translation in $\xpar$,
so the local estimate \eqref{eq:app-local} holds on unit half-balls with a
constant independent of the centre. Summing it over a locally finite cover of
$\Omex{a}$ by such half-balls upgrades $u \in \Ltwo{\Omex{a}}$ to
$u \in H^2(\Omex{a})$. Hence $u \in H^2(\Om) \cap \Ltwo{\Om}$, and $u$ satisfies
the boundary condition. By \eqref{eq:Hpert-def} it follows that
$u \in \Dom(\Hpert)$ and $\Hpert u = Eu$.
\end{proof}

The uniqueness statement that \S\ref{sec:lap-pert} actually consumes concerns
solutions on all of $\Om$, and for those the two energy regimes can now be
placed side by side.

\begin{corollary}[Uniqueness on $\Om$]
\label{cor:global-uniqueness}
Let $n \ge 3$, $\alpha > 0$, let $E > -\alpha^2$ with $E \ne 0$, and let
$u \in \Bstzd{\Om} \cap H^2_{\mathrm{loc}}(\overline{\Om})$ satisfy
\begin{equation}
  \label{eq:global-problem}
  (-\Delta - E)u = 0 \ \text{ in } \Om,
  \qquad
  \gamma_1 u - \alpha\,\gamma_0 u = 0 \ \text{ on } \Gpert .
\end{equation}
\begin{enumerate}
  \item[\textup{(a)}] If $E > 0$, then $u = 0$.
  \item[\textup{(b)}] If $E \in \Eband$, then $u \in \Dom(\Hpert)$ and
  $\Hpert u = Eu$. Accordingly $u = 0$ for every such $u$ if and only if
  $E \notin \spp(\Hpert)$.
\end{enumerate}
\end{corollary}

\begin{proof}
(a) Theorem~\ref{thm:rellich-positive}, applied with $R = \Rzero$, gives
$R_1 \ge \Rzero$ with $u = 0$ on $\Omex{R_1}$. Since $u$ solves
\eqref{eq:global-problem} throughout $\Om$ and $\Om$ is connected by
Assumption~\ref{assump:geometry}, Proposition~\ref{prop:ucp} of
Appendix~\ref{app:ucp} propagates the vanishing from the open set
$\Omex{R_1}$ to all of $\Om$.

(b) The first assertion is Theorem~\ref{thm:band-equivalence}. If
$E \notin \spp(\Hpert)$, then $\Hpert u = Eu$ forces $u = 0$. Conversely, let
$u \ne 0$ be an eigenfunction of $\Hpert$ at $E$. Then $u \in H^2(\Om)$ and $u$
satisfies the boundary condition by \eqref{eq:Hpert-def}, so $u$ satisfies
\eqref{eq:global-problem}. Moreover $u \in \Ltwo{\Om} \subset \Bstzd{\Om}$ by
Proposition~\ref{prop:AH-elementary}(b). Hence the vanishing fails for this $u$.
\end{proof}

\begin{remark}[Open and closed waveguides]
\label{rem:band-closed-waveguide}
Part~(b) of Corollary~\ref{cor:global-uniqueness} is weaker than part~(a) for a
structural reason. Decompose a solution in the transverse variable, as in
\eqref{eq:surface-band}--\eqref{eq:psizero}. The half-line Robin operator has
the simple eigenvalue $-\alpha^2$ with eigenfunction $\psizero$, together with
the continuum $[0,\infty)$. At energy $E$ the transverse mode of eigenvalue
$\lambda$ propagates in the tangential directions, with wavenumber
$(E-\lambda)^{1/2}$, if and only if $\lambda < E$. For $E > 0$ the propagating
modes are $\psizero$ and the continuum modes with $\lambda \in [0,E)$. For
$E \in \Eband$ only $\psizero$ remains.

This is the dichotomy between open and closed waveguides identified by Hazard
\cite{Hazard2015}. In an open waveguide the modal components of a trapped mode
extend analytically in the modal parameter and vanish on the interval of
propagating modes, hence identically. In a closed waveguide the propagating
modes form a finite set, and trapped modes do occur. For $E > 0$ our
configuration is open in this sense, and Theorem~\ref{thm:rellich-positive} is
the corresponding rigidity. On the band the bulk medium is below cutoff, that
is, $\mu(\xipar) = (|\xipar|^2+|E|)^{1/2}$ never vanishes. This is what makes
Proposition~\ref{prop:band-normal-decay} work, and it also removes the
continuum of propagating modes on which the rigidity for $E > 0$ rests.

We therefore do not expect the point spectrum of $\Hpert$ in $\Eband$ to be
empty in general, and we make no assumption to that effect. By
Theorem~\ref{thm:pert-sa-spec} its elements are embedded eigenvalues. They form
a discrete subset of $\Eband$ by Proposition~\ref{prop:exceptional-discrete}
below, a consequence of the exponential decay established here.
\end{remark}

\subsection{Outgoing Solutions and the Uniqueness Theorem}
\label{subsec:outgoing}

The limiting absorption principle is obtained in \S\ref{subsec:lap-pert-main} by
a compactness argument, and the only input to that argument which is not
elementary is a uniqueness statement: a solution of the homogeneous problem on
$\Om$ which is outgoing must vanish. Corollary~\ref{cor:global-uniqueness}
supplies such a statement, but only for solutions in $\Bstzd{\Om}$. There is no
evident reason for an outgoing solution to lie in that space. The outgoing
resolvent applied to a compactly supported source does not, its far field being
precisely what makes it outgoing. The content of this subsection is that a
solution of the homogeneous problem nevertheless does.

The mechanism is a flux computation. If $u$ solves the homogeneous problem on
all of $\Om$, the energy flux through every large sphere vanishes, because the
Robin boundary term is real. For an outgoing solution, on the other hand, the
flux through a large sphere is the square of its far field. The two statements
are joined by Part~I's formula for the shell average of the outgoing resolvent,
which is valid on both bands. In particular the argument does not pass through
the radiation condition of Part~I, whose formulation presupposes $E > 0$. The
perturbed radiation condition is taken up in \S\ref{sec:scattering-solution}.

\subsubsection*{Cut-offs}

Fix once and for all two radii with
\begin{equation}
  \label{eq:cutoff-radii}
  \Rzero < R_1 < R_2 .
\end{equation}
A function $\eta \in C^{\infty}(\Rn)$ is called an \emph{admissible cut-off} if
it is radial, $0 \le \eta \le 1$, and
\begin{equation}
  \label{eq:cutoff-admissible}
  \eta = 0 \ \text{ on } \Bfull{R_1},
  \qquad
  \eta = 1 \ \text{ on } \Rn \setminus \Bfull{R_2} .
\end{equation}
Radiality is needed: it gives
$\partial_{x_n}\eta = 0$ on $\{x_n = 0\}$, which is what makes multiplication by
$\eta$ preserve the Robin condition instead of turning it into an inhomogeneous
one. The same device is used in Appendix~\ref{app:regularity}. We write
\begin{equation}
  \label{eq:commutator}
  [-\Delta, \eta]\,u = -(\Delta\eta)\,u - 2\,\nabla\eta \cdot \nabla u ,
\end{equation}
a first-order operator whose coefficients are supported in the closed annulus
$\{R_1 \le |x| \le R_2\}$.

\begin{lemma}[Multiplication by an admissible cut-off]
\label{lem:gluing}
Let $\eta$ be an admissible cut-off and let
$u \in H^2_{\mathrm{loc}}(\overline{\Om})$ satisfy
$\gamma_1 u - \alpha\gamma_0 u = 0$ on $\Gpert$. Then the extension of $\eta u$
by zero, again denoted $\eta u$, belongs to
$H^2_{\mathrm{loc}}(\overline{\Rnp})$, satisfies
$\gamma_1(\eta u) - \alpha\gamma_0(\eta u) = 0$ on $\Gamma$, and
\begin{equation}
  \label{eq:gluing}
  (-\Delta - z)(\eta u)
    = \eta\,(-\Delta - z)u + [-\Delta,\eta]\,u
  \qquad \text{in } \Rnp,
\end{equation}
for every $z \in \mathbb{C}$. Moreover:
\begin{enumerate}
  \item[\textup{(a)}] $[-\Delta,\eta]u \in \Bsp$, with
  $\Bnorm{[-\Delta,\eta]u} \le C_\eta \norm{u}_{H^1(\Omin{R_2})}$.
  \item[\textup{(b)}] If $u \in \Bstd{\Om}$ then $\eta u \in \Bst$ with
  $\Bstnorm{\eta u} \le \BstnormO{u}$, and $u \in \Bstzd{\Om}$ if and only if
  $\eta u \in \Bstz$.
\end{enumerate}
\end{lemma}

\begin{proof}
By \eqref{eq:exterior-coincidence} the set $\Omex{\Rzero}$ coincides with
$\Rnp \setminus \overline{\Bfull{\Rzero}}$, and $\eta u$ vanishes on
$\Bfull{R_1} \supset \overline{\Bfull{\Rzero}}$ by \eqref{eq:cutoff-radii}. The
extension by zero is therefore well defined and lies in
$H^2_{\mathrm{loc}}(\overline{\Rnp})$. On $\Gamma \cap \{|x| > \Rzero\}$ the
outward normal of $\Om$ is $-e_n$ and
$\gamma_1(\eta u) = \eta\,\gamma_1 u + (\partial_\nu \eta)\,\gamma_0 u$. Writing
$\eta(x) = \tilde\eta(|x|)$ gives
$\partial_{x_n}\eta = \tilde\eta'(|x|)\,x_n/|x|$, which vanishes on
$\{x_n = 0\}$, so that
$\gamma_1(\eta u) = \eta\,\gamma_1 u = \alpha\,\gamma_0(\eta u)$. On the
remaining part of $\Gamma$ the function $\eta u$ vanishes identically near the
boundary. Identity \eqref{eq:gluing} is the Leibniz rule.

(a) The function $[-\Delta,\eta]u$ lies in $\Ltwo{\Om}$, is supported in
$\Omin{R_2}$, and its $\Ltwo{\Om}$ norm is at most a multiple of
$\norm{u}_{H^1(\Omin{R_2})}$ by \eqref{eq:commutator}. Let now
$v \in \Ltwo{\Om}$ be supported in $\Bfull{R_2}$. Choose $J$ with
$2^{J} \ge R_2$ and apply the Cauchy--Schwarz inequality to the at most
$J - \Jzero + 1$ nonzero terms of \eqref{eq:B-norms-pert}. This gives
$\BnormO{v} \le 2^{J/2}(J - \Jzero + 1)^{1/2}\norm{v}_{\Ltwo{\Om}}$. Finally
$\BnormO{v} = \Bnorm{v}$ by Lemma~\ref{lem:AH-exterior} once $v$ has been
extended by zero, since $v = [-\Delta,\eta]u$ vanishes on $\Omin{\Rzero}$.

(b) Since $\eta u$ vanishes on $\Omin{\Rzero}$, Lemma~\ref{lem:AH-exterior}
gives $\Bstnorm{\eta u} = \BstnormO{\eta u}$, and $\BstnormO{\eta u} \le
\BstnormO{u}$ because $0 \le \eta \le 1$. The same monotonicity gives one
implication for $\Bstz$. Conversely $\eta = 1$ on $\shellp{j}$ as soon as
$2^{j-1} \ge R_2$, so the two sequences whose limits appear in
\eqref{eq:Bzero-pert} and \eqref{eq:Bzero-flat} agree for all large $j$.
\end{proof}

\subsubsection*{Outgoing and incoming solutions}

\begin{definition}[Outgoing and incoming solutions on $\Om$]
\label{def:outgoing-pert}
Let $E \in \Ereg$ and let $u \in \Bstd{\Om} \cap
H^2_{\mathrm{loc}}(\overline{\Om})$ satisfy $\gamma_1 u - \alpha\gamma_0 u = 0$
on $\Gpert$ and $(-\Delta - E)u = 0$ in $\Omex{R_1}$. We call $u$
\emph{outgoing} (respectively \emph{incoming}) at energy $E$ if there is an
admissible cut-off $\eta$ with
\begin{equation}
  \label{eq:outgoing-pert}
  \eta u = \Ralpha{E + i0}\bigl( [-\Delta,\eta]u \bigr)
  \qquad
  \Bigl( \text{respectively } \eta u
    = \Ralpha{E - i0}\bigl( [-\Delta,\eta]u \bigr) \Bigr)
\end{equation}
as elements of $\Bst$.
\end{definition}

The right-hand side of \eqref{eq:outgoing-pert} is defined, because
$[-\Delta,\eta]u \in \Bsp$ by Lemma~\ref{lem:gluing}(a) and
$\Ralpha{E \pm i0} \in \Lop{\Bsp}{\Bst}$ by \cite[Theorem~2.9(b)]{PartI}. Both
sides of \eqref{eq:outgoing-pert} lie in $\Bst$ and solve the same equation on
$\Rnp$, namely $(-\Delta - E)\,\cdot\, = [-\Delta,\eta]u$ with the Robin
condition on $\Gamma$. This holds for the left-hand side by
Lemma~\ref{lem:gluing}, and for the right-hand side by
\cite[Theorem~2.9(e)]{PartI}. What
\eqref{eq:outgoing-pert} asserts is that $\eta u$ is the particular solution
selected by the limit from the upper (respectively lower) half-plane. The
solutions to which the definition will be applied in
\S\ref{subsec:lap-pert-main} arise as limits of $\Rpert{E \pm i\eps}f$, and for
them \eqref{eq:outgoing-pert} holds for every admissible cut-off at once. The
theorem below needs it for one.

\begin{lemma}[The flux through a large sphere vanishes]
\label{lem:zero-flux}
Let $E \in \mathbb{R}$ and let $u \in H^2_{\mathrm{loc}}(\overline{\Om})$
satisfy
\begin{equation*}
  (-\Delta - E)u = 0 \ \text{ in } \Om,
  \qquad
  \gamma_1 u - \alpha\,\gamma_0 u = 0 \ \text{ on } \Gpert .
\end{equation*}
Then
\begin{equation}
  \label{eq:zero-flux}
  \im \int_{\Om \cap \partial\Bfull{R}} \overline{u}\,\partial_r u \, dS = 0
  \qquad (R > \Rzero).
\end{equation}
\end{lemma}

\begin{proof}
The boundary of the bounded Lipschitz domain $\Omin{R}$ is, up to a set of
surface measure zero, the disjoint union of $\Gpert \cap \Bfull{R}$, on which
the outward normal is $\nu$, and $\Om \cap \partial\Bfull{R}$, on which it is
$\partial_r$. Green's first identity gives
\begin{equation*}
  \int_{\Omin{R}} (-\Delta u)\,\overline{u}\,dx
  = \int_{\Omin{R}} |\nabla u|^2\,dx
    - \int_{\Gpert \cap \Bfull{R}} (\gamma_1 u)\,\overline{\gamma_0 u}\,dS
    - \int_{\Om \cap \partial\Bfull{R}} \overline{u}\,\partial_r u\,dS ,
\end{equation*}
every term being finite because $u \in H^2$ on a neighbourhood of
$\overline{\Omin{R}}$. The left-hand side equals
$E\,\norm{u}_{\Ltwo{\Omin{R}}}^2$ and is real. On the right the first term is
real, and so is the second, which equals $\alpha\norm{\gamma_0
u}_{\Ltwo{\Gpert \cap \Bfull{R}}}^2$ by the boundary condition. Taking
imaginary parts gives \eqref{eq:zero-flux}.
\end{proof}

\begin{theorem}[Uniqueness of the outgoing solution]
\label{thm:outgoing-uniqueness}
Let $n \ge 3$, $\alpha > 0$ and $E \in \Ereg$, and let $u \in \Bstd{\Om} \cap
H^2_{\mathrm{loc}}(\overline{\Om})$ satisfy \eqref{eq:global-problem} and be
outgoing or incoming at energy $E$. Then $u \in \Bstzd{\Om}$. Consequently
the following hold.
\begin{enumerate}
  \item[\textup{(a)}] If $E > 0$, then $u = 0$.
  \item[\textup{(b)}] If $E \in \Eband$, then $u \in \Dom(\Hpert)$ and
  $\Hpert u = Eu$. In particular $u = 0$ unless $E \in \spp(\Hpert)$.
\end{enumerate}
\end{theorem}

\begin{proof}
Assume first that $u$ is outgoing, let $\eta$ be an admissible cut-off for
which \eqref{eq:outgoing-pert} holds, and put $g = [-\Delta,\eta]u \in \Bsp$, a
function supported in $\{R_1 \le |x| \le R_2\}$. By \eqref{eq:outgoing-pert}
and \cite[Proposition~2.22(i) and Theorem~2.35]{PartI},
\begin{equation}
  \label{eq:flux-vs-far-field}
  \im\,\bigl( \eta u, g \bigr)
  = \im\,\bigl( \Ralpha{E+i0}g, g \bigr)
  = \pi\,\bigl( \Adens{E}g, g \bigr)
  = \pi\,\bigl\| \Fzero(E)g \bigr\|^2_{\Echan{E}} .
\end{equation}
We show that the left-hand side vanishes.

Write $w = \eta u$, so that $w = u$ on $\{|x| > R_2\}$, and note that
\eqref{eq:gluing} with $z = E$ reads $g = -\Delta w - E w$ on $\Rnp$, since
$(-\Delta - E)u = 0$. Fix $R > R_2$. As $g$ is supported in $\Bfull{R_2}$,
\begin{equation*}
  (\eta u, g)
  = \int_{\Rnp \cap \Bfull{R}} w\,\overline{(-\Delta w)}\,dx
    - E \int_{\Rnp \cap \Bfull{R}} |w|^2\,dx ,
\end{equation*}
and the second term is real because $E$ is. For the first, Green's first
identity on $\Rnp \cap \Bfull{R}$ gives
\begin{equation*}
  \int_{\Rnp \cap \Bfull{R}} (-\Delta w)\,\overline{w}\,dx
  = \int_{\Rnp \cap \Bfull{R}} |\nabla w|^2\,dx
    - \alpha \int_{\Gamma \cap \Bfull{R}} |\gamma_0 w|^2\,dS
    - \int_{\Rnp \cap \partial\Bfull{R}} \overline{w}\,\partial_r w \, dS ,
\end{equation*}
where the boundary condition of Lemma~\ref{lem:gluing} has been used on
$\Gamma$. On the sphere $w = u$, and $\Rnp \cap \partial\Bfull{R}
= \Om \cap \partial\Bfull{R}$ by \eqref{eq:exterior-coincidence}, so the last
integral has vanishing imaginary part by Lemma~\ref{lem:zero-flux}. The first
two terms are real. Hence $\int (-\Delta w)\overline{w}$ is real, hence so is
its complex conjugate, and $\im(\eta u, g) = 0$.

By \eqref{eq:flux-vs-far-field} we conclude $\Fzero(E)g = 0$, and therefore
$\eta u = \Ralpha{E+i0}g \in \Bstz$ by \cite[Lemma~2.65]{PartI}. By
Lemma~\ref{lem:gluing}(b), $u \in \Bstzd{\Om}$.

If $u$ is incoming, observe that the differential expression, the boundary
condition and the coefficients of $[-\Delta,\eta]$ are all real, so that
$\overline{\Ralpha{E-i0}h} = \Ralpha{E+i0}\overline{h}$ for $h \in \Bsp$.
Consequently $\overline{u}$ satisfies \eqref{eq:global-problem} and is
outgoing, and $\overline{u} \in \Bstzd{\Om}$ by the case already treated.

Assertions (a) and (b) are now Corollary~\ref{cor:global-uniqueness}(a) and
(b).
\end{proof}

\subsubsection*{The exceptional set}

By Theorem~\ref{thm:outgoing-uniqueness} the uniqueness on which
\S\ref{subsec:lap-pert-main} rests holds at every $E > 0$ and at every
$E \in \Eband$ outside
\begin{equation}
  \label{eq:exceptional-set}
  \Sexc := \spp(\Hpert) \cap \Eband .
\end{equation}
Remark~\ref{rem:band-closed-waveguide} explains why this set need not be empty,
and no hypothesis is imposed to make it so. What is needed instead is that it
be small, and the following proposition supplies this. Its proof uses no operator theory beyond
Theorem~\ref{thm:pert-sa-spec}: the eigenfunctions on the band decay
exponentially at a rate which is uniform on compact subintervals, and a
uniformly decaying family which is bounded in $H^2(\Om)$ is precompact in
$\Ltwo{\Om}$.

\begin{proposition}[The exceptional set is discrete]
\label{prop:exceptional-discrete}
Let $n \ge 3$, $\alpha > 0$ and let $\Om$ satisfy
Assumption~\ref{assump:geometry}. Then $\Sexc \cap I$ is finite for every
compact interval $I \subset \Eband$, and every $E \in \Sexc$ is an eigenvalue
of $\Hpert$ of finite multiplicity. In particular $\Sexc$ is a closed discrete
subset of $\Eband$.
\end{proposition}

\begin{proof}
If the assertion fails, there are a compact interval $I \subset \Eband$ and an
orthonormal sequence $(u_k)_{k \ge 1}$ in $\Ltwo{\Om}$ with
$u_k \in \Dom(\Hpert)$ and $\Hpert u_k = E_k u_k$ for some $E_k \in I$: either
$\Sexc \cap I$ is infinite, in which case one eigenfunction is taken at each of
countably many distinct eigenvalues and orthogonality is automatic, or some
$E \in \Sexc \cap I$ has an infinite-dimensional eigenspace, in which case an
orthonormal basis of that eigenspace is taken. We show that $(u_k)$ has a
subsequence converging in $\Ltwo{\Om}$, contradicting
$\norm{u_k - u_l} = \sqrt{2}$ for $k \ne l$.

\emph{Step 1: uniform bounds.} By Proposition~\ref{prop:AH-elementary}(b),
$\BstnormO{u_k} \le \norm{u_k}_{\Ltwo{\Om}} = 1$. Since $I \subset \Eband$ we
have $|E_k| < \alpha^2$, so \eqref{eq:apriori-H2} gives
$\norm{u_k}_{H^2(\Om)} \le C(\alpha^2 + 1)$.

\emph{Step 2: the surface component vanishes outside $\Bfull{\Rzero}$.}
Restricted to $\Omex{\Rzero}$, each $u_k$ solves the exterior Robin problem at
energy $E_k$ outside $\Bfull{\Rzero}$ in the sense of
Definition~\ref{def:exterior-problem}. Pairing the equation with $\psizero$ and
integrating by parts in $x_n$, the boundary terms cancelling because
$\psizero'(0) + \alpha\psizero(0) = 0$, shows that the surface component
\eqref{eq:usurf-vanishes} of $u_k$ satisfies
$(-\Delta' - (E_k + \alpha^2))\usurf = 0$ on $\{|\xpar| > \Rzero\}$. Moreover
$\usurf \in \Ltwo{\Rnm}$ by the Cauchy--Schwarz inequality together with
$\norm{\psizero}_{\Ltwo{\Rhalf}} = 1$, so that in particular the shell average
of $\usurf$ over $\Rnm$ tends to $0$. Since $n - 1 \ge 2$ and
$E_k + \alpha^2 > 0$, the classical Rellich uniqueness theorem
\cite{Rellich1943, AgmonHormander1976} makes $\usurf$ vanish outside a compact
set, and unique continuation \cite{Aronszajn1957} on the connected set
$\{|\xpar| > \Rzero\}$, where $\usurf$ is real analytic, propagates this to
\begin{equation*}
  \usurf(\xpar) = 0 \qquad (|\xpar| > \Rzero) .
\end{equation*}
This is Steps~1--3 of \cite[Theorem~2.76]{PartI}, simplified by the fact that
$u_k$ is square integrable. The point is that the radius $\Rzero$ obtained here
does not depend on $k$.

\emph{Step 3: uniform exponential decay.} Step~2 is the additional hypothesis
of Theorem~\ref{thm:band-equivalence} with $R = \Rzero$, so
\eqref{eq:band-L2-quantitative} applies and yields $\beta > 0$ and $C$,
depending only on $n$, $\alpha$, $\Rzero$ and $I$, with
\begin{equation*}
  \int_{\Omex{\Rzero + 1}} e^{2\beta|x|}\,|u_k|^2\,dx
    \;\le\; C\,\BstnormO{u_k}^2 \;\le\; C
  \qquad (k \ge 1).
\end{equation*}

\emph{Step 4: compactness.} Let $\eps > 0$ and choose
$\rho \ge \Rzero + 1$ with $C e^{-2\beta\rho} < \eps^2$. Then
$\norm{u_k}_{\Ltwo{\Omex{\rho}}} < \eps$ for every $k$. By Step~1 the sequence
is bounded in $H^2(\Omin{\rho})$, a bounded Lipschitz domain, so the
Rellich--Kondrachov theorem provides a subsequence converging in
$\Ltwo{\Omin{\rho}}$. Applying this with $\eps = 1/m$ and extracting
diagonally produces a subsequence that is Cauchy in $\Ltwo{\Om}$.
This contradicts the choice of $u_k$.
Hence $\Sexc \cap I$ is finite for every compact $I \subset \Eband$, and
every point of $\Sexc$ is an eigenvalue of finite multiplicity. A subset of
$\Eband$ which meets every compact subinterval in finitely many points has
no accumulation point in $\Eband$, and is therefore discrete and closed in
$\Eband$.
\end{proof}

\begin{remark}[The thresholds]
\label{rem:exceptional-thresholds}
The two endpoints of $\Eband$ behave differently. At the lower endpoint the
proof of Proposition~\ref{prop:exceptional-discrete} gives more than its
statement. By the proof of Theorem~\ref{thm:band-equivalence}, the constants
of Step~3 depend on $E$ only through a lower bound for $|E|$, and Steps~1, 2
and~4 hold for every $E_k \in \Eband$. The argument therefore applies verbatim
to $I = (-\alpha^2, -\delta]$ for every $\delta \in (0, \alpha^2)$, so that
$\Sexc$ does not accumulate at $-\alpha^2$. At the upper endpoint the rate
$\beta$ of Step~3 degenerates as $E \to 0^{-}$ through
Proposition~\ref{prop:band-normal-decay}, and accumulation of eigenvalues of
$\Hpert$ at $0$ is not excluded here. The limiting absorption principle of
\S\ref{subsec:lap-pert-main} is asserted only on compact subsets of
$\Ereg \setminus \Sexc$, near $-\alpha^2$ as well, because the input
\cite[Theorem~2.9(a)]{PartI} is uniform only away from the thresholds. In
Part~I the constant of the resolvent estimate is explicit and grows at both
thresholds at the same rate $\dist(E,\Thr)^{-1/2}$.
\end{remark}

\begin{remark}[The parametrix approach]
\label{rem:no-fredholm}
The classical alternative for a geometric perturbation is a cut-off parametrix
$Q(z)$, built from $\Ralpha{z}$ in the exterior and from an interior
resolvent, with $(\Hpert - z)Q(z) = I + K(z)$ and $K(z)$ compact. The analytic
Fredholm theorem applies to $K(z)$ in each open half-plane, where it gives
nothing because $\Hpert$ is self-adjoint, and $K(z)$ is not analytic across the
real axis. It therefore does not yield the discreteness of $\Sexc$, which would
require a meromorphic continuation across $\Eband$.
Proposition~\ref{prop:exceptional-discrete} obtains it directly, and the
compactness argument of \S\ref{subsec:lap-pert-main} follows
\cite[Ch.~2, \S2.5]{IsozakiKurylev2014}.
\end{remark}

\subsection{The Limiting Absorption Principle for \texorpdfstring{$\Hpert$}{Hpert}}
\label{subsec:lap-pert-main}

We can now prove the principal analytic result of this section. The scheme is
the one of Isozaki and Kurylev \cite[Ch.~2, \S2.5]{IsozakiKurylev2014}: the uniform
bound is obtained by contradiction, a failing sequence of resolvents is shown to
have a convergent subsequence by elliptic regularity and the Rellich--Kondrachov
theorem, and the limit is an outgoing solution of the homogeneous problem, hence
zero by Theorem~\ref{thm:outgoing-uniqueness}. Only three inputs are used
besides that theorem: the uniform bound and the weak-$*$ continuity of the
unperturbed boundary values, from \cite[Theorem~2.9]{PartI}, and the local
elliptic estimate of Proposition~\ref{prop:local-regularity}. In particular no
continuity of $z \mapsto \Ralpha{E \pm i0}$ in operator norm is needed, and none
is available from Part~I.

The compactness is not visible in $\Bstd{\Om}$, whose bounded sets are not
relatively compact in any useful sense at infinity, so the argument is run in
the weighted space $L^{2,-s}(\Om)$ with $s > 1/2$, and the $\Bst$ bound is
recovered afterwards. Two elementary lemmas prepare for this.

\begin{lemma}[Tails and compactly supported elements]
\label{lem:AH-tail}
Let $s > 1/2$.
\begin{enumerate}
  \item[\textup{(a)}] For every $u \in \Bstd{\Om}$ and every $\rho \ge \Rzero$,
  \begin{equation}
    \label{eq:AH-tail}
    \norm{u}_{L^{2,-s}(\Omex{\rho})}
      \le C_s\, \rho^{1/2 - s}\, \BstnormO{u},
    \qquad
    C_s = \Bigl( \frac{2^{2s}}{1 - 2^{1-2s}} \Bigr)^{1/2} .
  \end{equation}
  \item[\textup{(b)}] The elements of $\Bspd{\Om}$ of compact support are dense
  in $\Bspd{\Om}$.
\end{enumerate}
\end{lemma}

\begin{proof}
(a) Since $\rho \ge \Rzero$, the set $\Omex{\rho}$ is contained in the union of
the shells $\shellp{j}$ with $2^{j} > \rho$, and on $\shell{j}$ one has
$(1+|x|)^{-2s} \le 2^{-2s(j-1)}$. Using
$\norm{u}_{\Ltwo{\shellp{j}}}^2 \le 2^{j}\BstnormO{u}^2$ and summing the
geometric series, whose exponent $1 - 2s$ is negative,
\begin{equation*}
  \int_{\Omex{\rho}} (1+|x|)^{-2s}|u|^2\,dx
    \le 2^{2s}\,\BstnormO{u}^2 \sum_{2^{j} > \rho} 2^{j(1-2s)}
    \le \frac{2^{2s}}{1 - 2^{1-2s}}\, \rho^{1-2s}\, \BstnormO{u}^2 ,
\end{equation*}
the last step because the smallest index in the sum satisfies
$2^{j} > \rho$.

(b) For $J > \Jzero$ put $u_J = \mathbf{1}_{\{|x| < 2^{J}\}}u$. Then
$\BnormO{u - u_J} = \sum_{j > J} 2^{j/2}\norm{u}_{\Ltwo{\shellp{j}}}$, a tail of
the convergent series \eqref{eq:B-norms-pert}, so $u_J \to u$.
\end{proof}

\begin{lemma}[From the weighted bound to the \texorpdfstring{$\Bst$}{B*} bound]
\label{lem:Bst-from-weighted}
Let $s > 1/2$ and let $\mathcal{K} \subset \Ereg$ be compact. There is
$C = C(n,\alpha,\Om,\mathcal{K},s)$ such that for every $f \in \Bspd{\Om}$ and
every $z \in \mathbb{C}\setminus\mathbb{R}$ with $\re z \in \mathcal{K}$ and
$|\im z| \le 1$, the function $u = \Rpert{z}f$ satisfies
\begin{equation}
  \label{eq:Bst-from-weighted}
  \BstnormO{u}
    \le C\bigl( \BnormO{f} + \norm{u}_{L^{2,-s}(\Om)} \bigr).
\end{equation}
\end{lemma}

\begin{proof}
Fix an admissible cut-off $\eta$. Since $u \in \Dom(\Hpert) \subset H^2(\Om)$,
Lemma~\ref{lem:gluing} shows that $\eta u \in \Dom(\Halpha)$ and that
$(\Halpha - z)(\eta u) = \eta f + [-\Delta,\eta]u$. As $\im z \ne 0$, this gives
\begin{equation}
  \label{eq:resolvent-gluing}
  \eta u = \Ralpha{z}\bigl( \eta f + [-\Delta,\eta]u \bigr).
\end{equation}
By \cite[Theorem~2.9(a)]{PartI} and Lemma~\ref{lem:gluing},
\begin{equation*}
  \Bstnorm{\eta u}
    \le C\bigl( \Bnorm{\eta f} + \Bnorm{[-\Delta,\eta]u} \bigr)
    \le C\bigl( \BnormO{f} + \norm{u}_{H^1(\Omin{R_2})} \bigr),
\end{equation*}
the constant being uniform for $\re z \in \mathcal{K}$. The local estimate
\eqref{eq:app-local}, applied with $\rho = R_2$, $\rho' = R_2 + 1$, $F = f$ and
$M = \sup_{\mathcal{K}}|\,\cdot\,| + 1$, bounds
$\norm{u}_{H^1(\Omin{R_2})}$ by
$C(\norm{f}_{\Ltwo{\Om}} + \norm{u}_{\Ltwo{\Omin{R_2+1}}})$, and here
$\norm{f}_{\Ltwo{\Om}} \le \BnormO{f}$ by
Proposition~\ref{prop:AH-elementary}(b) while
$\norm{u}_{\Ltwo{\Omin{R_2+1}}} \le (2 + R_2)^{s}\norm{u}_{L^{2,-s}(\Om)}$.
Finally $(1-\eta)u$ is supported in $\Omin{R_2}$, so
$\BstnormO{(1-\eta)u} \le \norm{u}_{\Ltwo{\Omin{R_2}}}$ by
Proposition~\ref{prop:AH-elementary}(b), and
$\BstnormO{u} \le \Bstnorm{\eta u} + \BstnormO{(1-\eta)u}$ by
Lemma~\ref{lem:gluing}(b).
\end{proof}

\subsubsection*{Extraction of a limit}

\begin{lemma}[Limits of resolvent sequences]
\label{lem:resolvent-limit}
Let $s > 1/2$, let $\mathcal{K} \subset \Ereg$ be compact, and let
$z_k \to E \in \mathcal{K}$ with $\im z_k > 0$ for every $k$ and
$\im z_k \to 0$. Let $f_k \to f$ in $\Bspd{\Om}$, put $u_k = \Rpert{z_k}f_k$,
and assume
\begin{equation*}
  \sup_{k} \norm{u_k}_{L^{2,-s}(\Om)} < \infty .
\end{equation*}
Then there are a subsequence and a function
$u \in \Bstd{\Om} \cap H^2_{\mathrm{loc}}(\overline{\Om})$ such that, along it,
$u_k \to u$ in $L^{2,-s}(\Om)$, in $H^1(\Omin{\rho})$ for every
$\rho \ge \Rzero$, and in the weak-$*$ topology of $\Bstd{\Om}$. Moreover $u$
satisfies
\begin{equation}
  \label{eq:limit-problem}
  (-\Delta - E)u = f \ \text{ in } \Om,
  \qquad
  \gamma_1 u - \alpha\,\gamma_0 u = 0 \ \text{ on } \Gpert,
\end{equation}
together with
\begin{equation}
  \label{eq:outgoing-source}
  \eta u = \Ralpha{E+i0}\bigl( \eta f + [-\Delta,\eta]u \bigr)
\end{equation}
for every admissible cut-off $\eta$. The same holds, with $\Ralpha{E-i0}$ in
\eqref{eq:outgoing-source}, if instead $\im z_k < 0$ for every $k$.
\end{lemma}

\begin{proof}
Write $A = \sup_k \norm{u_k}_{L^{2,-s}(\Om)}$ and $B = \sup_k \BnormO{f_k}$,
both finite. By Lemma~\ref{lem:Bst-from-weighted},
$\BstnormO{u_k} \le C(B + A)$ for all large $k$, and by
Proposition~\ref{prop:local-regularity}, applied with $\rho' = \rho + 1$ and
$M = \sup_{\mathcal{K}}|\,\cdot\,| + 1$,
\begin{equation}
  \label{eq:uk-H2-local}
  \norm{u_k}_{H^2(\Omin{\rho})}
    \le C_\rho \bigl( \norm{f_k}_{\Ltwo{\Om}}
      + \norm{u_k}_{\Ltwo{\Omin{\rho+1}}} \bigr)
    \le C_\rho'
  \qquad (\rho \ge \Rzero).
\end{equation}

\emph{Extraction.} Fix $\rho \ge \Rzero$. By \eqref{eq:uk-H2-local} and the
Rellich--Kondrachov theorem on the bounded Lipschitz domain $\Omin{\rho}$, a
subsequence of $(u_k)$ converges in $H^1(\Omin{\rho})$. Iterate over
$\rho \in \Rzero + \mathbb{N}$ and extract diagonally. This produces a
subsequence, not relabelled, and a function $u$ with
$u_k \to u$ in $H^1(\Omin{\rho})$ for every $\rho$. By
\eqref{eq:uk-H2-local} and weak compactness, $u \in H^2(\Omin{\rho})$ for every
$\rho$. Convergence in $L^{2,-s}(\Om)$ follows because
\begin{equation*}
  \norm{u_k - u_l}_{L^{2,-s}(\Om)}
    \le \norm{u_k - u_l}_{L^{2,-s}(\Omin{\rho})}
      + 2 C_s \rho^{1/2-s} \sup_j \BstnormO{u_j}
\end{equation*}
by Lemma~\ref{lem:AH-tail}(a), the second term being small uniformly in $k,l$
for $\rho$ large and the first tending to $0$ for fixed $\rho$. Since
$\norm{u}_{\Ltwo{\shellp{j}}} = \lim_k \norm{u_k}_{\Ltwo{\shellp{j}}}$ for every
$j$, we get $\BstnormO{u} \le \liminf_k \BstnormO{u_k} < \infty$, so
$u \in \Bstd{\Om}$. Weak-$*$ convergence in $\Bstd{\Om}$ holds because it
holds against compactly supported elements of $\Bspd{\Om}$, by convergence in
$\Ltwo{\Omin{\rho}}$. Those elements are dense by Lemma~\ref{lem:AH-tail}(b),
and the family is bounded in $\Bstd{\Om}$.

\emph{The equation.} Each $u_k$ satisfies \eqref{eq:local-weak} on
$\Omin{\rho'}$ with $F = f_k$ and $z = z_k$. Letting $k \to \infty$, using
$u_k \to u$ in $H^1(\Omin{\rho'})$, $f_k \to f$ in $\Ltwo{\Om}$ and
$z_k \to E$, we find that $u$ satisfies \eqref{eq:local-weak} with $F = f$ and
$z = E$. Proposition~\ref{prop:local-regularity} then gives the boundary
condition of \eqref{eq:limit-problem} on $\Gpert \cap \Bfull{\rho}$ for every
$\rho$, and the equation holds in the distributional sense. Hence
\eqref{eq:limit-problem}.

\emph{The outgoing representation.} Let $\eta$ be admissible and put
$g_k = \eta f_k + [-\Delta,\eta]u_k$ and $g = \eta f + [-\Delta,\eta]u$. Then
$g_k \to g$ in $\Bspd{\Om}$, because
\begin{equation*}
  \Bnorm{\eta(f_k - f)} \le \BnormO{f_k - f},
  \qquad
  \Bnorm{[-\Delta,\eta](u_k - u)}
    \le C \norm{u_k - u}_{H^1(\Omin{R_2})} ,
\end{equation*}
the second by Lemma~\ref{lem:gluing}(a), and both right-hand sides tend to
$0$. Consequently, for every $h \in \Bsp$,
\begin{equation*}
  \bigl( \Ralpha{z_k}g_k, h \bigr)
    = \bigl( \Ralpha{z_k}(g_k - g), h \bigr) + \bigl( \Ralpha{z_k}g, h \bigr)
    \longrightarrow \bigl( \Ralpha{E+i0}g, h \bigr),
\end{equation*}
the first term by the uniform bound of \cite[Theorem~2.9(a)]{PartI} and the
second by the weak-$*$ continuity of \cite[Theorem~2.9(c)]{PartI}. On the other
hand $\Ralpha{z_k}g_k = \eta u_k$ by \eqref{eq:resolvent-gluing}, and
$\eta u_k \to \eta u$ in the weak-$*$ topology of $\Bst$ by the extraction step.
Two weak-$*$ limits coincide, which is \eqref{eq:outgoing-source}.

For $\im z_k < 0$ apply the above to $\overline{u_k} = \Rpert{\overline{z_k}}
\overline{f_k}$, the operator $\Hpert$ having real coefficients and a real
boundary condition, and conjugate the conclusion.
\end{proof}

\subsubsection*{The uniform bound}

\begin{proposition}[Uniform resolvent bound]
\label{prop:uniform-bound}
Let $n \ge 3$, $\alpha > 0$, let $\Om$ satisfy
Assumption~\ref{assump:geometry}, let $s > 1/2$, and let
$\mathcal{K} \subset \Ereg \setminus \Sexc$ be compact. There is
$C = C(n,\alpha,\Om,\mathcal{K},s)$ such that
\begin{equation}
  \label{eq:uniform-weighted}
  \norm{\Rpert{z}f}_{L^{2,-s}(\Om)} \le C\,\BnormO{f}
\end{equation}
for every $f \in \Bspd{\Om}$ and every $z \in \mathbb{C}\setminus\mathbb{R}$
with $\re z \in \mathcal{K}$ and $|\im z| \le 1$.
\end{proposition}

\begin{proof}
Suppose not. Then there are $z_k$ with $\re z_k \in \mathcal{K}$,
$0 < |\im z_k| \le 1$, and $f_k \in \Bspd{\Om}$ with
\begin{equation}
  \label{eq:contradiction-normalization}
  \BnormO{f_k} \to 0,
  \qquad
  \norm{u_k}_{L^{2,-s}(\Om)} = 1,
  \qquad u_k = \Rpert{z_k}f_k ,
\end{equation}
obtained by normalizing a sequence along which the ratio of the two sides of
\eqref{eq:uniform-weighted} is unbounded. Passing to a subsequence we may assume
that all $\im z_k$ have the same sign, say $\im z_k > 0$, and that
$z_k \to z_\infty$ with $\re z_\infty \in \mathcal{K}$ and
$\im z_\infty = \tau \in [0,1]$.

If $\tau > 0$, then $|\im z_k| \ge \tau/2$ for large $k$ and the spectral
theorem gives
$\norm{u_k}_{\Ltwo{\Om}} \le 2\tau^{-1}\norm{f_k}_{\Ltwo{\Om}}
\le 2\tau^{-1}\BnormO{f_k} \to 0$
by Proposition~\ref{prop:AH-elementary}(b). Since
$\norm{\cdot}_{L^{2,-s}(\Om)} \le \norm{\cdot}_{\Ltwo{\Om}}$, this contradicts
\eqref{eq:contradiction-normalization}. Hence $\tau = 0$ and
$z_k \to E := \re z_\infty \in \mathcal{K}$.

Lemma~\ref{lem:resolvent-limit} applies with $f = 0$ and produces, along a
subsequence, a function $u \in \Bstd{\Om} \cap H^2_{\mathrm{loc}}(\overline{\Om})$
with $u_k \to u$ in $L^{2,-s}(\Om)$, satisfying \eqref{eq:global-problem} and,
by \eqref{eq:outgoing-source} with $f = 0$, outgoing at energy $E$ in the sense
of Definition~\ref{def:outgoing-pert}. By
Theorem~\ref{thm:outgoing-uniqueness}, $u = 0$ if $E > 0$, while if
$E \in \Eband$ then $u$ is an eigenfunction of $\Hpert$ at $E$, and $u = 0$
because $E \notin \Sexc$. Either way $u = 0$, contradicting
$\norm{u}_{L^{2,-s}(\Om)} = \lim_k \norm{u_k}_{L^{2,-s}(\Om)} = 1$.
\end{proof}

\subsubsection*{The theorem}

\begin{theorem}[Limiting absorption principle for $\Hpert$]
\label{thm:lap-pert}
Let $n \ge 3$, $\alpha > 0$, let $\Om$ satisfy
Assumption~\ref{assump:geometry}, and let $\Sexc$ be the discrete set
\eqref{eq:exceptional-set}.
\begin{enumerate}
  \item[\textup{(a)}] \textup{(Uniform bound.)} For every compact
  $\mathcal{K} \subset \Ereg \setminus \Sexc$ there is
  $C = C(n,\alpha,\Om,\mathcal{K})$ with
  \begin{equation}
    \label{eq:lap-pert-bound}
    \BstnormO{\Rpert{z}f} \le C\,\BnormO{f},
    \qquad f \in \Bspd{\Om},
  \end{equation}
  for all $z \in \mathbb{C}\setminus\mathbb{R}$ with $\re z \in \mathcal{K}$.
  \item[\textup{(b)}] \textup{(Boundary values.)} For every
  $E \in \Ereg \setminus \Sexc$ and every $f \in \Bspd{\Om}$ the limits
  \begin{equation*}
    \Rpert{E \pm i0}f := \lim_{\eps \downarrow 0} \Rpert{E \pm i\eps}f
  \end{equation*}
  exist in the weak-$*$ topology of $\Bstd{\Om}$. The operators so defined
  belong to $\Lop{\Bspd{\Om}}{\Bstd{\Om}}$ and satisfy
  \eqref{eq:lap-pert-bound} with the same constant.
  \item[\textup{(c)}] \textup{(Continuity.)} For all $f,g \in \Bspd{\Om}$ the
  maps $E \mapsto (\Rpert{E \pm i0}f, g)$ are continuous on
  $\Ereg \setminus \Sexc$.
  \item[\textup{(d)}] \textup{(Adjoints and positivity.)}
  $(\Rpert{E \pm i0}f, g) = (f, \Rpert{E \mp i0}g)$ for
  $f,g \in \Bspd{\Om}$, and $\im(\Rpert{E+i0}f,f) \ge 0$.
  \item[\textup{(e)}] \textup{(The equation.)} For $f \in \Bspd{\Om}$ and
  $E \in \Ereg \setminus \Sexc$ the function $u = \Rpert{E \pm i0}f$ belongs to
  $\Bstd{\Om} \cap H^2_{\mathrm{loc}}(\overline{\Om})$, solves
  \begin{equation*}
    (-\Delta - E)u = f \ \text{ in } \Om,
    \qquad
    \gamma_1 u - \alpha\gamma_0 u = 0 \ \text{ on } \Gpert,
  \end{equation*}
  and satisfies \eqref{eq:outgoing-source}, with $\Ralpha{E \pm i0}$
  respectively, for every admissible cut-off.
\end{enumerate}
\end{theorem}

\begin{proof}
(a) For $|\im z| \ge 1$ the operator $\Rpert{z}$ has $\Ltwo{\Om}$ norm at most
$1$, so \eqref{eq:lap-pert-bound} holds with $C = 1$ by
Proposition~\ref{prop:AH-elementary}(c). For $|\im z| \le 1$ combine
Proposition~\ref{prop:uniform-bound}, with any fixed $s > 1/2$, and
Lemma~\ref{lem:Bst-from-weighted}.

(b) and (e). Fix $E \in \Ereg \setminus \Sexc$, take $\mathcal{K} = \{E\}$ and
$f \in \Bspd{\Om}$, and let $\eps_k \downarrow 0$. By (a) and
Proposition~\ref{prop:uniform-bound} the sequence
$u_k = \Rpert{E+i\eps_k}f$ is bounded in $L^{2,-s}(\Om)$, so
Lemma~\ref{lem:resolvent-limit}, applied with $z_k = E + i\eps_k$ and
$f_k = f$, produces a subsequence converging in the weak-$*$ topology of
$\Bstd{\Om}$ to some $u$ with the properties listed in (e). If $u'$ is the
limit along another subsequence, or along another sequence $\eps_k' \downarrow
0$, then $w = u - u'$ lies in $\Bstd{\Om} \cap
H^2_{\mathrm{loc}}(\overline{\Om})$, satisfies \eqref{eq:global-problem}, and
satisfies \eqref{eq:outgoing-source} with $f$ replaced by $0$. That is, $w$
is outgoing in the sense of Definition~\ref{def:outgoing-pert}. By
Theorem~\ref{thm:outgoing-uniqueness} and $E \notin \Sexc$ we get $w = 0$.
The limit is therefore independent of the sequence, so
$\Rpert{E+i0}f$ is well defined and has the properties in (e). Linearity in $f$
is clear, and \eqref{eq:lap-pert-bound} passes to the weak-$*$ limit. The case
$E - i0$ is identical.

(c) Let $E_j \to E$ in $\Ereg \setminus \Sexc$ and let $\mathcal{K}$ be a
compact subset of $\Ereg \setminus \Sexc$ containing $E$ and all $E_j$, which
exists because $\Sexc$ is closed and discrete in $\Eband$
(Proposition~\ref{prop:exceptional-discrete}). Fix $f,g \in \Bspd{\Om}$. By (b)
we may pick $\eps_j \downarrow 0$ with
\begin{equation*}
  \bigl| \bigl( \Rpert{E_j + i\eps_j}f, g \bigr)
       - \bigl( \Rpert{E_j + i0}f, g \bigr) \bigr| < 1/j .
\end{equation*}
Put $z_j = E_j + i\eps_j \to E$. By (a) the sequence $\Rpert{z_j}f$ is bounded
in $L^{2,-s}(\Om)$, so every subsequence has, by
Lemma~\ref{lem:resolvent-limit}, a further subsequence converging weak-$*$ in
$\Bstd{\Om}$ to a function with the properties of (e) at energy $E$. By the
uniqueness just established, that function is $\Rpert{E+i0}f$. Hence
$(\Rpert{z_j}f,g) \to (\Rpert{E+i0}f,g)$, and the displayed inequality
transfers the convergence to $(\Rpert{E_j+i0}f,g)$.

(d) For $\im z \ne 0$ self-adjointness gives $\Rpert{z}^{*} = \Rpert{\bar z}$,
and
\begin{equation*}
  \im\bigl( \Rpert{z}f, f \bigr)
    = \im z \, \norm{\Rpert{\bar z}f}_{\Ltwo{\Om}}^2 ,
\end{equation*}
which is nonnegative for $\im z > 0$. Both statements pass to the limit
$\eps \downarrow 0$ by (b).
\end{proof}

\begin{corollary}[The outgoing solution]
\label{cor:outgoing-char}
Let $E \in \Ereg \setminus \Sexc$ and $f \in \Bspd{\Om}$. Then $u =
\Rpert{E+i0}f$ is the unique function in $\Bstd{\Om} \cap
H^2_{\mathrm{loc}}(\overline{\Om})$ that solves $(-\Delta - E)u = f$ in $\Om$
with $\gamma_1 u - \alpha\gamma_0 u = 0$ on $\Gpert$ and satisfies
\eqref{eq:outgoing-source} for some admissible cut-off. The corresponding
statement holds for $\Rpert{E-i0}f$ and $\Ralpha{E-i0}$.
\end{corollary}

\begin{proof}
Existence is Theorem~\ref{thm:lap-pert}(e). If $u$ and $u'$ both have the stated
properties then their difference is outgoing in the sense of
Definition~\ref{def:outgoing-pert} and solves \eqref{eq:global-problem}, hence
vanishes by Theorem~\ref{thm:outgoing-uniqueness}.
\end{proof}

Corollary~\ref{cor:outgoing-char} is the form in which the limiting absorption
principle will be used. It selects a solution by a condition imposed near
infinity alone, since \eqref{eq:outgoing-source} sees $u$ only through
$\eta u$ and $[-\Delta,\eta]u$. It is the condition itself, rather than
the resolvent, that is intrinsic. Recasting it as a radiation condition stated
directly on $u$, without reference to $\Ralpha{E \pm i0}$, is the subject of
\S\ref{sec:scattering-solution}. On the band this requires a formulation with a
single channel, the bulk channel being closed there, and Part~I's radiation
condition has no such version.

\section{Generalized Eigenfunctions and the Scattering Solution}
\label{sec:scattering-solution}
The limiting absorption principle of \S\ref{sec:lap-pert} selects a solution of
the inhomogeneous problem by the requirement \eqref{eq:outgoing-source}, and
that requirement mentions the free resolvent $\Ralpha{E \pm i0}$. As a
definition of what it means for a wave to be outgoing this is unsatisfactory on
two counts: it refers to an auxiliary operator attached to a different domain,
and it is not manifestly a statement about the behaviour of the solution near
infinity. The first subsection below removes both defects, by showing that
\eqref{eq:outgoing-source} is equivalent to a radiation condition imposed on the
far field of the solution itself. The remaining two subsections use it to
construct the scattering solution: a global solution of the perturbed problem
which, far from the perturbation, is a prescribed solution of the unperturbed
problem plus an outgoing wave.

Two features of the perturbed setting shape the formulation throughout. The
first is that an incident wave is a function on $\Rnp$ and not on $\Om$. The
two domains are not nested either way, since the perturbed boundary may dip
below the hyperplane $\Gamma$. An incident field can therefore enter the
perturbed problem only after being cut off. The resulting commutator is a
compactly supported source, and it is the source that drives the scattering. The second is that on
the surface band $\Eband = (-\alpha^2,0)$ the bulk channel is closed. There the
radiation condition involves a single wavenumber, incident waves are surface
waves only, and the exceptional set $\Sexc$ of \eqref{eq:exceptional-set} has to
be excluded.

\subsection{The Radiation Condition on the Perturbed Domain}
\label{subsec:radiation-pert}

Throughout this subsection $\Rzero < R_1 < R_2$ and the admissible cut-offs are
those of \eqref{eq:cutoff-radii}--\eqref{eq:cutoff-admissible}.

We first recall the two notions from Part~I on which the condition rests. For
$u, v \in \Bst$ one writes $u \simeq v$ when $u - v \in \Bstz$, that is when
\begin{equation}
  \label{eq:simeq-pert}
  \lim_{R \to \infty} \frac{1}{R}
    \int_{\Rnp \cap \Bfull{R}} |u(x) - v(x)|^2\,dx = 0,
\end{equation}
and one names the four elementary waves that can occur in the far field of a
two-channel solution: writing $x = (\xpar,x_n)$, $r = |x|$, $\theta = x/r$,
$\rho = |\xpar|$ and $\omega = \xpar/\rho$, one sets, for $A \in \Ltwo{\Sm}$ and
$C \in \Ltwo{\Ssurf}$,
\begin{equation}
  \label{eq:model-waves-pert}
  \Wb{\pm}{A}(x) = A(\theta)\,r^{-\frac{n-1}{2}}\,e^{\pm i\kbulk r},
  \qquad
  \Ws{\pm}{C}(x)
    = C(\omega)\,\rho^{-\frac{n-2}{2}}\,e^{\pm i\ksurf\rho}\,\psizero(x_n) .
\end{equation}
Both belong to $\Bst$, and the four of them are linearly independent modulo
$\Bstz$: if a linear combination of $\Wb{+}{A}$, $\Wb{-}{B}$, $\Ws{+}{C}$ and
$\Ws{-}{D}$ lies in $\Bstz$, then $A = B = C = D = 0$
\cite[Definition 2.36, equation (2.173) and Lemma 2.89]{PartI}. The exponents in
\eqref{eq:model-waves-pert} are the critical ones for $\Bst$ in the respective
dimensions $n$ and $n-1$, which is why a single space accommodates both
channels.

\begin{definition}[The radiation condition on $\Om$]
\label{def:radiation-pert}
Let $E \in \Ereg$ and $v \in \Bstd{\Om}$, and let $\eta$ be an admissible
cut-off. We say that $v$ \emph{satisfies the outgoing radiation condition} at
energy $E$ if there are $A \in \Ltwo{\Sm}$ and $C \in \Ltwo{\Ssurf}$ with
\begin{equation}
  \label{eq:radiation-out}
  \eta v \;\simeq\; \Wb{+}{A} + \Ws{+}{C}
  \qquad (E > 0),
  \qquad\qquad
  \eta v \;\simeq\; \Ws{+}{C}
  \qquad (E \in \Eband),
\end{equation}
the extension of $\eta v$ by zero being understood on the left. The
\emph{incoming} radiation condition is \eqref{eq:radiation-out} with
$\Wb{-}{A}$ and $\Ws{-}{C}$ in place of $\Wb{+}{A}$ and $\Ws{+}{C}$.
\end{definition}

\begin{lemma}[Elementary properties]
\label{lem:radiation-elementary}
Let $E \in \Ereg$.
\begin{enumerate}
  \item[\textup{(a)}] For $v \in \Bstd{\Om}$ and any admissible cut-off $\eta$
  the zero extension of $\eta v$ lies in $\Bst$, with
  $\Bstnorm{\eta v} \le \BstnormO{v}$.
  \item[\textup{(b)}] Condition \eqref{eq:radiation-out} does not depend on the
  choice of the admissible cut-off, and the amplitudes $A$ and $C$ occurring in
  it are uniquely determined by $v$.
  \item[\textup{(c)}] The solutions of the outgoing radiation condition form a
  linear subspace of $\Bstd{\Om}$ containing $\Bstzd{\Om}$, and on
  $\Bstzd{\Om}$ the amplitudes vanish.
\end{enumerate}
\end{lemma}

\begin{proof}
(a) is the computation in the proof of Lemma~\ref{lem:gluing}(b), which uses
only $0 \le \eta \le 1$, the vanishing of $\eta$ on $\Bfull{R_1}$, and
\eqref{eq:exterior-coincidence}. No equation and no boundary condition enter.

(b) If $\eta_1,\eta_2$ are admissible then $(\eta_1 - \eta_2)v$ is supported in
$\{R_1 \le |x| \le R_2\}$ and lies in $\Ltwo{\Rnp}$, hence in $\Bstz$. Therefore
$\eta_1 v \simeq \eta_2 v$. Uniqueness of $A$ and $C$ is the linear
independence modulo $\Bstz$ recalled above.

(c) Linearity is immediate from \eqref{eq:simeq-pert}. If $v \in \Bstzd{\Om}$
then $\eta v \in \Bstz$ by Lemma~\ref{lem:gluing}(b), so
\eqref{eq:radiation-out} holds with $A = 0$ and $C = 0$.
\end{proof}

Part~(b) is what makes Definition~\ref{def:radiation-pert} a condition on $v$
near infinity rather than on the pair $(v,\eta)$. Part~(c) will let a
compactly supported modification of the incident field be
absorbed without comment in \S\ref{subsec:scattering-sol}.

\begin{theorem}[The radiation condition and the outgoing resolvent]
\label{thm:radiation-equivalence}
Let $n \ge 3$, $\alpha > 0$ and $E \in \Ereg$. Let $f \in \Bspd{\Om}$ have
compact support and let $u \in \Bstd{\Om} \cap H^2_{\mathrm{loc}}(\overline{\Om})$
satisfy
\begin{equation}
  \label{eq:radiation-problem}
  (-\Delta - E)u = f \ \text{ in } \Om,
  \qquad
  \gamma_1 u - \alpha\,\gamma_0 u = 0 \ \text{ on } \Gpert .
\end{equation}
Then the following are equivalent.
\begin{enumerate}
  \item[\textup{(a)}] Identity \eqref{eq:outgoing-source} holds for
  \emph{some} admissible cut-off $\eta$.
  \item[\textup{(b)}] Identity \eqref{eq:outgoing-source} holds for
  \emph{every} admissible cut-off $\eta$.
  \item[\textup{(c)}] $u$ satisfies the outgoing radiation condition of
  Definition~\ref{def:radiation-pert}.
\end{enumerate}
The corresponding statement holds with $\Ralpha{E-i0}$ in
\eqref{eq:outgoing-source} and the incoming radiation condition.
\end{theorem}

\begin{proof}
Fix an admissible cut-off $\eta$ and put
\begin{equation}
  \label{eq:transplanted-source}
  h_\eta = \eta f + [-\Delta,\eta]\,u .
\end{equation}
We first record that $h_\eta \in \Bsp$ and has compact support. Indeed
$[-\Delta,\eta]u \in \Bsp$ is supported in $\{R_1 \le |x| \le R_2\}$ by
Lemma~\ref{lem:gluing}(a). Moreover $\eta f$ vanishes on $\Omin{\Rzero}$, so its
extension by zero satisfies $\Bnorm{\eta f} = \BnormO{\eta f} \le \BnormO{f}$
by Lemma~\ref{lem:AH-exterior} and $0 \le \eta \le 1$, while its support is
contained in that of $f$. Moreover, by Lemma~\ref{lem:gluing} the zero
extension of $\eta u$ lies in $\Bst \cap H^2_{\mathrm{loc}}(\overline{\Rnp})$
and satisfies
\begin{equation}
  \label{eq:transplanted-equation}
  (-\Delta - E)(\eta u) = h_\eta \ \text{ in } \Rnp,
  \qquad
  \gamma_1(\eta u) - \alpha\,\gamma_0(\eta u) = 0 \ \text{ on } \Gamma .
\end{equation}

\emph{(b) $\Rightarrow$ (a)} is trivial.

\emph{(a) $\Rightarrow$ (c).} Let $\eta$ be a cut-off for which
\eqref{eq:outgoing-source} holds, so that $\eta u = \Ralpha{E+i0}h_\eta$. If
$E > 0$, then \cite[Theorem 2.49]{PartI} gives
$\Ralpha{E+i0}h_\eta \simeq \Wb{+}{A} + \Ws{+}{C}$ with amplitudes
proportional to $\Fbulk^{(+)}(E)h_\eta$ and $\Fsurf(E)h_\eta$. If
$E \in \Eband$, the same conclusion in its single-channel form
$\Ralpha{E+i0}h_\eta \simeq \Ws{+}{C}$ is
\cite[Proposition 2.64]{PartI}. Either way \eqref{eq:radiation-out} holds.

\emph{(c) $\Rightarrow$ (b).} Let $\eta$ be an arbitrary admissible cut-off and
put $u_0 = \Ralpha{E+i0}h_\eta$, which is defined because $h_\eta \in \Bsp$,
and set $w = \eta u - u_0$. By \eqref{eq:transplanted-equation} and
\cite[Theorem 2.9(e)]{PartI} both $\eta u$ and $u_0$ solve the same
inhomogeneous problem on $\Rnp$, so
\begin{equation*}
  w \in \Bst \cap H^2_{\mathrm{loc}}(\overline{\Rnp}),
  \qquad
  (-\Delta - E)w = 0 \ \text{ in } \Rnp,
  \qquad
  \gamma_1 w - \alpha\gamma_0 w = 0 \ \text{ on } \Gamma ,
\end{equation*}
that is, $w \in \Nsol{E}$ in the notation of \cite[Theorem 2.68]{PartI}. Now
$\eta u$ satisfies \eqref{eq:radiation-out} by hypothesis and
Lemma~\ref{lem:radiation-elementary}(b), while $u_0$ satisfies it by the
implication already proved. Subtracting, and using that $\simeq$ is compatible
with differences, we get
\begin{equation}
  \label{eq:w-purely-outgoing}
  w \;\simeq\; \Wb{+}{A'} + \Ws{+}{C'}
\end{equation}
for some $A' \in \Ltwo{\Sm}$ and $C' \in \Ltwo{\Ssurf}$, the first summand
being absent on the band. On the other hand, by
\cite[Corollary 2.70]{PartI} every element of $\Nsol{E}$ has uniquely
determined incoming and outgoing far-field data
$\yin, \yout \in \Echan{E}$, its far field being the superposition of the four
waves \eqref{eq:model-waves-pert} with amplitudes that are nonzero multiples of
$\yout$ and $\yin$, and the map $w \mapsto \yin$ is a bijection of $\Nsol{E}$
onto $\Echan{E}$. Comparing that expansion with
\eqref{eq:w-purely-outgoing} and invoking the linear independence modulo
$\Bstz$ recalled before Definition~\ref{def:radiation-pert}, we find
$\yin = 0$, whence $w = 0$. Thus $\eta u = \Ralpha{E+i0}h_\eta$, which is
\eqref{eq:outgoing-source} for the cut-off $\eta$. And $\eta$ was arbitrary.

For the incoming case, the differential expression, the boundary condition and
the coefficients of $[-\Delta,\eta]$ are real, so
$\overline{\Ralpha{E-i0}g} = \Ralpha{E+i0}\overline{g}$ for $g \in \Bsp$ and
$\overline{\Wb{-}{A}} = \Wb{+}{\overline A}$,
$\overline{\Ws{-}{C}} = \Ws{+}{\overline C}$. The statement follows by applying
the case already treated to $\overline u$ and $\overline f$.
\end{proof}

We record two consequences. The first closes a gap deliberately left
open in Definition~\ref{def:outgoing-pert}, which required
\eqref{eq:outgoing-pert} for one admissible cut-off only.

\begin{corollary}[Intrinsic form of the results of \S\ref{sec:lap-pert}]
\label{cor:outgoing-intrinsic}
Let $n \ge 3$, $\alpha > 0$ and $E \in \Ereg$.
\begin{enumerate}
  \item[\textup{(a)}] A function $u \in \Bstd{\Om} \cap
  H^2_{\mathrm{loc}}(\overline{\Om})$ satisfying \eqref{eq:global-problem} is
  outgoing in the sense of Definition~\ref{def:outgoing-pert} if and only if it
  satisfies the outgoing radiation condition of
  Definition~\ref{def:radiation-pert}. In that case
  \eqref{eq:outgoing-pert} holds for every admissible cut-off.
  \item[\textup{(b)}] Let in addition $E \notin \Sexc$ and let
  $f \in \Bspd{\Om}$ have compact support. Then $\Rpert{E+i0}f$ is the unique
  $u \in \Bstd{\Om} \cap H^2_{\mathrm{loc}}(\overline{\Om})$ that solves
  \eqref{eq:radiation-problem} and satisfies the outgoing radiation condition.
\end{enumerate}
\end{corollary}

\begin{proof}
(a) is Theorem~\ref{thm:radiation-equivalence} with $f = 0$, for which
\eqref{eq:outgoing-source} reduces to \eqref{eq:outgoing-pert}. For (b),
existence is Theorem~\ref{thm:lap-pert}(e) together with
Theorem~\ref{thm:radiation-equivalence}. If $u$ and $u'$ both have the stated
properties, their difference lies in $\Bstd{\Om}$, satisfies
\eqref{eq:global-problem}, and satisfies the outgoing radiation condition by
Lemma~\ref{lem:radiation-elementary}(c), hence is outgoing in the sense of
Definition~\ref{def:outgoing-pert} by (a) and vanishes by
Theorem~\ref{thm:outgoing-uniqueness} and $E \notin \Sexc$.
\end{proof}

\begin{remark}[The Sommerfeld form of the condition]
\label{rem:sommerfeld}
For $E > 0$ Part~I also formulates the radiation condition in the classical
manner, through the radial derivative. Following
\cite{DuranMugaNedelec2009}, one fixes an exponent $\alpha_1 \in (0,1)$,
unrelated to the impedance $\alpha$, splits the shell $\{R < |x| < 2R\}$ by
the interface $x_n = r^{\alpha_1}$ into
\begin{equation}
  \label{eq:sommerfeld-regions}
  \Sigma_1(R) = \{R<|x|<2R\} \cap \{x_n > r^{\alpha_1}\},
  \qquad
  \Sigma_2(R) = \{R<|x|<2R\} \cap \{0 \le x_n \le r^{\alpha_1}\},
\end{equation}
and requires
\begin{equation}
  \label{eq:sommerfeld}
  \lim_{R \to \infty} \frac1R \int_{\Sigma_1(R)}
    \bigl| \partial_r u - i\,\kbulk\,u \bigr|^2 dx = 0,
  \qquad
  \lim_{R \to \infty} \frac1R \int_{\Sigma_2(R)}
    \bigl| \partial_r u - i\,\ksurf\,u \bigr|^2 dx = 0 .
\end{equation}
The interface separates the two channels, the surface wave being confined to a
layer of thickness $O(1/\alpha)$ along $\Gamma$. For $R > \Rzero$ the regions
\eqref{eq:sommerfeld-regions} lie in $\Omex{\Rzero}$, where $\Om$ and $\Rnp$
coincide, so \eqref{eq:sommerfeld} is imposed on $u$ itself without a cut-off.
There $\partial_r u$ is classically defined, by the elliptic regularity up to
the flat boundary recalled in \cite[Remark 2.73]{PartI}.

We do not adopt \eqref{eq:sommerfeld} as the definition, for two reasons
inherited from Part~I. First, it has no counterpart on the band, where there is
no second wavenumber, whereas Definition~\ref{def:radiation-pert} merely drops
one summand. Secondly, Part~I verifies \eqref{eq:sommerfeld} for the outgoing
resolvent only for sources in $C_c^\infty(\Rnp)$
\cite[Proposition 2.85]{PartI}, because the proof uses a Hankel expansion and a
Paley--Wiener argument, and a condition on $\partial_r u$ is not stable under
$\Bst$-limits \cite[Remark 2.86]{PartI}. The source
\eqref{eq:transplanted-source} is compactly supported, but its support meets
$\Gamma$, since the transplantation annulus must meet the boundary along which
the surface channel propagates. We therefore do not assert the equivalence of
\eqref{eq:sommerfeld} with Definition~\ref{def:radiation-pert} for the
solutions constructed below. By \cite[Theorem 2.90]{PartI} the two agree for
the solutions covered by the verification of Part~I.
\end{remark}

\subsection{Incident Waves}
\label{subsec:incident}

The incident field is a solution of the \emph{unperturbed} problem, and the
class of such solutions of $\Bst$ growth was determined in Part~I. For
$E \in \Ereg$ put
\begin{equation}
  \label{eq:Nsol-recall}
  \Nsol{E} = \Bigl\{ \phi \in \Bst \cap H^2_{\mathrm{loc}}(\overline{\Rnp}) :
    (-\Delta - E)\phi = 0 \text{ in } \Rnp,\
    \gamma_1\phi - \alpha\gamma_0\phi = 0 \text{ on } \Gamma \Bigr\} .
\end{equation}
By \cite[Theorem 2.68 and Corollary 2.70]{PartI}, $\Nsol{E}$ is the range of
$\Fzero(E)^{*}$. Every $\phi \in \Nsol{E}$ possesses uniquely determined
incoming and outgoing far-field data $\yin[\phi], \yout[\phi] \in \Echan{E}$.
The map $\phi \mapsto \yin[\phi]$ is a bijection of $\Nsol{E}$ onto
$\Echan{E}$. The norms are comparable,
\begin{equation}
  \label{eq:incident-norms}
  c_1\,\norm{\yin[\phi]}_{\Echan{E}}
  \le \Bstnorm{\phi}
  \le c_2\,\norm{\yin[\phi]}_{\Echan{E}} ,
\end{equation}
with constants depending only on $n$, $\alpha$ and $E$. The channel space is
$\Echan{E} = \Ltwo{\Sm} \oplus \Ltwo{\Ssurf}$ for $E > 0$ and
$\Echan{E} = \Ltwo{\Ssurf}$ for $E \in \Eband$.

\begin{definition}[Incident wave]
\label{def:incident}
Let $E \in \Ereg$. An \emph{incident wave} at energy $E$ is an element
$\phi \in \Nsol{E}$. It is prescribed by its incoming data
$\yin[\phi] \in \Echan{E}$.
\end{definition}

Two comments on the shape of this definition. It is a superposition over the
energy shell, and not an individual generalized eigenfunction: the functions
$\Psib(\cdot,\kbulk\theta)$ and $\Psis(\cdot,\ksurf\omega)$ of
\cite[Definition 2.23]{PartI} are bounded but do not lie in $\Bst$, their shell
averages growing like $r^{n-1}$ and $r^{n-2}$ respectively. Since the whole
construction takes place in $\Bsp$--$\Bst$, the incident field has to lie in
$\Bst$ as well, and $\Nsol{E}$ is exactly the class of unperturbed solutions
that does. Individual directions reappear later, as kernels rather than as
functions. Secondly, on the band the only incident waves are surface waves,
$\Echan{E}$ having lost its bulk summand. This reflects the fact that below the
bulk threshold nothing arrives from, or escapes to, the interior of the
half-space.

Transplanting an incident wave into $\Om$ requires a cut-off, because $\phi$ is a function on $\Rnp$ while the perturbed problem is posed
on $\Om$, and neither domain contains the other. The next lemma is the
counterpart of Lemma~\ref{lem:gluing}, read in the opposite direction.

\begin{lemma}[Transplantation of an incident wave]
\label{lem:incident-transplant}
Let $E \in \Ereg$, let $\phi \in \Nsol{E}$ and let $\eta$ be an admissible
cut-off. Then the extension of $\eta\phi$ by zero belongs to
$\Bstd{\Om} \cap H^2_{\mathrm{loc}}(\overline{\Om})$, satisfies
$\gamma_1(\eta\phi) - \alpha\gamma_0(\eta\phi) = 0$ on $\Gpert$, and
\begin{equation}
  \label{eq:incident-source}
  (-\Delta - E)(\eta\phi) = [-\Delta,\eta]\,\phi
  \qquad \text{in } \Om .
\end{equation}
Moreover $[-\Delta,\eta]\phi \in \Bspd{\Om}$ is supported in
$\{R_1 \le |x| \le R_2\}$ and
\begin{equation}
  \label{eq:incident-source-bound}
  \BnormO{[-\Delta,\eta]\phi} \le C_\eta\,\Bstnorm{\phi},
\end{equation}
with $C_\eta$ depending only on $n$, $E$, $\alpha$ and $\eta$.
\end{lemma}

\begin{proof}
By \eqref{eq:cutoff-admissible} the function $\eta\phi$ vanishes on
$\Bfull{R_1}$, and $\Rzero < R_1$. Since $\Omex{\Rzero}$ and
$\Rnp \setminus \overline{\Bfull{\Rzero}}$ coincide by
\eqref{eq:exterior-coincidence}, the extension by zero is a well-defined
element of $H^2_{\mathrm{loc}}(\overline{\Om})$, and it lies in $\Bstd{\Om}$
with $\BstnormO{\eta\phi} = \Bstnorm{\eta\phi} \le \Bstnorm{\phi}$ by
Lemma~\ref{lem:AH-exterior}. On $\Gpert \cap \{|x| > \Rzero\}$ the boundary is
the flat $\Gamma$ and the argument of Lemma~\ref{lem:gluing} applies verbatim:
writing $\eta(x) = \tilde\eta(|x|)$ gives
$\partial_{x_n}\eta = \tilde\eta'(|x|)x_n/|x|$, which vanishes on
$\{x_n = 0\}$, so that
$\gamma_1(\eta\phi) = \eta\,\gamma_1\phi = \alpha\,\gamma_0(\eta\phi)$. On the
remaining part of $\Gpert$, which lies in $\Bfull{\Rzero}$, the function
$\eta\phi$ vanishes identically near the boundary. Identity
\eqref{eq:incident-source} is the Leibniz rule together with
$(-\Delta - E)\phi = 0$.

For \eqref{eq:incident-source-bound}, the coefficients of $[-\Delta,\eta]$ are
supported in $\{R_1 \le |x| \le R_2\}$, so by \eqref{eq:commutator}
\begin{equation*}
  \norm{[-\Delta,\eta]\phi}_{\Ltwo{\Om}}
  \le C \sup_{R_1 \le |x| \le R_2} \bigl( |\phi(x)| + |\nabla\phi(x)| \bigr)
      \,\bigl| \{ R_1 \le |x| \le R_2 \} \bigr|^{1/2} .
\end{equation*}
The supremum is bounded by $C\,R_2^{1/2}\,\Bstnorm{\phi}$. Indeed the pointwise
growth bound
\begin{equation*}
  |\phi(x)| + |\nabla\phi(x)| \le C\,\Bstnorm{\phi}\,|x|^{1/2}
\end{equation*}
of \cite[Lemma 2.75]{PartI} requires only that $\phi$ solve the Robin problem
outside some radius $\varrho$, and then holds on $\{|x| > \varrho + 2\}$. Here
$\phi$ solves it outside every positive radius, so the bound is available on
$\{|x| > 3\}$. The annulus lies there, because $R_1 > \Rzero \ge 4$, the
normalization \eqref{eq:R0-dyadic} allowing $\Rzero$ to be enlarged at will.
Passing from the $\Ltwo{}$ norm to the $\Bspd{\Om}$ norm on a set contained in
$\Bfull{R_2}$ is the computation in the proof of Lemma~\ref{lem:gluing}(a).
\end{proof}

\subsection{The Scattering Solution}
\label{subsec:scattering-sol}

We can now state the main result of this section. The total field is the
transplanted incident wave plus a correction which is outgoing. The correction
is not small, and no smallness of the perturbation is assumed anywhere.

\begin{theorem}[Existence and uniqueness of the scattering solution]
\label{thm:scattering-solution}
Let $n \ge 3$, $\alpha > 0$, let $\Om$ satisfy
Assumption~\ref{assump:geometry}, let $E \in \Ereg \setminus \Sexc$ and let
$\phi \in \Nsol{E}$. Then there is exactly one function
$u \in H^2_{\mathrm{loc}}(\overline{\Om})$ such that
\begin{enumerate}
  \item[\textup{(i)}] $(-\Delta - E)u = 0$ in $\Om$ and
  $\gamma_1 u - \alpha\gamma_0 u = 0$ on $\Gpert$.
  \item[\textup{(ii)}] For some admissible cut-off $\eta$ the difference
  $u - \eta\phi$ belongs to $\Bstd{\Om}$ and satisfies the outgoing radiation
  condition of Definition~\ref{def:radiation-pert}.
\end{enumerate}
It is given by
\begin{equation}
  \label{eq:scattering-solution}
  u = \eta\phi - \Rpert{E+i0}\bigl( [-\Delta,\eta]\phi \bigr),
\end{equation}
the right-hand side being independent of the admissible cut-off $\eta$. The
condition in \textup{(ii)} then holds for every admissible cut-off, and the
scattered field $v = u - \eta\phi$ obeys
\begin{equation}
  \label{eq:scattered-bound}
  \BstnormO{v} \le C\,\Bstnorm{\phi},
  \qquad C = C(n,\alpha,\Om,E) .
\end{equation}
The corresponding statements hold with $\Rpert{E-i0}$ and the incoming
radiation condition.
\end{theorem}

\begin{proof}
\emph{Existence.} Fix an admissible cut-off $\eta$ and put
$g = [-\Delta,\eta]\phi$, which by Lemma~\ref{lem:incident-transplant} lies in
$\Bspd{\Om}$ and is supported in $\{R_1 \le |x| \le R_2\}$. Since
$E \in \Ereg \setminus \Sexc$, Theorem~\ref{thm:lap-pert}(e) provides
$v := -\Rpert{E+i0}g \in \Bstd{\Om} \cap H^2_{\mathrm{loc}}(\overline{\Om})$
with $(-\Delta - E)v = -g$ in $\Om$, $\gamma_1 v - \alpha\gamma_0 v = 0$ on
$\Gpert$, and \eqref{eq:outgoing-source} for every admissible cut-off. By
Theorem~\ref{thm:radiation-equivalence}, applied with the compactly supported
source $-g$, the function $v$ satisfies the outgoing radiation condition. Put
$u = \eta\phi + v$. Then $u \in H^2_{\mathrm{loc}}(\overline{\Om})$, and adding
\eqref{eq:incident-source} to the equation for $v$ gives $(-\Delta - E)u = 0$
in $\Om$, while both summands satisfy the Robin condition on $\Gpert$. Thus
$u$ satisfies (i) and (ii), and it is of the form
\eqref{eq:scattering-solution}.

\emph{Uniqueness.} Let $u_1,u_2$ both satisfy (i) and (ii), with admissible
cut-offs $\eta_1,\eta_2$, and put $w = u_1 - u_2$. Then
\begin{equation*}
  w = \bigl( u_1 - \eta_1\phi \bigr) - \bigl( u_2 - \eta_2\phi \bigr)
      + (\eta_1 - \eta_2)\phi .
\end{equation*}
The first two brackets lie in $\Bstd{\Om}$ and satisfy the outgoing radiation
condition by hypothesis. The last term is supported in
$\{R_1 \le |x| \le R_2\}$ and lies in $\Ltwo{\Om}$, hence in $\Bstzd{\Om}$ by
Proposition~\ref{prop:AH-elementary}(b), and therefore also satisfies the
outgoing radiation condition, with vanishing amplitudes. By
Lemma~\ref{lem:radiation-elementary}(c) the sum $w$ does too, and
$w \in \Bstd{\Om} \cap H^2_{\mathrm{loc}}(\overline{\Om})$ satisfies
\eqref{eq:global-problem} by (i). Corollary~\ref{cor:outgoing-intrinsic}(a)
makes $w$ outgoing in the sense of Definition~\ref{def:outgoing-pert}, and
Theorem~\ref{thm:outgoing-uniqueness} together with $E \notin \Sexc$ gives
$w = 0$.

\emph{The remaining assertions.} Independence of $\eta$ in
\eqref{eq:scattering-solution} follows from the uniqueness just proved, since
for each admissible $\eta$ the right-hand side satisfies (i) and (ii). Once
$u$ is known to be independent of $\eta$, so is the validity of (ii) for a
given cut-off, by Lemma~\ref{lem:radiation-elementary}(b) and (c). Finally
\eqref{eq:scattered-bound} follows from \eqref{eq:lap-pert-bound} with
$\mathcal{K} = \{E\}$ and \eqref{eq:incident-source-bound}. The incoming case
is obtained by replacing $\Rpert{E+i0}$ with $\Rpert{E-i0}$ throughout, or by
conjugation as in the proof of Theorem~\ref{thm:radiation-equivalence}.
\end{proof}

\begin{definition}[The perturbed generalized eigenfunctions]
\label{def:Phipert}
Let $E \in \Ereg \setminus \Sexc$. For $\phi \in \Nsol{E}$ let
$\Phipert{\pm}{E}\phi$ denote the solution furnished by
Theorem~\ref{thm:scattering-solution} for the outgoing, respectively incoming,
radiation condition.
\end{definition}

\begin{corollary}[Parametrization by the incoming data]
\label{cor:scattering-parametrization}
Let $E \in \Ereg \setminus \Sexc$. For every $g \in \Echan{E}$ there is
exactly one solution $u$ of the perturbed problem \textup{(i)} of
Theorem~\ref{thm:scattering-solution} whose incident wave has incoming data
$g$, namely $u = \Phipert{+}{E}\phi$ with $\phi \in \Nsol{E}$ the unique
element with $\yin[\phi] = g$. The map $g \mapsto u - \eta\phi$ is linear and
bounded from $\Echan{E}$ into $\Bstd{\Om}$, and $\Phipert{+}{E}$ is injective
on $\Nsol{E}$.
\end{corollary}

\begin{proof}
The first assertion combines the bijectivity of $\phi \mapsto \yin[\phi]$
recalled after \eqref{eq:Nsol-recall} with
Theorem~\ref{thm:scattering-solution}. Boundedness is
\eqref{eq:scattered-bound} together with the right-hand inequality in
\eqref{eq:incident-norms}. For injectivity, suppose $\Phipert{+}{E}\phi = 0$.
Then the scattered field is $v = -\eta\phi$, so $\eta\phi$ satisfies the
outgoing radiation condition by Theorem~\ref{thm:scattering-solution}, and
$\eta\phi \simeq \phi$ because $(\eta - 1)\phi$ is compactly supported and lies
in $\Ltwo{\Rnp} \subset \Bstz$. Hence $\phi \simeq \Wb{+}{A} + \Ws{+}{C}$ for
some amplitudes $A$, $C$. But $\phi \in \Nsol{E}$, so by
\cite[Corollary 2.70]{PartI} its far field is the superposition of the four
waves \eqref{eq:model-waves-pert} whose incoming amplitudes are nonzero
multiples of $\yin[\phi]$. Comparing the two expansions and using the linear
independence modulo $\Bstz$ gives $\yin[\phi] = 0$, whence $\phi = 0$.
\end{proof}

The scattered field $v = \Phipert{+}{E}\phi - \eta\phi$ therefore has, by
Definition~\ref{def:radiation-pert}, a well-defined pair of outgoing amplitudes
$(A,C) \in \Ltwo{\Sm} \oplus \Ltwo{\Ssurf}$, reducing to $C$ alone on the band.
These amplitudes depend linearly and boundedly on the incoming data of the
incident wave, and the resulting map $\Echan{E} \to \Echan{E}$ is the scattering
amplitude of the perturbed problem. Its identification with a generalized
Fourier transform of $\Hpert$, and with it the block structure of the
scattering matrix, is the subject of \S\ref{sec:gft-pert} and
\S\ref{sec:smatrix}.

\begin{remark}[The exceptional set]
\label{rem:exceptional-cost}
Theorem~\ref{thm:scattering-solution} excludes $E \in \Sexc$, and the exclusion
cannot be removed. At an eigenvalue of $\Hpert$ embedded in the band the
homogeneous problem has a nontrivial solution in $\Bstzd{\Om}$, by
Corollary~\ref{cor:global-uniqueness}(b). It satisfies the outgoing radiation
condition with vanishing amplitudes by
Lemma~\ref{lem:radiation-elementary}(c), so uniqueness fails there. The same
eigenfunctions have vanishing far-field data, and at $E \in \Sexc$ they also
obstruct the characterization of Theorem~\ref{thm:characterization-pert}. For
$E > 0$ no such energies occur, by Corollary~\ref{cor:global-uniqueness}(a),
and on the band they form a discrete set by
Proposition~\ref{prop:exceptional-discrete}. The source of this asymmetry is
the one described in Remark~\ref{rem:band-closed-waveguide}.
\end{remark}

\section{The Generalized Fourier Transform and the Far Field}
\label{sec:gft-pert}
Section~\ref{sec:scattering-solution} attached to each incident wave a unique
scattering solution, and thereby a map $\Phipert{\pm}{E}$ defined on the space
$\Nsol{E}$ of unperturbed solutions of $\Bst$ growth. That map is the perturbed
analogue of the eigenoperator $\Fzero(E)^{*}$ of Part~I, and the present
section draws the two consequences that follow from it: a spectral
representation of $\Hpert$ on $\Ereg \setminus \Sexc$, and a formula
identifying the far field of the outgoing resolvent.

The two are the same statement seen from opposite ends, and the object that
joins them is the transform defined in \S\ref{subsec:gft-pert-def}. Its
usefulness rests on a single algebraic identity, proved in
\S\ref{subsec:transplant}: the cut-off used to transplant functions between
$\Rnp$ and $\Om$ enters through a commutator which is \emph{antisymmetric} with
respect to the pairing, so that moving it from one factor to the other costs
only a sign. Everything in this section, including the completeness of the
expansion, is obtained by stationary means. No wave operator is constructed
anywhere in this paper, and no time-dependent argument is used at all.

\subsection{The Perturbed Generalized Fourier Transform}
\label{subsec:gft-pert-def}

We first recall what Part~I supplies. For $E \in \Ereg$ the channel space is
$\Echan{E} = \Ltwo{\Sm} \oplus \Ltwo{\Ssurf}$ for $E > 0$ and
$\Echan{E} = \Ltwo{\Ssurf}$ for $E \in \Eband$, and the generalized Fourier
transform of the unperturbed operator is a bounded operator
$\Fzero(E) : \Bsp \to \Echan{E}$ whose adjoint
$\Fzero(E)^{*} \in \Lop{\Echan{E}}{\Bst}$ is an eigenoperator: for every
$g \in \Echan{E}$ the function $\Fzero(E)^{*}g$ solves the unperturbed problem
on $\Rnp$ \cite[Definition 2.29 and Proposition 2.34]{PartI}. By
\cite[Theorem 2.68]{PartI} the map
\begin{equation}
  \label{eq:Fzero-bijection}
  \Fzero(E)^{*} : \Echan{E} \longrightarrow \Nsol{E}
\end{equation}
is a bijection, bicontinuous with
\begin{equation}
  \label{eq:Fzero-two-sided}
  \frac{1}{\sqrt{8\pi\ksurf}}\,\norm{g}_{\Echan{E}}
  \;\le\; \Bstnorm{\Fzero(E)^{*}g}
  \;\le\; \norm{\Fzero(E)^{*}}\,\norm{g}_{\Echan{E}} .
\end{equation}
In particular $\Fzero(E)^{*}$ is injective, and its range is exactly the space
$\Nsol{E}$ of \eqref{eq:Nsol-recall} on which $\Phipert{\pm}{E}$ was defined in
Definition~\ref{def:Phipert}. The two maps may therefore be composed.

\begin{definition}[The perturbed generalized Fourier transform]
\label{def:Fpert}
Let $E \in \Ereg \setminus \Sexc$. Define
\begin{equation}
  \label{eq:Fpert-adjoint-def}
  \Fpert{\pm}{E}^{*} \;:=\; \Phipert{\pm}{E} \circ \Fzero(E)^{*}
  \;:\; \Echan{E} \longrightarrow \Bstd{\Om} ,
\end{equation}
and let $\Fpert{\pm}{E} : \Bspd{\Om} \to \Echan{E}$ be determined by
\begin{equation}
  \label{eq:Fpert-def}
  \bigl( \Fpert{\pm}{E}f, g \bigr)_{\Echan{E}}
  = \bigl( f, \Fpert{\pm}{E}^{*}g \bigr),
  \qquad f \in \Bspd{\Om},\ g \in \Echan{E},
\end{equation}
the right-hand side being the pairing of $\Bspd{\Om}$ with $\Bstd{\Om}$ of
Proposition~\ref{prop:AH-elementary}(a). We write
$\Fpert{\pm}{E} = \bigl( \Fbulkpert{\pm}{E}, \Fsurfpert{\pm}{E} \bigr)$ for
$E > 0$ and $\Fpert{\pm}{E} = \Fsurfpert{\pm}{E}$ for $E \in \Eband$,
according to the splitting of $\Echan{E}$.
\end{definition}

The definition is stated on the adjoint side because that is the side on which
the objects of \S\ref{sec:scattering-solution} live: $\Phipert{\pm}{E}$ produces
solutions, not coefficients. The next lemma records that
\eqref{eq:Fpert-def} does determine an operator, and collects what follows at
once.

\begin{lemma}[Elementary properties]
\label{lem:Fpert-elementary}
Let $E \in \Ereg \setminus \Sexc$.
\begin{enumerate}
  \item[\textup{(a)}] $\Fpert{\pm}{E}^{*} \in \Lop{\Echan{E}}{\Bstd{\Om}}$, and
  \eqref{eq:Fpert-def} defines a unique
  $\Fpert{\pm}{E} \in \Lop{\Bspd{\Om}}{\Echan{E}}$ with
  $\norm{\Fpert{\pm}{E}} = \norm{\Fpert{\pm}{E}^{*}}$.
  \item[\textup{(b)}] \textup{(Eigenoperator property.)} For every
  $g \in \Echan{E}$ the function $u = \Fpert{\pm}{E}^{*}g$ lies in
  $\Bstd{\Om} \cap H^2_{\mathrm{loc}}(\overline{\Om})$ and satisfies
  \eqref{eq:global-problem}. Moreover $u - \eta\,\Fzero(E)^{*}g$ satisfies the
  outgoing, respectively incoming, radiation condition of
  Definition~\ref{def:radiation-pert} for every admissible cut-off $\eta$.
  \item[\textup{(c)}] $\Fpert{\pm}{E}^{*}$ is injective.
  \item[\textup{(d)}] $\Fpert{\pm}{E}\bigl( C_c^{\infty}(\Om) \bigr)$ is dense
  in $\Echan{E}$.
\end{enumerate}
\end{lemma}

\begin{proof}
(a) Boundedness of the composition \eqref{eq:Fpert-adjoint-def} follows from
that of $\Fzero(E)^{*}$, recalled in \eqref{eq:Fzero-two-sided}, and from
Theorem~\ref{thm:scattering-solution}: for $\phi \in \Nsol{E}$ and an
admissible cut-off $\eta$ one has
$\Phipert{\pm}{E}\phi = \eta\phi + v$ with
$\BstnormO{v} \le C\Bstnorm{\phi}$ by \eqref{eq:scattered-bound}, while
$\BstnormO{\eta\phi} \le \Bstnorm{\phi}$ by
Lemma~\ref{lem:incident-transplant}. Let now $f \in \Bspd{\Om}$. The map
$g \mapsto \overline{( f, \Fpert{\pm}{E}^{*}g )}$ is a bounded linear
functional on $\Echan{E}$, so the Riesz representation theorem produces a
unique $\Fpert{\pm}{E}f \in \Echan{E}$ satisfying \eqref{eq:Fpert-def}. The
resulting map is linear and bounded. The equality of norms is the usual
duality computation, the pairing of Proposition~\ref{prop:AH-elementary}(a)
having constant $1$.

(b) Put $\phi = \Fzero(E)^{*}g \in \Nsol{E}$. Then $u = \Phipert{\pm}{E}\phi$
is the solution furnished by Theorem~\ref{thm:scattering-solution}, so it has
the stated regularity and solves \eqref{eq:global-problem}. The last assertion
is part (ii) of that theorem together with the statement, proved there, that it
then holds for every admissible cut-off.

(c) If $\Fpert{\pm}{E}^{*}g = 0$ then $\Phipert{\pm}{E}\phi = 0$ with
$\phi = \Fzero(E)^{*}g$, whence $\phi = 0$ by the injectivity of
$\Phipert{\pm}{E}$ (Corollary~\ref{cor:scattering-parametrization}), and then
$g = 0$ by \eqref{eq:Fzero-two-sided}.

(d) Let $g \in \Echan{E}$ be orthogonal to
$\Fpert{\pm}{E}\bigl( C_c^{\infty}(\Om) \bigr)$. By \eqref{eq:Fpert-def} the
function $\Fpert{\pm}{E}^{*}g$, which lies in
$\Bstd{\Om} \subset L^2_{\mathrm{loc}}(\Om)$, is orthogonal to every
$f \in C_c^{\infty}(\Om)$ and hence vanishes. Now apply (c).
\end{proof}

\begin{remark}[The crossing of the signs]
\label{rem:Fpert-crossing}
Definition~\ref{def:Fpert} attaches $\Fpert{+}{E}$ to the \emph{outgoing}
scattering solutions. That is the transparent convention, but it is not the one
that reads off the far field of the outgoing resolvent. That far field is given
by $\Fpert{-}{E}$, as Theorem~\ref{thm:Fpert-explicit} below makes explicit
and \S\ref{subsec:far-field} exploits. The same crossing occurs in potential
scattering, where
$\Phi^{(+)} = \phi_0 - R(E+i0)V\phi_0$ and the amplitude of $R(E+i0)f$ is the
pairing of $f$ against $\Phi^{(-)}$. Its source is the single algebraic fact
that the adjoint of $\Rpert{E+i0}$ is $\Rpert{E-i0}$
(Theorem~\ref{thm:lap-pert}(d)). On the unperturbed half-space the two
transforms coincide, both reducing to $\Fzero(E)$, so no discrepancy with
Part~I arises. The distinction is created by the perturbation, and it is what
the scattering matrix of \S\ref{sec:smatrix} measures.
\end{remark}

\subsection{Transplantation Identities}
\label{subsec:transplant}

Throughout this subsection $\eta$ is an admissible cut-off in the sense of
\eqref{eq:cutoff-radii}--\eqref{eq:cutoff-admissible}, so that $\eta$ is radial,
$0 \le \eta \le 1$, $\eta = 0$ on $\Bfull{R_1}$ and $\eta = 1$ outside
$\Bfull{R_2}$, with $\Rzero < R_1 < R_2$. The coefficients of the commutator
$[-\Delta,\eta] = -(\Delta\eta) - 2\nabla\eta \cdot \nabla$ are then supported
in the compact annulus
\begin{equation}
  \label{eq:annulus}
  \mathcal{A} = \{ R_1 \le |x| \le R_2 \} \cap \Rnp
              = \{ R_1 \le |x| \le R_2 \} \cap \Om ,
\end{equation}
the two descriptions agreeing by \eqref{eq:exterior-coincidence}. Pairings of
functions supported in $\mathcal{A}$ are therefore unambiguous, and we write
them without indicating the domain.

The following identity is the mechanism of the whole section. It uses neither
the equation nor the boundary condition, and the two functions to which it is
applied need not be defined on the same domain.

\begin{lemma}[Antisymmetry of the transplantation commutator]
\label{lem:commutator-antisym}
Let $u \in H^2_{\mathrm{loc}}(\overline{\Om})$ and
$\psi \in H^2_{\mathrm{loc}}(\overline{\Rnp})$. Then
\begin{equation}
  \label{eq:commutator-antisym}
  \bigl( [-\Delta,\eta]u, \psi \bigr)
  = - \bigl( u, [-\Delta,\eta]\psi \bigr) .
\end{equation}
\end{lemma}

\begin{proof}
Fix $R_0'$ and $R$ with $\Rzero < R_0' < R_1$ and $R > R_2$, and put
$D = \{ R_0' < |x| < R \} \cap \Rnp$. By \eqref{eq:exterior-coincidence} we have
$\overline{D} \subset \overline{\Om}$ as well, so both $u$ and $\psi$ belong to
$H^2(D)$, and $\mathcal{A} \subset D$ carries the supports of both integrands
in \eqref{eq:commutator-antisym}. It therefore suffices to prove the identity
with both pairings taken over $D$.

Write $[-\Delta,\eta]u = -\Delta(\eta u) + \eta\,\Delta u$ and treat the two
summands separately. Green's second identity on $D$, applied to the pair
$(\eta u, \psi)$ and to the pair $(u, \eta\psi)$, gives
\begin{align*}
  \int_D \bigl( -\Delta(\eta u) \bigr)\overline{\psi}\,dx
  &= - \int_D \eta u\,\overline{\Delta\psi}\,dx
    - \int_{\partial D}
        \Bigl( \partial_\nu(\eta u)\,\overline{\psi}
             - \eta u\,\partial_\nu\overline{\psi} \Bigr) dS,
  \\
  \int_D \eta\,\Delta u\,\overline{\psi}\,dx
  &= \int_D u\,\overline{\Delta(\eta\psi)}\,dx
    + \int_{\partial D}
        \Bigl( \partial_\nu u\,\overline{\eta\psi}
             - u\,\partial_\nu\overline{\eta\psi} \Bigr) dS,
\end{align*}
where $\nu$ is the outward unit normal of $D$ and $\eta$ is real. Adding, the
volume terms combine into
\begin{equation*}
  \int_D u\,\overline{\bigl( \Delta(\eta\psi) - \eta\,\Delta\psi \bigr)}\,dx
  = - \int_D u\,\overline{[-\Delta,\eta]\psi}\,dx ,
\end{equation*}
which is the right-hand side of \eqref{eq:commutator-antisym}. It remains to
see that the boundary terms cancel. The boundary $\partial D$ consists of the
two spherical caps $\{|x| = R_0'\} \cap \Rnp$ and $\{|x| = R\} \cap \Rnp$ and
of the flat piece $\Gamma \cap \{ R_0' < |x| < R \}$.

On the inner cap $\eta$ vanishes identically in a neighbourhood, because
$R_0' < R_1$, hence so do $\eta u$, $\eta\psi$ and all their derivatives, and
both boundary integrands vanish. On the outer cap $\eta \equiv 1$ and
$\nabla\eta = 0$, because $R > R_2$. The two integrands then reduce to
$\partial_\nu u\,\overline{\psi} - u\,\partial_\nu\overline{\psi}$ and its
negative, and cancel. On $\Gamma$ one has $\nu = -e_n$, and $\eta$ being radial
gives $\partial_{x_n}\eta(x) = \tilde\eta'(|x|)\,x_n/|x| = 0$ at $x_n = 0$.
Consequently $\partial_\nu(\eta u) = \eta\,\partial_\nu u$ and
$\partial_\nu(\eta\psi) = \eta\,\partial_\nu\psi$ there, and the two integrands
again cancel.
\end{proof}

We can now identify $\Fpert{\pm}{E}$ in terms of the unperturbed transform. The
source appearing on the right is the one produced by transplantation in
\eqref{eq:transplanted-source}.

\begin{theorem}[The transplanted source formula]
\label{thm:Fpert-explicit}
Let $E \in \Ereg \setminus \Sexc$, let $f \in \Bspd{\Om}$ and let $\eta$ be an
admissible cut-off. Then
\begin{equation}
  \label{eq:Fpert-explicit}
  \Fpert{\mp}{E}f
  = \Fzero(E)\Bigl( \eta f + [-\Delta,\eta]\,\Rpert{E \pm i0}f \Bigr),
\end{equation}
the argument of $\Fzero(E)$ being understood as extended by zero to $\Rnp$. In
particular the right-hand side does not depend on the choice of $\eta$.
\end{theorem}

\begin{proof}
We treat the upper sign. The lower one is obtained by exchanging $E+i0$ and
$E-i0$ throughout. Put $u = \Rpert{E+i0}f$ and
\begin{equation}
  \label{eq:heta-pert}
  h_\eta = \eta f + [-\Delta,\eta]u .
\end{equation}
That $h_\eta$ belongs to $\Bsp$ and is compactly supported was checked in the
proof of Theorem~\ref{thm:radiation-equivalence}: the zero extension of
$\eta f$ satisfies $\Bnorm{\eta f} \le \BnormO{f}$ by
Lemma~\ref{lem:AH-exterior}, and $[-\Delta,\eta]u \in \Bsp$ is supported in
$\mathcal{A}$ by Lemma~\ref{lem:gluing}(a).

Let $g \in \Echan{E}$ and put $\psi = \Fzero(E)^{*}g \in \Nsol{E}$. By
Theorem~\ref{thm:scattering-solution}, in its incoming form,
\begin{equation*}
  \Fpert{-}{E}^{*}g = \Phipert{-}{E}\psi
  = \eta\psi - \Rpert{E-i0}\bigl( [-\Delta,\eta]\psi \bigr),
\end{equation*}
and therefore, by \eqref{eq:Fpert-def},
\begin{equation}
  \label{eq:Fpert-expand}
  \bigl( \Fpert{-}{E}f, g \bigr)_{\Echan{E}}
  = \bigl( f, \eta\psi \bigr)
    - \bigl( f, \Rpert{E-i0}([-\Delta,\eta]\psi) \bigr) .
\end{equation}
The first term equals $(\eta f, \psi)$, the function $\eta f$ being supported
in $\Omex{R_1} = \Rnp \setminus \overline{\Bfull{R_1}}$ where the two domains
agree. For the second, Theorem~\ref{thm:lap-pert}(d) gives
$\bigl( f, \Rpert{E-i0}k \bigr) = \bigl( \Rpert{E+i0}f, k \bigr)$ for
$k \in \Bspd{\Om}$, so that the second term equals
$-\bigl( u, [-\Delta,\eta]\psi \bigr)$, which by
Lemma~\ref{lem:commutator-antisym} is $\bigl( [-\Delta,\eta]u, \psi \bigr)$.
Both terms of \eqref{eq:Fpert-expand} are thus paired against the same $\psi$,
and
\begin{equation*}
  \bigl( \Fpert{-}{E}f, g \bigr)_{\Echan{E}}
  = \bigl( h_\eta, \psi \bigr)
  = \bigl( h_\eta, \Fzero(E)^{*}g \bigr)
  = \bigl( \Fzero(E)h_\eta, g \bigr)_{\Echan{E}} .
\end{equation*}
Since $g \in \Echan{E}$ was arbitrary, \eqref{eq:Fpert-explicit} follows. The
left-hand side involves no cut-off, so neither does the right.
\end{proof}

Formula \eqref{eq:Fpert-explicit} is what makes the perturbed transform
computable: it expresses it through the unperturbed transform of a compactly
supported source, and every property of $\Fzero(E)$ established in Part~I is
thereby available. The first application is the regularity in the spectral
parameter that the completeness argument of \S\ref{subsec:stone-pert} will
require. The difficulty there is that
Theorem~\ref{thm:lap-pert}(c) gives only weak-$*$ continuity of
$E \mapsto \Rpert{E+i0}f$, which is not enough to control the commutator. The
strong local convergence has to be recovered, and
Lemma~\ref{lem:resolvent-limit} provides it.

\begin{corollary}[Continuity in the spectral parameter]
\label{cor:Fpert-continuity}
Let $f \in \Bspd{\Om}$. Then $E \mapsto \Fpert{\pm}{E}f$ is continuous from
$\Ereg \setminus \Sexc$ into the channel spaces, in the sense that it is
continuous on each of $\Eband \setminus \Sexc$ and $(0,\infty)$ with values in
$\Ltwo{\Ssurf}$ and in $\Ltwo{\Sm} \oplus \Ltwo{\Ssurf}$ respectively.
Moreover $\norm{\Fpert{\pm}{E}}$ is bounded uniformly on compact subsets of
$\Ereg \setminus \Sexc$.
\end{corollary}

\begin{proof}
Again we treat $\Fpertg{-}$. Let $E_j \to E$ within $\Ereg \setminus \Sexc$,
all lying in one of the two intervals, and fix an admissible cut-off $\eta$.
The set $\Sexc$ is closed and discrete in $\Eband$
(Proposition~\ref{prop:exceptional-discrete}), so some compact
$\mathcal{K} \subset \Ereg \setminus \Sexc$ contains $E$ together with $E_j$
for every large $j$. Discard the finitely many remaining terms.

Fix $j$. Lemma~\ref{lem:resolvent-limit}, applied to an arbitrary sequence
$\eps \downarrow 0$, produces a subsequence converging in $H^1(\Omin{R_2})$ to
a function which the uniqueness established in the proof of
Theorem~\ref{thm:lap-pert}(b) identifies as $\Rpert{E_j+i0}f$. Hence
$\Rpert{E_j + i\eps}f \to \Rpert{E_j+i0}f$ in $H^1(\Omin{R_2})$ as
$\eps \downarrow 0$, and we may pick $\eps_j \in (0,1/j)$ with
\begin{equation*}
  \norm{\Rpert{E_j + i\eps_j}f - \Rpert{E_j+i0}f}_{H^1(\Omin{R_2})}
  < 1/j .
\end{equation*}

Put $z_j = E_j + i\eps_j \to E$. By Theorem~\ref{thm:lap-pert}(a) and
Proposition~\ref{prop:uniform-bound} the sequence $\Rpert{z_j}f$ is bounded in
$L^{2,-s}(\Om)$ for a fixed $s > 1/2$, so every subsequence has, by
Lemma~\ref{lem:resolvent-limit}, a further subsequence converging in
$H^1(\Omin{R_2})$ to a function satisfying \eqref{eq:limit-problem} and
\eqref{eq:outgoing-source} at energy $E$. That function is $\Rpert{E+i0}f$ by
Theorem~\ref{thm:lap-pert}(b) and Corollary~\ref{cor:outgoing-char}. Hence the
full sequence converges, and with the displayed inequality,
\begin{equation*}
  \Rpert{E_j+i0}f \longrightarrow \Rpert{E+i0}f
  \qquad \text{in } H^1(\Omin{R_2}) .
\end{equation*}
The commutator $[-\Delta,\eta]$ is a first-order operator with bounded
coefficients supported in $\mathcal{A} \subset \Omin{R_2}$, so
$[-\Delta,\eta]\Rpert{E_j+i0}f \to [-\Delta,\eta]\Rpert{E+i0}f$ in
$\Ltwo{\Rnp}$, and hence in $\Bsp$, all the functions being supported in the
fixed compact set $\mathcal{A}$. Write $h_\eta(E')$ for the transplanted source
\eqref{eq:heta-pert} built at energy $E'$. Consequently
$h_\eta(E_j) \to h_\eta(E)$ in $\Bsp$. Now
\begin{equation*}
  \Fzero(E_j)h_\eta(E_j) - \Fzero(E)h_\eta(E)
  = \Fzero(E_j)\bigl( h_\eta(E_j) - h_\eta(E) \bigr)
    + \bigl( \Fzero(E_j) - \Fzero(E) \bigr)h_\eta(E),
\end{equation*}
and both terms tend to $0$: the first because $\norm{\Fzero(E_j)}$ is bounded
on $\mathcal{K}$ and the second by the continuity of
$E \mapsto \Fzero(E)h$ for fixed $h \in \Bsp$, both being part of
\cite[Proposition 2.33]{PartI}. By \eqref{eq:Fpert-explicit} this is the
assertion. The uniform bound on $\norm{\Fpert{-}{E}}$ follows from the same
formula, from \cite[Proposition 2.33]{PartI} and from
\eqref{eq:lap-pert-bound}.
\end{proof}

The last identity of this subsection transplants the energy flux. In Part~I the
imaginary part of the resolvent at a boundary value is the squared norm of the
scattering data \cite[equation (2.76)]{PartI}. For $\Hpert$ the corresponding
quantity is equal to the one computed on $\Rnp$ for the transplanted problem.
The reason is geometric: the flux is carried by a sphere of large radius, and out there
the cut-off is identically $1$ and the two domains coincide.

\begin{lemma}[Transplantation of the flux]
\label{lem:flux-transplant}
Let $E \in \Ereg \setminus \Sexc$, let $f \in \Bspd{\Om}$, put
$u = \Rpert{E+i0}f$ and let $h_\eta$ be as in \eqref{eq:heta-pert}. Then
\begin{equation}
  \label{eq:flux-transplant}
  \im\bigl( \Rpert{E+i0}f, f \bigr)_{\Ltwo{\Om}}
  = \im\bigl( \Ralpha{E+i0}h_\eta, h_\eta \bigr)_{\Ltwo{\Rnp}} ,
\end{equation}
both pairings being those of $\Bst$ with $\Bsp$. The same holds with $E - i0$
throughout.
\end{lemma}

\begin{proof}
Let $R > R_2$. Green's second identity on $\Omin{R}$, applied to $u$ and
$\overline{u}$ and combined with $(-\Delta - E)u = f$ and $E \in \mathbb{R}$,
gives
\begin{equation*}
  2i \im \int_{\Omin{R}} u\,\overline{f}\,dx
  = \int_{\partial\Omin{R}}
      \bigl( \partial_\nu u\,\overline{u} - u\,\partial_\nu\overline{u}
      \bigr) dS .
\end{equation*}
On $\Gpert \cap \Bfull{R}$ the Robin condition
$\gamma_1 u = \alpha\gamma_0 u$ and the reality of $\alpha$ make the integrand
$\alpha|\gamma_0 u|^2 - \alpha|\gamma_0 u|^2 = 0$. Since $R > R_2 > \Rzero$,
the remaining part of $\partial\Omin{R}$ is the cap
$\Sigma_R := \{|x| = R\} \cap \Rnp$ by \eqref{eq:exterior-coincidence}, and
therefore
\begin{equation}
  \label{eq:flux-pert-sphere}
  \im \int_{\Omin{R}} u\,\overline{f}\,dx
  = \im \int_{\Sigma_R} \partial_r u\,\overline{u}\,dS .
\end{equation}
Now repeat the computation on $\Rnp$ for the transplanted solution. By
Lemma~\ref{lem:gluing} the zero extension of $\eta u$ lies in
$\Bst \cap H^2_{\mathrm{loc}}(\overline{\Rnp})$ and satisfies
$(-\Delta - E)(\eta u) = h_\eta$ in $\Rnp$ together with the Robin condition on
$\Gamma$. The same two steps give
\begin{equation}
  \label{eq:flux-flat-sphere}
  \im \int_{\Rnp \cap \Bfull{R}} \eta u\,\overline{h_\eta}\,dx
  = \im \int_{\Sigma_R} \partial_r(\eta u)\,\overline{\eta u}\,dS
  = \im \int_{\Sigma_R} \partial_r u\,\overline{u}\,dS ,
\end{equation}
the last equality because $\eta \equiv 1$ and $\nabla\eta = 0$ on $\Sigma_R$.
The right-hand sides of \eqref{eq:flux-pert-sphere} and
\eqref{eq:flux-flat-sphere} are literally the same, so the left-hand sides
agree for every $R > R_2$. Now let $R \to \infty$. The left-hand side of
\eqref{eq:flux-pert-sphere} converges to $\im( u, f )$, because
$u \in \Bstd{\Om}$ and $f \in \Bspd{\Om}$. That of
\eqref{eq:flux-flat-sphere} converges to $\im( \eta u, h_\eta )$, because
$\eta u \in \Bst$ and $h_\eta \in \Bsp$. Finally
$\eta u = \Ralpha{E+i0}h_\eta$ by Theorem~\ref{thm:lap-pert}(e), which is
\eqref{eq:flux-transplant}. For $E - i0$ the argument is unchanged.
\end{proof}

Lemma~\ref{lem:flux-transplant}, the flux identity of Part~I and
Theorem~\ref{thm:Fpert-explicit} now compose into a single chain, and that
chain is the quadratic form of Stone's formula for $\Hpert$. The next
subsection begins there.

\subsection{Stone's Formula and the Spectral Representation}
\label{subsec:stone-pert}

Throughout this subsection $P(\,\cdot\,)$ denotes the spectral measure of
$\Hpert$, and we write
\begin{equation}
  \label{eq:Adens-pert-def}
  \Adenspert{E} := \frac{1}{2\pi i}
    \bigl( \Rpert{E+i0} - \Rpert{E-i0} \bigr),
  \qquad E \in \Ereg \setminus \Sexc ,
\end{equation}
the boundary values being those of Theorem~\ref{thm:lap-pert}(b).

\begin{theorem}[Stone's formula for $\Hpert$]
\label{thm:stone-pert}
Let $E \in \Ereg \setminus \Sexc$. Then
\begin{equation}
  \label{eq:stone-pert}
  \Adenspert{E} = \Fpert{+}{E}^{*}\,\Fpert{+}{E}
                = \Fpert{-}{E}^{*}\,\Fpert{-}{E}
  \qquad \text{in } \Lop{\Bspd{\Om}}{\Bstd{\Om}} .
\end{equation}
In particular $\Adenspert{E} \ge 0$, and
\begin{equation}
  \label{eq:flux-pert}
  \im\bigl( \Rpert{E+i0}f, f \bigr)
  = \pi\,\norm{\Fpert{-}{E}f}_{\Echan{E}}^2
  = \pi\,\norm{\Fpert{+}{E}f}_{\Echan{E}}^2 ,
  \qquad f \in \Bspd{\Om} .
\end{equation}
\end{theorem}

\begin{proof}
Fix $f \in \Bspd{\Om}$ and an admissible cut-off $\eta$, and let $h_\eta$ be the
transplanted source \eqref{eq:heta-pert} built from $u = \Rpert{E+i0}f$. By
Lemma~\ref{lem:flux-transplant}, by the flux identity
\cite[equation (2.76)]{PartI} applied to $h_\eta \in \Bsp$, and by
Theorem~\ref{thm:Fpert-explicit},
\begin{equation}
  \label{eq:flux-chain}
  \im\bigl( \Rpert{E+i0}f, f \bigr)
  = \im\bigl( \Ralpha{E+i0}h_\eta, h_\eta \bigr)
  = \pi\,\norm{\Fzero(E)h_\eta}_{\Echan{E}}^2
  = \pi\,\norm{\Fpert{-}{E}f}_{\Echan{E}}^2 .
\end{equation}
Now run the same three steps at $E - i0$, with the transplanted source built
from $\Rpert{E-i0}f$. On the flat side $\Ralpha{E-i0}$ is the adjoint of
$\Ralpha{E+i0}$, so the conjugate of \cite[equation (2.76)]{PartI} reads
$\im( \Ralpha{E-i0}h, h ) = -\pi\norm{\Fzero(E)h}^2$. The chain therefore
gives $\im( \Rpert{E-i0}f, f ) = -\pi\norm{\Fpert{+}{E}f}^2$. On the other
hand Theorem~\ref{thm:lap-pert}(d) gives
$( f, \Rpert{E+i0}f ) = ( \Rpert{E-i0}f, f )$, so that
$\im( \Rpert{E-i0}f, f ) = -\im( \Rpert{E+i0}f, f )$. The two norms in
\eqref{eq:flux-pert} therefore agree, and \eqref{eq:flux-pert} is proved.

Finally, again by Theorem~\ref{thm:lap-pert}(d),
\begin{equation*}
  \bigl( \Adenspert{E}f, f \bigr)
  = \frac{1}{2\pi i}\Bigl[ \bigl( \Rpert{E+i0}f, f \bigr)
      - \overline{\bigl( \Rpert{E+i0}f, f \bigr)} \Bigr]
  = \frac{1}{\pi}\,\im\bigl( \Rpert{E+i0}f, f \bigr) ,
\end{equation*}
so \eqref{eq:flux-pert} says that the two sides of \eqref{eq:stone-pert} have
the same quadratic form on $\Bspd{\Om}$. Both are sesquilinear, so the
polarization identity upgrades this to equality of the forms, and hence of the
operators, the pairing of Proposition~\ref{prop:AH-elementary}(a) being
nondegenerate.
\end{proof}

Identity \eqref{eq:flux-pert} says that the total flux
radiated by a source $f$ is the same whether it is measured by
$\Fpert{-}{E}f$ or by $\Fpert{+}{E}f$. That symmetry is the germ of the
unitarity of the scattering matrix, and \S\ref{sec:smatrix} makes the
intertwining operator explicit.

\begin{proposition}[The spectral measure of $\Hpert$]
\label{prop:spectral-measure-pert}
Let $J = (a,b)$ be an interval with $\overline{J} \subset \Ereg \setminus \Sexc$.
Then for all $f,g \in \Bspd{\Om}$
\begin{equation}
  \label{eq:stone-measure-pert}
  \bigl( P(J)f, g \bigr)_{\Ltwo{\Om}}
  = \int_a^b \bigl( \Fpert{\pm}{E}f, \Fpert{\pm}{E}g \bigr)_{\Echan{E}}\,dE ,
\end{equation}
with either choice of sign.
\end{proposition}

\begin{proof}
Stone's formula for the spectral measure of the self-adjoint operator $\Hpert$
gives
\begin{equation*}
  \tfrac12\bigl( P(\overline{J}) + P(J) \bigr)
  = \operatorname*{s-lim}_{\eps \downarrow 0} \frac{1}{2\pi i}
    \int_a^b \bigl( \Rpert{E+i\eps} - \Rpert{E-i\eps} \bigr) dE .
\end{equation*}
The two projections on the left coincide, because neither $a$ nor $b$ is an
eigenvalue of $\Hpert$: for $a > 0$ this is
Corollary~\ref{cor:global-uniqueness}(a), and for $a \in \Eband$ it is the
definition \eqref{eq:exceptional-set} of $\Sexc$ together with
$a \notin \Sexc$. The same applies to $b$. By Theorem~\ref{thm:lap-pert}(b) the
integrand converges pointwise in $E$ to $2\pi i\,\Adenspert{E}$ in the weak
sense of the pairing, and by Theorem~\ref{thm:lap-pert}(a) it is dominated,
uniformly in $\eps \in (0,1]$, by the constant $2C\BnormO{f}\BnormO{g}$
attached to the compact set $\overline{J} \subset \Ereg \setminus \Sexc$.
Dominated convergence gives
$( P(J)f, g ) = \int_a^b ( \Adenspert{E}f, g )\,dE$, and
Theorem~\ref{thm:stone-pert} rewrites the integrand as
$( \Fpert{\pm}{E}f, \Fpert{\pm}{E}g )_{\Echan{E}}$.
\end{proof}

\begin{proposition}[The absolutely continuous subspace]
\label{prop:ac-subspace-pert}
Let $n \ge 3$, $\alpha > 0$ and let $\Om$ satisfy
Assumption~\ref{assump:geometry}. Then
\begin{equation}
  \label{eq:ac-subspace}
  \Hac(\Hpert) = \Ran P\bigl( \Ereg \setminus \Sexc \bigr),
\end{equation}
and $\ssc(\Hpert) = \varnothing$. Moreover
$\spp(\Hpert) \cap \Ereg = \Sexc$, and
$\sigma(\Hpert) \setminus \bigl( \Ereg \setminus \Sexc \bigr)$ is countable.
\end{proposition}

\begin{proof}
Write $\mathcal{G} = \Ereg \setminus \Sexc$ and
$\mathcal{Z} = \mathbb{R} \setminus \mathcal{G}$.

\emph{$\Ran P(\mathcal{G}) \subset \Hac(\Hpert)$.} Let $I$ be an interval with
$\overline{I} \subset \mathcal{G}$ and let $f \in \Bspd{\Om}$. By
\eqref{eq:stone-measure-pert} the measure $I' \mapsto \norm{P(I')f}^2$ has, on
$I$, the density $E \mapsto \norm{\Fpert{-}{E}f}^2$, and is therefore
absolutely continuous there. Hence $P(I)f \in \Hac(\Hpert)$. Since
$\Bspd{\Om} \supset C_c^{\infty}(\Om)$ is dense in $\Ltwo{\Om}$ and $P(I)$ is
bounded, $P(I)\Ltwo{\Om} \subset \Hac(\Hpert)$, the latter being closed.
Finally $\mathcal{G}$ is a countable disjoint union of such intervals, because
$\Sexc$ is closed and discrete in $\Eband$
(Proposition~\ref{prop:exceptional-discrete}) and
$\Ereg = \Eband \cup (0,\infty)$. Summing gives the inclusion.

\emph{$\mathcal{Z} \cap \sigma(\Hpert)$ is countable.} One has
$\mathcal{Z} = (-\infty,-\alpha^2] \cup \{0\} \cup \Sexc$. By
\eqref{eq:sess} the part of $\sigma(\Hpert)$ in $(-\infty,-\alpha^2)$ consists
of at most countably many eigenvalues, by
Theorem~\ref{thm:pert-sa-spec}(c). The set $\Sexc$ is countable as well, being
discrete in $\Eband$.

\emph{$\Hac(\Hpert) \subset \Ran P(\mathcal{G})$.} Let
$f \in \Hac(\Hpert)$. The measure $\norm{P(\,\cdot\,)f}^2$ is absolutely
continuous and supported in $\sigma(\Hpert)$, hence gives no mass to the
countable set $\mathcal{Z} \cap \sigma(\Hpert)$. Therefore
$P(\mathcal{Z})f = 0$ and $f = P(\mathcal{G})f$. This proves
\eqref{eq:ac-subspace}.

\emph{The remaining assertions.} By \eqref{eq:ac-subspace} the singular
subspace of $\Hpert$ is $\Ran P(\mathcal{Z})$, on which the spectral measure
lives on a countable set and is therefore purely atomic. Hence
$\ssc(\Hpert) = \varnothing$. Finally
$\spp(\Hpert) \cap \Ereg = \Sexc$ is Corollary~\ref{cor:global-uniqueness}: no
$E > 0$ is an eigenvalue by (a), and in $\Eband$ the eigenvalues are exactly
$\Sexc$ by \eqref{eq:exceptional-set}.
\end{proof}

We can now state the spectral representation. Let
\begin{equation}
  \label{eq:direct-integral}
  \Hdirect = \int_{\Ereg}^{\oplus} \Echan{E}\,dE
\end{equation}
be the direct integral of the channel spaces, exactly the space in which
Part~I's transform $\Fzero$ takes its values \cite[Theorem 2.30]{PartI}. Since
$\Sexc$ is countable it is a Lebesgue null set, so nothing is lost by leaving
$\Fpert{\pm}{E}$ undefined on it.

\begin{theorem}[Spectral representation and completeness]
\label{thm:spectral-rep-pert}
Let $n \ge 3$, $\alpha > 0$ and let $\Om$ satisfy
Assumption~\ref{assump:geometry}. Then, with either choice of sign, the map
$f \mapsto \bigl( E \mapsto \Fpert{\pm}{E}f \bigr)$, defined for
$f \in \Bspd{\Om}$, extends uniquely to a bounded operator
\begin{equation*}
  \Fpertg{\pm} : \Ltwo{\Om} \longrightarrow \Hdirect ,
\end{equation*}
and this operator is a partial isometry with initial space $\Hac(\Hpert)$ and
final space all of $\Hdirect$. It diagonalizes $\Hpert$: for every bounded
Borel function $\varphi$ on $\mathbb{R}$ and every $f \in \Ltwo{\Om}$,
\begin{equation}
  \label{eq:diagonalization}
  \bigl( \Fpertg{\pm}\,\varphi(\Hpert)f \bigr)(E)
  = \varphi(E)\,\bigl( \Fpertg{\pm}f \bigr)(E)
  \qquad \text{for a.e. } E \in \Ereg .
\end{equation}
In particular $\Fpertg{\pm}$ restricts to a unitary operator from
$\Hac(\Hpert)$ onto $\Hdirect$, and $\Hpert$ restricted to $\Hac(\Hpert)$ is
unitarily equivalent to $\Halpha$.
\end{theorem}

\begin{proof}
We argue for $\Fpertg{-}$. The other case is identical.

\emph{Step 1: the map is well defined and contractive.} Let
$f \in \Bspd{\Om}$. By Corollary~\ref{cor:Fpert-continuity} the section
$E \mapsto \Fpert{-}{E}f$ is continuous on
$\mathcal{G} = \Ereg \setminus \Sexc$, hence measurable on $\Ereg$, the set
$\Sexc$ being a Lebesgue null set. We denote this section by $\Fpertg{-}f$, so
that $(\Fpertg{-}f)(E) = \Fpert{-}{E}f$ for $E \in \mathcal{G}$. Exhausting
$\mathcal{G}$ by intervals $J$ with $\overline{J} \subset \mathcal{G}$ and
using \eqref{eq:stone-measure-pert} together with monotone convergence,
\begin{equation}
  \label{eq:global-parseval}
  \norm{\Fpertg{-}f}_{\Hdirect}^2
  = \int_{\Ereg} \norm{\Fpert{-}{E}f}_{\Echan{E}}^2\,dE
  = \norm{P(\mathcal{G})f}_{\Ltwo{\Om}}^2
  \le \norm{f}_{\Ltwo{\Om}}^2 .
\end{equation}
In particular $\Fpertg{-}f \in \Hdirect$. Since $\Bspd{\Om}$ contains
$C_c^{\infty}(\Om)$ and is therefore dense in $\Ltwo{\Om}$, the map
$f \mapsto \Fpertg{-}f$ extends uniquely to a contraction
$\Fpertg{-} : \Ltwo{\Om} \to \Hdirect$, and the identity
$\norm{\Fpertg{-}f}_{\Hdirect} = \norm{P(\mathcal{G})f}_{\Ltwo{\Om}}$ of
\eqref{eq:global-parseval} passes to the limit. By
Proposition~\ref{prop:ac-subspace-pert}, $P(\mathcal{G})$ is the orthogonal
projection onto $\Hac(\Hpert)$, so $\Fpertg{-}$ is isometric on
$\Hac(\Hpert)$ and vanishes on its orthogonal complement. It is therefore a
partial isometry with initial space $\Hac(\Hpert)$. Its range
$\mathcal{M} := \Ran\Fpertg{-} = \Fpertg{-}\Hac(\Hpert)$ is the isometric
image of a complete space, hence complete, and therefore closed in $\Hdirect$.

\emph{Step 2: diagonalization.} We first extend
\eqref{eq:stone-measure-pert} to all of $\Ltwo{\Om}$ and to all Borel sets.
Let $I$ be an open interval with $\overline{I} \subset \mathcal{G}$. For
$f,g \in \Bspd{\Om}$, \eqref{eq:stone-measure-pert} reads
\begin{equation*}
  \bigl( P(I)f, g \bigr)
  = \int_I \bigl( (\Fpertg{-}f)(E), (\Fpertg{-}g)(E) \bigr)_{\Echan{E}}\,dE .
\end{equation*}
Both sides are continuous sesquilinear forms on $\Ltwo{\Om}$, bounded by
$\norm{f}\norm{g}$. For the right-hand side this follows from the
Cauchy--Schwarz inequality, first in $\Echan{E}$ and then in $E$, and from
\eqref{eq:global-parseval}. Since $\Bspd{\Om}$ is dense in $\Ltwo{\Om}$, the
identity holds for all $f,g \in \Ltwo{\Om}$.

Now fix $f,g \in \Ltwo{\Om}$, and for Borel sets $B \subset \mathbb{R}$ put
\begin{equation*}
  \mu_1(B) = \bigl( P(B \cap \mathcal{G})f, g \bigr),
  \qquad
  \mu_2(B) = \int_{B \cap \Ereg}
    \bigl( (\Fpertg{-}f)(E), (\Fpertg{-}g)(E) \bigr)_{\Echan{E}}\,dE .
\end{equation*}
Both are complex measures, the density of $\mu_2$ being integrable by the
bound just used. They agree on every open interval $I \subset \mathbb{R}$.
Indeed, $I \cap \mathcal{G}$ is open, hence a countable disjoint union of open
intervals. Each of these is the union of an increasing sequence of open
intervals whose closures lie in $\mathcal{G}$, and on those the two measures
agree. Countable additivity transfers the agreement to $I \cap \mathcal{G}$,
and $\mu_2(I) = \mu_2(I \cap \mathcal{G})$ because $\Sexc$ is a Lebesgue null
set. The open intervals, $\mathbb{R}$ itself and the empty set included, form
a family closed under finite intersections which generates the Borel sets.
Dynkin's
$\pi$--$\lambda$ theorem, applied to the real and imaginary parts, therefore
gives $\mu_1 = \mu_2$. Thus, for all
$f,g \in \Ltwo{\Om}$ and every Borel set $B \subset \mathbb{R}$,
\begin{equation*}
  \tag{$\ast$}
  \label{eq:spec-measure-all}
  \bigl( P(B \cap \mathcal{G})f, g \bigr)
  = \int_{B \cap \Ereg}
    \bigl( (\Fpertg{-}f)(E), (\Fpertg{-}g)(E) \bigr)_{\Echan{E}}\,dE .
\end{equation*}
Let $\psi$ be a bounded Borel function on $\mathbb{R}$. By the functional
calculus, $( \psi(\Hpert)P(\mathcal{G})f, g ) = \int_{\mathbb{R}} \psi\,d\mu$,
where $\mu(B) = ( P(B)P(\mathcal{G})f, g ) = ( P(B \cap \mathcal{G})f, g )$.
By \eqref{eq:spec-measure-all} the measure $\mu$ has the density
$E \mapsto ( (\Fpertg{-}f)(E), (\Fpertg{-}g)(E) )_{\Echan{E}}$ on $\Ereg$ and
vanishes outside $\Ereg$. Hence
\begin{equation*}
  \tag{$\ast\ast$}
  \label{eq:fc-all}
  \bigl( \psi(\Hpert)P(\mathcal{G})f, g \bigr)
  = \int_{\Ereg} \psi(E)\,
    \bigl( (\Fpertg{-}f)(E), (\Fpertg{-}g)(E) \bigr)_{\Echan{E}}\,dE .
\end{equation*}

Let now $\varphi$ be a bounded Borel function on $\mathbb{R}$ and
$f \in \Ltwo{\Om}$. Write $\varphi\,\Fpertg{-}f$ for the element
$E \mapsto \varphi(E)(\Fpertg{-}f)(E)$ of $\Hdirect$, and put
\begin{equation*}
  D = \Fpertg{-}\varphi(\Hpert)f - \varphi\,\Fpertg{-}f ,
  \qquad
  I_\varphi = \int_{\Ereg} |\varphi(E)|^2\,
    \norm{(\Fpertg{-}f)(E)}_{\Echan{E}}^2\,dE .
\end{equation*}
Then
\begin{equation*}
  \norm{D}_{\Hdirect}^2
  = \norm{\Fpertg{-}\varphi(\Hpert)f}_{\Hdirect}^2
    - 2\re\bigl( \Fpertg{-}\varphi(\Hpert)f,
                 \varphi\,\Fpertg{-}f \bigr)_{\Hdirect}
    + \norm{\varphi\,\Fpertg{-}f}_{\Hdirect}^2 ,
\end{equation*}
and we show that each of the three terms equals $I_\varphi$. For the last one
this is the definition of the norm of $\Hdirect$. For the first, Step~1 gives
$\norm{\Fpertg{-}h} = \norm{P(\mathcal{G})h}$ for $h = \varphi(\Hpert)f$.
Since $\varphi(\Hpert)^{*} = \overline{\varphi}(\Hpert)$ and bounded Borel
functions of $\Hpert$ commute with each other,
\begin{equation*}
  \norm{\Fpertg{-}\varphi(\Hpert)f}_{\Hdirect}^2
  = \bigl( \overline{\varphi}(\Hpert)P(\mathcal{G})\varphi(\Hpert)f,
           f \bigr)
  = \bigl( |\varphi|^2(\Hpert)P(\mathcal{G})f, f \bigr)
  = I_\varphi ,
\end{equation*}
the last step being \eqref{eq:fc-all} with $\psi = |\varphi|^2$ and $g = f$.
For the middle term, the definition of the inner product of $\Hdirect$ and
\eqref{eq:fc-all} with $\psi = \overline{\varphi}$, with $\varphi(\Hpert)f$ in
place of $f$ and with $g = f$, give
\begin{equation*}
  \begin{aligned}
    \bigl( \Fpertg{-}\varphi(\Hpert)f, \varphi\,\Fpertg{-}f \bigr)_{\Hdirect}
    &= \int_{\Ereg} \overline{\varphi(E)}\,
      \bigl( (\Fpertg{-}\varphi(\Hpert)f)(E), (\Fpertg{-}f)(E)
      \bigr)_{\Echan{E}}\,dE
    \\
    &= \bigl( \overline{\varphi}(\Hpert)P(\mathcal{G})\varphi(\Hpert)f,
             f \bigr) ,
  \end{aligned}
\end{equation*}
which is $I_\varphi$ by the preceding computation. Hence
$\norm{D}_{\Hdirect}^2 = I_\varphi - 2I_\varphi + I_\varphi = 0$, that is,
$\Fpertg{-}\varphi(\Hpert)f = \varphi\,\Fpertg{-}f$ in $\Hdirect$. This is
\eqref{eq:diagonalization}.

\emph{Step 3: the range is everything.} Let $v \in \Hdirect$ be orthogonal to
$\mathcal{M}$, and fix a representative $E \mapsto v(E)$ of $v$. For
$f \in \Bspd{\Om}$ put
\begin{equation*}
  a_f(E) = \bigl( \Fpert{-}{E}f, v(E) \bigr)_{\Echan{E}} ,
  \qquad E \in \mathcal{G} .
\end{equation*}
The function $a_f$ is measurable, as the inner product of two measurable
sections, and the Cauchy--Schwarz inequality, first in $\Echan{E}$ and then in
$E$, gives $\int_{\mathcal{G}} |a_f|\,dE \le \norm{\Fpertg{-}f}\norm{v}$.
Choose a bounded Borel function $\varphi$ on $\mathbb{R}$, vanishing outside
$\mathcal{G}$, such that $\varphi\,a_f = |a_f|$ almost everywhere on
$\mathcal{G}$. For instance, take a Borel function which agrees almost
everywhere with $\overline{a_f}/|a_f|$ on the set where $a_f \ne 0$ and
vanishes elsewhere. By Step~2,
$\varphi\,\Fpertg{-}f = \Fpertg{-}\varphi(\Hpert)f$ lies in $\mathcal{M}$, and
therefore
\begin{equation*}
  0 = \bigl( \varphi\,\Fpertg{-}f, v \bigr)_{\Hdirect}
    = \int_{\mathcal{G}} \varphi(E)\,a_f(E)\,dE
    = \int_{\mathcal{G}} |a_f(E)|\,dE .
\end{equation*}
Hence there is a null set $N_f$ such that $a_f(E) = 0$ for every
$E \in \mathcal{G} \setminus N_f$.

The null set $N_f$ depends on $f$. Lemma~\ref{lem:Fpert-elementary}(d), on the
other hand, concerns all $f \in C_c^{\infty}(\Om)$ at one and the same energy,
and an uncountable union of null sets need not be null. We therefore pass to a
countable family, which is possible for two reasons. First, for fixed
$E \in \mathcal{G}$ the map $f \mapsto a_f(E)$ is a bounded linear functional
on $\Bspd{\Om}$, by Lemma~\ref{lem:Fpert-elementary}(a). The set of
$f \in \Bspd{\Om}$ with $a_f(E) = 0$ is therefore closed in $\Bspd{\Om}$.
Second, $C_c^{\infty}(\Om)$ contains a countable subset which is dense in it
for the norm $\BnormO{\,\cdot\,}$. We now construct such a subset.

We begin with a bound of the norm of $\Bspd{\Om}$ on bounded sets. Let $m$ be
an integer with $m \ge \Rzero$, and let $u \in \Ltwo{\Om}$ vanish outside
$\Omin{m}$. A shell $\shellp{j}$ meets $\Omin{m}$ only if $2^{j-1} < m$. Let
$j_m$ be the largest such $j$, so that $2^{j_m} < 2m$. If $j_m \le \Jzero$,
only the first term of \eqref{eq:B-norms-pert} is present, and
$\BnormO{u} = \Rzero^{1/2}\norm{u}_{\Ltwo{\Om}}
\le m^{1/2}\norm{u}_{\Ltwo{\Om}}$.
Otherwise the Cauchy--Schwarz inequality for the finitely many terms of
\eqref{eq:B-norms-pert}, together with \eqref{eq:shells-split}, gives
\begin{equation*}
  \BnormO{u}
  \le \Bigl( \Rzero + \sum_{j = \Jzero + 1}^{j_m} 2^{j} \Bigr)^{1/2}
      \Bigl( \norm{u}_{\Ltwo{\Omin{\Rzero}}}^2
        + \sum_{j = \Jzero + 1}^{j_m} \norm{u}_{\Ltwo{\shellp{j}}}^2
      \Bigr)^{1/2}
  = \bigl( 2^{j_m + 1} - \Rzero \bigr)^{1/2}\,\norm{u}_{\Ltwo{\Om}} .
\end{equation*}
In either case
\begin{equation*}
  \BnormO{u} \le 2m^{1/2}\,\norm{u}_{\Ltwo{\Om}} .
\end{equation*}

For each integer $m \ge \Rzero$ put
\begin{equation*}
  K_m = \bigl\{\, x \in \Om : |x| \le m,\ \dist(x,\Gpert) \ge 1/m \,\bigr\} .
\end{equation*}
This set is compact, since a limit of points of $K_m$ lies in
$\overline{\Om}$ at distance at least $1/m$ from $\Gpert$, hence in $\Om$.
Every compact subset of $\Om$ is bounded and has positive distance from
$\Gpert$, so it lies in $K_m$ for some $m$. Fix
$\chi_m \in C_c^{\infty}(\Om)$ with $0 \le \chi_m \le 1$, with $\chi_m = 1$ on
$K_m$ and with support in $\Omin{m+1}$. Let $\{p_k\}_{k \ge 1}$ be an
enumeration of the polynomials on $\Rn$ with coefficients in
$\mathbb{Q} + i\mathbb{Q}$, and put
\begin{equation*}
  f_{m,k} = \chi_m\,p_k \in C_c^{\infty}(\Om) ,
  \qquad m \ge \Rzero,\ k \ge 1 .
\end{equation*}
The family $\{f_{m,k}\}$ is countable, and the polynomials $p_k$ are the same
for every $m$. Let $g \in C_c^{\infty}(\Om)$ and $\delta > 0$, and choose $m$
with $\operatorname{supp} g \subset K_m$. Extended by zero, $g$ is continuous
on $\Rn$. The Weierstrass approximation theorem, applied to its real and
imaginary parts on the compact set $\overline{\Bfull{m+1}}$, provides a
polynomial within $\delta/2$ of $g$ there. Moving its finitely many
coefficients to nearby points of $\mathbb{Q} + i\mathbb{Q}$ changes it by
less than $\delta/2$ on $\overline{\Bfull{m+1}}$, so that
$\sup_{\overline{\Bfull{m+1}}} |g - p_k| < \delta$ for some $k$. Since
$\chi_m g = g$, the difference $g - f_{m,k} = \chi_m(g - p_k)$ vanishes
outside $\Omin{m+1}$, and the bound above with $m + 1$ in place of $m$ gives
\begin{equation*}
  \BnormO{g - f_{m,k}}
  \le 2(m+1)^{1/2}\,\norm{\chi_m(g - p_k)}_{\Ltwo{\Om}}
  \le 2(m+1)^{1/2}\,|\Bfull{m+1}|^{1/2}\,\delta .
\end{equation*}
As $\delta > 0$ is arbitrary, $\{f_{m,k}\}$ is dense in $C_c^{\infty}(\Om)$
for the norm $\BnormO{\,\cdot\,}$.

We can now conclude. The set
\begin{equation*}
  N = \Sexc \cup \bigcup_{m \ge \Rzero,\,k \ge 1} N_{f_{m,k}}
\end{equation*}
is a countable union of null sets, hence null. Let $E \in \Ereg \setminus N$.
Then $a_{f_{m,k}}(E) = 0$ for all $m$ and $k$. By the two facts established
above, $a_f(E) = 0$ for every $f \in C_c^{\infty}(\Om)$, that is, $v(E)$ is
orthogonal to $\Fpert{-}{E}C_c^{\infty}(\Om)$. Therefore $v(E) = 0$ by
Lemma~\ref{lem:Fpert-elementary}(d). Consequently
\begin{equation*}
  \norm{v}_{\Hdirect}^2 = \int_{\Ereg} \norm{v(E)}_{\Echan{E}}^2\,dE = 0 ,
\end{equation*}
so $v = 0$. Thus $\mathcal{M}^{\perp} = \{0\}$, and since $\mathcal{M}$ is
closed by Step~1, $\mathcal{M} = (\mathcal{M}^{\perp})^{\perp} = \Hdirect$.

\emph{Step 4: conclusion.} Steps~1 and~3 make $\Fpertg{-}$ a partial isometry
with initial space $\Hac(\Hpert)$ and final space $\Hdirect$, and Step~2 is
\eqref{eq:diagonalization}. Since $\Fzero$ is a unitary operator from
$\Ltwo{\Rnp}$ onto $\Hdirect$ diagonalizing $\Halpha$
\cite[Theorem 2.30]{PartI}, the operator
$\Fzero^{*}\,\Fpertg{-}\big|_{\Hac(\Hpert)}$ is unitary from $\Hac(\Hpert)$
onto $\Ltwo{\Rnp}$ and intertwines $\Hpert$ with $\Halpha$.
\end{proof}

\begin{corollary}
\label{cor:ac-spectrum-pert}
Under the hypotheses of Theorem~\ref{thm:spectral-rep-pert},
\begin{equation}
  \label{eq:ac-spectrum-pert}
  \sac(\Hpert) = [-\alpha^2,\infty),
\end{equation}
and the absolutely continuous part of $\Hpert$ has the same spectral
multiplicity function as $\Halpha$: multiplicity $\dim\Echan{E}$ at $E$, with
the surface channel alone contributing on $\Eband$ and both channels
contributing on $(0,\infty)$.
\end{corollary}

\begin{proof}
Immediate from Theorem~\ref{thm:spectral-rep-pert}, the direct integral
\eqref{eq:direct-integral} being taken over
$\Ereg = \Eband \cup (0,\infty)$, whose closure is $[-\alpha^2,\infty)$.
\end{proof}

\begin{remark}[The stationary route]
\label{rem:stationary-completeness}
Theorem~\ref{thm:spectral-rep-pert} asserts that the final space is all of
$\Hdirect$, and its proof constructs no wave operator. Three ingredients carry
it, and all of them come from Parts~I and~II. The injectivity of
$\Phipert{\pm}{E}$ (Corollary~\ref{cor:scattering-parametrization}), that is,
the uniqueness of the scattering solution, makes the range of
$\Fpert{\pm}{E}$ dense in $\Echan{E}$ at each energy. The diagonalization
\eqref{eq:diagonalization} carries this density from the fibres to the direct
integral. Finally, each transform $\Fpert{\pm}{E}$ is an operator defined at
every $E \in \Ereg \setminus \Sexc$ and bounded on $\Bspd{\Om}$
(Lemma~\ref{lem:Fpert-elementary}(a)), not merely an object determined for
almost every $E$. The almost everywhere statements of the direct integral can
therefore be tested at a single energy against a countable dense family of
sources. Corollary~\ref{cor:Fpert-continuity} enters only through the
measurability of the sections $E \mapsto \Fpert{\pm}{E}f$. Step~3 of the proof
carries out directly the part of a general theorem that the argument needs.
That theorem identifies the commutant of the diagonalizable operators in a
direct integral with the algebra of decomposable operators
\cite[Part II, Chapter 2]{Dixmier1981} (see also
\cite[\S XIII.16]{ReedSimonIV}). From that point of view the range of
$\Fpertg{\pm}$ is a decomposable subspace whose fibres contain the ranges of
the transforms $\Fpert{\pm}{E}$, and it is therefore all of $\Hdirect$.
No hypothesis excluding embedded eigenvalues in the band is needed,
$\Sexc$ being countable and hence Lebesgue null. The comparison with the routes
of \cite{Isozaki2001,IsozakiKurylev2014} is made in
\S\ref{subsec:intro-novelty}. The scattering matrix of \S\ref{sec:smatrix} is
likewise obtained at a fixed energy from the two transforms, and we do not
claim that the operators implicit in \eqref{eq:intertwining} agree with
time-dependent wave operators.
\end{remark}

\subsection{The Far Field and the Scattering Amplitude}
\label{subsec:far-field}

The far field of the outgoing resolvent is now a corollary of what precedes.
Theorem~\ref{thm:radiation-equivalence} already established that
$\Rpert{E+i0}f$ satisfies the outgoing radiation condition. What was missing
was the identification of the two amplitudes.
Theorem~\ref{thm:Fpert-explicit} supplies it, the amplitudes of the
transplanted problem on $\Rnp$ being those of Part~I.

\begin{theorem}[Far field of the outgoing resolvent]
\label{thm:far-field-resolvent}
Let $E \in \Ereg \setminus \Sexc$ and $f \in \Bspd{\Om}$, and let $\eta$ be an
admissible cut-off.
\begin{enumerate}
  \item[\textup{(a)}] If $E > 0$, then
  $u = \Rpert{E+i0}f$ satisfies the outgoing radiation condition with
  \begin{equation}
    \label{eq:far-field-resolvent}
    \eta u \;\simeq\; \Wb{+}{A} + \Ws{+}{C},
    \qquad
    \begin{aligned}
      A &= c_n\,\cph{\kbulk\theta_n}\,\bigl(\Fbulkpert{-}{E}f\bigr),
      \\
      C &= c_n'\,\bigl(\Fsurfpert{-}{E}f\bigr),
    \end{aligned}
  \end{equation}
  where $c_n$ and $c_n'$ are the constants of
  \cite[equations (2.102)--(2.103)]{PartI}, with
  $|c_n|^2 = \pi/\kbulk$ and $|c_n'|^2 = \pi/\ksurf$.
  \item[\textup{(b)}] If $E \in \Eband \setminus \Sexc$, then
  \begin{equation}
    \label{eq:far-field-resolvent-band}
    \eta\,\Rpert{E+i0}f \;\simeq\; \Ws{+}{C},
    \qquad
    C = c_n'\,\bigl(\Fsurfpert{-}{E}f\bigr) .
  \end{equation}
\end{enumerate}
In both cases the amplitudes are independent of the admissible cut-off $\eta$.
\end{theorem}

\begin{proof}
Put $u = \Rpert{E+i0}f$ and let $h_\eta$ be the transplanted source
\eqref{eq:heta-pert}. By Theorem~\ref{thm:lap-pert}(e),
$\eta u = \Ralpha{E+i0}h_\eta$.

(a) The zero extension of $\eta f$ lies in $\Bsp$, and $[-\Delta,\eta]u$ is
supported in $\mathcal{A}$. Hence $h_\eta \in \Bsp$, and
\cite[Theorem 2.49]{PartI} applies, giving
\begin{equation*}
  \Ralpha{E+i0}h_\eta
  \;\simeq\;
  c_n\,\bigl(\Fbulk^{(+)}(E)h_\eta\bigr)(\theta)\,
    r^{-\frac{n-1}{2}}e^{i\kbulk r}
  + c_n'\,\bigl(\Fsurf(E)h_\eta\bigr)(\omega)\,
    \rho^{-\frac{n-2}{2}}e^{i\ksurf\rho}\,\psizero(x_n),
\end{equation*}
which in the notation \eqref{eq:model-waves-pert} is
$\Wb{+}{c_n\Fbulk^{(+)}(E)h_\eta} + \Ws{+}{c_n'\Fsurf(E)h_\eta}$. By
\cite[Definition 2.43]{PartI} one has
$\Fbulk^{(+)}(E) = \cph{\kbulk\theta_n}\,\Fbulk(E)$ and
$\Fsurf^{(+)}(E) = \Fsurf(E)$, while Theorem~\ref{thm:Fpert-explicit} gives
$\Fzero(E)h_\eta = \Fpert{-}{E}f$ in both channels. This is
\eqref{eq:far-field-resolvent}.

(b) On the band the bulk channel is closed and
\cite[Proposition 2.64]{PartI} gives
$\Ralpha{E+i0}h_\eta \simeq \Ws{+}{c_n'\Fsurf(E)h_\eta}$. The identification
of the amplitude is as before.

Independence of $\eta$ is Lemma~\ref{lem:radiation-elementary}(b). It is also
visible directly, the right-hand sides of
\eqref{eq:far-field-resolvent}--\eqref{eq:far-field-resolvent-band} containing
no cut-off.
\end{proof}

Combining the two amplitudes with Stone's formula turns \eqref{eq:flux-pert}
into a statement about radiated energy, in which the channel wavenumbers appear
as the weights.

\begin{corollary}[The radiated energy]
\label{cor:flux-far-field}
Let $E > 0$, let $f \in \Bspd{\Om}$, and let $A$ and $C$ be the amplitudes
\eqref{eq:far-field-resolvent}. Then
\begin{equation}
  \label{eq:radiated-energy}
  \im\bigl( \Rpert{E+i0}f, f \bigr)
  = \kbulk\,\norm{A}_{\Ltwo{\Sm}}^2
    + \ksurf\,\norm{C}_{\Ltwo{\Ssurf}}^2 .
\end{equation}
On $\Eband \setminus \Sexc$ the same holds for every $f \in \Bspd{\Om}$ with
the first term absent.
\end{corollary}

\begin{proof}
Since $|\cph{\,\cdot\,}| = 1$ pointwise, \eqref{eq:far-field-resolvent} and
\cite[equation (2.104)]{PartI} give
\begin{equation*}
  \begin{aligned}
    \norm{A}^2 &= |c_n|^2\,\norm{\Fbulkpert{-}{E}f}^2
                = \frac{\pi}{\kbulk}\,\norm{\Fbulkpert{-}{E}f}^2,
    \\[1mm]
    \norm{C}^2 &= |c_n'|^2\,\norm{\Fsurfpert{-}{E}f}^2
                = \frac{\pi}{\ksurf}\,\norm{\Fsurfpert{-}{E}f}^2 .
  \end{aligned}
\end{equation*}
Substituting into
\begin{equation*}
  \norm{\Fpert{-}{E}f}_{\Echan{E}}^2
  = \norm{\Fbulkpert{-}{E}f}_{\Ltwo{\Sm}}^2
    + \norm{\Fsurfpert{-}{E}f}_{\Ltwo{\Ssurf}}^2
\end{equation*}
and using \eqref{eq:flux-pert} yields \eqref{eq:radiated-energy}.
\end{proof}

\begin{remark}[The phase factor and the normalization of the transform]
\label{rem:phase-normalization}
The factor $\cph{\kbulk\theta_n}$ in \eqref{eq:far-field-resolvent} is the one
that distinguishes the outgoing transform $\Fbulk^{(+)}(E)$ from the
standing-wave transform $\Fbulk(E)$ in Part~I
\cite[Definition 2.43]{PartI}, and it appears here only because
Definition~\ref{def:Fpert} was built on the latter. It is unimodular, so it
affects neither the norms in Corollary~\ref{cor:flux-far-field} nor the
scattering amplitude introduced below, and multiplication by it is a unitary
automorphism of $\Echan{E}$. The surface channel carries no such factor, since
$\Fsurf^{(\pm)}(E) = \Fsurf(E)$. Which of the two normalizations to fix is a
question about the scattering matrix rather than about the far field, and it is
settled in \S\ref{sec:smatrix}.
\end{remark}

\subsubsection*{The scattering amplitude}

We can now define the object promised at the end of
\S\ref{sec:scattering-solution}. Let $E \in \Ereg \setminus \Sexc$ and
$\phi \in \Nsol{E}$, and recall from Theorem~\ref{thm:scattering-solution} that
the outgoing scattering solution is
\begin{equation}
  \label{eq:scattering-solution-recall}
  \Phipert{+}{E}\phi = \eta\phi + v,
  \qquad
  v = -\,\Rpert{E+i0}\bigl( [-\Delta,\eta]\phi \bigr) .
\end{equation}
The scattered field $v$ is the outgoing resolvent applied to a compactly
supported element of $\Bspd{\Om}$, so Theorem~\ref{thm:far-field-resolvent}
applies to it verbatim. Its amplitudes are what an experiment measures.

\begin{definition}[The scattering amplitude]
\label{def:scattering-amplitude}
Let $E \in \Ereg \setminus \Sexc$. The \emph{scattering amplitude} of $\Hpert$
at energy $E$ is the operator
\begin{equation}
  \label{eq:scattering-amplitude}
  \Ampl{E} : \Echan{E} \longrightarrow \Echan{E},
  \qquad
  \Ampl{E}g
  := -\,\Fpert{-}{E}\Bigl( [-\Delta,\eta]\,\Fzero(E)^{*}g \Bigr),
\end{equation}
where $\eta$ is an admissible cut-off. Part~(a) of the next theorem shows that
the right-hand side does not depend on which one is taken.
\end{definition}

\begin{theorem}[Properties of the scattering amplitude]
\label{thm:scattering-amplitude}
Let $E \in \Ereg \setminus \Sexc$.
\begin{enumerate}
  \item[\textup{(a)}] Definition~\ref{def:scattering-amplitude} does not depend
  on the choice of the admissible cut-off $\eta$, and $\Ampl{E}$ is a bounded
  operator on $\Echan{E}$ whose norm is bounded uniformly on compact subsets of
  $\Ereg \setminus \Sexc$. For each fixed $g$ the map $E \mapsto \Ampl{E}g$ is
  continuous.
  \item[\textup{(b)}] Let $g \in \Echan{E}$, put $\phi = \Fzero(E)^{*}g$ and
  let $\yout,\yin \in \Echan{E}$ be the far-field data of $\phi$ furnished by
  \cite[Corollary 2.70]{PartI}. Then the scattering solution
  $\Fpert{+}{E}^{*}g = \Phipert{+}{E}\phi$ has the far field
  \begin{equation}
    \label{eq:far-field-scattering-solution}
    \eta\,\Fpert{+}{E}^{*}g
    \;\simeq\;
    \Wb{+}{A} + \Wb{-}{B} + \Ws{+}{C} + \Ws{-}{D},
  \end{equation}
  where the incoming amplitudes are those of the incident wave alone,
  \begin{equation}
    \label{eq:incoming-amplitudes}
    B = \overline{\gambulk(E)}\,\yin_{\mathrm{bulk}},
    \qquad
    D = \overline{\gamsurf(E)}\,\yin_{\mathrm{surf}},
  \end{equation}
  while the outgoing ones carry the scattering amplitude,
  \begin{equation}
    \label{eq:outgoing-amplitudes}
    \begin{aligned}
      A &= \gambulk(E)\,\yout_{\mathrm{bulk}}
           + c_n\,\cph{\kbulk\theta_n}\,\bigl(\Ampl{E}g\bigr)_{\mathrm{bulk}},
      \\
      C &= \gamsurf(E)\,\yout_{\mathrm{surf}}
           + c_n'\,\bigl(\Ampl{E}g\bigr)_{\mathrm{surf}} .
    \end{aligned}
  \end{equation}
  The amplitudes $A$ and $B$ are absent when $E \in \Eband$, and
  $\gambulk(E)$, $\gamsurf(E)$ are the constants of
  \cite[equation (2.122)]{PartI}.
  \item[\textup{(c)}] With $g$ and $\phi$ as in \textup{(b)}, one has
  $\Ampl{E}g = 0$ if and only if $\Phipert{+}{E}\phi$ and $\phi$ have the same
  far field. In particular $\Ampl{E} = 0$ for every $E$ when $\Om = \Rnp$.
\end{enumerate}
\end{theorem}

\begin{proof}
(a) Let $\eta$ and $\tilde\eta$ be admissible cut-offs and let
$v, \tilde v$ be the corresponding scattered fields
\eqref{eq:scattering-solution-recall}. By
Theorem~\ref{thm:scattering-solution} the scattering solution
$\Phipert{+}{E}\phi$ does not depend on the cut-off, so
$v - \tilde v = (\tilde\eta - \eta)\phi$, which is supported in a compact
subset of $\Om$ and therefore lies in $\Bstzd{\Om}$ by
Proposition~\ref{prop:AH-elementary}(b). By
Lemma~\ref{lem:radiation-elementary}(c) the amplitudes of $v$ and $\tilde v$
coincide, and by Theorem~\ref{thm:far-field-resolvent} those amplitudes are
nonzero multiples of the two channel components of
$-\Fpert{-}{E}([-\Delta,\eta]\phi)$ and of
$-\Fpert{-}{E}([-\Delta,\tilde\eta]\phi)$ respectively, the constants $c_n$,
$c_n'$ and $\cph{\,\cdot\,}$ being the same in both cases. Hence the two
definitions agree.

Boundedness follows from \eqref{eq:scattering-amplitude}. The map
$g \mapsto [-\Delta,\eta]\Fzero(E)^{*}g$ is bounded from $\Echan{E}$ to
$\Bspd{\Om}$, because $\Fzero(E)^{*}$ is bounded into $\Bst$ by
\eqref{eq:Fzero-two-sided} and $[-\Delta,\eta]$ maps
$\Bst \cap H^2_{\mathrm{loc}}$ boundedly into the compactly supported elements
of $\Bsp$ by Lemma~\ref{lem:gluing}(a). Moreover $\Fpert{-}{E}$ is bounded by
Lemma~\ref{lem:Fpert-elementary}(a). Both bounds are uniform on compact subsets
of $\Ereg \setminus \Sexc$, the first by \cite[Proposition 2.33]{PartI} and
the second by Corollary~\ref{cor:Fpert-continuity}.

There remains the continuity in $E$. Let $E_j \to E$ within
$\Ereg \setminus \Sexc$, all lying in one of the two intervals, so that the
channel space is the same for all of them. Put $\phi_j = \Fzero(E_j)^{*}g$ and
$k_j = [-\Delta,\eta]\phi_j$, and let $\phi$ and $k$ be the corresponding
objects at $E$. Then
\begin{equation*}
  \Ampl{E_j}g - \Ampl{E}g
  = -\,\Fpert{-}{E_j}\bigl( k_j - k \bigr)
    - \bigl( \Fpert{-}{E_j} - \Fpert{-}{E} \bigr)k .
\end{equation*}
The second term tends to $0$ by Corollary~\ref{cor:Fpert-continuity}, and the
norms $\norm{\Fpert{-}{E_j}}$ are bounded by the same corollary. It therefore
suffices to prove that $k_j \to k$ in $\Bspd{\Om}$. All the $k_j$ are
supported in the fixed compact set $\mathcal{A}$, so it suffices in turn to
prove that $\phi_j \to \phi$ in $H^1$ of a neighbourhood of $\mathcal{A}$.

Fix radii $R_2 < \rho < \rho'$ and put $D_\rho = \Bfull{\rho} \cap \Rnp$.
Four facts give the convergence. First, $\Bstnorm{\phi_j}$ is bounded, by
\eqref{eq:Fzero-two-sided} and the uniform bound of
\cite[Proposition 2.33]{PartI}, and hence so is $\norm{\phi_j}_{\Ltwo{D_{\rho'}}}$.
Second, each $\phi_j$ solves the unperturbed problem at energy $E_j$, so
Proposition~\ref{prop:local-regularity}, applied on $\Rnp$ with $z = E_j$ and
$F = 0$, bounds $\norm{\phi_j}_{H^2(D_\rho)}$ by
$C\norm{\phi_j}_{\Ltwo{D_{\rho'}}}$ with $C$ independent of $j$. Third, the
embedding $H^2(D_\rho) \hookrightarrow H^1(D_\rho)$ is compact, so the sequence
has $H^1(D_\rho)$ convergent subsequences. Fourth, for
$f \in C_c^{\infty}(\Rnp) \subset \Bsp$ one has
$( f, \phi_j ) = ( \Fzero(E_j)f, g )_{\Echan{E}}$, which converges to
$( f, \phi )$ by \cite[Proposition 2.33]{PartI}. The fourth fact identifies
every subsequential limit of the third as $\phi$, so the full sequence
converges in $H^1(D_\rho)$. Since $\mathcal{A} \subset D_\rho$, this gives
$k_j \to k$ in $\Ltwo{D_\rho}$ and therefore in $\Bspd{\Om}$.

(b) By \eqref{eq:scattering-solution-recall},
$\Fpert{+}{E}^{*}g = \eta\phi + v$. The two summands are treated separately.
The zero extension of $\eta\phi$ differs from $\phi$ by a compactly supported
function, hence $\eta\phi \simeq \phi$, and the far field of $\phi \in \Nsol{E}$
is \cite[equation (2.149)]{PartI}, which in the notation
\eqref{eq:model-waves-pert} reads
\begin{equation*}
  \phi \;\simeq\; \Wb{+}{\gambulk(E)\yout_{\mathrm{bulk}}}
    + \Wb{-}{\overline{\gambulk(E)}\yin_{\mathrm{bulk}}}
    + \Ws{+}{\gamsurf(E)\yout_{\mathrm{surf}}}
    + \Ws{-}{\overline{\gamsurf(E)}\yin_{\mathrm{surf}}} .
\end{equation*}
For $v$, apply Theorem~\ref{thm:far-field-resolvent} with
$f = -[-\Delta,\eta]\phi$, which is a compactly supported element of
$\Bspd{\Om}$. By \eqref{eq:scattering-amplitude} its amplitudes are
$c_n\cph{\kbulk\theta_n}(\Ampl{E}g)_{\mathrm{bulk}}$ and
$c_n'(\Ampl{E}g)_{\mathrm{surf}}$. Adding the two expansions gives
\eqref{eq:far-field-scattering-solution}, the incoming amplitudes receiving no
contribution from $v$ because $v$ satisfies the outgoing radiation condition.

(c) By (b) the far fields of $\Phipert{+}{E}\phi$ and of $\phi$ differ exactly
by the outgoing waves with amplitudes
$c_n\cph{\kbulk\theta_n}(\Ampl{E}g)_{\mathrm{bulk}}$ and
$c_n'(\Ampl{E}g)_{\mathrm{surf}}$. These vanish if and only if
$\Ampl{E}g = 0$, the constants being nonzero and $\cph{\,\cdot\,}$ unimodular.
The amplitudes are determined by the solution, by
Lemma~\ref{lem:radiation-elementary}(b).

Suppose now that $\Om = \Rnp$. Then $\Rpert{E+i0} = \Ralpha{E+i0}$, and
$\phi \in \Nsol{E}$ gives $[-\Delta,\eta]\phi = (-\Delta - E)(\eta\phi)$, so
that $v = -\Ralpha{E+i0}(-\Delta-E)(\eta\phi) = -\eta\phi + w$ for some
$w \in \Nsol{E}$. The difference $w - \phi$ equals $v + (\eta\phi - \phi)$,
whose second summand is compactly supported and hence lies in $\Bstz$ by
Proposition~\ref{prop:AH-elementary}(b). By
Lemma~\ref{lem:radiation-elementary}(c) the difference is therefore outgoing,
and it solves the unperturbed problem, so
Theorem~\ref{thm:outgoing-uniqueness} places it in $\Bstz$. Consequently
$v = (w - \phi) + (\phi - \eta\phi) \in \Bstz$, and $\Ampl{E}g = 0$ for
every $g$.
\end{proof}

Since $\yout = \Smat(E)\Jchan\,\yin$ for the unperturbed operator
\cite[equation (2.150)]{PartI}, with $\Jchan$ the involution of
\cite[Definition 2.57]{PartI}, formula \eqref{eq:outgoing-amplitudes} shows
that the map from incoming to outgoing data of the perturbed problem is
$\Smat(E)\Jchan$ plus a term built from $\Ampl{E}$. Section~\ref{sec:smatrix}
identifies this map as the scattering matrix of $\Hpert$.

\section{The Scattering Matrix}
\label{sec:smatrix}
Section~\ref{sec:gft-pert} produced two objects: a spectral representation of
$\Hpert$ carried by the transforms $\Fpert{\pm}{E}$, and a scattering amplitude
$\Ampl{E}$ reading the far field of the scattered wave. The present section
shows that these are two descriptions of a single unitary operator on the
channel space, the scattering matrix, and identifies it with the map that
sends the incoming far-field data of a solution to its outgoing data.

Two structural differences from the flat problem appear along the way. The
first concerns the off-diagonal blocks. In the flat half-space they vanish
identically \cite[equation (2.151)]{PartI}, and here nothing forces them to.
The second concerns the surface band. There Part~I had the identity and
nothing else, whereas the perturbed matrix is an arbitrary unitary operator on
$\Ltwo{\Ssurf}$ as far as the present section can tell. Both differences are
read off the single formula \eqref{eq:smatrix-amplitude} below. Neither is a
claim that the perturbation actually produces the effect in question, and
Remark~\ref{rem:band-single-block} says what is and is not decided.

Everything here is stationary. Wave operators are not constructed, and the
scattering matrix is defined directly, either through the far field or through
the transforms. The comparison with the time-dependent formulation is recorded,
without being used, in Remark~\ref{rem:stationary-completeness}.

\subsection{The Space of Perturbed Solutions}
\label{subsec:solution-space-pert}

The far-field data of Part~I are attached to elements of $\Nsol{E}$, and
$\Nsol{E}$ was characterized there as the range of the eigenoperator
\cite[Theorem 2.68]{PartI}. We begin by proving the perturbed counterpart,
which is what makes the scattering matrix of \S\ref{subsec:smatrix-def} a map
between solutions rather than between parameters.

\begin{definition}[Perturbed solutions of $\Bst$ growth]
\label{def:Nsolpert}
For $E \in \Ereg$ put
\begin{equation}
  \label{eq:Nsolpert-def}
  \Nsolpert{E} = \Bigl\{ u \in \Bstd{\Om} \cap H^2_{\mathrm{loc}}(\overline{\Om})
    : (-\Delta - E)u = 0 \text{ in } \Om,\
    \gamma_1 u - \alpha\gamma_0 u = 0 \text{ on } \Gpert \Bigr\} .
\end{equation}
\end{definition}

No radiation condition is imposed in \eqref{eq:Nsolpert-def}, and no decay
beyond membership in $\Bstd{\Om}$. The space is the exact analogue of
$\Nsol{E}$ of \eqref{eq:Nsol-recall}, with $\Om$ in place of $\Rnp$.

\begin{theorem}[Characterization of perturbed $\Bst$-solutions]
\label{thm:characterization-pert}
Let $E \in \Ereg \setminus \Sexc$. Then
\begin{equation}
  \label{eq:Phipert-bijection}
  \Phipert{\pm}{E} : \Nsol{E} \longrightarrow \Nsolpert{E}
\end{equation}
is a bijection. In particular every $u \in \Nsolpert{E}$ is the outgoing
scattering solution of exactly one incident wave, and also the incoming
scattering solution of exactly one incident wave.
\end{theorem}

\begin{proof}
That $\Phipert{\pm}{E}$ maps $\Nsol{E}$ into $\Nsolpert{E}$ is part of
Theorem~\ref{thm:scattering-solution}, and injectivity is
Corollary~\ref{cor:scattering-parametrization}. Only surjectivity has to be
proved, and we treat the outgoing case. The incoming one follows by complex
conjugation. Indeed $\Om$, the equation and the boundary condition are all
real, so $u \mapsto \overline{u}$ maps $\Nsolpert{E}$ and $\Nsol{E}$ onto
themselves and interchanges $\Phipert{+}{E}$ with $\Phipert{-}{E}$, the four
model waves \eqref{eq:model-waves-pert} being interchanged in pairs by
$\overline{\Wb{\pm}{A}} = \Wb{\mp}{\overline{A}}$ and
$\overline{\Ws{\pm}{C}} = \Ws{\mp}{\overline{C}}$.

Let $u \in \Nsolpert{E}$ and fix an admissible cut-off $\eta$. By
Lemma~\ref{lem:gluing} the zero extension of $\eta u$ belongs to
$\Bst \cap H^2_{\mathrm{loc}}(\overline{\Rnp})$ and satisfies
\begin{equation}
  \label{eq:glued-solution}
  (-\Delta - E)(\eta u) = h_\eta := [-\Delta,\eta]u \quad\text{in } \Rnp,
  \qquad \gamma_1(\eta u) - \alpha\gamma_0(\eta u) = 0 \quad\text{on } \Gamma,
\end{equation}
with $h_\eta \in \Bsp$ supported in the compact annulus $\mathcal{A}$ of
\eqref{eq:annulus}. Put
\begin{equation}
  \label{eq:phi-from-u}
  \phi := \eta u - \Ralpha{E+i0}h_\eta .
\end{equation}
Then $\phi \in \Bst$ solves the unperturbed problem on $\Rnp$, so
$\phi \in \Nsol{E}$.

We claim that $d := u - \Phipert{+}{E}\phi \in \Nsolpert{E}$ satisfies the
outgoing radiation condition of Definition~\ref{def:radiation-pert}. Here
\cite[Theorem 2.49]{PartI} applies when $E > 0$, and
\cite[Proposition 2.64]{PartI} when $E \in \Eband$. In either case $\Ralpha{E+i0}h_\eta$ is
asymptotically equal to a sum of the outgoing model waves alone, so that
\eqref{eq:phi-from-u} gives
\begin{equation}
  \label{eq:etau-split}
  \eta u \;\simeq\; \phi + \Wb{+}{A_0} + \Ws{+}{C_0}
\end{equation}
for some $A_0, C_0$, the first of the two waves being absent on the band. On
the other hand $\Phipert{+}{E}\phi - \eta\phi$ satisfies the outgoing radiation
condition by Lemma~\ref{lem:Fpert-elementary}(b), and
$(1-\eta^2)\phi$ is compactly supported, hence lies in $\Bstz$. Therefore
\begin{equation}
  \label{eq:etaPhi-split}
  \eta\,\Phipert{+}{E}\phi \;\simeq\; \phi + \Wb{+}{A_1} + \Ws{+}{C_1} .
\end{equation}
Subtracting \eqref{eq:etaPhi-split} from \eqref{eq:etau-split} leaves
$\eta d \simeq \Wb{+}{A_0 - A_1} + \Ws{+}{C_0 - C_1}$, which is the claim.

By Theorem~\ref{thm:radiation-equivalence}, applied with vanishing source, $d$
is outgoing in the sense of Definition~\ref{def:outgoing-pert}, so
Theorem~\ref{thm:outgoing-uniqueness} applies. If $E > 0$ it gives $d = 0$
directly. If $E \in \Eband$ it gives that $d$ is an eigenfunction of $\Hpert$
with eigenvalue $E$, which is excluded because $E \notin \Sexc$. So $d = 0$
again. Hence $u = \Phipert{+}{E}\phi$.
\end{proof}

The far-field data of a perturbed solution now follow, together with the fact
that each of the two data determines the solution.

\begin{corollary}[Far field of a perturbed solution]
\label{cor:far-field-data-pert}
Let $E \in \Ereg \setminus \Sexc$ and $u \in \Nsolpert{E}$, and let $\eta$ be an
admissible cut-off. Then there are unique $\yin[u], \yout[u] \in \Echan{E}$,
independent of $\eta$, such that
\begin{equation}
  \label{eq:far-field-pert}
  \begin{aligned}
  \eta u \;\simeq\;
  & \Wb{+}{\gambulk(E)\,\yout_{\mathrm{bulk}}[u]}
    + \Wb{-}{\overline{\gambulk(E)}\,\yin_{\mathrm{bulk}}[u]}
  \\[1mm]
  + \;& \Ws{+}{\gamsurf(E)\,\yout_{\mathrm{surf}}[u]}
    + \Ws{-}{\overline{\gamsurf(E)}\,\yin_{\mathrm{surf}}[u]} ,
  \end{aligned}
\end{equation}
the first line being absent when $E \in \Eband$, where $\gambulk(E)$ and
$\gamsurf(E)$ are the constants of \cite[equation (2.122)]{PartI}. Moreover,
writing $\phi = (\Phipert{+}{E})^{-1}u$ and $\psi = (\Phipert{-}{E})^{-1}u$ for
the two incident waves furnished by
Theorem~\ref{thm:characterization-pert},
\begin{equation}
  \label{eq:data-inherited}
  \yin[u] = \yin[\phi],
  \qquad
  \yout[u] = \yout[\psi] ,
\end{equation}
and each of the maps $u \mapsto \yin[u]$ and $u \mapsto \yout[u]$ is a
bijection of $\Nsolpert{E}$ onto $\Echan{E}$.
\end{corollary}

\begin{proof}
Existence of the expansion is \eqref{eq:etau-split} of the preceding proof
combined with the far field \cite[equation (2.149)]{PartI} of
$\phi \in \Nsol{E}$: the resolvent term contributes only to the outgoing
waves, so the incoming amplitudes of $\eta u$ are those of $\phi$, which is
the first identity in \eqref{eq:data-inherited}. The second follows by
conjugation as in Theorem~\ref{thm:characterization-pert}. Uniqueness of the
four amplitudes is the rigidity \cite[Lemma 2.89]{PartI}, and the constants
$\gambulk(E)$, $\gamsurf(E)$ are nonzero by \cite[equation (2.124)]{PartI}.
Independence of $\eta$ holds because two admissible cut-offs differ on a
compact set, so that $\eta u - \tilde\eta u \in \Bstz$.

For the last assertion, $\phi \mapsto \yin[\phi]$ is a bijection of $\Nsol{E}$
onto $\Echan{E}$ by \cite[Corollary 2.70]{PartI}, and
$u \mapsto \phi$ is a bijection of $\Nsolpert{E}$ onto $\Nsol{E}$ by
Theorem~\ref{thm:characterization-pert}. Now use \eqref{eq:data-inherited}.
The argument for $\yout$ is the same with $\psi$ in place of $\phi$.
\end{proof}

\subsection{The Scattering Matrix}
\label{subsec:smatrix-def}

Two normalizations of the unperturbed transform coexist in Part~I, the
standing-wave $\Fzero(E)$ and the outgoing $\Fzero^{(+)}(E)$, and they differ
by a unimodular multiplier in the bulk channel
\cite[Definition 2.43]{PartI}. Definition~\ref{def:Fpert} was built on the
former, and Remark~\ref{rem:phase-normalization} deferred to the present
subsection the question of which to fix. The answer is that the multiplier
should be kept in evidence, because it is a square root of the unperturbed
scattering matrix.

\begin{definition}[The phase multiplier]
\label{def:phase-multiplier}
For $E \in \Ereg$ let $\Phase{E}$ be the operator on $\Echan{E}$ acting as
multiplication by $\cph{\kbulk\theta_n}$ on $\Ltwo{\Sm}$ and as the identity on
$\Ltwo{\Ssurf}$, where $\cph{\cdot}$ is the phase factor of
\cite[equation (2.50)]{PartI}. For $E \in \Eband$ this reads
$\Phase{E} = \mathrm{Id}_{\Ltwo{\Ssurf}}$.
\end{definition}

\begin{lemma}[Properties of the phase multiplier]
\label{lem:phase-multiplier}
Let $E \in \Ereg$.
\begin{enumerate}
  \item[\textup{(a)}] $\Phase{E}$ is unitary on $\Echan{E}$, block diagonal for
  the splitting $\Echan{E} = \Ltwo{\Sm} \oplus \Ltwo{\Ssurf}$, and commutes
  with the involution $\Jchan$ of \cite[Definition 2.57]{PartI}.
  \item[\textup{(b)}] $\Phase{E}^2 = \Smat(E)$, the unperturbed scattering
  matrix of \cite[equation (2.151)]{PartI}.
  \item[\textup{(c)}] If $\phi = \Fzero(E)^{*}g$ with $g \in \Echan{E}$, then
  $g = \Phase{E}\,\Jchan\,\yin[\phi]$.
\end{enumerate}
\end{lemma}

\begin{proof}
(a) The multiplier has modulus $1$ by \cite[equation (2.50)]{PartI}, and it
depends on $\theta$ only through $\theta_n$, whereas $\Jbulk$ reverses only the
tangential components. The surface block is the identity and commutes with
$\Jsurf$.

(b) On the bulk block, \cite[equation (2.50)]{PartI} states
$\refl{k} = \cph{k}/\overline{\cph{k}}$, and $|\cph{k}| = 1$ gives
$\refl{k} = \cph{k}^2$. Multiplication by $\refl{\kbulk\theta_n}$ is the bulk
block of $\Smat(E)$ by \cite[equation (2.151)]{PartI}. On the surface block
both sides are the identity, by the same equation.

(c) This is \cite[equation (2.152)]{PartI}.
\end{proof}

The next lemma is the point at which the antisymmetry of the transplantation
commutator is used for the last time. It says that the two ways of inserting
$[-\Delta,\eta]$ between the unperturbed and the perturbed transform differ by
passing to the adjoint.

\begin{lemma}[The amplitude as a transplanted composition]
\label{lem:amplitude-adjoint}
Let $E \in \Ereg \setminus \Sexc$ and let $\eta$ be an admissible cut-off.
Then
\begin{equation}
  \label{eq:amplitude-adjoint}
  \Fzero(E)\Bigl( [-\Delta,\eta]\,\Fpert{-}{E}^{*}g \Bigr) = \Ampl{E}^{*}g,
  \qquad g \in \Echan{E} .
\end{equation}
\end{lemma}

\begin{proof}
By Lemma~\ref{lem:Fpert-elementary}(b) the function $\Fpert{-}{E}^{*}g$ lies in
$H^2_{\mathrm{loc}}(\overline{\Om})$, so $[-\Delta,\eta]\Fpert{-}{E}^{*}g$ is
supported in the annulus $\mathcal{A}$ and its zero extension lies in $\Bsp$.
The left-hand side of \eqref{eq:amplitude-adjoint} is therefore defined. Let
$h \in \Echan{E}$. Using in turn the adjoint relation for $\Fzero(E)$,
Lemma~\ref{lem:commutator-antisym} with $u = \Fpert{-}{E}^{*}g$ and
$\psi = \Fzero(E)^{*}h$, the defining relation \eqref{eq:Fpert-def}, and
Definition~\ref{def:scattering-amplitude}, we obtain
\begin{equation*}
  \begin{aligned}
  \Bigl( \Fzero(E)\bigl( [-\Delta,\eta]\Fpert{-}{E}^{*}g \bigr), h
  \Bigr)_{\Echan{E}}
  &= \bigl( [-\Delta,\eta]\Fpert{-}{E}^{*}g, \Fzero(E)^{*}h \bigr)
  \\
  &= - \bigl( \Fpert{-}{E}^{*}g, [-\Delta,\eta]\Fzero(E)^{*}h \bigr)
  \\
  &= - \overline{\bigl( [-\Delta,\eta]\Fzero(E)^{*}h, \Fpert{-}{E}^{*}g \bigr)}
  \\
  &= - \overline{\bigl( \Fpert{-}{E}\bigl( [-\Delta,\eta]\Fzero(E)^{*}h \bigr),
       g \bigr)_{\Echan{E}}}
  \\
  &= \overline{\bigl( \Ampl{E}h, g \bigr)_{\Echan{E}}}
   = \bigl( \Ampl{E}^{*}g, h \bigr)_{\Echan{E}} .
  \end{aligned}
\end{equation*}
Since $h$ was arbitrary, \eqref{eq:amplitude-adjoint} follows.
\end{proof}

\begin{theorem}[The two transforms are intertwined by the amplitude]
\label{thm:intertwining}
Let $E \in \Ereg \setminus \Sexc$ and put
\begin{equation}
  \label{eq:smatrix-rel-def}
  \Smatrel{E} := \mathrm{Id}_{\Echan{E}} + 2\pi i\,\Ampl{E} .
\end{equation}
Then
\begin{equation}
  \label{eq:intertwining}
  \Fpert{+}{E} = \Smatrel{E}^{*}\,\Fpert{-}{E},
  \qquad
  \Fpert{-}{E} = \Smatrel{E}\,\Fpert{+}{E}
  \qquad\text{on } \Bspd{\Om} ,
\end{equation}
and $\Smatrel{E}$ is unitary on $\Echan{E}$. Moreover $E \mapsto \Smatrel{E}g$
is continuous on each of $\Eband \setminus \Sexc$ and $(0,\infty)$, for every
$g$ in the channel space of that interval. Finally
$\Smatrel{E} = \mathrm{Id}$ for every $E$ when $\Om = \Rnp$.
\end{theorem}

\begin{proof}
Let $f \in \Bspd{\Om}$. Subtracting the two instances of
Theorem~\ref{thm:Fpert-explicit} and using \eqref{eq:Adens-pert-def},
\begin{equation*}
  \begin{aligned}
  \Fpert{-}{E}f - \Fpert{+}{E}f
  &= \Fzero(E)\Bigl( [-\Delta,\eta]
       \bigl( \Rpert{E+i0} - \Rpert{E-i0} \bigr) f \Bigr)
  \\
  &= 2\pi i\,\Fzero(E)\Bigl( [-\Delta,\eta]\,\Adenspert{E}f \Bigr) .
  \end{aligned}
\end{equation*}
By Theorem~\ref{thm:stone-pert} we may write
$\Adenspert{E}f = \Fpert{-}{E}^{*}\bigl( \Fpert{-}{E}f \bigr)$, and then
Lemma~\ref{lem:amplitude-adjoint} applied with $g = \Fpert{-}{E}f$ gives
\begin{equation}
  \label{eq:transform-difference}
  \Fpert{-}{E}f - \Fpert{+}{E}f = 2\pi i\,\Ampl{E}^{*}\,\Fpert{-}{E}f ,
\end{equation}
that is, $\Fpert{+}{E} = T\,\Fpert{-}{E}$ with
$T := \mathrm{Id} - 2\pi i\Ampl{E}^{*} = \Smatrel{E}^{*}$.

We show that $T$ is unitary. By \eqref{eq:flux-pert} one has
$\norm{T\Fpert{-}{E}f} = \norm{\Fpert{+}{E}f} = \norm{\Fpert{-}{E}f}$ for every
$f \in \Bspd{\Om}$, so $T$ is isometric on the range of $\Fpert{-}{E}$, which
is dense in $\Echan{E}$ by Lemma~\ref{lem:Fpert-elementary}(d). Being bounded,
$T$ is then isometric on all of $\Echan{E}$. Its range is therefore closed, and
it contains $T\bigl( \Ran\Fpert{-}{E} \bigr) = \Ran\Fpert{+}{E}$, which is
dense by the same lemma. Hence $T$ is onto. Thus $T$ is unitary, and so is
$\Smatrel{E} = T^{*}$. Being unitary, $T$ satisfies $T^{-1} = T^{*} =
\Smatrel{E}$, so applying $T^{-1}$ to $\Fpert{+}{E} = T\Fpert{-}{E}$ yields the
second identity in \eqref{eq:intertwining}.

The continuity in $E$ is Theorem~\ref{thm:scattering-amplitude}(a), and the
last assertion is Theorem~\ref{thm:scattering-amplitude}(c).
\end{proof}

We call $\Smatrel{E}$ the \emph{relative scattering matrix} of $\Hpert$ at
energy $E$. The reference is the unperturbed operator: $\Smatrel{E}$ compares
$\Hpert$ with $\Halpha$, and accordingly reduces to the identity when
$\Om = \Rnp$. The scattering matrix proper, which acts on far-field data, is
obtained from it by a conjugation.

\begin{definition}[The scattering matrix]
\label{def:smatrix-pert}
For $E \in \Ereg \setminus \Sexc$ the \emph{scattering matrix} of $\Hpert$ at
energy $E$ is the operator
\begin{equation}
  \label{eq:smatrix-pert-def}
  \Smatpert{E} := \Phase{E}\,\Smatrel{E}\,\Phase{E}
  \;:\; \Echan{E} \longrightarrow \Echan{E} .
\end{equation}
\end{definition}

\begin{theorem}[The scattering matrix and the far field]
\label{thm:smatrix-far-field}
Let $E \in \Ereg \setminus \Sexc$.
\begin{enumerate}
  \item[\textup{(a)}] $\Smatpert{E}$ is unitary on $\Echan{E}$, and
  \begin{equation}
    \label{eq:smatrix-amplitude}
    \Smatpert{E} = \Smat(E) + 2\pi i\,\Phase{E}\,\Ampl{E}\,\Phase{E} ,
  \end{equation}
  where $\Smat(E)$ is the unperturbed scattering matrix of
  \cite[equation (2.151)]{PartI}. In particular
  $\Smatpert{E} = \Smat(E)$ for every $E$ when $\Om = \Rnp$.
  \item[\textup{(b)}] For every $u \in \Nsolpert{E}$ the far-field data of
  Corollary~\ref{cor:far-field-data-pert} satisfy
  \begin{equation}
    \label{eq:smatrix-far-field}
    \yout[u] = \Smatpert{E}\,\Jchan\,\yin[u],
    \qquad
    \norm{\yout[u]}_{\Echan{E}} = \norm{\yin[u]}_{\Echan{E}} .
  \end{equation}
  \item[\textup{(c)}] $E \mapsto \Smatpert{E}g$ is continuous on each of
  $\Eband \setminus \Sexc$ and $(0,\infty)$, for every $g$ in the channel
  space of that interval.
\end{enumerate}
\end{theorem}

\begin{proof}
(a) Unitarity is that of $\Smatrel{E}$ (Theorem~\ref{thm:intertwining})
together with that of $\Phase{E}$ (Lemma~\ref{lem:phase-multiplier}(a)), and
\eqref{eq:smatrix-amplitude} follows from \eqref{eq:smatrix-rel-def} and
$\Phase{E}^2 = \Smat(E)$ (Lemma~\ref{lem:phase-multiplier}(b)). The last
sentence is the corresponding statement in Theorem~\ref{thm:intertwining}.

(b) Let $\phi \in \Nsol{E}$ be the incident wave with
$u = \Phipert{+}{E}\phi$, and let $g \in \Echan{E}$ be its channel datum, so
that $\phi = \Fzero(E)^{*}g$ and $u = \Fpert{+}{E}^{*}g$. By
\eqref{eq:data-inherited} and Lemma~\ref{lem:phase-multiplier}(c),
\begin{equation}
  \label{eq:datum-from-incoming}
  \yin[u] = \yin[\phi],
  \qquad
  g = \Phase{E}\,\Jchan\,\yin[u] .
\end{equation}
Comparing the outgoing amplitudes in
Theorem~\ref{thm:scattering-amplitude}(b) with those in
\eqref{eq:far-field-pert} gives
\begin{equation*}
  \begin{aligned}
  \yout_{\mathrm{bulk}}[u]
  &= \yout_{\mathrm{bulk}}[\phi]
     + \frac{c_n}{\gambulk(E)}\,\cph{\kbulk\theta_n}\,
       \bigl( \Ampl{E}g \bigr)_{\mathrm{bulk}},
  \\[1mm]
  \yout_{\mathrm{surf}}[u]
  &= \yout_{\mathrm{surf}}[\phi]
     + \frac{c_n'}{\gamsurf(E)}\,
       \bigl( \Ampl{E}g \bigr)_{\mathrm{surf}} .
  \end{aligned}
\end{equation*}
The two quotients coincide. Indeed, by
\cite[equations (2.102), (2.103) and (2.122)]{PartI},
\begin{equation}
  \label{eq:constant-quotients}
  \frac{c_n}{\gambulk(E)}
  = 2\pi\,e^{i\frac{-(n-3)+(n-1)}{4}\pi} = 2\pi i,
  \qquad
  \frac{c_n'}{\gamsurf(E)}
  = 2\pi\,e^{i\frac{-(n-4)+(n-2)}{4}\pi} = 2\pi i,
\end{equation}
the moduli being $2\pi$ in both cases because
$|c_n| = \sqrt{\pi/\kbulk}$ and $|\gambulk(E)| = (4\pi\kbulk)^{-1/2}$, and
likewise one dimension lower for the surface channel. Hence, by
Definition~\ref{def:phase-multiplier},
\begin{equation}
  \label{eq:outgoing-data-split}
  \yout[u] = \yout[\phi] + 2\pi i\,\Phase{E}\,\Ampl{E}\,g .
\end{equation}
By \cite[equation (2.150)]{PartI} one has
$\yout[\phi] = \Smat(E)\Jchan\yin[\phi]$, so that
\eqref{eq:datum-from-incoming} and
Lemma~\ref{lem:phase-multiplier}(b) turn this into
$\yout[\phi] = \Phase{E}^{2}\Jchan\yin[u]$. Substituting it, together with
$g = \Phase{E}\Jchan\yin[u]$, into \eqref{eq:outgoing-data-split} gives
\begin{equation*}
  \yout[u]
  = \Bigl( \Phase{E}^{2} + 2\pi i\,\Phase{E}\Ampl{E}\Phase{E} \Bigr)
      \Jchan\,\yin[u]
  = \Smatpert{E}\,\Jchan\,\yin[u] ,
\end{equation*}
which is the first half of \eqref{eq:smatrix-far-field}. The second half
follows because $\Smatpert{E}$ and $\Jchan$ are unitary.

(c) Immediate from Theorem~\ref{thm:intertwining}, the operator $\Phase{E}$
depending continuously on $E$ in the strong sense because
$\cph{\kbulk\theta_n}$ is continuous in $E$ and bounded by $1$.
\end{proof}

\begin{remark}[The choice of normalization]
\label{rem:normalization-fixed}
Remark~\ref{rem:phase-normalization} left open whether to build
$\Fpert{\pm}{E}$ on $\Fzero(E)$ or on $\Fzero^{(+)}(E)$. By
Lemma~\ref{lem:phase-multiplier}(b) the factor separating the two is a square
root of $\Smat(E)$. Absorbing it into the transform would therefore not remove
$\Smat(E)$ from the formulas, and keeping $\Fzero(E)$ costs only the
conjugation in \eqref{eq:smatrix-pert-def}. In return $\Smatrel{E}$ is the
identity in the unperturbed case, while $\Smatpert{E}$ reduces to $\Smat(E)$.
Since $\Phase{E}$ is block diagonal, the two matrices have the same pattern of
vanishing blocks, and \S\ref{subsec:optical} uses whichever is convenient.
The matrix $\Smatrel{E}$ measures the perturbed system against $\Halpha$, and
$\Smatpert{E}$ measures it against no scattering at all.
\end{remark}

\subsection{Channel Coupling and the Optical Theorem}
\label{subsec:optical}

For $E > 0$ the channel space splits, and every bounded operator on it is a
$2 \times 2$ operator matrix. Applying this to
Definition~\ref{def:smatrix-pert} turns the unitarity of
Theorem~\ref{thm:smatrix-far-field}(a) into a set of relations between the four
blocks, and those relations are what carries the physical content of the
section.

\begin{proposition}[The block relations]
\label{prop:block-relations}
Let $E \in \Ereg \setminus \Sexc$ with $E > 0$, and write
\begin{equation}
  \label{eq:smatrix-blocks}
  \Smatpert{E}
  = \begin{pmatrix}
      \Sbb & \Sbs \\[1mm]
      \Ssb & \Sss
    \end{pmatrix}
  \quad\text{on } \Echan{E} = \Ltwo{\Sm} \oplus \Ltwo{\Ssurf},
\end{equation}
the energy being suppressed from the four blocks, so that
$\Sbs : \Ltwo{\Ssurf} \to \Ltwo{\Sm}$ converts an incident surface wave into a
radiated bulk wave and $\Ssb : \Ltwo{\Sm} \to \Ltwo{\Ssurf}$ does the converse.
Then
\begin{align}
  \label{eq:block-bb}
  \Sbb^{*}\Sbb + \Ssb^{*}\Ssb &= \mathrm{Id}_{\Ltwo{\Sm}},
  \\
  \label{eq:block-ss}
  \Sss^{*}\Sss + \Sbs^{*}\Sbs &= \mathrm{Id}_{\Ltwo{\Ssurf}},
  \\
  \label{eq:block-cross}
  \Sbb^{*}\Sbs + \Ssb^{*}\Sss &= 0 .
\end{align}
Moreover $\Sbs = 0$ if and only if the corresponding block of $\Ampl{E}$
vanishes, and likewise for $\Ssb$. Finally $\Smatpert{E} = \Smat(E)$ if and
only if $\Ampl{E} = 0$.
\end{proposition}

\begin{proof}
The three identities are the entries of $\Smatpert{E}^{*}\Smatpert{E} =
\mathrm{Id}$, which holds by Theorem~\ref{thm:smatrix-far-field}(a). For the
last two assertions, $\Phase{E}$ is block diagonal and unitary by
Lemma~\ref{lem:phase-multiplier}(a), and $\Smat(E)$ is block diagonal by
\cite[equation (2.151)]{PartI}. So \eqref{eq:smatrix-amplitude} shows that
$\Smatpert{E} - \Smat(E)$ and $\Ampl{E}$ have the same blocks up to
composition with unitary multipliers on either side.
\end{proof}

Identity \eqref{eq:block-bb} is a conservation law with a trade-off built in.
Reading it on a purely bulk incident datum $F \in \Ltwo{\Sm}$, the two terms
$\norm{\Sbb F}^2$ and $\norm{\Ssb F}^2$ are the energies leaving the
perturbation as a bulk wave and as a surface wave, and their sum is
$\norm{F}^2$: whatever the geometry inside $\Bfull{\Rzero}$, energy diverted
into the boundary channel is energy removed from the radiated one. Nothing in
the argument computes either term, and neither is claimed to be nonzero. What
is claimed is that they cannot both be small.

The same computation performed on $\Smatrel{E}$ rather than on $\Smatpert{E}$
gives the optical theorem, in the form in which the scattering amplitude enters
directly.

\begin{theorem}[Optical theorem]
\label{thm:optical}
Let $E \in \Ereg \setminus \Sexc$. Then
\begin{equation}
  \label{eq:optical-operator}
  \frac{1}{2i}\bigl( \Ampl{E} - \Ampl{E}^{*} \bigr) = \pi\,\Ampl{E}^{*}\Ampl{E}
  \qquad\text{on } \Echan{E} .
\end{equation}
Equivalently, let $g \in \Echan{E}$, let $\phi = \Fzero(E)^{*}g$ be the
corresponding incident wave and let $A$ and $C$ be the outgoing amplitudes of
the scattered field $v = \Fpert{+}{E}^{*}g - \eta\phi$, as in
\eqref{eq:outgoing-amplitudes}. Then
\begin{equation}
  \label{eq:optical-channels}
  \im\bigl( \Ampl{E}g, g \bigr)_{\Echan{E}}
  = \kbulk\,\norm{A}_{\Ltwo{\Sm}}^2
    + \ksurf\,\norm{C}_{\Ltwo{\Ssurf}}^2 ,
\end{equation}
the first term on the right being absent when $E \in \Eband$.
\end{theorem}

\begin{proof}
By Theorem~\ref{thm:intertwining} the operator $\Smatrel{E}$ is unitary, so
\begin{equation*}
  \begin{aligned}
  \mathrm{Id}
  = \Smatrel{E}^{*}\Smatrel{E}
  &= \bigl( \mathrm{Id} - 2\pi i\Ampl{E}^{*} \bigr)
     \bigl( \mathrm{Id} + 2\pi i\Ampl{E} \bigr)
  \\
  &= \mathrm{Id} + 2\pi i\bigl( \Ampl{E} - \Ampl{E}^{*} \bigr)
     + 4\pi^2\,\Ampl{E}^{*}\Ampl{E},
  \end{aligned}
\end{equation*}
so that $\Ampl{E} - \Ampl{E}^{*} = 2\pi i\,\Ampl{E}^{*}\Ampl{E}$. Dividing by
$2i$ gives \eqref{eq:optical-operator}.

Taking the quadratic form of \eqref{eq:optical-operator} at $g$ turns the left
side into $\im( \Ampl{E}g, g )$ and the right one into
$\pi\norm{\Ampl{E}g}^2$. By Theorem~\ref{thm:scattering-amplitude}(b) the
amplitudes of $v$ are
\begin{equation*}
  A = c_n\,\cph{\kbulk\theta_n}\bigl( \Ampl{E}g \bigr)_{\mathrm{bulk}},
  \qquad
  C = c_n'\,\bigl( \Ampl{E}g \bigr)_{\mathrm{surf}} ,
\end{equation*}
so that $\norm{A}^2 = |c_n|^2\norm{(\Ampl{E}g)_{\mathrm{bulk}}}^2$ and
$\norm{C}^2 = |c_n'|^2\norm{(\Ampl{E}g)_{\mathrm{surf}}}^2$, the factor
$\cph{\,\cdot\,}$ being unimodular. Now \cite[equation (2.104)]{PartI} gives
$|c_n|^2 = \pi/\kbulk$ and $|c_n'|^2 = \pi/\ksurf$, whence
\begin{equation*}
  \pi\,\norm{\Ampl{E}g}_{\Echan{E}}^2
  = \kbulk\,\norm{A}_{\Ltwo{\Sm}}^2 + \ksurf\,\norm{C}_{\Ltwo{\Ssurf}}^2 ,
\end{equation*}
which is \eqref{eq:optical-channels}. On the band $\Echan{E} = \Ltwo{\Ssurf}$
and the bulk term is absent throughout.
\end{proof}

\begin{remark}[Reading the optical theorem]
\label{rem:optical-reading}
The right-hand side of \eqref{eq:optical-channels} is the radiated energy of
Corollary~\ref{cor:flux-far-field}, computed there as
$\im( \Rpert{E+i0}f, f )$ with $f = -[-\Delta,\eta]\phi$. Thus
\eqref{eq:optical-channels} says that the imaginary part of the forward
quantity $( \Ampl{E}g, g )$ accounts for the entire energy removed from the
incident wave, whether it leaves as a bulk wave or travels along the boundary.
This is the two-channel form of the classical optical theorem. The weights are
the channel wavenumbers, because the factors $\pi$ in
$|c_n|^2 = \pi/\kbulk$ and $|c_n'|^2 = \pi/\ksurf$ cancel against the $\pi$
produced by Stone's formula. In particular both sides of
\eqref{eq:optical-channels} are nonnegative, which is a consistency check on
the signs.
\end{remark}

\begin{remark}[Comparison with the flat half-space]
\label{rem:band-single-block}
Compare \eqref{eq:block-ss} with the corresponding statement for the flat
half-space. There the whole matrix is diagonal, and its surface block
is the identity for every $E$ and every $\alpha$
\cite[Proposition 3.3(i) and (ii)]{PartI}. The inverse problem of Part~I
consequently draws nothing at all from that block, and recovers $\alpha$ from
the bulk block alone by the explicit formula of
\cite[Corollary 3.4]{PartI}. Two things change here. First, on the
band the bulk channel is closed, $\Echan{E} = \Ltwo{\Ssurf}$, and
Theorem~\ref{thm:smatrix-far-field} produces a unitary operator on
$\Ltwo{\Ssurf}$ where Part~I had the identity: the surface channel alone
constitutes a complete one-channel scattering problem, with data available at
every $E \in \Eband \setminus \Sexc$. Second, for $E > 0$ nothing in
\eqref{eq:block-bb}--\eqref{eq:block-ss} forces the off-diagonal blocks of
\eqref{eq:smatrix-blocks} to vanish, whereas in the flat case they vanish
identically.

Proposition~\ref{prop:block-relations} does not assert that any of this
actually happens for a given $\Om$: it gives the equivalence
$\Smatpert{E} = \Smat(E) \iff \Ampl{E} = 0$, and nothing stronger. Which
perturbations are invisible at a given energy is a question about the inverse
problem, and it is settled in \S\ref{sec:inverse}: by
Corollary~\ref{cor:no-invisible} none are, at any positive energy, so that the
left-hand side of the equivalence never occurs for a nontrivial $\Om$. That
argument locates no individual block, however, and whether the off-diagonal
ones can vanish while $\Smatpert{E} \ne \Smat(E)$ is not decided here.
\end{remark}

\begin{remark}[The kernel of the amplitude]
\label{rem:amplitude-kernel}
The classical optical theorem is pointwise. It is stated for the kernels
$A_{\mu\nu}(E;\cdot\,,\cdot)$, $\mu,\nu \in \{\mathrm{b},\mathrm{s}\}$,
that $\Ampl{E}$ would have as an integral operator, and it evaluates
$\im A_{\mathrm{bb}}(E;\theta_0,\theta_0)$ against the integrals of
$|A_{\mathrm{bb}}(E;\theta,\theta_0)|^2$ over $\Sm$ and of
$|A_{\mathrm{sb}}(E;\omega,\theta_0)|^2$ over $\Ssurf$. Such kernels require
the generalized eigenfunctions $\Psib(\cdot,\theta)$ of
\cite[Definition 2.23]{PartI}, which do not belong to $\Bst$, and no statement
of the present section needs them. They are constructed in \S\ref{sec:inverse}
through the sources $[-\Delta,\eta]\Psich{\mu}(\cdot,\theta)$, which do lie
in $\Bspd{\Om}$.
\end{remark}

\section{The Inverse Problem}
\label{sec:inverse}
The direct theory is now complete. The scattering amplitude $\Ampl{E}$ of
Definition~\ref{def:scattering-amplitude} is defined for every
$E \in \Ereg \setminus \Sexc$, it determines and is determined by the
scattering matrix of Definition~\ref{def:smatrix-pert} through
\eqref{eq:smatrix-amplitude}, and Theorem~\ref{thm:optical} relates it to the
energy carried away in the two channels. The inverse problem asks what the
perturbation itself can be recovered from these data.

We proceed in three steps. The first, in \S\ref{subsec:inverse-kernels}, is a
change of viewpoint rather than a new result. The transforms of
\S\ref{sec:gft-pert} were deliberately defined as operators on $\Nsol{E}$
rather than through individual incident directions, and for the inverse
problem the individual directions are what one measures. We construct the
corresponding kernels and show that they obey a reciprocity relation. This
discharges the promise made in Remark~\ref{rem:amplitude-kernel} and yields the
pointwise form of the optical theorem. The second step,
\S\ref{subsec:inverse-exterior}, shows that equality of the scattering
amplitudes at a single positive energy forces the two scattering solutions of
every common incident wave to coincide identically outside the ball in which
the perturbations are confined. The third, \S\ref{subsec:schiffer}, propagates this
into the region where the two configurations differ and derives a contradiction
there. What the surface channel can and cannot do on its own is discussed
separately in \S\ref{subsec:surface-only}.

Throughout \S\ref{subsec:inverse-kernels} a single perturbed domain is fixed
and the notation of \S\S\ref{sec:setting}--\ref{sec:smatrix} is retained
unchanged. Two configurations are compared only from
\S\ref{subsec:inverse-exterior} onwards, and the notation for that comparison
is fixed in Definition~\ref{def:admissible-pair}.

\subsection{Directional Eigenfunctions and the Scattering Kernels}
\label{subsec:inverse-kernels}

Fix $E \in \Ereg \setminus \Sexc$ and an admissible cut-off $\eta$. We recall
from \cite[Definition 2.23]{PartI} the generalized eigenfunctions $\Psib$ and
$\Psis$ of the unperturbed operator, and from \cite[Definition 2.29]{PartI}
that the two components of $\Fzero(E)$ are integral transforms with those
eigenfunctions as kernels. It is convenient to index the two channels
uniformly. We let $\mu$ range over $\{\mathrm{b},\mathrm{s}\}$ and put
\begin{equation}
  \label{eq:channel-index}
  \Sch{\mathrm{b}} = \Sm,
  \quad
  \Sch{\mathrm{s}} = \Ssurf,
  \qquad
  \Psich{\mathrm{b}}(x,\theta) = \Psib(x,\kbulk\,\theta),
  \quad
  \Psich{\mathrm{s}}(x,\omega) = \Psis(x,\ksurf\,\omega),
\end{equation}
together with the positive constants
\begin{equation}
  \label{eq:kappa-constants}
  \kapch{\mathrm{b}}{E} = 2^{-1/2}\,E^{\frac{n-2}{4}},
  \qquad
  \kapch{\mathrm{s}}{E} = 2^{-1/2}\,(E+\alpha^2)^{\frac{n-3}{4}} ,
\end{equation}
so that \cite[Definition 2.29]{PartI} reads, in either channel,
\begin{equation}
  \label{eq:Fzero-kernels}
  \bigl( \Fzero(E)h \bigr)_\mu(\theta)
  = \kapch{\mu}{E} \int_{\Rnp}
      h(x)\,\overline{\Psich{\mu}(x,\theta)}\,dx,
  \qquad \theta \in \Sch{\mu} .
\end{equation}
When $E \in \Eband$ the index $\mu$ takes only the value $\mathrm{s}$, and
every statement below is to be read with the bulk channel deleted.

The perturbed counterpart of $\Psich{\mu}(\cdot,\theta)$ is produced by the
same transplantation that gave the scattering solution in
Theorem~\ref{thm:scattering-solution}, but one point has to be watched.
The function $\Psich{\mu}(\cdot,\theta)$ is \emph{not} an element of $\Bst$: it
is bounded and does not decay, so its shell averages grow like a power of the
radius, and this already on $\Rnp$. It is therefore not an incident wave in the
sense of Definition~\ref{def:incident}, and $\Phipert{\pm}{E}$ cannot be
applied to it. What does survive the passage to a single direction is the
transplanted source.

\begin{definition}[Directional eigenfunctions]
\label{def:directional}
Let $E \in \Ereg \setminus \Sexc$, let $\eta$ be an admissible cut-off, let
$\mu \in \{\mathrm{b},\mathrm{s}\}$ and $\theta \in \Sch{\mu}$. The
\emph{outgoing} and \emph{incoming directional eigenfunctions} of $\Hpert$ are
\begin{equation}
  \label{eq:directional-def}
  \Phidir{\pm}{\mu}(\cdot,\theta)
    := \eta\,\Psich{\mu}(\cdot,\theta)
       - \Rpert{E \pm i0}
         \bigl( [-\Delta,\eta]\,\Psich{\mu}(\cdot,\theta) \bigr) ,
\end{equation}
the first term being restricted to $\Om$ and the second understood through
Theorem~\ref{thm:lap-pert}(b).
\end{definition}

\begin{lemma}[Elementary properties of the directional eigenfunctions]
\label{lem:directional}
Let $E \in \Ereg \setminus \Sexc$ and $\mu \in \{\mathrm{b},\mathrm{s}\}$.
\begin{enumerate}
  \item[\textup{(a)}] For every $\theta \in \Sch{\mu}$ the source
  $[-\Delta,\eta]\Psich{\mu}(\cdot,\theta)$ is smooth, supported in the compact
  annulus $\mathcal{A}$ of \eqref{eq:annulus}, and belongs to $\Bspd{\Om}$.
  Hence \eqref{eq:directional-def} is well defined and
  \begin{equation}
    \label{eq:directional-affine}
    \Phidir{\pm}{\mu}(\cdot,\theta) - \eta\,\Psich{\mu}(\cdot,\theta)
      \;\in\; \Bstd{\Om} \cap H^2_{\mathrm{loc}}(\overline{\Om}) .
  \end{equation}
  \item[\textup{(b)}] $\Phidir{\pm}{\mu}(\cdot,\theta)$ belongs to
  $H^2_{\mathrm{loc}}(\overline{\Om})$ and solves
  \begin{equation}
    \label{eq:directional-equation}
    (-\Delta - E)\,\Phidir{\pm}{\mu}(\cdot,\theta) = 0 \ \text{ in } \Om,
    \qquad
    \gamma_1 \Phidir{\pm}{\mu} - \alpha\,\gamma_0\Phidir{\pm}{\mu} = 0
    \ \text{ on } \Gpert .
  \end{equation}
  \item[\textup{(c)}] $\Phidir{\pm}{\mu}(\cdot,\theta)$ does not depend on the
  choice of the admissible cut-off $\eta$.
  \item[\textup{(d)}] The maps
  $\theta \mapsto [-\Delta,\eta]\Psich{\mu}(\cdot,\theta)$, with values in
  $\Ltwo{\mathcal{A}}$ and hence in $\Bspd{\Om}$, and
  $\theta \mapsto \Phidir{\pm}{\mu}(\cdot,\theta)
   - \eta\,\Psich{\mu}(\cdot,\theta) \in \Bstd{\Om}$
  are continuous on $\Sch{\mu}$. Consequently, for every compactly supported
  $f \in \Ltwo{\Om}$ the function
  $\theta \mapsto \bigl( f, \Phidir{\pm}{\mu}(\cdot,\theta) \bigr)$ is
  continuous and bounded on $\Sch{\mu}$.
\end{enumerate}
\end{lemma}

\begin{proof}
(a) By \cite[equation (2.51)]{PartI} the function $\Psib(x,\xi)$ is a linear
combination, with coefficients depending on $\xi$ alone and of modulus at most
one, of the two exponentials $e^{ix\cdot\xi}$ and $e^{ix\cdot\xi^{*}}$ with
$\xi^{*}$ the reflection of $\xi$ in $\Gamma$. In the surface channel
$\Psis(x,\xipar) = (2\pi)^{-\frac{n-1}{2}}\,e^{i\xpar\cdot\xipar}\,\psizero(x_n)$ by
\cite[Definition 2.23]{PartI}. In both cases $\Psich{\mu}$ is smooth in $x$ and
bounded together with all its $x$-derivatives on $\Rnp$, uniformly in
$\theta \in \Sch{\mu}$. By \eqref{eq:commutator} the coefficients of
$[-\Delta,\eta]$ are supported in $\mathcal{A}$, which gives the first two
assertions. A bounded function supported in a compact subset of
$\overline{\Om}$ lies in $\Bspd{\Om}$, only finitely many dyadic blocks
contributing to the norm of Definition~\ref{def:AH-pert}. Membership of the
resolvent term in $\Bstd{\Om} \cap H^2_{\mathrm{loc}}(\overline{\Om})$ is then
Theorem~\ref{thm:lap-pert}(b) and (e).

(b) The cut-off $\eta$ is radial, so by Lemma~\ref{lem:gluing} the function
$\eta\,\Psich{\mu}(\cdot,\theta)$ satisfies the boundary condition on $\Gpert$
and
$(-\Delta - E)\bigl(\eta\Psich{\mu}\bigr) = [-\Delta,\eta]\Psich{\mu}$ in
$\Om$, the transplantation being valid because $\nabla\eta$ is supported
in $\mathcal{A} \subset \Omex{\Rzero}$, where $\Om$ and $\Rnp$ agree by
\eqref{eq:exterior-coincidence}. Subtracting Theorem~\ref{thm:lap-pert}(e),
applied to that same source, cancels the right-hand side and preserves the
boundary condition, which is \eqref{eq:directional-equation}.

(c) Let $\eta,\tilde\eta$ be admissible cut-offs, write $\Psi =
\Psich{\mu}(\cdot,\theta)$, and let $d$ be the difference of the two functions
\eqref{eq:directional-def}. Then
\begin{equation*}
  d = (\eta - \tilde\eta)\Psi
      - \Rpert{E \pm i0}\bigl( [-\Delta, \eta - \tilde\eta]\Psi \bigr),
\end{equation*}
and $d$ solves \eqref{eq:directional-equation} on all of $\Om$ by (b). Both
terms satisfy the outgoing, respectively incoming, radiation condition of
Definition~\ref{def:radiation-pert}. The first does because
$\eta - \tilde\eta$ vanishes outside a compact set, so that
$(\eta-\tilde\eta)\Psi \in \Bstzd{\Om}$, and that space is contained in the
class by Lemma~\ref{lem:radiation-elementary}(c). The second does by
Corollary~\ref{cor:outgoing-intrinsic}(b), the source
$[-\Delta,\eta-\tilde\eta]\Psi$ being compactly supported. The class is a
linear subspace, so $d$ belongs to it, and
Theorem~\ref{thm:radiation-equivalence} makes $d$ outgoing, respectively
incoming, in the sense of Definition~\ref{def:outgoing-pert}. Now
Theorem~\ref{thm:outgoing-uniqueness} gives $d \in \Bstzd{\Om}$ and then
$d = 0$, either directly when $E > 0$ or, when $E \in \Eband$, because $d$
would otherwise be an eigenfunction of $\Hpert$ at an energy
$E \notin \Sexc$.

(d) By the explicit formulas quoted in (a), $\Psich{\mu}(x,\theta)$ and its
first $x$-derivatives depend continuously on $\theta$, uniformly for $x$ in the
compact set $\mathcal{A}$. Hence
$\theta \mapsto [-\Delta,\eta]\Psich{\mu}(\cdot,\theta)$ is continuous into
$\Ltwo{\mathcal{A}}$, and therefore into $\Bspd{\Om}$, because on functions
supported in a fixed compact set the two norms are comparable, by the argument
of (a). The second map is that one composed with $-\Rpert{E \pm i0}$, which is
bounded from $\Bspd{\Om}$ to $\Bstd{\Om}$ by Theorem~\ref{thm:lap-pert}(b). For
the last assertion, split the pairing according to
\eqref{eq:directional-affine}. The term
$(f, \eta\,\Psich{\mu}(\cdot,\theta))$ is continuous in $\theta$ by dominated
convergence, $f$ being integrable and $\Psich{\mu}$ bounded. The remaining term
is bounded by $\BnormO{f}$ times the $\Bstd{\Om}$ norm of
\eqref{eq:directional-affine} through \eqref{eq:AH-pairing}, and is continuous
because that norm is.
\end{proof}

We can now identify the kernels of the perturbed transforms. The computation
uses only the transplanted source formula, the antisymmetry of the
transplantation commutator, and the adjoint relation between the two boundary
values of the resolvent.

\begin{theorem}[Kernels of the perturbed transforms and reciprocity]
\label{thm:kernel-reciprocity}
Let $E \in \Ereg \setminus \Sexc$, let $\mu \in \{\mathrm{b},\mathrm{s}\}$ and
let $f \in \Ltwo{\Om}$ have compact support. Then
\begin{equation}
  \label{eq:kernel-reciprocity}
  \bigl( \Fpert{\pm}{E}f \bigr)_\mu(\theta)
  = \kapch{\mu}{E}
    \int_{\Om} f(x)\,\overline{\Phidir{\pm}{\mu}(x,\theta)}\,dx ,
  \qquad \theta \in \Sch{\mu} ,
\end{equation}
and both sides are continuous functions of $\theta$. In particular
$\Fpert{-}{E}f = 0$ in $\Echan{E}$ if and only if $f$ is orthogonal in
$\Ltwo{\Om}$ to $\Phidir{-}{\mu}(\cdot,\theta)$ for every channel $\mu$ and
every direction $\theta \in \Sch{\mu}$.
\end{theorem}

\begin{proof}
We prove the statement with the upper signs. The lower ones are obtained by
exchanging $E + i0$ and $E - i0$ throughout. Write
$\Psi = \Psich{\mu}(\cdot,\theta)$. Since $f$ is compactly supported it lies in
$\Bspd{\Om}$ by \eqref{eq:AH-inclusions-pert}, so
Theorem~\ref{thm:Fpert-explicit} applies and gives, in the notation
\eqref{eq:Fzero-kernels},
\begin{equation*}
  \bigl( \Fpert{+}{E}f \bigr)_\mu(\theta)
  = \kapch{\mu}{E}\,
    \bigl( \eta f + [-\Delta,\eta]\,\Rpert{E-i0}f,\ \Psi
    \bigr)_{\Ltwo{\Rnp}} ,
\end{equation*}
the argument of $\Fzero(E)$ being extended by zero to $\Rnp$. This
extension is justified because both $\eta f$ and $[-\Delta,\eta]\Rpert{E-i0}f$ are
supported in $\Om \setminus \Bfull{R_1}$, a subset of $\Rnp$ by
\eqref{eq:exterior-coincidence}. For the same reason the first term may be read
on $\Om$,
\begin{equation*}
  \bigl( \eta f, \Psi \bigr)_{\Ltwo{\Rnp}}
  = \bigl( f, \eta\Psi \bigr)_{\Ltwo{\Om}} ,
\end{equation*}
the cut-off being real valued. For the second term,
$\Rpert{E-i0}f \in H^2_{\mathrm{loc}}(\overline{\Om})$ by
Theorem~\ref{thm:lap-pert}(e) and
$\Psi \in H^2_{\mathrm{loc}}(\overline{\Rnp})$, so
Lemma~\ref{lem:commutator-antisym} applies and, followed by
Theorem~\ref{thm:lap-pert}(d), yields
\begin{equation*}
  \bigl( [-\Delta,\eta]\,\Rpert{E-i0}f,\ \Psi \bigr)
  = -\bigl( \Rpert{E-i0}f,\ [-\Delta,\eta]\Psi \bigr)
  = -\bigl( f,\ \Rpert{E+i0}\bigl([-\Delta,\eta]\Psi\bigr) \bigr) ,
\end{equation*}
the last step being applicable because
$[-\Delta,\eta]\Psi \in \Bspd{\Om}$ by Lemma~\ref{lem:directional}(a). Adding
the two contributions gives
\begin{equation*}
  \bigl( \Fpert{+}{E}f \bigr)_\mu(\theta)
  = \kapch{\mu}{E}\,
    \Bigl( f,\ \eta\Psi
      - \Rpert{E+i0}\bigl([-\Delta,\eta]\Psi\bigr) \Bigr)_{\Ltwo{\Om}} ,
\end{equation*}
which is \eqref{eq:kernel-reciprocity} by Definition~\ref{def:directional}.

Both sides are continuous in $\theta$: the right-hand side by
Lemma~\ref{lem:directional}(d), and the left-hand side because
$\eta f + [-\Delta,\eta]\Rpert{E-i0}f$ is one fixed compactly supported
element of $\Bsp$, to which \eqref{eq:Fzero-kernels} applies with the same
dominated convergence argument. For the final assertion, an element of
$\Echan{E}$ represented by a continuous function vanishes in $\Echan{E}$ if and
only if that function vanishes at every point.
\end{proof}

The scattering amplitude now acquires an integral kernel, and the four blocks
announced in Remark~\ref{rem:amplitude-kernel} appear.

\begin{corollary}[The kernel of the scattering amplitude]
\label{cor:amplitude-kernel-formula}
Let $E \in \Ereg \setminus \Sexc$, let $\mu,\nu \in \{\mathrm{b},\mathrm{s}\}$
and let $\eta$ be an admissible cut-off. Put
\begin{equation}
  \label{eq:amplitude-kernel}
  \Akern{\mu\nu}{E}{\theta}{\theta'}
  := -\,\kapch{\mu}{E}\,\kapch{\nu}{E}
     \int_{\Om}
       \bigl( [-\Delta,\eta]\,\Psich{\nu} \bigr)(x,\theta')\,
       \overline{\Phidir{-}{\mu}(x,\theta)}\,dx
\end{equation}
for $(\theta,\theta') \in \Sch{\mu} \times \Sch{\nu}$. Then
$\Akern{\mu\nu}{E}{\cdot}{\cdot}$ is continuous, independent of the choice of
$\eta$, and is the integral kernel of the corresponding block of the scattering
amplitude: for every $g \in \Echan{E}$,
\begin{equation}
  \label{eq:amplitude-kernel-action}
  \bigl( \Ampl{E}g \bigr)_\mu(\theta)
  = \sum_{\nu} \int_{\Sch{\nu}}
      \Akern{\mu\nu}{E}{\theta}{\theta'}\,g_\nu(\theta')\,d\theta' .
\end{equation}
In particular $\Ampl{E}$ is a Hilbert--Schmidt operator on $\Echan{E}$, and the
kernel of its adjoint is
\begin{equation}
  \label{eq:amplitude-kernel-adjoint}
  \bigl( \Ampl{E}^{*} \bigr)_{\mu\nu}(\theta,\theta')
  = \overline{\Akern{\nu\mu}{E}{\theta'}{\theta}} .
\end{equation}
\end{corollary}

\begin{proof}
By \cite[Proposition 2.34]{PartI} the adjoint of $\Fzero(E)$ acts by
\begin{equation}
  \label{eq:Fzero-adjoint-kernel}
  \Fzero(E)^{*}g
  = \sum_\nu \kapch{\nu}{E} \int_{\Sch{\nu}}
      g_\nu(\theta')\,\Psich{\nu}(\cdot,\theta')\,d\theta' .
\end{equation}
Apply $[-\Delta,\eta]$, whose coefficients are supported in $\mathcal{A}$. The
integrand is continuous in $\theta'$ with values in $\Ltwo{\mathcal{A}}$ by
Lemma~\ref{lem:directional}(d), and $\Sch{\nu}$ is compact, so
$[-\Delta,\eta]\Fzero(E)^{*}g$ is a compactly supported element of
$\Ltwo{\Om}$ given by an absolutely convergent $\Ltwo{\mathcal{A}}$-valued
integral. Theorem~\ref{thm:kernel-reciprocity} therefore applies to it, and by
Definition~\ref{def:scattering-amplitude}
\begin{equation*}
  \begin{aligned}
  \bigl( \Ampl{E}g \bigr)_\mu(\theta)
  &= -\bigl( \Fpert{-}{E}
      \bigl([-\Delta,\eta]\Fzero(E)^{*}g\bigr) \bigr)_\mu(\theta)
  \\
  &= -\,\kapch{\mu}{E}
    \Bigl( [-\Delta,\eta]\Fzero(E)^{*}g,\
      \Phidir{-}{\mu}(\cdot,\theta) \Bigr)_{\Ltwo{\Om}} ,
  \end{aligned}
\end{equation*}
and inserting \eqref{eq:Fzero-adjoint-kernel} and interchanging the integral
with the pairing gives \eqref{eq:amplitude-kernel-action} with the kernel
\eqref{eq:amplitude-kernel}.

Continuity of the kernel in $\theta$ is Lemma~\ref{lem:directional}(d) and
continuity in $\theta'$ is the continuity, established there as well, of
$\theta' \mapsto [-\Delta,\eta]\Psich{\nu}(\cdot,\theta')$ in
$\Ltwo{\mathcal{A}}$. The pairing is bounded and bilinear in the two factors,
so it is jointly continuous. Independence of $\eta$ holds because the
left-hand side of \eqref{eq:amplitude-kernel-action} is independent of $\eta$
by Theorem~\ref{thm:scattering-amplitude}(a), and a continuous kernel is
determined by the operator it represents. The sets $\Sch{\mu}$ are compact, so a continuous kernel is square integrable and $\Ampl{E}$ is
Hilbert--Schmidt. Finally \eqref{eq:amplitude-kernel-adjoint} is the usual
formula for the kernel of the adjoint of an integral operator.
\end{proof}

\begin{theorem}[The optical theorem in pointwise form]
\label{thm:optical-pointwise}
Let $E \in \Ereg \setminus \Sexc$, let $\nu \in \{\mathrm{b},\mathrm{s}\}$ and
let $\theta_0 \in \Sch{\nu}$. Then
\begin{equation}
  \label{eq:optical-pointwise}
  \im \Akern{\nu\nu}{E}{\theta_0}{\theta_0}
  = \pi \sum_{\mu} \int_{\Sch{\mu}}
      \bigl| \Akern{\mu\nu}{E}{\theta}{\theta_0} \bigr|^2 \, d\theta .
\end{equation}
For $E > 0$ and an incident bulk direction this reads
\begin{equation}
  \label{eq:optical-pointwise-bulk}
  \im \Akern{\mathrm{bb}}{E}{\theta_0}{\theta_0}
  = \pi \int_{\Sm}
      \bigl| \Akern{\mathrm{bb}}{E}{\theta}{\theta_0} \bigr|^2 d\theta
    + \pi \int_{\Ssurf}
      \bigl| \Akern{\mathrm{sb}}{E}{\omega}{\theta_0} \bigr|^2 d\omega .
\end{equation}
\end{theorem}

\begin{proof}
Theorem~\ref{thm:optical} gives the operator identity
\eqref{eq:optical-operator} on $\Echan{E}$. By
Corollary~\ref{cor:amplitude-kernel-formula} both of its sides are integral
operators with continuous kernels. On the left the kernel is
$\tfrac{1}{2i}\bigl( \Akern{\mu\nu}{E}{\theta}{\theta'}
- \overline{\Akern{\nu\mu}{E}{\theta'}{\theta}} \bigr)$ by
\eqref{eq:amplitude-kernel-adjoint}. On the right, a product of two
Hilbert--Schmidt operators with continuous kernels on a compact space again has
a continuous kernel, namely
\begin{equation*}
  \pi \sum_\lambda \int_{\Sch{\lambda}}
    \overline{\Akern{\lambda\mu}{E}{\theta''}{\theta}}\,
    \Akern{\lambda\nu}{E}{\theta''}{\theta'}\,d\theta'' .
\end{equation*}
Two integral operators with continuous kernels coincide only if their kernels
coincide at every point. Setting $\mu = \nu$ and $\theta = \theta' = \theta_0$
turns the left-hand kernel into $\im \Akern{\nu\nu}{E}{\theta_0}{\theta_0}$ and
the right-hand one into
$\pi\sum_\lambda \int_{\Sch{\lambda}}
|\Akern{\lambda\nu}{E}{\theta''}{\theta_0}|^2 d\theta''$, which is
\eqref{eq:optical-pointwise}. On the band the bulk channel is absent, and for
$E > 0$ and $\nu = \mathrm{b}$ the sum has the two terms displayed in
\eqref{eq:optical-pointwise-bulk}.
\end{proof}

\begin{remark}[The scope of the kernels]
\label{rem:kernel-scope}
Relation \eqref{eq:kernel-reciprocity} is the analogue, for the impedance
half-space, of the mixed reciprocity relation of obstacle scattering, in which
the far field of the Green function, taken in the observation variable, is
proportional to the total field for plane-wave incidence evaluated at the
source point \cite[\S3.3]{ColtonKress2013}. The classical proof applies Green's
second identity on a large sphere. The proof above uses instead
Lemma~\ref{lem:commutator-antisym}, which involves neither the equation nor
the boundary condition, so both channels are handled by one argument and no
asymptotic expansion of a fundamental solution is needed.

The restriction to compactly supported $f$ in
Theorem~\ref{thm:kernel-reciprocity} cannot be lifted, since the kernel grows
like $\eta\,\Psich{\mu}(\cdot,\theta) \notin \Bst$. The uniqueness argument
of \S\S\ref{subsec:inverse-exterior}--\ref{subsec:schiffer} does not use the
double kernel \eqref{eq:amplitude-kernel}. It uses
Definition~\ref{def:directional} and the superposition formula
\eqref{eq:superposition} below. The kernel is constructed because it is what
an experiment records and because \eqref{eq:optical-pointwise} is the form in
which the optical theorem is usually stated. The incoming eigenfunction
$\Phidir{-}{\mu}$ appears in \eqref{eq:amplitude-kernel} for the reason given
in Remark~\ref{rem:Fpert-crossing}.
\end{remark}

\subsection{From the Scattering Data to the Common Exterior}
\label{subsec:inverse-exterior}

We now compare two configurations. The following convention fixes the notation
for the rest of the section.

\begin{definition}[Admissible pair]
\label{def:admissible-pair}
An \emph{admissible pair} is a pair $\Om_1, \Om_2$ of domains, each satisfying
Assumption~\ref{assump:geometry}, with one and the same radius $\Rzero$ in
(G2), and with the same impedance $\alpha > 0$ and the same dimension
$n \ge 3$. Every object attached to $\Om_j$ carries the index $j$: thus
$\Hpertj{j}$ is the perturbed Robin Laplacian of $\Om_j$, $\Rpertj{j}{z}$ its
resolvent, $\Phipertj{j}{\pm}{E}$ the maps of Definition~\ref{def:Phipert},
$\Nsolpertj{j}{E}$ the space of Definition~\ref{def:Nsolpert}, $\Sexcj{j}$ the
exceptional set of \eqref{eq:exceptional-set}, and $\Amplj{j}{E}$ the
scattering amplitude of Definition~\ref{def:scattering-amplitude}.
\end{definition}

Requiring a common $\Rzero$ is no restriction. Hypothesis~(G2) persists when
$\Rzero$ is enlarged, so one replaces the two radii by their maximum. With a
common $\Rzero$ the two domains agree outside $\Bfull{\Rzero}$,
\begin{equation}
  \label{eq:common-exterior}
  \Om_1 \setminus \overline{\Bfull{\Rzero}}
  = \Om_2 \setminus \overline{\Bfull{\Rzero}}
  = \Rnp \setminus \overline{\Bfull{\Rzero}} ,
\end{equation}
so that the notation $\Omex{R}$, $R \ge \Rzero$, is unambiguous and we use it
without an index. What matters more than the equality of sets is that the whole
apparatus of Part~I lives on $\Omex{\Rzero}$ for both problems at once: the
channel space $\Echan{E}$, the space $\Nsol{E}$ of incident waves of
Definition~\ref{def:incident}, the four model waves
\eqref{eq:model-waves-pert}, and the far-field data $\yin$, $\yout$ of
\cite[Corollary 2.70]{PartI} are literally the same objects in the two
problems. It is therefore meaningful to compare $\Amplj{1}{E}$ with
$\Amplj{2}{E}$, and to send one and the same incident wave into both.

\begin{theorem}[Equality of the fields outside the perturbation]
\label{thm:exterior-equality}
Let $\Om_1,\Om_2$ be an admissible pair, let $E > 0$, and suppose that
\begin{equation}
  \label{eq:equal-data}
  \Amplj{1}{E} = \Amplj{2}{E}
  \qquad \text{on } \Echan{E} .
\end{equation}
Then, for every incident wave $\phi \in \Nsol{E}$,
\begin{equation}
  \label{eq:exterior-equality}
  \Phipertj{1}{+}{E}\phi = \Phipertj{2}{+}{E}\phi
  \qquad \text{in } \Omex{\Rzero} ,
\end{equation}
and the same holds with $\Phipertj{j}{-}{E}$ in place of
$\Phipertj{j}{+}{E}$.
\end{theorem}

\begin{proof}
Both exceptional sets are contained in $\Eband$ by
\eqref{eq:exceptional-set}, so at the energy $E > 0$ every result of
\S\S\ref{sec:lap-pert}--\ref{sec:smatrix} is available for either domain with
no set to avoid.

Let $\phi \in \Nsol{E}$, put $g = (\Fzero(E)^{*})^{-1}\phi \in \Echan{E}$,
which is well defined by \eqref{eq:Fzero-bijection}, and set
$u_j = \Phipertj{j}{+}{E}\phi$ and $w = u_1 - u_2$, a function on
$\Omex{\Rzero}$. Fix an admissible cut-off $\eta$. By
Theorem~\ref{thm:scattering-amplitude}(b) each $\eta u_j$ is asymptotically
equal to a sum of the four model waves. Its incoming amplitudes
\eqref{eq:incoming-amplitudes} are determined by $\phi$ alone, hence are the
same for $j = 1$ and $j = 2$. Its outgoing amplitudes
\eqref{eq:outgoing-amplitudes} are determined by $\phi$ together with
$\Amplj{j}{E}g$, and by \eqref{eq:equal-data} these agree as well. All four
amplitudes of $\eta u_1$ and of $\eta u_2$ therefore coincide, and by
Lemma~\ref{lem:radiation-elementary}(b),
\begin{equation}
  \label{eq:w-Bstz}
  \eta\,w \;\simeq\; 0,
  \qquad\text{that is,}\qquad
  \eta\,w \in \Bstzd{\Om_1} ,
\end{equation}
with $\eta w$ extended by zero to $\Om_1$. By Lemma~\ref{lem:AH-exterior} the
same statement holds verbatim with $\Om_2$ in place of $\Om_1$, the two norms
agreeing on functions supported outside $\Bfull{\Rzero}$.

Since $\eta \equiv 1$ outside $\Bfull{R_2}$ and both $u_j$ solve the same
equation with the same boundary condition on $\Omex{\Rzero}$ by
Theorem~\ref{thm:scattering-solution}, the function $\eta w$ solves the
exterior Robin problem at energy $E$ outside $\Bfull{R_2}$ in the sense of
Definition~\ref{def:exterior-problem}. As $E > 0$,
Theorem~\ref{thm:rellich-positive} yields a radius $R' \ge R_2$ with
$w = 0$ in $\Omex{R'}$.

It remains to remove $R'$. By Theorem~\ref{thm:scattering-solution},
$w \in H^2_{\mathrm{loc}}(\Omex{\Rzero})$ and $(-\Delta - E)w = 0$ in
$\Omex{\Rzero}$. The set $\Omex{\Rzero}$ is open and connected, and $w$
vanishes on its nonempty open subset $\Omex{R'}$.
Proposition~\ref{prop:ucp} therefore gives $w = 0$ on all of $\Omex{\Rzero}$,
which is \eqref{eq:exterior-equality}. The incoming case follows by complex
conjugation, exactly as in the proof of
Theorem~\ref{thm:characterization-pert}.
\end{proof}

\begin{corollary}[The Cauchy data on an artificial boundary]
\label{cor:cauchy-data}
Under the hypotheses of Theorem~\ref{thm:exterior-equality}, let $R \ge \Rzero$
and let $\Sart{R} = \{ x \in \Rnp : |x| = R \}$ be the corresponding artificial
boundary. For $j = 1,2$ put
\begin{equation}
  \label{eq:cauchy-set}
  \Cdata_j(E)
  = \bigl\{\, \bigl( u|_{\Sart{R}},\ \partial_r u|_{\Sart{R}} \bigr)
     \;:\; u \in \Nsolpertj{j}{E} \,\bigr\}
  \;\subset\; H^{3/2}(\Sart{R}) \times H^{1/2}(\Sart{R}) .
\end{equation}
Then $\Cdata_1(E) = \Cdata_2(E)$. In particular the two sets of Dirichlet
traces coincide, and a given Dirichlet trace is accompanied by the same set of
Neumann traces in the two problems.
\end{corollary}

\begin{proof}
Let $u \in \Nsolpertj{1}{E}$. By Theorem~\ref{thm:characterization-pert} there
is exactly one $\phi \in \Nsol{E}$ with $u = \Phipertj{1}{+}{E}\phi$, and
Theorem~\ref{thm:exterior-equality} gives $u = \Phipertj{2}{+}{E}\phi$ on
$\Omex{\Rzero}$, the right-hand side being an element of $\Nsolpertj{2}{E}$.
The sphere $\Sart{R}$ lies in $\overline{\Omex{\Rzero}}$ for $R \ge \Rzero$,
and both functions belong to
$H^2_{\mathrm{loc}}(\overline{\Omex{\Rzero}})$. Their traces and radial
derivatives on $\Sart{R}$, taken from the exterior side, therefore agree, and
the stated regularity is the trace theorem applied on the collar
$\{R \le |x| < R+\delta\} \cap \Rnp$ for small $\delta > 0$. We take the
collar on one side because at $R = \Rzero$ the two functions need not be
defined on the same set inside $\Bfull{\Rzero}$. Hence
$\Cdata_1(E) \subset \Cdata_2(E)$, and the reverse inclusion is the same
argument with the indices exchanged.
\end{proof}

\begin{remark}[The Dirichlet-to-Neumann map]
\label{rem:dn-map}
The first stage of an inverse scattering argument is often phrased as the
statement that the far-field data determine the Dirichlet-to-Neumann map
$\DN_j$ of the region enclosed by an artificial boundary.
Corollary~\ref{cor:cauchy-data} stops short of that. To define $\DN_j$ one
needs the interior problem on $\Om_j \cap \Bfull{R}$, with the Robin condition
on $\partial\Om_j \cap \Bfull{R}$ and Dirichlet data on $\Sart{R}$, to be
uniquely solvable at the energy $E$, which is a condition on $E$ and on the
unknown geometry. Granting it, Corollary~\ref{cor:cauchy-data} says that
$\DN_1$ and $\DN_2$ agree on the Dirichlet traces of perturbed solutions, and
$\DN_1 = \DN_2$ would require these traces to be dense in
$H^{1/2}(\Sart{R})$. We neither prove nor use that density. The argument of
\S\ref{subsec:schiffer} uses Theorem~\ref{thm:exterior-equality} directly and
needs no hypothesis on $E$.
\end{remark}

\subsection{Uniqueness of the Perturbation}
\label{subsec:schiffer}

The argument of this subsection is the one introduced by Schiffer for the
sound-soft obstacle. If two configurations produce the same data, the region in
which they differ carries too many solutions of an overdetermined boundary
value problem. Two features of the present setting have to be accommodated.
First, the boundary condition is of impedance type. The overdetermined problem
on the difference region is therefore not a Dirichlet problem but a mixed one,
and the impedance $\alpha$ appears in it with \emph{opposite signs} on the two
parts of the boundary. Secondly, the perturbation is not an obstacle but an
arbitrary deformation of the domain. The difference region may accordingly be a
dent in the boundary, a detached obstacle, or a combination of the two. Neither
feature requires a separate treatment. Indeed, the argument never asks what
kind of perturbation produced the difference region.

Two hypotheses on the pair are needed, and we isolate them.

\begin{assumption}
\label{assump:pair}
The admissible pair $\Om_1,\Om_2$ of Definition~\ref{def:admissible-pair}
satisfies
\begin{enumerate}
  \item[\textup{(U1)}] the open set $\Omcap := \Om_1 \cap \Om_2$ is connected,
  and
  \item[\textup{(U2)}] every connected component of
  $\Om_1 \setminus \overline{\Om_2}$ and of $\Om_2 \setminus \overline{\Om_1}$
  has Lipschitz boundary.
\end{enumerate}
\end{assumption}

Both are conditions on the pair rather than on either domain separately, which
is unusual in a uniqueness theorem. Remark~\ref{rem:hypotheses} discusses what
they exclude and where each is used. For the moment we record only that
$\Omcap$ contains $\Omex{\Rzero}$ by \eqref{eq:common-exterior}. Hence
\textup{(U1)} says exactly that the region common to the two configurations is
reachable, in one piece, from infinity.

One convention removes an apparent restriction on the energy. Every energy
considered in this subsection is positive, and both exceptional sets are
contained in $\Eband$ by \eqref{eq:exceptional-set}. Every $E > 0$ therefore lies in
$\Ereg \setminus \Sexcj{j}$ for $j = 1,2$, so that the directional
eigenfunctions of Definition~\ref{def:directional} are available in both
configurations with no further hypothesis.

The first step transports Theorem~\ref{thm:exterior-equality} from
$\Omex{\Rzero}$ to all of $\Omcap$, and at the same time from the scattering
solutions to the individual directions.

\begin{lemma}[Transfer to the common region]
\label{lem:common-region}
Let $\Om_1,\Om_2$ be an admissible pair satisfying \textup{(U1)}, let $E > 0$
and assume \eqref{eq:equal-data}. Write $\Phidirj{j}{-}{\mu}$ for the incoming
directional eigenfunction \eqref{eq:directional-def} of $\Hpertj{j}$. Then
\begin{equation}
  \label{eq:directional-equality}
  \Phidirj{1}{-}{\mu}(x,\theta) = \Phidirj{2}{-}{\mu}(x,\theta)
  \qquad \text{for all } x \in \Omcap,\
  \mu \in \{\mathrm{b},\mathrm{s}\},\ \theta \in \Sch{\mu} .
\end{equation}
\end{lemma}

\begin{proof}
Let $\phi \in \Nsol{E}$ and put $w = \Phipertj{1}{-}{E}\phi -
\Phipertj{2}{-}{E}\phi$, a function on $\Omcap$. Since $\Omcap \subset \Om_j$
for both $j$, Theorem~\ref{thm:scattering-solution} gives
$w \in H^2_{\mathrm{loc}}(\Omcap)$ and $(-\Delta - E)w = 0$ in $\Omcap$. The
set $\Omcap$ is open and connected by \textup{(U1)} and contains the nonempty
open set $\Omex{\Rzero}$ by \eqref{eq:common-exterior}, and $w$ vanishes there
by Theorem~\ref{thm:exterior-equality}. Proposition~\ref{prop:ucp} therefore
gives
\begin{equation}
  \label{eq:w-vanishes-common}
  \Phipertj{1}{-}{E}\phi = \Phipertj{2}{-}{E}\phi
  \qquad \text{in } \Omcap, \quad \phi \in \Nsol{E} .
\end{equation}

It remains to pass from the scattering solutions to the individual
directions. Let $g \in \Echan{E}$ and put $\phi = \Fzero(E)^{*}g$, which by
\eqref{eq:Fzero-bijection} is the general element of $\Nsol{E}$. We claim
that, for $j = 1,2$ and every $R > 0$,
\begin{equation}
  \label{eq:superposition}
  \Phipertj{j}{-}{E}\Fzero(E)^{*}g
  = \sum_\nu \kapch{\nu}{E} \int_{\Sch{\nu}}
      g_\nu(\theta')\,\Phidirj{j}{-}{\nu}(\cdot,\theta')\,d\theta'
  \qquad \text{in } \Ltwo{\Om_j \cap \Bfull{R}} ,
\end{equation}
the integrals being Bochner integrals with values in
$\Ltwo{\Om_j \cap \Bfull{R}}$.

It suffices to take $R > R_2$. By the formulas quoted in the proof of
Lemma~\ref{lem:directional}(a), the map
$\theta' \mapsto \Psich{\nu}(\cdot,\theta')$ is continuous and bounded on
$\Sch{\nu}$ with values in $H^1(\Rnp \cap \Bfull{R})$. Since $\Sch{\nu}$ has
finite measure, $g_\nu \in \Ltwo{\Sch{\nu}} \subset L^1(\Sch{\nu})$, so
\eqref{eq:Fzero-adjoint-kernel} is a Bochner integral in
$H^1(\Rnp \cap \Bfull{R})$. The operators
\begin{equation*}
  \eta : \Ltwo{\Rnp \cap \Bfull{R}} \to \Ltwo{\Om_j \cap \Bfull{R}},
  \qquad
  [-\Delta,\eta] : H^1(\Rnp \cap \Bfull{R}) \to \Bspd{\Om_j},
\end{equation*}
$\Rpertj{j}{E-i0} : \Bspd{\Om_j} \to \Bstd{\Om_j}$ and the restriction
$\Bstd{\Om_j} \to \Ltwo{\Om_j \cap \Bfull{R}}$ are bounded: the first because
$\eta$ vanishes on $\Bfull{R_1}$ and $\Om_j$ coincides with $\Rnp$ outside
$\Bfull{\Rzero}$ by \eqref{eq:common-exterior}, the second because its
coefficients are supported in $\mathcal{A} \subset \Bfull{R}$ (see the proof of
Lemma~\ref{lem:directional}(a)), the third by Theorem~\ref{thm:lap-pert}(b),
and the fourth by Definition~\ref{def:AH-pert}. Bounded operators commute with
Bochner integrals, so applying them to \eqref{eq:Fzero-adjoint-kernel} and
inserting the result into the incoming counterpart of
\eqref{eq:scattering-solution} gives \eqref{eq:superposition}.

Fix $\mu \in \{\mathrm{b},\mathrm{s}\}$ and take $g$ with $g_\nu = 0$ for
$\nu \neq \mu$. Put
$D(\theta') = \bigl[ \Phidirj{1}{-}{\mu} - \Phidirj{2}{-}{\mu}
\bigr](\cdot,\theta')$ on $\Omcap$. Subtracting the two identities
\eqref{eq:superposition} and using \eqref{eq:w-vanishes-common} gives
\begin{equation*}
  \kapch{\mu}{E} \int_{\Sch{\mu}} g_\mu(\theta')\,D(\theta')\,d\theta' = 0
  \qquad \text{in } \Ltwo{\Omcap \cap \Bfull{R}}
\end{equation*}
for every $g_\mu \in \Ltwo{\Sch{\mu}}$ and every $R > 0$. The radius $\Rzero$
is common to the two domains, so one admissible cut-off serves both, and the
terms $\eta\Psich{\mu}(\cdot,\theta')$ of Definition~\ref{def:directional}
cancel in $D(\theta')$. Hence $D$ is continuous on $\Sch{\mu}$ with values in
$\Ltwo{\Omcap \cap \Bfull{R}}$ by Lemma~\ref{lem:directional}(d) and
Definition~\ref{def:AH-pert}. Let $h \in \Ltwo{\Omcap \cap \Bfull{R}}$ and
choose $g_\mu(\theta') = \overline{( D(\theta'), h )}$. Extended by zero, $h$ is
a compactly supported element of $\Ltwo{\Om_j}$ for both $j$, so $g_\mu$ is
continuous and bounded on $\Sch{\mu}$ by the last assertion of
Lemma~\ref{lem:directional}(d), and in particular
$g_\mu \in \Ltwo{\Sch{\mu}}$. Taking the inner product of the last identity with $h$ gives
\begin{equation*}
  \kapch{\mu}{E} \int_{\Sch{\mu}}
    \bigl| ( D(\theta'), h ) \bigr|^2\,d\theta' = 0 .
\end{equation*}
Since $\kapch{\mu}{E} > 0$ and the integrand is continuous,
$( D(\theta'), h ) = 0$ for every $\theta' \in \Sch{\mu}$, and $h$ being
arbitrary, $D(\theta') = 0$ in $\Ltwo{\Omcap \cap \Bfull{R}}$. As $R$ is
arbitrary, $\Phidirj{1}{-}{\mu}(\cdot,\theta) = \Phidirj{2}{-}{\mu}(\cdot,\theta)$
almost everywhere in $\Omcap$ for every $\theta \in \Sch{\mu}$. Both functions
belong to $H^2_{\mathrm{loc}}(\Omcap)$ and solve $(-\Delta - E)u = 0$ in
$\Omcap$ by Lemma~\ref{lem:directional}(b), so they are smooth in $\Omcap$ by
interior elliptic regularity, and \eqref{eq:directional-equality} holds at
every point.
\end{proof}

The second step is the topological one. It says that if the two domains differ
at all, then they differ somewhere that can be seen from $\Omcap$. The proof is
the observation that a subset of a connected set which is both open and closed
in it is everything.

\begin{lemma}[Accessibility of a difference]
\label{lem:accessibility}
Let $\Om_1,\Om_2$ be an admissible pair satisfying \textup{(U1)} and suppose
$\Om_1 \ne \Om_2$. Then, after exchanging the two indices if necessary, there
are a point $p \in \Om_1 \cap \partial\Om_2$ and a radius $\delta > 0$ such
that $\Bfull{\delta}(p) \subset \Om_1$, such that
$\Bfull{\delta}(p) \setminus \partial\Om_2$ has exactly the two components
\begin{equation}
  \label{eq:local-splitting}
  \Bfull{\delta}(p) \cap \Om_2 \;\subset\; \Omcap,
  \qquad
  \Bfull{\delta}(p) \setminus \overline{\Om_2}
    \;\subset\; \Om_1 \setminus \overline{\Om_2} ,
\end{equation}
and such that $\Bfull{\delta}(p) \cap \partial\Om_2$ is a relatively open
$C^2$ piece of $\partial\Om_2$.
\end{lemma}

\begin{proof}
Suppose first that $\partial\Omcap \subset \partial\Om_1$. Since $\Om_1$ is
open, $\partial\Om_1 \cap \Om_1 = \emptyset$, so $\partial\Omcap \cap \Om_1 =
\emptyset$ and the nonempty set $\Omcap$ is both open and relatively closed in
$\Om_1$. The domain $\Om_1$ is connected by Assumption~\ref{assump:geometry},
so this forces $\Omcap = \Om_1$. The same argument with the indices exchanged
shows that $\partial\Omcap \subset \partial\Om_2$ forces $\Omcap = \Om_2$.
Hence, if $\partial\Omcap \subset \partial\Om_1 \cap \partial\Om_2$, then
$\Om_1 = \Omcap = \Om_2$, contrary to hypothesis.

There is therefore a point $p \in \partial\Omcap$ lying outside
$\partial\Om_1 \cap \partial\Om_2$. Now
$\partial\Omcap \subset \overline{\Om_1} \cap \overline{\Om_2}$ and
$p \notin \Omcap$, so $p$ lies in $\partial\Om_1 \cup \partial\Om_2$. Exchange
the indices if necessary so that $p \notin \partial\Om_1$. Then
$p \in \overline{\Om_1} \setminus \partial\Om_1 = \Om_1$, and consequently
$p \in \partial\Om_2$. Thus $p \in \Om_1 \cap \partial\Om_2$.

Choose $\delta > 0$ so small that $\Bfull{\delta}(p) \subset \Om_1$ and that,
by hypothesis (G1) of Assumption~\ref{assump:geometry},
$\Bfull{\delta}(p) \cap \partial\Om_2$ is a connected relatively open $C^2$
piece of $\partial\Om_2$ splitting $\Bfull{\delta}(p)$ into the two connected
open sets $\Bfull{\delta}(p) \cap \Om_2$ and
$\Bfull{\delta}(p) \setminus \overline{\Om_2}$, both nonempty. The first is
contained in $\Om_1 \cap \Om_2 = \Omcap$ and the second in
$\Om_1 \setminus \overline{\Om_2}$, which is \eqref{eq:local-splitting}.
\end{proof}

The third step is where the impedance condition enters. On one part of the
boundary of the difference region, the field of the first configuration
satisfies its own boundary condition. On the other part it satisfies the
boundary condition of the \emph{second} configuration, because the two fields
agree on $\Omcap$. There the normal is reversed, and with it the sign of
$\alpha$.

\begin{lemma}[The overdetermined problem on a difference region]
\label{lem:mixed-problem}
Let the hypotheses of Lemma~\ref{lem:common-region} hold, assume
$\Om_1 \ne \Om_2$, and let $p,\delta$ be as in
Lemma~\ref{lem:accessibility}. Let $\Odiff$ be the connected component of
$\Om_1 \setminus \overline{\Om_2}$ containing
$\Bfull{\delta}(p) \setminus \overline{\Om_2}$, and put
\begin{equation}
  \label{eq:Gfar-Gacc}
  \Gfar = \partial\Odiff \cap \partial\Om_1,
  \qquad
  \Gacc = \partial\Odiff \setminus \partial\Om_1 .
\end{equation}
Then $\Odiff$ is a nonempty bounded open set with
$\partial\Odiff = \Gfar \cup \Gacc$. The part $\Gacc$ is contained in
$\partial\Om_2 \cap \Om_1$ and contains the nonempty relatively open set
$\Bfull{\delta}(p) \cap \partial\Om_2$. Moreover, for every
$\mu \in \{\mathrm{b},\mathrm{s}\}$ and $\theta \in \Sch{\mu}$ the restriction
$u := \Phidirj{1}{-}{\mu}(\cdot,\theta)|_{\Odiff}$ belongs to $H^2(\Odiff)$ and
satisfies
\begin{equation}
  \label{eq:mixed-problem}
  (-\Delta - E)u = 0 \ \text{ in } \Odiff,
  \qquad
  \partial_{\nu}u - \alpha u = 0 \ \text{ on } \Gfar,
  \qquad
  \partial_{\nu}u + \alpha u = 0 \ \text{ on } \Gacc ,
\end{equation}
$\nu$ denoting the outward unit normal of $\Odiff$.
\end{lemma}

\begin{proof}
The set $\Odiff$ is nonempty by construction and bounded, because
$\Om_1$ and $\Om_2$ agree outside $\Bfull{\Rzero}$ by
\eqref{eq:common-exterior}, so that
$\Om_1 \setminus \overline{\Om_2} \subset \Bfull{\Rzero}$. Since
$\Odiff \subset \Om_1$ and $\Om_1$ is open, no point of $\partial\Odiff$
interior to $\Om_1$ can lie outside $\partial\Om_2$. Hence
$\partial\Odiff \subset \partial\Om_1 \cup \partial\Om_2$, and
$\Gacc \subset \partial\Om_2 \cap \Om_1$ by the definition
\eqref{eq:Gfar-Gacc}. That $\Gacc$ contains
$\Bfull{\delta}(p) \cap \partial\Om_2$ is \eqref{eq:local-splitting}.

Write $u_j = \Phidirj{j}{-}{\mu}(\cdot,\theta)$. By
Lemma~\ref{lem:directional}(b) the function $u_1$ lies in
$H^2_{\mathrm{loc}}(\overline{\Om_1})$ and solves the equation and the boundary
condition on $\Om_1$. The closure $\overline{\Odiff}$ is compact and contained
in $\overline{\Om_1}$, so $u_1 \in H^2(\Odiff)$, and the equation in
\eqref{eq:mixed-problem} holds.

On $\Gfar$ the set $\Odiff$ lies on the $\Om_1$ side of $\partial\Om_1$, so the
outward normal of $\Odiff$ agrees with the outward normal of $\Om_1$. The
boundary condition $\gamma_1 u_1 - \alpha\gamma_0 u_1 = 0$ of
\eqref{eq:directional-equation} is therefore the second identity in
\eqref{eq:mixed-problem}.

On $\Gacc$ the situation is reversed. Let $q \in \Gacc$. Then $q \in \Om_1$, so
a small ball around $q$ lies in $\Om_1$ and is split by the $C^2$ hypersurface
$\partial\Om_2$ into a part contained in $\Om_1 \cap \Om_2 = \Omcap$ and a part
contained in $\Odiff$. Because $u_1 \in H^2_{\mathrm{loc}}(\Om_1)$ near $q$,
its trace and normal derivative on $\partial\Om_2$ are the same computed from
either side. Computed from the $\Omcap$ side they coincide with those of $u_2$,
by Lemma~\ref{lem:common-region}, and $u_2$ satisfies
$\gamma_1 u_2 - \alpha\gamma_0 u_2 = 0$ on $\partial\Om_2$ with respect to the
outward normal $\nu_2$ of $\Om_2$. That normal points out of $\Om_2$, hence
into $\Odiff$, so $\nu = -\nu_2$ on $\Gacc$ and
\begin{equation*}
  \partial_\nu u_1 = -\partial_{\nu_2}u_1 = -\partial_{\nu_2}u_2
    = -\alpha u_2 = -\alpha u_1
  \qquad \text{on } \Gacc ,
\end{equation*}
which is the third identity in \eqref{eq:mixed-problem}.
\end{proof}

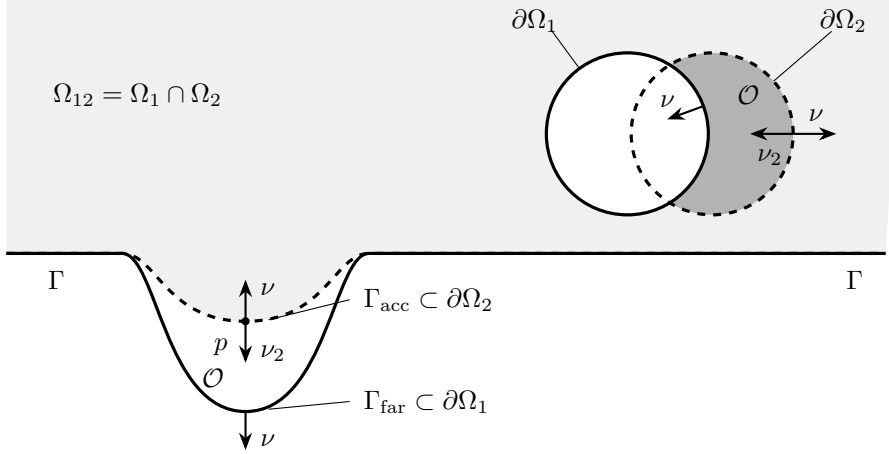
\begin{figure}[htbp]
\centering
\begin{tikzpicture}[scale=1.02, >=Stealth, line join=round]

\def\bdryone{
  (-5.7,0) -- (-4.20,0)
  .. controls (-3.75,0) and (-3.65,-2.05) .. (-2.60,-2.05)
  .. controls (-1.55,-2.05) and (-1.45,0) .. (-1.00,0)
  -- (5.7,0)
}
\def\bdrytwo{
  (-5.7,0) -- (-4.20,0)
  .. controls (-3.75,0) and (-3.65,-0.88) .. (-2.60,-0.88)
  .. controls (-1.55,-0.88) and (-1.45,0) .. (-1.00,0)
  -- (5.7,0)
}

\begin{scope}
  \clip (-5.7,-2.60) rectangle (5.9,3.30);
  \fill[black!6] \bdrytwo -- (5.9,3.30) -- (-5.7,3.30) -- cycle;
  \fill[white] (2.35,1.55) circle (1.05);
  \fill[white] (3.45,1.55) circle (1.05);
\end{scope}

\begin{scope}
  \clip \bdryone -- (5.9,-2.60) -- (-5.7,-2.60) -- cycle;
  \fill[black!30] \bdrytwo -- (5.9,3.30) -- (-5.7,3.30) -- cycle;
\end{scope}

\begin{scope}
  \clip (3.45,1.55) circle (1.05);
  \fill[black!30] (-1,-1) rectangle (10,10);
  \fill[white] (2.35,1.55) circle (1.05);
\end{scope}

\draw[very thick,dashed] \bdrytwo;
\draw[very thick] \bdryone;
\draw[very thick] (2.35,1.55) circle (1.05);
\draw[very thick,dashed] (3.45,1.55) circle (1.05);

\node[font=\small] at (-4.00,2.05) {$\Omcap = \Om_1 \cap \Om_2$};
\node[font=\small] at (-3.02,-1.62) {$\Odiff$};
\node[font=\small] at (3.92,2.06) {$\Odiff$};
\node[font=\small,anchor=north] at (-5.05,-0.10) {$\Gamma$};
\node[font=\small,anchor=north] at (5.30,-0.10) {$\Gamma$};
\node[font=\small] at (1.15,3.02) {$\partial\Om_1$};
\draw[thin] (1.36,2.90) -- (1.729,2.397);
\node[font=\small] at (5.15,3.02) {$\partial\Om_2$};
\draw[thin] (5.00,2.88) -- (4.247,2.234);

\draw[thin] (-2.30,-2.02) -- (-1.25,-1.92);
\node[font=\small,anchor=west] at (-1.20,-1.92) {$\Gfar \subset \partial\Om_1$};
\draw[thin] (-2.30,-0.86) -- (-1.25,-0.60);
\node[font=\small,anchor=west] at (-1.20,-0.60) {$\Gacc \subset \partial\Om_2$};

\draw[->,thick] (-2.60,-2.05) -- (-2.60,-2.55);
\node[font=\small,anchor=west] at (-2.54,-2.42) {$\nu$};
\draw[->,thick] (-2.60,-0.88) -- (-2.60,-0.34);
\node[font=\small,anchor=west] at (-2.54,-0.46) {$\nu$};
\draw[->,thick] (-2.60,-0.88) -- (-2.60,-1.42);
\node[font=\small,anchor=west] at (-2.54,-1.30) {$\nu_2$};
\fill (-2.60,-0.88) circle (1.5pt);
\node[font=\small,anchor=north east] at (-2.70,-0.97) {$p$};

\draw[->,thick] (3.336,1.909) -- (2.866,1.738);
\node[font=\small,anchor=south east] at (3.10,1.76) {$\nu$};
\draw[->,thick] (4.50,1.55) -- (5.06,1.55);
\node[font=\small,anchor=south] at (4.80,1.59) {$\nu$};
\draw[->,thick] (4.50,1.55) -- (3.94,1.55);
\node[font=\small,anchor=north] at (4.20,1.48) {$\nu_2$};

\end{tikzpicture}
\caption{The geometry of Lemma~\ref{lem:mixed-problem}. The boundary of $\Om_1$
is drawn solid and that of $\Om_2$ dashed. The lightly shaded region is the
common domain $\Omcap$, and the two darker ones are components of
$\Om_1 \setminus \overline{\Om_2}$ arising respectively from a difference of
the boundary deformations and from a difference of the obstacles. The argument
does not distinguish the two cases. On $\Gfar$ the outward normal of $\Odiff$
is the outward normal of $\Om_1$, and the field of the first configuration
satisfies its own impedance condition there. On $\Gacc$ it is the
\emph{opposite} of the outward normal $\nu_2$ of $\Om_2$, and the field
inherits the impedance condition of the second configuration through
Lemma~\ref{lem:common-region}. This reversal is the origin of the two signs in
\eqref{eq:mixed-problem}.}
\label{fig:schiffer}
\end{figure}

The fourth step bounds the solutions of \eqref{eq:mixed-problem}. Here the sign
reversal is harmless. The quadratic form is Hermitian for either sign of
$\alpha$, and the boundary term is subordinate to the Dirichlet integral
whatever its sign.

\begin{lemma}[Finite dimensionality]
\label{lem:finite-dim}
Let $\Odiff \subset \Rn$ be a bounded open set with Lipschitz boundary, let
$\partial\Odiff = \Gfar \cup \Gacc$ be a measurable partition, let
$\alpha > 0$ and $E \in \mathbb{R}$, and let
\begin{equation}
  \label{eq:Vmix}
  \Vmix{\Odiff}
  = \bigl\{\, u \in H^2(\Odiff) \;:\; u \text{ satisfies }
    \eqref{eq:mixed-problem} \,\bigr\} .
\end{equation}
Then $\dim \Vmix{\Odiff} < \infty$.
\end{lemma}

\begin{proof}
On $H^1(\Odiff)$ define the Hermitian form
\begin{equation*}
  \mathfrak{a}(u,v)
  = \int_{\Odiff} \nabla u \cdot \overline{\nabla v}\,dx
    - \alpha \int_{\Gfar} u\overline{v}\,dS
    + \alpha \int_{\Gacc} u\overline{v}\,dS .
\end{equation*}
Since $\partial\Odiff$ is Lipschitz and $\Odiff$ is bounded, the trace map
$H^1(\Odiff) \to \Ltwo{\partial\Odiff}$ is compact, so for every
$\varepsilon > 0$ there is $C_\varepsilon$ with
$\norm{u}^2_{\Ltwo{\partial\Odiff}} \le \varepsilon\norm{\nabla u}^2 +
C_\varepsilon \norm{u}^2$ for $u \in H^1(\Odiff)$. Taking
$\varepsilon = 1/(2\alpha)$ gives the Gårding inequality
\begin{equation*}
  \mathfrak{a}(u,u) + \lambda \norm{u}^2_{\Ltwo{\Odiff}}
  \;\ge\; \tfrac12 \norm{u}^2_{H^1(\Odiff)},
  \qquad \lambda = \alpha C_\varepsilon + \tfrac12 .
\end{equation*}
Hence $\mathfrak{a}$ is closed, semibounded and Hermitian on $H^1(\Odiff)$, and
the associated self-adjoint operator $A$ on $\Ltwo{\Odiff}$ has resolvent
factoring through the embedding $H^1(\Odiff) \hookrightarrow \Ltwo{\Odiff}$,
which is compact for a bounded Lipschitz domain. Therefore $A$ has compact
resolvent and all its eigenspaces are finite dimensional.

It remains to see that $\Vmix{\Odiff} \subset \ker(A - E)$. Let
$u \in \Vmix{\Odiff}$ and $v \in H^1(\Odiff)$. As $u \in H^2(\Odiff)$ and
$\Odiff$ is Lipschitz, Green's formula gives
\begin{equation*}
  \int_{\Odiff} \nabla u \cdot \overline{\nabla v}\,dx
  = \int_{\Odiff} (-\Delta u)\overline{v}\,dx
    + \int_{\partial\Odiff} (\partial_\nu u)\overline{v}\,dS ,
\end{equation*}
the boundary integral being absolutely convergent because
$\partial_\nu u \in H^{1/2}(\partial\Odiff)$. Substituting
\eqref{eq:mixed-problem} in both terms turns this into
$\mathfrak{a}(u,v) = E(u,v)_{\Ltwo{\Odiff}}$, so $u \in \Dom(A)$ and
$Au = Eu$.
\end{proof}

The last step supplies the infinitely many solutions. It is also where the
failure of $\Psich{\mu}(\cdot,\theta)$ to lie in $\Bst$, an obstruction in
\S\ref{subsec:inverse-kernels}, becomes the mechanism of the proof. Either
channel alone would supply them, and the lemma treats both at once.

\begin{lemma}[Linear independence of the directional eigenfunctions]
\label{lem:independence}
Let $n \ge 3$, let $E \in \Ereg \setminus \Sexc$ and let $U \subset \Om$ be a
nonempty open set. Let $(\mu_1,\theta_1),\dots,(\mu_m,\theta_m)$ be distinct
channel--direction pairs. If $c_1,\dots,c_m \in \mathbb{C}$ satisfy
\begin{equation}
  \label{eq:independence-hyp}
  \sum_{j=1}^{m} c_j\,\Phidir{-}{\mu_j}(x,\theta_j) = 0
  \qquad \text{for almost every } x \in U ,
\end{equation}
then $c_1 = \dots = c_m = 0$. In particular, if $U$ is bounded, the
restrictions to $U$ of the functions $\Phidir{-}{\mu}(\cdot,\theta)$ span an
infinite dimensional subspace of $\Ltwo{U}$.
\end{lemma}

\begin{proof}
Write $v$ for the left-hand side of \eqref{eq:independence-hyp}, considered on
all of $\Om$. By Lemma~\ref{lem:directional}(b), $v \in
H^2_{\mathrm{loc}}(\overline{\Om})$ and $(-\Delta - E)v = 0$ in $\Om$. The set
$\Om$ is open and connected by Assumption~\ref{assump:geometry}, and $v$
vanishes almost everywhere on its nonempty open subset $U$.
Proposition~\ref{prop:ucp} therefore gives $v = 0$ in $\Om$. By
Definition~\ref{def:directional} this means $\eta P = Q$ on $\Om$, where
\begin{equation*}
  P := \sum_j c_j\,\Psich{\mu_j}(\cdot,\theta_j) ,
  \qquad
  Q := \sum_{j} c_j\,
    \Rpert{E-i0}\bigl( [-\Delta,\eta]\Psich{\mu_j}(\cdot,\theta_j) \bigr),
\end{equation*}
and $Q \in \Bstd{\Om}$ by Theorem~\ref{thm:lap-pert}(b). Write
$D_R = \Bfull{R} \cap \Rnp$. For $R > R_2$ the set $D_R \setminus \Bfull{R_2}$
is contained in $\Om$, and $\eta \equiv 1$ there, so that $P = Q$ on it.
Summing the dyadic blocks of Definition~\ref{def:AH-pert} that meet
$\Bfull{R}$ gives
$\norm{Q}^2_{\Ltwo{\Om \cap \Bfull{R}}} \le C R\,\BstnormO{Q}^2$, while $P$ is
bounded on $D_{R_2}$. Hence
\begin{equation}
  \label{eq:P-growth}
  \norm{P}^2_{\Ltwo{D_R}} \le C R ,
  \qquad R > R_2 ,
\end{equation}
with $C$ independent of $R$. We show that this forces all $c_j$ to vanish.

By \cite[equation (2.51)]{PartI} and \cite[Definition 2.23]{PartI},
\begin{align*}
  \Psich{\mathrm{b}}(x,\theta)
  &= (2\pi)^{-\frac n2}\,\overline{\cph{\zeta_n}}\,
    \bigl( e^{ix\cdot\zeta^{*}} + \refl{\zeta_n}\,e^{ix\cdot\zeta} \bigr),
  &&\zeta = \kbulk\,\theta,
  \\
  \Psich{\mathrm{s}}(x,\omega)
  &= (2\pi)^{-\frac{n-1}{2}}\,e^{i\xpar\cdot\zeta'}\psizero(x_n),
  &&\zeta' = \ksurf\,\omega ,
\end{align*}
where $\zeta^{*}$ is the reflection of $\zeta$ in $\Gamma$,
$|\cph{\zeta_n}| = |\refl{\zeta_n}| = 1$, and
$\norm{\psizero}_{\Ltwo{(0,\infty)}} = 1$. Split
$P = P_{\mathrm{b}} + P_{\mathrm{s}}$ into its bulk and surface parts, and
suppose that some bulk coefficient is nonzero. Write
$\zeta_j = \kbulk\,\theta_j$ for the $j$ with $\mu_j = \mathrm{b}$. Then
\begin{equation*}
  P_{\mathrm{b}} = \sum_{k} d_k\,e^{ix\cdot\xi_k},
\end{equation*}
where $\xi_k$ runs through the vectors $\zeta_j$ and $\zeta_j^{*}$, and the
coefficient of each of the two frequencies attached to $\theta_j$ has modulus
$|d_k| = (2\pi)^{-n/2}|c_j|$. The $\xi_k$ are pairwise distinct. Indeed the
$\theta_j$ are distinct and reflection in $\Gamma$ is injective, and no
$\zeta_j^{*}$ equals a $\zeta_k$, because $\theta_j \in \Sm$ gives
$\zeta_{j,n} > 0$, so that the last component of $\zeta_j^{*}$ is negative and
that of $\zeta_k$ positive. Expanding the square,
\begin{equation*}
  \norm{P_{\mathrm{b}}}^2_{\Ltwo{D_R}}
  = |D_R| \sum_k |d_k|^2
    + \sum_{k \ne l} d_k\overline{d_l}
      \int_{D_R} e^{ix\cdot(\xi_k - \xi_l)}\,dx ,
  \qquad
  |D_R| = \tfrac12 |B_1|\,R^n ,
\end{equation*}
where $|B_1|$ is the volume of the unit ball of $\Rn$. The substitution
$x = Ry$ gives
\begin{equation*}
  \int_{D_R} e^{ix\cdot\xi}\,dx = R^n \int_{D_1} e^{iRy\cdot\xi}\,dy ,
\end{equation*}
and for $\xi \ne 0$ the last integral tends to zero as $R \to \infty$ by the
Riemann--Lebesgue lemma, the indicator of $D_1$ being integrable. Hence
\begin{equation*}
  \lim_{R \to \infty} R^{-n}\norm{P_{\mathrm{b}}}^2_{\Ltwo{D_R}}
  = \tfrac12 |B_1| \sum_k |d_k|^2
  = (2\pi)^{-n}\,|B_1| \sum_{\mu_j = \mathrm{b}} |c_j|^2 > 0 ,
\end{equation*}
each $j$ contributing two frequencies. On the other hand
$|P_{\mathrm{s}}(x)| \le C\psizero(x_n)$, so that
\begin{equation*}
  \norm{P_{\mathrm{s}}}^2_{\Ltwo{D_R}}
  \le C^2 \int_{|\xpar| < R} d\xpar \int_0^\infty \psizero(x_n)^2\,dx_n
  = O(R^{n-1}),
\end{equation*}
and the cross term
$|(P_{\mathrm{b}},P_{\mathrm{s}})_{\Ltwo{D_R}}|$ is $O(R^{n-1/2})$ by the
Cauchy--Schwarz inequality. Therefore
$R^{-n}\norm{P}^2_{\Ltwo{D_R}}$ has the same positive limit, whereas
\eqref{eq:P-growth} gives $R^{-n}\norm{P}^2_{\Ltwo{D_R}} = O(R^{1-n}) \to 0$,
as $n > 1$. This contradiction shows that all bulk coefficients vanish.

Now $P = P_{\mathrm{s}}$, and
\begin{equation*}
  P(x) = (2\pi)^{-\frac{n-1}{2}}\,\psizero(x_n)\,S(\xpar),
  \qquad
  S(\xpar) := \sum_{\mu_j = \mathrm{s}} c_j\,e^{i\xpar\cdot\zeta'_j},
  \qquad
  \zeta'_j = \ksurf\,\theta_j ,
\end{equation*}
the $\zeta'_j \in \Rnm$ being pairwise distinct. The cylinder
$\{|\xpar| < R/2\} \times (0,R/2)$ is contained in $D_R$, and on it the two
variables separate, so that
\begin{equation*}
  \norm{P}^2_{\Ltwo{D_R}}
  \ge (2\pi)^{-(n-1)} \int_0^{R/2} \psizero(x_n)^2\,dx_n
      \int_{|\xpar| < R/2} |S(\xpar)|^2\,d\xpar .
\end{equation*}
The first integral tends to $\norm{\psizero}^2 = 1$. The computation just made
for $P_{\mathrm{b}}$, carried out in $\Rnm$ with the ball $\{|\xpar| < \rho\}$
in place of $D_R$, gives
$\int_{|\xpar| < \rho} |S|^2\,d\xpar
= |B'_1|\,\rho^{n-1} \sum_{\mu_j = \mathrm{s}} |c_j|^2 + o(\rho^{n-1})$ as
$\rho \to \infty$, where $|B'_1|$ is the volume of the unit ball of $\Rnm$.
Taking $\rho = R/2$,
\begin{equation*}
  \liminf_{R \to \infty} R^{-(n-1)}\norm{P}^2_{\Ltwo{D_R}}
  \;\ge\; (4\pi)^{-(n-1)}\,|B'_1| \sum_{\mu_j = \mathrm{s}} |c_j|^2 ,
\end{equation*}
whereas \eqref{eq:P-growth} gives
$R^{-(n-1)}\norm{P}^2_{\Ltwo{D_R}} = O(R^{2-n}) \to 0$, as $n \ge 3$. Hence
the remaining $c_j$ vanish as well.

For the final assertion let $U$ be bounded. By
Lemma~\ref{lem:directional}(b) the restriction to $U$ of each
$\Phidir{-}{\mu}(\cdot,\theta)$ lies in $\Ltwo{U}$. Given $m \in \mathbb{N}$,
choose distinct $\omega_1,\dots,\omega_m \in \Ssurf$, which is possible because
$\Ssurf$ is an infinite set for $n \ge 3$. If
$\sum_k c_k \Phidir{-}{\mathrm{s}}(\cdot,\omega_k) = 0$ in $\Ltwo{U}$, the sum
vanishes almost everywhere in $U$, and the first part of the lemma gives
$c_1 = \dots = c_m = 0$. The span therefore contains $m$ linearly independent
elements for every $m$, and is infinite dimensional.
\end{proof}

\begin{theorem}[Uniqueness of the perturbation]
\label{thm:main-uniqueness}
Let $n \ge 3$ and $\alpha > 0$, let $\Om_1,\Om_2$ be an admissible pair
satisfying Assumption~\ref{assump:pair}, and let $E > 0$. If
\begin{equation}
  \label{eq:main-hypothesis}
  \Amplj{1}{E} = \Amplj{2}{E} ,
\end{equation}
then $\Om_1 = \Om_2$.
\end{theorem}

\begin{proof}
Suppose $\Om_1 \ne \Om_2$. By Lemma~\ref{lem:accessibility}, after exchanging
the indices if necessary, there are $p$ and $\delta$ as in that lemma.
Exchanging the indices leaves the hypotheses \eqref{eq:main-hypothesis} and
Assumption~\ref{assump:pair} intact. Let $\Odiff$, $\Gfar$ and $\Gacc$ be as in
Lemma~\ref{lem:mixed-problem}. That lemma gives
\begin{equation*}
  \Phidirj{1}{-}{\mu}(\cdot,\theta)\big|_{\Odiff} \in \Vmix{\Odiff}
  \qquad \text{for every } \mu \in \{\mathrm{b},\mathrm{s}\}
  \text{ and } \theta \in \Sch{\mu} ,
\end{equation*}
while $\Odiff$ is a bounded open set with Lipschitz boundary by \textup{(U2)},
so that $\Vmix{\Odiff}$ is finite dimensional by
Lemma~\ref{lem:finite-dim}. Let $W \subset \Ltwo{\Odiff}$ be the linear span of
the restrictions $\Phidirj{1}{-}{\mu}(\cdot,\theta)|_{\Odiff}$ with
$\mu \in \{\mathrm{b},\mathrm{s}\}$ and $\theta \in \Sch{\mu}$. The space
$\Vmix{\Odiff}$ is a linear subspace of $H^2(\Odiff) \subset \Ltwo{\Odiff}$
containing each of these restrictions, so $W \subset \Vmix{\Odiff}$ and
\begin{equation*}
  \dim W \le \dim \Vmix{\Odiff} < \infty .
\end{equation*}
On the other hand $\Odiff$ is a nonempty bounded open subset of $\Om_1$, so
Lemma~\ref{lem:independence}, applied to $\Om_1$ with $U = \Odiff$, gives
$\dim W = \infty$. This contradiction proves $\Om_1 = \Om_2$.
\end{proof}

Taking one of the two configurations to be the unperturbed half-space costs
nothing, and answers a question which is usually asked on its own.

\begin{corollary}[No nontrivial perturbation is invisible]
\label{cor:no-invisible}
Let $\Om$ satisfy Assumption~\ref{assump:geometry} and suppose that the pair
$\Om_1 = \Om$, $\Om_2 = \Rnp$ satisfies Assumption~\ref{assump:pair}. If
$\Om \ne \Rnp$, then
\[
  \Ampl{E} \ne 0
  \qquad\text{and equivalently}\qquad
  \Smatpert{E} \ne \Smat(E)
  \qquad\text{for every } E > 0 .
\]
\end{corollary}

\begin{proof}
The half-space itself satisfies Assumption~\ref{assump:geometry}. It is
connected, its boundary is the hyperplane $\Gamma$ and so is a $C^2$
hypersurface, and \eqref{eq:exterior-coincidence} holds trivially for the
radius $\Rzero$ of $\Om$. Hence $\Om_1,\Om_2$ is an admissible pair in the
sense of Definition~\ref{def:admissible-pair}. By the last assertion of
Theorem~\ref{thm:smatrix-far-field}(a), the scattering matrix of
$\Om_2 = \Rnp$ equals $\Smat(E)$ for every $E$, so
\eqref{eq:smatrix-amplitude}, applied to $\Om_2$, gives
$\Phase{E}\Amplj{2}{E}\Phase{E} = 0$, and $\Amplj{2}{E} = 0$ because
$\Phase{E}$ is unitary (Lemma~\ref{lem:phase-multiplier}(a)). If now
$\Amplj{1}{E} = 0$ for some $E > 0$ then $\Amplj{1}{E} = \Amplj{2}{E}$, and
Theorem~\ref{thm:main-uniqueness} forces $\Om = \Rnp$. The equivalence of the
two displayed statements is \eqref{eq:smatrix-amplitude} together with the
unitarity of $\Phase{E}$.
\end{proof}

We note two features of this result. The first is that the hypotheses
degenerate here into conditions on $\Om$ alone, and into recognizable ones. In
Theorem~\ref{thm:main-uniqueness} they are conditions on a pair of unknowns.
Here $\Rnp \setminus \overline{\Om}$ is what the perturbation removes from the
half-space, and $\Om \setminus \overline{\Rnp}$ is what it adds. Hypothesis
\textup{(U2)} therefore asks that each piece added or removed have Lipschitz
boundary, and \textup{(U1)} that what remains of the half-space not be
disconnected by the perturbation. By Remark~\ref{rem:hypotheses},
\textup{(U2)} can fail only where $\Gpert$ meets $\Gamma$ tangentially. The
second feature is that the conclusion holds at \emph{every} positive energy, so
that a perturbation is visible at any single frequency of the experimenter's
choosing.

We record what the corollary does not say. It asserts that
$\Ampl{E} \ne 0$, that is, that at least one of the four blocks of
$\Smatpert{E} - \Smat(E)$ is nonzero. It does not locate which. In particular
it does not exclude the configuration $\Sbs = \Ssb = 0$ with
$\Sbb \ne \refl{\kbulk\theta_n}$, in which the perturbation scatters within
each channel without transferring energy between them. Nothing in this paper
decides whether such a configuration occurs. The question of which blocks are
nonzero, as opposed to which matrices differ, remains open, and
Remark~\ref{rem:band-single-block} bears on it.

\begin{remark}[Where each hypothesis is used]
\label{rem:hypotheses}
The proof consumes four hypotheses.

Hypothesis \textup{(U1)} is used twice, and both times to run a unique
continuation argument on a connected set. It enters in
Lemma~\ref{lem:common-region}, to propagate the equality of the fields from
$\Omex{\Rzero}$ into $\Omcap$, and in Lemma~\ref{lem:accessibility}, where the
connectedness of $\Omcap$ is what makes the clopen argument produce an
accessible point. It holds automatically in the two configurations one has in
mind. If the perturbations are deformations of the boundary given by graphs
$x_n = \zeta_j(\xpar)$, then $\Omcap = \{x_n > \max(\zeta_1,\zeta_2)\}$ is
connected. If they are obstacles, \textup{(U1)} is the standard hypothesis that
the complement of the union of the two obstacles is connected
\cite[\S5.1]{ColtonKress2013}.

Hypothesis \textup{(U2)} is used once, in Lemma~\ref{lem:finite-dim}, and only
through the compactness of the embedding $H^1(\Odiff) \hookrightarrow
\Ltwo{\Odiff}$ and of the trace $H^1(\Odiff) \to \Ltwo{\partial\Odiff}$. It is
the one place where the impedance condition costs something that a Dirichlet
condition would not. For the sound-soft problem the corresponding space is
contained in $H^1_0(\Odiff)$, whose embedding into $\Ltwo{\Odiff}$ is compact
for \emph{every} bounded open set, and no regularity of $\partial\Odiff$ is
needed. Here the boundary term in $\mathfrak{a}$ has to be controlled, and a
trace is unavoidable. Since $\partial\Odiff$ is contained in the union of two
$C^2$ hypersurfaces, \textup{(U2)} can fail only where $\partial\Om_1$ and
$\partial\Om_2$ touch tangentially, and it holds in particular whenever the two
boundaries are transversal wherever they meet.

That each $\Om_j$ is separately connected, which is part of
Assumption~\ref{assump:geometry}, is used in Lemma~\ref{lem:accessibility} and
again in Lemma~\ref{lem:independence}, there through
Proposition~\ref{prop:ucp}. The positivity of $E$ enters only
through Theorem~\ref{thm:exterior-equality}, as
\S\ref{subsec:surface-only} explains.
\end{remark}

The relation of Theorem~\ref{thm:main-uniqueness} to the uniqueness theorem of
Zhao, Hu and Yan \cite{ZhaoHuYan2021} for the Dirichlet and Neumann half-spaces,
and the role of the reflection principle there, are discussed in
\S\ref{subsec:intro-novelty}.

\subsection{Determination by the Surface Channel Alone}
\label{subsec:surface-only}

Theorem~\ref{thm:main-uniqueness} uses the scattering amplitude at one positive
energy, where both channels are open and the data consist of the four blocks of
\eqref{eq:amplitude-kernel}. It is natural to ask how much of the perturbation
is already determined by the surface channel alone. On the band
$\Eband = (-\alpha^2,0)$ the question is not a restriction of the data but a
description of them. There is no bulk channel, $\Echan{E} = \Ltwo{\Ssurf}$, and
$\Ampl{E}$ consists of the single block $A_{\mathrm{ss}}$, as recorded in
Remark~\ref{rem:band-single-block}.

The proof of Theorem~\ref{thm:main-uniqueness} is almost insensitive to this.
Inspecting it, the positivity of $E$ enters at exactly one place, namely
Theorem~\ref{thm:exterior-equality}. Everything downstream is valid for every
$E \in \Ereg$, and that includes the transfer to $\Omcap$, the accessibility
argument, the mixed problem, the finite dimensionality and the linear
independence. In particular Lemma~\ref{lem:independence} needs only that
infinitely many distinct directions are available, and $\Ssurf$ is infinite for
$n \ge 3$. The surface family is a smaller family than the bulk one, being
parametrized by $n-2$ rather than $n-1$ angles. That costs nothing here, since
the argument asks for infinitely many members and not for a family of a
prescribed size. We may therefore state what is true.

\begin{theorem}[Conditional uniqueness on the band]
\label{thm:band-conditional}
Let $n \ge 3$ and $\alpha > 0$, let $\Om_1,\Om_2$ be an admissible pair
satisfying Assumption~\ref{assump:pair}, and let
$E \in \Eband \setminus (\Sexcj{1} \cup \Sexcj{2})$ satisfy
$\Amplj{1}{E} = \Amplj{2}{E}$. Assume in addition that
\begin{equation}
  \label{eq:band-hypothesis}
  \Phipertj{1}{-}{E}\phi = \Phipertj{2}{-}{E}\phi
  \quad \text{in } \Omex{\Rzero},
  \qquad \phi \in \Nsol{E} .
\end{equation}
Then $\Om_1 = \Om_2$.
\end{theorem}

\begin{proof}
Hypothesis \eqref{eq:band-hypothesis} is the conclusion of
Theorem~\ref{thm:exterior-equality}, and that theorem is the only point at
which $E > 0$ was used. The proof of Lemma~\ref{lem:common-region} appeals to
it in its first sentence and nowhere else, so \eqref{eq:band-hypothesis} may be
substituted there. The proofs of Lemmas~\ref{lem:accessibility},
\ref{lem:mixed-problem}, \ref{lem:finite-dim} and \ref{lem:independence} and of
Theorem~\ref{thm:main-uniqueness} then apply verbatim with
$E \in \Eband \setminus (\Sexcj{1} \cup \Sexcj{2})$, which is a subset of
$\Ereg \setminus \Sexcj{j}$ for both $j$, the bulk channel being absent
throughout.
\end{proof}

What is missing is therefore a single implication. On the band, equality of the
amplitudes gives a difference $w$ on $\Omex{\Rzero}$ which lies in
$\Bstzd{\Om}$ and solves the exterior Robin problem, and one wants to conclude
that $w$ vanishes. For $E > 0$ this is Theorem~\ref{thm:rellich-positive}. On the band no such
statement is available, an exterior problem there behaving like a closed
waveguide (Remark~\ref{rem:band-closed-waveguide}), and
Theorem~\ref{thm:band-equivalence} gives only exponential decay. The substitute available in
\S\ref{sec:lap-pert}, namely Corollary~\ref{cor:global-uniqueness}(b), requires
its argument to be a solution on the whole of a perturbed domain. The
difference of two scattering solutions is not such a function. It satisfies the
Robin condition of $\Om_1$ on $\partial\Om_1 \cap \partial\Om_2$ only, and on
$\partial\Omcap \setminus (\partial\Om_1 \cap \partial\Om_2)$ one of the two
fields is interior and satisfies no boundary condition at all. We record the
gap as a problem.

\begin{problem}
\label{prob:band-rellich}
Let $\Om_1,\Om_2$ be an admissible pair, let
$E \in \Eband \setminus (\Sexcj{1} \cup \Sexcj{2})$ and let
$\Amplj{1}{E} = \Amplj{2}{E}$. Does \eqref{eq:band-hypothesis} hold?
Equivalently, if a function on $\Omex{\Rzero}$ which solves the exterior Robin
problem at a band energy and belongs to $\Bstzd{\Om}$ arises as the difference
of the scattering solutions of two admissible configurations, must it vanish?
\end{problem}

A positive answer would, by Theorem~\ref{thm:band-conditional}, determine the
whole perturbation, boundary deformation and obstacles alike, from the single
surface block $A_{\mathrm{ss}}$ at one band energy. Those are data which can be
collected without ever launching a wave into the bulk. We are not able to
decide the question. Four comments delimit what is at stake.

First, the exponential decay furnished by Theorem~\ref{thm:band-equivalence} is
already a strong constraint. The difference is invisible not only at leading
order but at every order of the far-field expansion, so a counterexample would
be a pair of configurations which are exponentially indistinguishable in the
band and yet distinct. Whether such pairs exist is precisely
Problem~\ref{prob:band-rellich}.

Secondly, the difficulty is one about \emph{differences} of two problems, not about a single one. For a single configuration
Corollary~\ref{cor:global-uniqueness}(b) settles the matter up to the
exceptional set, and the exceptional set is harmless there because it is
intrinsic to the given domain and discrete by
Proposition~\ref{prop:exceptional-discrete}. Here no single domain is available
to attach it to.

Thirdly, the flat case shows that the surface channel can be blind. On the band
the scattering matrix of $\Rnp$ is the identity for every $\alpha$ and every
band energy \cite[Proposition 3.3(ii)]{PartI}, so the unperturbed band data are
one fixed operator, and Part~I recovers the impedance from the bulk block
instead \cite[Corollary 3.4]{PartI}. Remark~\ref{rem:band-single-block} reads
that degeneracy as a baseline against which a perturbation can be seen. Read the
other way, it says that no argument can extract information from the surface
block without using the geometry, there being none to extract when the geometry
is flat. That is consistent with Theorem~\ref{thm:band-conditional} being
conditional.

Finally, one comment points outside the stationary framework. In the
time-dependent theory of Lax--Phillips type, Kawashita and Soga
\cite{KawashitaSoga2010} study the elastic wave equation in a half-space of
$\mathbb{R}^3$ whose boundary carries a bounded undulation, a perturbation of
the same local character as ours. They decompose the scattering kernel into
channels and analyse the one carried by the Rayleigh wave, which is the
analogue of the surface block here. Their result locates its singularities. The
singular support is contained in the half-line determined by the maximal
sojourn time of the flow generated on the perturbed boundary by the Hamiltonian
governing the surface wave, and, when that maximum is attained at a single
point, the endpoint belongs to it. Their hypotheses are not ours, the boundary
flow being assumed non-trapping, and their conclusion is one about
singularities, hence about all high frequencies rather than about a single
energy. Nor is an inverse problem formulated there. The mechanism is
nevertheless suggestive for Problem~\ref{prob:band-rellich}. Because $\alpha$
is constant, the transverse profile $\psizero$ is the same at every point of
$\Gamma$, so that what a perturbation changes is the metric induced on
$\Gpert$. The recovery of a metric from travel times of this kind is a question
of Riemannian geometry, and it would then remain to pass from the intrinsic
metric to its embedding in $\Rn$. We record this for two reasons. It indicates
that the surface block constrains the geometry of $\Gpert$ along a route
independent of \S\ref{subsec:schiffer}. Moreover, the non-trapping hypothesis
of \cite{KawashitaSoga2010} and the exceptional set of
Proposition~\ref{prop:exceptional-discrete} are the time-dependent and the
stationary form of one difficulty, that of trapping.

\appendix

\section{Elliptic Regularity and Traces on the Perturbed Domain}
\label{app:regularity}
Throughout this appendix $\Om$ satisfies Assumption~\ref{assump:geometry} and
$\alpha > 0$ is fixed. We prove part~(a) of
Theorem~\ref{thm:pert-sa-spec}: that the form $\hpert$ of
\eqref{eq:form-pert} is closed and semibounded on $H^1(\Om)$, that the
self-adjoint operator it generates is the operator $\Hpert$ of
\eqref{eq:Hpert-def}, and that the a priori estimate \eqref{eq:apriori-H2}
holds. Parts~(b) and~(c) are proved in \S\ref{subsec:pert-sa-spec} and require
nothing from here beyond \eqref{eq:apriori-H2}.

The only step that is not immediate is the regularity assertion
$\Dom(\Hpert) \subset H^2(\Om)$, and it is separated off as
Proposition~\ref{prop:regularity}.

\subsection{The Form and the Operator It Generates}
\label{app:form}

\begin{lemma}
\label{lem:form-closed}
Let $n \ge 2$. The form $\hpert$ with form domain $H^1(\Om)$ is symmetric,
semibounded and closed, and
\begin{equation}
  \label{eq:form-equivalence}
  \tfrac12 \norm{u}_{H^1(\Om)}^2 - c_0 \norm{u}_{\Ltwo{\Om}}^2
  \;\le\; \hpert[u]
  \;\le\; \norm{u}_{H^1(\Om)}^2 ,
  \qquad
  c_0 := 2\alpha^2 + \alpha\Cdiv + \tfrac12 .
\end{equation}
\end{lemma}

\begin{proof}
Symmetry is clear from \eqref{eq:form-pert}. Apply
Lemma~\ref{lem:trace-pert} with $\varepsilon = 1/(2\alpha)$. For
$u \in H^1(\Om)$ this gives
\begin{equation*}
  \alpha \norm{\gamma_0 u}_{\Ltwo{\Gpert}}^2
    \le \tfrac12 \norm{\nabla u}_{\Ltwo{\Om}}^2
      + \bigl(2\alpha^2 + \alpha\Cdiv\bigr) \norm{u}_{\Ltwo{\Om}}^2 ,
\end{equation*}
whence
$\hpert[u] \ge \tfrac12\norm{\nabla u}^2
- (2\alpha^2 + \alpha\Cdiv)\norm{u}^2$, which is
\eqref{eq:lower-bound-pert}. Adding and subtracting $\tfrac12\norm{u}^2$ turns
that bound into the left inequality of \eqref{eq:form-equivalence}. The right
one is $\hpert[u] \le \norm{\nabla u}^2$, obtained by discarding the negative
boundary term. Thus the form norm
$\bigl(\hpert[u] + (c_0+1)\norm{u}^2\bigr)^{1/2}$ is equivalent to the
$H^1(\Om)$ norm, and $H^1(\Om)$ is complete, so $\hpert$ is closed.
\end{proof}

By the representation theorem for closed semibounded forms
\cite[Thm.~VI.2.1 and Thm.~VI.2.6]{Kato1976} there is a unique self-adjoint
operator $H_{\mathfrak{h}}$ in $\Ltwo{\Om}$ associated with $\hpert$, namely
\begin{equation}
  \label{eq:app-weak-form}
  \Dom(H_{\mathfrak{h}})
    = \bigl\{\, u \in H^1(\Om) \;\bigm|\;
        \exists f \in \Ltwo{\Om} \ \text{ such that } \
        \hpert[u,v] = (f, v)_{\Ltwo{\Om}} \ \ \forall v \in H^1(\Om) \,\bigr\},
\end{equation}
with $H_{\mathfrak{h}} u = f$. It remains to identify $H_{\mathfrak{h}}$ with
$\Hpert$.

One inclusion is Green's formula. If $u \in H^2(\Om)$ satisfies
$\gamma_1 u = \alpha\gamma_0 u$, then for every $v \in H^1(\Om)$
\begin{equation}
  \label{eq:green}
  (-\Delta u, v)_{\Ltwo{\Om}}
    = (\nabla u, \nabla v)_{\Ltwo{\Om}}
      - (\gamma_1 u, \gamma_0 v)_{\Ltwo{\Gpert}}
    = \hpert[u, v],
\end{equation}
so that $u \in \Dom(H_{\mathfrak{h}})$ with $H_{\mathfrak{h}}u = -\Delta u$.
That is, $\Hpert \subset H_{\mathfrak{h}}$. Note the sign. The trace $\gamma_1$
is taken with respect to the \emph{outward} normal $\nu$, which is why the
boundary term in \eqref{eq:green} is subtracted, and it is the substitution
$\gamma_1 u = +\alpha\gamma_0 u$ that produces the negative boundary term of
\eqref{eq:form-pert}.

The reverse inclusion is the content of the next proposition. Granting it, let
$u \in \Dom(H_{\mathfrak{h}})$, so that $u \in H^2(\Om)$. Taking
$v \in C_0^\infty(\Om)$ in \eqref{eq:app-weak-form} gives $-\Delta u = f$ in
the distributional and hence in the $\Ltwo{\Om}$ sense. Substituting this back
and letting $v$ range over $H^1(\Om)$ leaves
$(\gamma_1 u - \alpha\gamma_0 u, \gamma_0 v)_{\Ltwo{\Gpert}} = 0$ for every
such $v$. The trace map $\gamma_0$ carries $H^1(\Om)$ onto $H^{1/2}(\Gpert)$,
which is dense in $\Ltwo{\Gpert}$, so the boundary condition holds. Hence
$H_{\mathfrak{h}} \subset \Hpert$, and the two operators coincide. Together
with Lemma~\ref{lem:form-closed} this proves
Theorem~\ref{thm:pert-sa-spec}(a), apart from \eqref{eq:apriori-H2}, which is
part of Proposition~\ref{prop:regularity}.

\subsection{Elliptic Regularity up to the Boundary}
\label{app:regularity-proper}

\begin{proposition}
\label{prop:regularity}
Let $n \ge 2$, let $\Om$ satisfy Assumption~\ref{assump:geometry}, let
$f \in \Ltwo{\Om}$, and let $u \in H^1(\Om)$ satisfy
\begin{equation}
  \label{eq:weak-solution}
  \hpert[u, v] = (f, v)_{\Ltwo{\Om}}
  \qquad \text{for every } v \in H^1(\Om).
\end{equation}
Then $u \in H^2(\Om)$, and there is a constant $C = C(\Om, \alpha, n)$,
independent of $u$ and $f$, with
\begin{equation}
  \label{eq:app-apriori}
  \norm{u}_{H^2(\Om)}
    \le C \bigl( \norm{f}_{\Ltwo{\Om}} + \norm{u}_{\Ltwo{\Om}} \bigr).
\end{equation}
\end{proposition}

\begin{proof}
Fix a radial cut-off $\cut \in C^\infty(\Rn)$ with $\cut = 1$ on
$\Bfull{\Rzero}$ and $\cut = 0$ outside $\Bfull{\Rzero + 1}$, and split $u = \cut u + (1 - \cut)u$. Radiality is needed: it gives $\partial_{x_n}\cut = 0$ on $\{x_n = 0\}$, hence
$\gamma_1(\cut u) = \cut\,\gamma_1 u$ on the flat part of the boundary, so that
each piece continues to satisfy a Robin condition rather than an inhomogeneous
one. Throughout, $\norm{u}_{H^1(\Om)} \le C(\norm{f} + \norm{u})$, which
follows from \eqref{eq:weak-solution} with $v = u$ and
Lemma~\ref{lem:form-closed}. We use it without further comment.

\emph{Step 1: the compact piece.} Let $D$ be a bounded open set of class $C^2$
with
\begin{equation*}
  \Om \cap \Bfull{\Rzero + \tfrac32} \;\subset\; D
  \;\subset\; \Om \cap \Bfull{\Rzero + 2},
  \qquad
  \partial D \cap \Bfull{\Rzero + \tfrac32}
    = \Gpert \cap \Bfull{\Rzero + \tfrac32} .
\end{equation*}
Such a set exists. On the annulus $\Bfull{\Rzero+2} \setminus \Bfull{\Rzero}$
the boundary of $\Om$ is the hyperplane $\Gamma$ by
\eqref{eq:exterior-coincidence}, so the edge along which $\Gamma$ meets the
sphere $\partial\Bfull{\Rzero+2}$ can be rounded off there. The truncation
$\Omin{\Rzero+2}$ itself would not serve, its boundary having an edge along
that sphere and the regularity theory quoted below tolerating none.

The function $w = \cut u$ vanishes on a neighbourhood of
$\partial D \setminus \Gpert$, and on $\Gpert$ it satisfies
$\gamma_1 w = \alpha \gamma_0 w$. Testing \eqref{eq:weak-solution} against
$\cut v$ and commuting $\cut$ past the gradient shows that $w \in H^1(D)$ is a
weak solution of
\begin{equation*}
  -\Delta w = \cut f - 2\nabla\cut \cdot \nabla u - (\Delta\cut) u
  \quad \text{in } D,
  \qquad
  \gamma_1 w - \alpha \gamma_0 w = 0 \quad \text{on } \partial D,
\end{equation*}
whose right-hand side
$\Phi := \cut f - 2\nabla\cut \cdot \nabla u - (\Delta\cut) u$ lies in
$\Ltwo{D}$ with norm bounded by $C(\norm{f} + \norm{u}_{H^1(\Om)})$.

The regularity is quoted from \cite[Thm.~2.4.2.7]{Grisvard1985}, but not before
the sign of the boundary condition has been attended to. That theorem asserts
existence and uniqueness, and it requires the zeroth-order coefficient of the
boundary operator to have the dissipative sign, which \eqref{eq:robin-pert}
does not have. For our sign the form \eqref{eq:form-pert} is not coercive,
which is the same fact that gives $\Hpert$ negative spectrum and, on the half-space, the surface wave
\eqref{eq:psizero}. The remedy is to move the offending term into the data. For
$v \in H^1(D)$ the weak formulation of the problem above reads
\begin{equation*}
  (\nabla w, \nabla v)_{\Ltwo{D}}
    - \alpha\,(\gamma_0 w, \gamma_0 v)_{\Ltwo{\partial D}}
  = (\Phi, v)_{\Ltwo{D}} ,
\end{equation*}
and adding
$(w,v)_{\Ltwo{D}} + 2\alpha(\gamma_0 w, \gamma_0 v)_{\Ltwo{\partial D}}$ to
both sides exhibits $w$ as a variational solution of
\begin{equation*}
  -\Delta w + w = \Phi + w \ \text{ in } D,
  \qquad
  \partial_\nu w + \alpha w = 2\alpha\,\gamma_0 w \ \text{ on } \partial D .
\end{equation*}

Every hypothesis of \cite[Thm.~2.4.2.7]{Grisvard1985} is met by this
problem. The boundary $\partial D$ is of class $C^2$, hence of
class $C^{1,1}$. The zeroth-order coefficient in the interior is $1$, which is
bounded below by a positive constant. The coefficient vector of the boundary
operator is the outward normal $\nu$, which is Lipschitz on $\overline{D}$
because $\partial D$ is $C^{1,1}$, and $\nu \cdot \nu = 1$ does not vanish.
The zeroth-order coefficient on the boundary is $\alpha$, and the product of
the two is $\alpha > 0$. Finally $\Phi + w \in \Ltwo{D}$, and
$2\alpha\gamma_0 w \in H^{1/2}(\partial D)$ because $w \in H^1(D)$. The
theorem therefore produces a solution of it in $H^2(D)$.
That solution and $w$ are both variational solutions in $H^1(D)$ of the same
problem, whose form is coercive there, so they coincide. Hence
$w \in H^2(D)$. The solution operator is a bijection between Banach spaces, so
it is bounded, and
\begin{equation*}
  \norm{w}_{H^2(D)}
    \le C\bigl( \norm{\Phi + w}_{\Ltwo{D}}
        + \norm{\gamma_0 w}_{H^{1/2}(\partial D)} \bigr)
    \le C\bigl( \norm{f}_{\Ltwo{\Om}} + \norm{u}_{\Ltwo{\Om}} \bigr),
\end{equation*}
the last step by the $H^1$ bound recorded at the start of the proof.

\emph{Step 2: the flat piece.} Put $\tilde{w} = (1 - \cut)u$, extended by zero,
so that $\tilde{w} \in H^1(\Rnp)$ is supported in
$\Rnp \setminus \Bfull{\Rzero}$, where by \eqref{eq:exterior-coincidence} the
domain is the unperturbed half-space. As in Step~1, $\tilde{w}$ is a weak
solution of $-\Delta \tilde{w} = g$ on $\Rnp$ with
$g = (1-\cut)f + 2\nabla\cut \cdot \nabla u + (\Delta\cut)u \in \Ltwo{\Rnp}$
and $\gamma_1\tilde{w} - \alpha\gamma_0\tilde{w} = 0$ on $\Gamma$. Here the
constant-coefficient argument of \cite[proof of Theorem~2.2(a)]{PartI} applies
verbatim, and we recall it in one line. The half-space is invariant under
translations in $\xpar$, so the difference quotients $D_h^j \tilde{w}$,
$1 \le j \le n-1$, are admissible test functions and remain in the form domain.
Testing the weak formulation against $D_h^{-j}D_h^j\tilde{w}$ together with
Lemma~\ref{lem:trace-pert} gives a bound on $\nabla'\nabla\tilde{w}$ in
$\Ltwo{\Rnp}$ uniform in $h$. Hence $\nabla'\nabla\tilde{w} \in \Ltwo{\Rnp}$,
and the remaining second derivative is recovered from the equation itself,
\begin{equation*}
  \partial_{x_n}^2 \tilde{w} = -g - \Delta' \tilde{w} \in \Ltwo{\Rnp},
\end{equation*}
so that $\tilde{w} \in H^2(\Rnp)$ with
$\norm{\tilde{w}}_{H^2(\Rnp)} \le C(\norm{g} + \norm{\tilde{w}})
\le C(\norm{f}_{\Ltwo{\Om}} + \norm{u}_{\Ltwo{\Om}})$.

\emph{Step 3.} Both pieces lie in $H^2$ and are supported where $\Om$ and the
respective model domain agree, so $u = w + \tilde{w} \in H^2(\Om)$, and adding
the two estimates gives \eqref{eq:app-apriori}.
\end{proof}

Applying Proposition~\ref{prop:regularity} to $u \in \Dom(\Hpert)$ and
$f = \Hpert u$ yields \eqref{eq:apriori-H2}, which completes the proof of
Theorem~\ref{thm:pert-sa-spec}(a).

The limiting absorption principle of \S\ref{subsec:lap-pert-main} is proved by
extracting a convergent subsequence from a family of resolvents, and for that
the global estimate \eqref{eq:app-apriori} is not enough. The members of the
family satisfy the equation, but the limit is only known to satisfy it after
the passage to the limit has been performed. What is needed is a local version,
on a truncation of $\Om$, with a constant that does not depend on the spectral
parameter beyond its modulus.

\begin{proposition}[Local a priori estimate]
\label{prop:local-regularity}
Let $n \ge 2$, let $\Om$ satisfy Assumption~\ref{assump:geometry}, and let
\begin{equation*}
  \Rzero \le \rho < \rho' , \qquad M > 0 .
\end{equation*}
There is a constant $C = C(\Om, \alpha, n, \rho, \rho', M)$ with the following
property. Let $z \in \mathbb{C}$ with $|z| \le M$, let
$F \in \Ltwo{\Omin{\rho'}}$, and let $u \in H^1(\Omin{\rho'})$ satisfy
\begin{equation}
  \label{eq:local-weak}
  (\nabla u, \nabla v)_{\Ltwo{\Omin{\rho'}}}
    - \alpha\,(\gamma_0 u, \gamma_0 v)_{\Ltwo{\Gpert \cap \Bfull{\rho'}}}
    - z\,(u,v)_{\Ltwo{\Omin{\rho'}}}
  = (F, v)_{\Ltwo{\Omin{\rho'}}}
\end{equation}
for every $v \in H^1(\Omin{\rho'})$ that vanishes in a neighbourhood of
$\Om \cap \partial\Bfull{\rho'}$. Then $u \in H^2(\Omin{\rho})$, it satisfies
$\gamma_1 u - \alpha\gamma_0 u = 0$ on $\Gpert \cap \Bfull{\rho}$, and
\begin{equation}
  \label{eq:app-local}
  \norm{u}_{H^2(\Omin{\rho})}
    \le C\bigl( \norm{F}_{\Ltwo{\Omin{\rho'}}}
               + \norm{u}_{\Ltwo{\Omin{\rho'}}} \bigr).
\end{equation}
\end{proposition}

\begin{proof}
Fix $\rho < \rho_1 < \rho_2 < \rho'$ and radial cut-offs
$\cut_1, \cut_2 \in C^\infty(\Rn)$ with $\cut_1 = 1$ on $\Bfull{\rho}$ and
$\operatorname{supp}\cut_1 \subset \Bfull{\rho_1}$, and $\cut_2 = 1$ on
$\Bfull{\rho_1}$ with $\operatorname{supp}\cut_2 \subset \Bfull{\rho_2}$. Since
$\rho \ge \Rzero$, the gradients $\nabla\cut_1$ and $\nabla\cut_2$ are
supported in $\{|x| \ge \Rzero\}$, where the boundary is the hyperplane
$\Gamma$ by \eqref{eq:exterior-coincidence}. Radiality therefore gives
$\partial_\nu \cut_i = -\partial_{x_n}\cut_i = 0$ at every boundary point of
$\operatorname{supp}\nabla\cut_i$, so that multiplication by $\cut_i$ preserves
the homogeneous Robin condition.

\emph{Step 1: a Caccioppoli inequality.} Take $v = \cut_2^2 u$ in
\eqref{eq:local-weak}, which is admissible, and take real parts. Since
$\gamma_0(\cut_2 u) = \cut_2\gamma_0 u$, the boundary term is
$-\alpha\norm{\gamma_0(\cut_2 u)}^2$, and
$\re(\nabla u, \nabla(\cut_2^2 u))
= \norm{\cut_2\nabla u}^2 + 2\re(\cut_2\nabla u, u \nabla\cut_2)$. Bound the
last term by
$\tfrac12\norm{\cut_2\nabla u}^2 + 2\norm{|\nabla\cut_2|u}^2$ and apply
Lemma~\ref{lem:trace-pert} to $\cut_2 u \in H^1(\Om)$ with
$\varepsilon = 1/(8\alpha)$. Since
$\norm{\nabla(\cut_2 u)}^2 \le 2\norm{\cut_2\nabla u}^2 +
2\norm{|\nabla\cut_2|u}^2$, the trace term absorbs a further
$\tfrac14\norm{\cut_2\nabla u}^2$, and what is left is
\begin{equation*}
  \tfrac14 \norm{\cut_2 \nabla u}_{\Ltwo{\Om}}^2
  \le C\bigl( \norm{F}^2 + (1 + M) \norm{u}^2 \bigr),
\end{equation*}
all norms on the right being those of $\Ltwo{\Omin{\rho'}}$. In particular
\begin{equation}
  \label{eq:app-caccioppoli}
  \norm{u}_{H^1(\Omin{\rho_1})}
    \le C\bigl( \norm{F}_{\Ltwo{\Omin{\rho'}}}
               + \norm{u}_{\Ltwo{\Omin{\rho'}}} \bigr).
\end{equation}

\emph{Step 2: second derivatives.} Let $D$ be a bounded open set of class $C^2$
with
\begin{equation*}
  \Om \cap \Bfull{\rho_1} \;\subset\; D \;\subset\; \Om \cap \Bfull{\rho_2},
  \qquad
  \partial D \cap \Bfull{\rho_1} = \Gpert \cap \Bfull{\rho_1} ,
\end{equation*}
which exists for the same reason as the set $D$ of Step~1 in the proof of
Proposition~\ref{prop:regularity}, the boundary of $\Om$ being the hyperplane
$\Gamma$ on $\Bfull{\rho_2} \setminus \Bfull{\rho_1}$. Put $w = \cut_1 u$, a
function vanishing on a neighbourhood of $\partial D \setminus \Gpert$. Testing
\eqref{eq:local-weak} against $\cut_1 v$ and commuting $\cut_1$ past the
gradient shows that $w \in H^1(D)$ is a weak solution of
\begin{equation*}
  -\Delta w
    = \cut_1 F + z\,\cut_1 u - 2\nabla\cut_1 \cdot \nabla u
      - (\Delta\cut_1)u
  \quad \text{in } D,
  \qquad
  \gamma_1 w - \alpha\gamma_0 w = 0 \quad \text{on } \partial D,
\end{equation*}
the right-hand side lying in $\Ltwo{D}$ with norm at most
$C(\norm{F} + (1+M)\norm{u}_{H^1(\Omin{\rho_1})})$, hence, by
\eqref{eq:app-caccioppoli}, at most the right-hand side of
\eqref{eq:app-local}. Flip the sign of the boundary term as in Step~1 of the
proof of Proposition~\ref{prop:regularity} and apply
\cite[Thm.~2.4.2.7]{Grisvard1985} to the resulting problem. This gives
$w \in H^2(D)$ with $\norm{w}_{H^2(D)}$ bounded by the same quantity together
with $\norm{w}_{H^1(D)}$, which \eqref{eq:app-caccioppoli} controls as well.
Since $u = w$ on $\Omin{\rho}$, this is \eqref{eq:app-local}. The Robin
condition for $w$ on $\Gpert \cap \Bfull{\rho}$, where $\cut_1 = 1$, is that
for $u$.
\end{proof}

\section{Unique Continuation}
\label{app:ucp}
The uniqueness arguments of \S\ref{sec:lap-pert} and \S\ref{subsec:schiffer}
need one classical fact, which we isolate here. What is required is unique
continuation \emph{inside} the domain. A solution of the homogeneous equation which vanishes on a nonempty
open subset vanishes identically. No continuation across $\Gpert$ is needed,
and the boundary condition plays no role at all.

\begin{proposition}[Unique continuation from an open subset]
\label{prop:ucp}
Let $D \subset \Rn$ be open and connected, let $E \in \mathbb{C}$, and let
$u \in H^2_{\mathrm{loc}}(D)$ satisfy $(-\Delta - E)u = 0$ in $D$. If $u$
vanishes almost everywhere on some nonempty open subset $\omega \subset D$, then
$u = 0$ in $D$.
\end{proposition}

\begin{proof}
Interior elliptic regularity gives $u \in C^{\infty}(D)$. The operator
$-\Delta - E$ is elliptic with constant, hence real-analytic, coefficients, so
$u$ is real-analytic in $D$ \cite{MorreyNirenberg1957}. Let $A$ be the interior
of the set where $u$ vanishes. It is open and contains $\omega$, so it is
nonempty. It is also relatively closed in $D$. Indeed, let $x_0 \in D$ lie in
the closure of $A$. Every partial derivative of $u$ vanishes at $x_0$, being
continuous and identically zero on $A$, so the Taylor series of $u$ at $x_0$
vanishes identically. That series represents $u$ on a ball about $x_0$
contained in $D$, whence $x_0 \in A$. Since $D$ is connected, $A = D$.
\end{proof}

The proposition is applied four times, and to four different sets $D$. In the
proof of Corollary~\ref{cor:global-uniqueness}(a) it is applied with
$D = \Om$, which is connected by Assumption~\ref{assump:geometry}, and with
$\omega = \Omex{R_1}$ for the radius $R_1$ produced by
Theorem~\ref{thm:rellich-positive}. In the proof of
Theorem~\ref{thm:exterior-equality} it is applied with $D = \Omex{\Rzero}$, in
that of Lemma~\ref{lem:common-region} with $D = \Omcap$, and in that of
Lemma~\ref{lem:independence} again with $D = \Om$.

\begin{remark}[What the proposition uses]
\label{rem:ucp-scope}
Proposition~\ref{prop:ucp} uses neither the boundary condition nor any
regularity of $\Gpert$, and no continuation across $\Gpert$ or from Cauchy data
on it is needed. Because the coefficients are constant, real analyticity
replaces the Carleman estimates required for variable coefficients
\cite{Aronszajn1957,KochTataru2001}. A nonconstant impedance
$\alpha = \alpha(x)$ would leave the proposition applicable, whereas a
perturbation of the metric or a potential would not. The connectedness of $\Om$
is essential. If $\Om$ had a bounded component, that component would carry the
eigenfunctions of a Robin problem on a bounded domain, and the conclusion of
Corollary~\ref{cor:global-uniqueness}(a) would fail at the corresponding
energies.
\end{remark}

\section*{Statements and Declarations}

\noindent\textbf{Funding.}
The author did not receive support from any organization for the submitted work.

\smallskip
\noindent\textbf{Competing interests.}
The author has no competing interests to declare that are relevant to the
content of this article.

\smallskip
\noindent\textbf{Data availability.}
Data sharing is not applicable to this article, as no datasets were generated
or analysed during the current study.

\smallskip
\noindent\textbf{Use of AI tools.}
The author used an AI assistant (Claude, Anthropic) for language editing, for
drafting parts of the LaTeX source under the author's direction, and for
checking cross-references and citations. All mathematical content, proofs and
references were verified by the author, who takes full responsibility for the
manuscript.
\bibliographystyle{amsplain}
\bibliography{references}

\end{document}